\documentclass[10pt, twoside]{amsbook} 
\pdfmapfile{=cm-super-t1.map}
\pdfmapfile{=boondox.map}
\usepackage[paperwidth=7in, paperheight=10in, textwidth=4.75in, textheight=8in, includehead=true, bindingoffset=.5in]{geometry}
\makeatletter
\renewcommand{\descriptionlabel}[1]{\hspace\labelsep \upshape\bfseries #1} %CHANGE THIS LINE
\renewenvironment{description}{\list{}{%
  \advance\leftmargini6\p@ \itemindent-12\p@
  \labelwidth\z@ }%
}{
  \endlist
}
\makeatother
\makeindex
\usepackage[parfill]{parskip}
\usepackage{amsfonts}
\usepackage{amssymb}
\usepackage{amsthm}
\usepackage{textcomp}
\usepackage{blkarray}
\usepackage[scaled=.85]{beramono}% a typewriter font must be defined
\usepackage[leqno]{amsmath} %use eq no on left
\usepackage{url}
\usepackage[bb=boondox]{mathalfa}
\usepackage{bm}% load after all math to give access to bold math
\def\Mid{\,\vert\,}
\numberwithin{equation}{chapter}

\def\nmid{\not\mid}
\let\varmathbb\mathbb
\newcommand{\bF}{{\mathbb{F}}}
\newcommand{\bC}{{\mathbb{C}}}
\newcommand{\bR}{{\mathbb{R}}}
\newcommand{\bP}{{\mathbb{P}}}
\newcommand{\bN}{{\mathbb{N}}}
\newcommand{\bZ}{{\mathbb{Z}}}%multi-character has problems
\newcommand{\bK}{{\mathbb{K}}}
\newcommand{\bQ}{{\mathbb{Q}}}
\newcommand{\bL}{{\mathbb{L}}}

\usepackage{graphicx}
\usepackage{calc}
\usepackage{epstopdf}
\usepackage{microtype}

\newcommand{\beq}{\begin{equation}}
\newcommand{\eeq}{\end{equation}}
\newtheorem{thm}[equation]{Theorem}
\newtheorem{theorem}[equation]{Theorem}
\newtheorem{lem}[equation]{Lemma}
\newtheorem{lemma}[equation]{Lemma}
\newtheorem{cor}[equation]{Corollary}
\newtheorem{corollary}[equation]{Corollary}

\theoremstyle{remark}
\newtheorem{rem}[equation]{Remark}
\newtheorem{remark}[equation]{Remark}

\theoremstyle{definition}
\newtheorem{defn}[equation]{Definition}
\newtheorem{definition}[equation]{Definition}
\newtheorem{ex}[equation]{Exercise}
\newtheorem{exercise}[equation]{Exercise}
\newtheorem{examp}[equation]{Example}

\newtheorem{review}[equation]{Review}
\title{\textsf{\Huge Matrix Theory}\\[15pt]{\em proof techniques and intuition}}
\author{\textsf{\huge S. Gill Williamson}\\[95pt]\textsf{\small\textcopyright\ S. Gill Williamson}}
\date{}                                           % Activate to display a given date or no date
\usepackage{setspace}
\usepackage[pdftex,colorlinks=true, pdfstartview=FitV, linkcolor=blue, citecolor=blue, urlcolor=blue,bookmarks,bookmarksnumbered]{hyperref}
\hypersetup{
pdftitle={Matrix Theory},
pdfauthor={S. Gill Williamson}%{\textcopyright\ S. Gill Williamson}
}
\begin{document}
\thispagestyle{empty}
\begin{center}
\vspace*{2in}
\textsf{\Huge Matrix Canonical Forms}\\[.2in]
\textsf{\Large notational skills and proof techniques}\\[.7in]

\textsf{\huge S. Gill Williamson}
\vfill
\textcopyright  \textsf{S. Gill Williamson 2012. All rights reserved.}
\end{center}

\newpage
\thispagestyle{empty}
\hspace{1 pt}
\newpage
%PREFACE
\index{preface}
\begin{center}
{\huge Preface}
\end{center}
\pagestyle{plain}
This material is a rewriting of notes handed out by me to beginning graduate students in seminars in combinatorial mathematics (Department of Mathematics, University of California San Diego).  Topics covered in this seminar were in algebraic and algorithmic combinatorics.  Solid skills in linear and multilinear algebra were required of students in these seminars  - especially in algebraic combinatorics.    
I developed these notes to review the students' undergraduate linear algebra and improve their proof skills.  We focused on a careful development of the general matrix canonical forms as a training ground. 

I would like to thank Dr. Tony Trojanowski for a careful reading of this material and numerous corrections and helpful suggestions.  I would also like  to thank Professor Mike Sharpe, UCSD Department of Mathematics, for considerable LaTeX typesetting assistance and for his Linux Libertine font options to the newtxmath package.

S. Gill Williamson, 2012\\
\url{http://cseweb.ucsd.edu/~gill}
%${\rm cseweb.ucsd.edu\slash\sim gill\slash}$

\newpage
{\color{white}x}\newpage
\centerline{\Large CONTENTS}

\hspace*{1in}
%\section*[Contents]
\renewcommand\contentsname{} 
\tableofcontents
%\index{contents}
\hyperlink{index}{}
%\newpage
%\thispagestyle{empty}

%Begin usage examples
%\iffalse 
%\section{Preliminaries}
%\begin{defn}[\bfseries Differential] $dx$ means a little bit of $x$.\label{def:diffs}\end{defn}
%\beq \mbox{This equation left intentionally blank.}\eeq
%\begin{thm} [Differentiating a power of $x$] \label{def:diffpower}$\frac{d}{dx} x^n=n x^{n-1}$.\end{thm}
%\begin{proof} Obvious.\end{proof}
%For example,
%\beq \frac{d}{dx} x^3=3 x^2.\eeq
%\end{ex}
%\begin{ex} Calculate $\frac{d}{dx}x^4$.\end{ex}
%\fi
%End usage examples

%START OF CHAPTER 1
\chapter{Functions and Permutations}
\section*{Algebraic terminology}
In this first section,  we summarize for reference certain basic concepts in algebra.  These concepts are useful for the material we develop  here  and are essential for reading related online sources (e.g., Wikipedia).
\index{set notation!${\bN}=\{1,2, \ldots, \}$}
\index{set notation!${\bN}_0=\{0,1,2, \ldots, \}$}
\index{set notation!${\bZ}=\{0, \pm 1, \pm 2, \ldots\}$}
\index{set notation!$\underline{n}=\{1,\ldots, n\}$}
\index{set notation!$R^D$ functions from $D$ to $R$}
\index{delta!$\delta({\rm statement})$ $1$ true, $0$ false}
\begin{remark}[\bfseries Basic sets and notation]
\label{rem:basicsets}
We use the notation ${\bN}=\{1,2, \ldots \}$ for the positive integers.
Let ${\bN}_0=\{0,1,2, \ldots \}$ denote the nonnegative integers, and
let ${\bZ}=\{0, \pm 1, \pm 2, \ldots\}$ denote the set of all integers.
Let $\times^n S$ ($n$-fold Cartesian product of $S$) be the set of $n$-tuples 
 from a nonempty set $S$.  We also use $S^n$ for this cartesian product.  A slightly more general notation is to write $S^{\underline{n}}$ for this product where
 $\underline{n}=\{1,\ldots, n\}$ and the exponential notation $R^D$ denotes all $n$-tuples (i.e., functions~\ref{def:function})
 from $D$ to $R$ ( We use {\em delta notation}: $\delta({\rm Statement})=1$ if Statement is true, $0$ if Statement is false.)
 \end{remark}
\[
\framebox{\textbf{\emph{Semigroup$\rightarrow$Monoid$\rightarrow$Group: sets with one binary operation}}}
\]
\index{binary operation!associative}
\index{binary operation!commutative}
A function $w:S^2\rightarrow S$ is called a {\em binary operation}.
It is sometimes useful to
write $w(x,y)$ in a simpler form such as $x\,w\,y$ or simply $x\cdot y$ or even just $x\,y$.  To tie the binary operation $w$ to $S$ explicitly, we write $(S, w)$ or $(S, \cdot)$.
%semigroup
\index{semigroup!definition}
\index{semigroup!commutative (abelian)}
\index{semigroup!nonempty always}
\begin{defn}[\bfseries Semigroup]
\label{def:semigroup}
Let $(S,\cdot)$ be a {\em nonempty} set $S$ with a binary operation ``$\cdot$'' .  If $(x\cdot y)\cdot z = x\cdot (y\cdot z)$, for all $x$, $y$, $z\in S$, then the binary operation ``$\cdot$'' is called {\em associative} and $(S, \cdot)$ is 
called a {\em semigroup}.  If two elements, $s,\, t\in S$ satisfy $s\cdot t = t\cdot s$ then we say $s$ and $t$ {\em commute}.  If for all $x,\,y\in S$ we have $x\cdot y = y\cdot x$ then
$(S, \cdot)$ is a {\em commutative} (or {\em abelian}) semigroup.
\end{defn}

\begin{rem}[\bfseries Semigroup]
\label{rem:semigroup}
Let $S={\bf M}_{2,2}({\bZ}_{\rm e})$ be the set of $2\times 2$ matrices with entries in ${\bZ}_{\rm e}=\{0, \pm 2,\pm 4, \ldots \}$, the
set of even integers. Define $w(X,Y)=XY$ to be 
the standard multiplication of matrices (which is associative).  Then  $(S, w)$ is a semigroup. This semigroup is not commutative (alternatively, it is a noncommutative semigroup or a semigroup with non-commuting elements).  The semigroup of even integers, 
$({\bZ}_{\rm e}, \cdot$), 
where ``$\cdot$'' denotes multiplication of integers, is commutative.
\end{rem}
\[
\framebox{\textbf{\emph{Associative$\;\;+\;\;$Identity$\;\;=\;\;$ Monoid}}}
\]
%monoid
\index{monoid!definition}
\index{monoid!inverse element}
\begin{defn}[\bfseries Monoid]
\label{def:monoid}
Let $(S,\cdot)$ be a semigroup.  If there exists an element $e\in S$ such that for all $x\in S$,
$e\cdot x= x\cdot e=x$, then $e$ is called an {\em identity} for the semigroup.  A semigroup
with an identity is called a {\em monoid}.   If $x \in S$ and there is a $y\in S$ such that 
$x\cdot y=y\cdot x=e$ then $y$ is called an {\em inverse} of $x$. 
 \end{defn}
%\hspace*{1in}
%remark monoid
\index{monoid!uniqueness of inverse}
\begin{rem}[\bfseries Monoid]
\label{rem:monoid}
The identity is unique (i.e., if $e$ and $e'$ are both identities then $e=e\cdot e'=e'$).  
Likewise, if $y$ and $y'$ are inverses of $x$, 
then $y'=y'\cdot e=y' \cdot(x\cdot y) = (y' \cdot x)\cdot y = e\cdot y = y$  so the inverse of $x$ is unique.
Note that this last computation shows that if $y'$ satisfies $y'x=e$ ($y'$ is a ``left inverse'') and
$y$ satisfies $xy=e$ ($y$ is a ``right inverse'') then $y'=y$.
The $2\times 2$ matrices, ${\bf M}_{2,2}({\bZ})$, with matrix multiplication form
a monoid  (identity $I_2$, the $2\times 2$ identity matrix).
\end{rem}
\[
\framebox{\textbf{\emph{Associative$\;\;+\;\;$Identity$\;\;+\;\;$
Inverses$\;\;=\;\;$ Group}}}
\]
%group
\index{group!definition}
\begin{defn}[\bfseries Group]
\label{def:group}
Let $(S,\cdot)$ be a monoid with identity $e$ and let $x\in S$.  If there is a $y\in S$ such that 
$x\cdot y=y\cdot x=e$ then $y$ is called an {\em inverse} of $x$ (see \ref{def:monoid}).
A monoid in which every element has an inverse is a {\em group}.
\end{defn}
%remark group
\begin{rem}[\bfseries Group]
\index{group!commutative(abelian)}
Commutative groups, 
$x\cdot y = y\cdot x$ 
for all $x$ and $y$, play an important role in group theory.  They are also called {\em abelian} groups.  
Note that the inverse of an element $x$ in a group is unique: if $y$ and $y'$ are inverses of $x$, 
then $y'=y'\cdot e=y' \cdot(x\cdot y) = (y' \cdot x)\cdot y = e\cdot y = y$ 
(see \ref{rem:monoid}).
\end{rem}
\[
\framebox{\textbf{\emph{Ring: one set with two intertwined binary operations }}}
\]
%ring
\index{ring!definition}
\index{ring!definition}
\begin{defn}[\bfseries Ring and Field]
\label{def:ring}
A {\em ring}, $R=(S, +, \cdot)$, is a set with two binary operations such that $(S,+)$ is an abelian group with identity denoted by $0$ (``$+$'' is called ``addition'') and
$(S-\{0\},\cdot)$ is a semigroup (``$\cdot$'' is called ``multiplication''). 
The two operations are related by distributive rules which state that for all $x$, $y$, $z$ in 
$S$:
\[ 
{\bf (left)}\;\;
z\cdot (x+y)= z\cdot x + z\cdot y\;\;{\rm and}\;\;(x+y)\cdot z = x\cdot z + y\cdot z
\;\;{\bf (right)}.
\]
\end{defn}
%END DEFINITION
%REMARK
\index{ring!group of units}
\index{ring!commutative}
\index{ring!skew-field}
\index{ring!field}
\index{skew-field!also division ring}
\index{ring!invertible element (unit)}
\index{ring!with identity}
\begin{rem}[\bfseries Notation, special rings and group of units]
\label{rem:specialrings}
The definition of a ring assumes only that $(S,\cdot)$ is a semigroup  (\ref{def:semigroup}). Thus, a ring may not have a multplicative identity. We refer to such a structure as a {\em ring without an identity}.
Let $R=(S, +, \cdot)$ be a ring.
The identity of the abelian group $(S,+)$ is denoted by $0_R$  (or $0$)  and 
is called the {\em zero} of the ring $(S, +, \cdot)$. 
If $r\in S$ then the inverse of $r$ in $(S,+)$ is denoted by $-r$ so $r+(-r)=(-r)+r=0$.
Suppose $(S-\{0\},\cdot)$ is a monoid with identity $1_R$ 
(we say ``$R$ is a {\em ring with identity} $1_R$'');
its invertible elements (or {\em units}), $U(R)$, form a group, $(U(R), \cdot)$, with $1_R$ as the group identity.   
The group $(U(R), \cdot)$, or simply $U(R)$, is the {\em group of units} of the ring $R$.  
\index{division ring!also skew-field}
If $(S-\{0\},\cdot)$ is commutative then $(S, +, \cdot)$ is {\em a commutative ring}.
If $(S-\{0\},\cdot)$ is a group (i.e., $(U(R), \cdot) = (S-\{0\},\cdot)$) then the ring $R$ is called a {\em skew-field} or {\em division ring}.  
If this group is {\em abelian} then the ring is called a {\em field}.
\end{rem}
%END REMARK
%REMARK
\begin{rem}[\bfseries Basic ring identities]
\label{rem:basicidents}
If $r,s,t$ are in a ring  $(S, +, \cdot)$ then the following basic identities (in braces, plus hints for proof) hold:
%DESCRIPTION ITEMS
\index{ring!basic identities}
\begin{description}
\item[${\bf (1)}\;\{r\cdot 0 = 0\cdot r=0\}$] If $x+x=x$ then $x=0$.  
Take $x=r\cdot0$ and $x=0\cdot r$.
\item[${\bf (2)}\;\{(-r)\cdot s = r\cdot (-s)=-(r\cdot s)\}$] $r\cdot s +  (-r)\cdot s = 0 \implies
(-r)\cdot s = -(r\cdot s)$.
\item[${\bf (3)}\;\{(-r)\cdot (-s) = r\cdot s\}$] Replace $r$ by $-r$ in {\bf (2)}.  
Note that $-(-r)=r$.
\end{description}
In particular, if 
$(S-\{0\},\cdot)$ has identity $1_R\equiv 1$, then
$(-1)\cdot a = -a$ for any $a\in S$ and, taking $a=-1$, $(-1)\cdot(-1)  = 1$.
It is convenient to define $r-s = r+(-s)$.  
Then we have $(r-s)\cdot t=r\cdot t - s\cdot t$ and $t\cdot(r-s)=t\cdot r - t\cdot s$:
\[t\cdot(r-s)=t\cdot( r+(-s))=t\cdot r + t\cdot(-s)) = (t\cdot r + -(t\cdot s))= t\cdot r - t\cdot s.\]
\end{rem}
\[
\framebox{\textbf{\emph{A field is a ring $(S, +, \cdot)$ where  $(S-\{0\},\cdot)$ is an abelian group}}}
\]
%remark ring
\index{ring!examples}
\begin{rem}[\bfseries Rings, fields, identities and units]
\label{rem:ringfield}
The  $2\times 2$ matrices over the even integers, ${\bf M}_{2,2}({\bZ}_{\rm e})$, with the usual multiplication and addition of matrices, is a noncommutative ring without an identity.  
The matrices, ${\bf M}_{2,2}({\bZ})$, over all integers, is a noncommutative ring with an identity. 
The group of units, $U({\bf M}_{2,2}({\bZ}))$, is all invertible $2\times 2$ integral matrices.
The  matrix  $P=\left(\begin{array}{cc} +1&-1\\-2& +3 \end{array}\right)$ is  a unit in
${\bf M}_{2,2}({\bZ})$ with $P^{-1}=\left(\begin{array}{cc} 3&1\\2& 1 \end{array}\right).$
 $U({\bf M}_{2,2}({\bZ}))$ is usually denoted by $\rm{GL}(2,\bZ)$ and is called a 
 {\em general linear group}. 
The ring  of $2\times 2$ matrices of the form 
$
\left(
\begin{array}{rc}
x & y\\
-\overline{y}&\overline{x}
\end{array}
\right)
$
\index{ring!skew-field!not field}
where $x$ and $y$ are complex numbers is a skew-field but not a field. This skew-field is equivalent to (i.e, a ``matrix representation of'') the {\em skew field of quaternions} (see Wikipedia article on quaternions).  The most important fields for us will be the fields of real and complex numbers.
\end{rem}
%Definition Ideal
\index{subring!definition}
\index{ideal!definition}
\begin{defn}[\bfseries Ideal]
\label{def:subringideal}
Let $(R, +, \cdot)$ be a ring and let $A\subseteq R$ be a subset of $R$.  
If  $(A,+)$ is a subgroup of $(R,+)$ then $A$ is a {\em left ideal} if  for every $x\in R$ and $y\in A$, $xy\in A$.  
A {\em right ideal} is similarly defined.  If $A$ is both a left and right ideal then it is a 
{\em two-sided ideal} or, simply, an {\em ideal}. 
If $(R, +, \cdot)$ is commutative then all ideals are two sided.
\end{defn}
%Remark Ideal   
\index{ideal!examples}
\begin{rem}[\bfseries Ideal]
\label{rem:subringideal}
The set $A$ of all matrices  
$
a=
\left(
\begin{array}{cc}
x & y\\
0 & 0
\end{array}
\right)
$
forms a subgroup $(A,+)$ of $({\bf M}_{2,2}({\bZ}),+)$.  
The subset $A$ is a right ideal but not a left ideal of the ring   ${\bf M}_{2,2}({\bZ}).$
Note that $R=(A, +, \cdot)$ is itself a ring.  This ring has pairs of zero divisors -
pairs of elements $(a,b)$ where $a\neq 0$ and $b\neq 0$ such that $a\cdot b = 0$.   
For example, take the pair $(a,b)$ to be
$
a=\left(
\begin{array}{cc}
0 & 1\\
0 &0
\end{array}
\right)
$
and
$
b=\left(
\begin{array}{cc}
1& 1\\
0 &0
\end{array}
\right).
$
\index{zero divisor pairs}
The pair $(a,b)$ is, of course, also a pair of zero divisors in ${\bf M}_{2,2}({\bZ})$.

Another example of an ideal is 
the set of even integers, ${\bZ_e}$, which is a subset of the integers, ${\bZ}$ (which, it is worth noting, forms a ring with {\em no zero divisor pairs}).  
The subset ${\bZ_e}$, is an ideal (two-sided) in  ${\bZ}$. 
Given any integer $n\neq 0$, the set $\{k\cdot n\Mid k\in {\bZ}\}$ of multiples of $n$, is an ideal of the ring ${\bZ}$ which we denote by $(n)=n{\bZ}= {\bZ}n$.  
\index{ideal!principal}
Such an ideal (i.e., generated by a single element, $n$) in ${\bZ}$ is called a {\em principal ideal}.  
It is easy to see that all ideals $A$ in ${\bZ}$ are principal ideals, $(n)$, where $|n|>0$ is  minimal over the set $A$.  Another nice property of integers is that they {\em uniquely factor} into primes (up to order and sign).
\end{rem} 
\index{ring!characteristic of}
\begin{defn}[\bfseries Characteristic of a ring]
\label{def:characteristic}
 Let $R$ be a ring. Given $a\in R$ and an integer $n>0$, define $na\equiv a+ a +\cdots a$ where there are $n$ terms in the sum.  If there is an integer $n>0$ such that $na=0$ for all $a\in R$ then the {\em characteristic} of $R$ is the least such $n$.  If no such $n$ exists, then $R$ has {\em characteristic zero} 
 (see Wikipedia article ``Characteristic (algebra)'' for generalizations). 
\end{defn}

Algebraists have defined several important abstractions of the ring of integers, 
${\bZ}$.  We next discuss four such abstractions: integral domains, principal ideal domains (PID), unique factorization domains (UFD), and Euclidean domains - each more restrictive than the other.
\[
\framebox{\textbf{\emph{Euclidean Domain $\implies$ PID $\implies$ UFD }}}
\]
%DEFINITION
\index{ring!commutative!integral domain}
\index{integral domain!definition}
\begin{defn}[\bfseries {\bf Integral domain}] 
\label{def:integraldomain}
An {\em integral domain} is a {\em commutative} ring with identity, $(R, +, \cdot)$, in which there are no  {\em zero divisor pairs}: pairs of nonzero elements $(a,b)$ where $ab=0.$ (See \ref{rem:subringideal} for ring {\em with} pairs of zero divisors.)  
\end{defn}   
%END DEFINITION
%DEFINITION

%END DEFINITION
%REMARK
\index{ring!commutative!divides relation}
\index{ring!commutative!zero divisor pairs}
\index{ring!commutative!irreducible element}
\index{ring!commutative!prime element}
\index{ring!commutative!associates}
\begin{rem}[\bfseries Divisors, units, associates and primes] 
\label{rem:unitassoc} 
For noncommutative rings, an element $a$ is a left zero divisor if there exists $x\neq 0$ such that $ax\neq 0$ (right zero divisors similarly defined).
Let $R$ be a commutative ring.  
If $a\neq 0$ and $b$ are elements of a $R$, we say that 
$a$ is a {\em divisor} of $b$ (or $a$ {\em divides} $b$),
$a\Mid b$, if there exists $c$ such that $b =ac$.  
Otherwise, $a$ does not divide $b$, $a\nmid b$.
Note that if $b=0$ and $a\neq 0$ then $a\Mid b$ because $b=a\,0$ ($c=0$). Thus $a\Mid 0$ or $a$ divides $0$.
(The term {\em zero divisors} of \ref{def:integraldomain} refers to pairs $(a,b)$ of nonzero
elements and is not the same as ``$a$ is a divisor of $0$'' or ``$a$ divides $0$''.)
An element $u$ in $R-\{0\}$ is an {\em invertible element} or a {\em unit of $R$} if $u$ has an inverse in $(R-\{0\},\cdot)$. The units form a group, $U(R)$ (\ref{rem:specialrings}).
For commutative rings, $ab=u$, $u$ a unit, implies that both $a$ and $b$ are units: $ab=u$ implies $a(bu^{-1})= e$ and $(au^{-1})b = e$ so both $a$ and $b$ are units.
Two elements $a$ and $b$ of $R$ are {\em associates in $R$} if $a=bu$ where $u$ is a unit.
An element $p$ in $R-\{0\}$ is  {\em irreducible} if $p=a b$ implies that either $a$ or $b$ is a unit and {\em prime} if $p\Mid a b$ implies $p\Mid a$ or $p\Mid b$.  
For unique factorization domains (\ref{def:ufd}), $p$ is irreducible if and only if it is prime.
In the ring ${\bZ}$, the only invertible elements are $\{+1, - 1\}$.  
The only associates of an integer $n\neq 0$ are $+n$ and $-n$.  
The integer $12=3\cdot 4$ is the product of two non-units so $12$ is not irreducible (i.e., {\em is} reducible) or, equivalently in this case, not a prime.  
The integer $13$  is a prime with the two associates $+13$ and $-13$.  
A field is an integral domain in which every nonzero element is a unit.  
In a field, if $0\neq p = ab$ then {\em both} $a$ and $b$ are nonzero and hence both are units (so at least one is a unit) and thus every nonzero element  in a field is irreducible (and prime).\\
\end{rem}
%END REMARK
%DEFINITION
\index{ring!commutative!UFD}
\begin{defn}[\bfseries Unique factorization domain] 
\label{def:ufd}
An integral domain $R$ is a {\em unique factorization domain} (UFD) if 
\begin{description}
\item[(1)] Every $0\neq a\in R$ can be factored into a finite (perhaps empty) product of primes and a unit: $a=up_1\cdots p_r.$ %(\ref{rem:unitassoc}).  
An empty product is defined as $1_R$.
 \item[(2)] If $a=up_1\cdots p_r$ and $a=wq_1\cdots q_s$ are two such factorizations
 then either both products of primes are empty (and $u=w$) or
 $r=s$ and the $q_i$ can be reindexed so that $p_i$ and $q_i$ are associates for
 $i=1, \ldots, s$.\\
\end{description}
\end{defn} 
%END DEFINITION
%REMARK
\index{ring!commutative!UFD examples}
\begin{remark}[\bfseries Unique factorization domains]
\label{rem:ufd}
The integers, ${\bZ}$, are a unique factorization domain.  Every field is also a unique factorization domain because every nonzero element is a unit times the empty product.  If $R$ is a UFD then so are the polynomial rings $R[x]$ and 
$R[x_1, \ldots, x_n].$  
If $a_1, \ldots, a_n$ are nonzero elements of a UFD, then there exists  a greatest common divisor $d={\rm gcd}(a_1, \ldots, a_n)$ which is unique up to multiplication by units.
The divisor $d$ is {\em greatest} in the sense that any element $\hat{d}$ such that
$\hat{d}\Mid a_i$, $i=1,\ldots n$, also divides $d$ (i.e., $\hat{d}\Mid d$).\\ 
\end{remark}
%END REMARK
%DEFINITION
\index{ring!commutative!PID}
\begin{defn}[\bfseries Principal ideal domain] 
\label{def:pid}
An integral domain $R$ is a {\em principal ideal domain} (PID) if every ideal in $R$ is a 
principal ideal (\ref{rem:subringideal}).\\
\end{defn} 
%END DEFINITION
%REMARK
\begin{remark} [\bfseries Principal ideal domains] 
\label{rem:pid}
We noted in Remark~\ref{rem:subringideal} that every ideal in ${\bZ}$ is a principal ideal.
If $(F, +, \cdot)$ is a field, then any ideal $A\neq (0)$ contains a nonzero and hence invertible element $a$.  The ideal $(a) = F$.  
There is only one nontrivial ideal in a field and that is a principal ideal that equals $F$.  Thus, any field $F$ is a PID.
Let $a_1, \ldots, a_n$  be nonzero elements of a PID, $R$.  
It can be shown that if  $d={\rm gcd}(a_1, \ldots, a_n)$ in $R$ then there exists 
$r_1, \ldots, r_n$ in $R$ such that $r_1a_1 + \cdots + r_na_n = d$.   
The ring of polynomials, $F[x_1,\ldots x_n]$, in $n\geq 2$ variables over 
field $F$ is not a PID.  
Also, the ring of polynomials with integral coefficients, $\bZ[x]$, is not a PID.
For example, the ideal $<2,x>=\{2a(x)+xb(x)\Mid a, b \in\bZ[x]\}$ is not a principal ideal $(p(x))$, $p\in\bZ[x]$. 
Otherwise, $2=q(x)p(x)$ for some $q\in \bZ[x]$ which implies $p=\pm1$ or
$p=\pm 2$, either case leading to a contradiction.

\end{remark}
%END REMARK
%DEFINITION
\begin{defn}[\bfseries Euclidean valuation] 
\label{def:euclidean}
Let $R$ be an integral domain and let $\nu: R-\{0\} \rightarrow \bN_0$ (nonnegative integers).  $\nu$ is a  {\em valuation} on $R$ if
\begin{description}
\item[(1)] For all $a, b \in R$ with $b\neq 0$, there exist $q$ and $r$ in $R$ such that
$a=b\cdot q + r$ where either $r=0$ or $\nu(r) < \nu( b)$.
\item[(2)] For all $a, b \in R$ with $a\neq 0$ and $b\neq 0$, $\nu(a)\leq \nu(a\cdot b)$.\\
 \end{description}
\end{defn} 
%END DEFINITION
%DEFINITION
\index{ring!commutative!Euclidean}
\begin{defn}[\bfseries Euclidean domain] 
\label{def:eud}
An integral domain $R$ is a {\em Euclidean domain} if there exists a Euclidean valuation on $R$ (see \ref{def:euclidean}).\\
\end{defn} 
%END DEFINITION
%REMARK
\begin{remark}[\bfseries Euclidean domains] 
\label{rem:eud}
It can be shown that every Euclidian domain is a principal ideal domain and every principal ideal domain is a unique factorization domain.  The integers ${\bZ}$ are a Euclidean domain with $\nu(n) = |n|$.  The units of $\bZ$ are $\{-1, +1\}.$  In general in a Euclidean domain, if $a$ and $b$ are nonzero and $b$ is not a unit then $\nu(a) < \nu(ab)$  (check this out for $\bZ$).
Intuitively, the units have minimal valuations over all elements of the Euclidean domain and multiplying any element by a non-unit increases the valuation.
Any field $(F, +, \cdot)$ is a Euclidean domain with $\nu(x)=1$ for all nonzero $x$.
The polynomials, $\bF[x]$, with coefficients in a field $\bF$ form a Euclidean domain with $\nu(p(x))$ the degree of $p(x)$.  The units of $\bF[x]$ are all nonzero constant polynomials (degree zero).
The ring of polynomials with integral coefficients, $\bZ[x]$, is  not a PID (\ref{rem:pid}) 
and thus not a Euclidean domain.
 Likewise, the ring of polynomials in $n$ variables, $n>1$, over a field $\bF$, 
 $\bF[x_1, \ldots, x_n]$, is not a PID (\ref{rem:pid}) and hence not a Euclidean domain.  
Rings that are PIDs but not Euclidean domains are rarely discussed (the ring $\bZ[\alpha] =\{a + b\alpha\Mid a, b \in \bZ,\; \alpha = (1+ (19)^{1/2} i)\}$ is a standard example).
\end{remark}
%END REMARK
\[
\framebox{\textbf{\emph{We now combine a ring with an abelian group to get a module.}}}
\]
%DEFINITION
\index{module!definition}
\begin{defn}[\bfseries Module]
\label{def:module}
Let $(R, +, \cdot)$ be a ring with identity $1_R$.
 Let $(M, \oplus )$ be an abelian group.
We define an operation with domain $R\times M$ and range $M$ which for each 
$r\in R$ and $x\in M$ takes $(r,x)$ to $rx$ (juxtaposition of $r$ and $x$).  
This operation, called {\em scalar multiplication}, defines a {\em left} $R$-{\em module} $M$ if the following hold for every 
$r,s \in R$ and $x,y \in M$:
\[
(1)\;r(x\oplus y) = rx\oplus ry\;\;(2)\;(r+s)x = rx\oplus sx\;\;(3)\;(r\cdot s)x=r(sx)\;\;(4)\;1_Rx=x.
\]
 We sometimes use ``$+$'' for the addition in both abelian groups and replace ``$\cdot$'' with juxtaposition.  
Thus, we have: $(2)\;(r+s)x = rx+sx\;\;(3)\;(rs)x=r(sx).$ 
\end{defn}
%END DEFINITION
Sometimes a module is defined without the assumption of the identity $1_R.$ 
In that case, what we call a module  is called a {\em unitary} module.
%REMARK
\index{module!examples}
\begin{rem}[\bfseries Module]
\label{rem:modulerank}
Let $R$ be the ring of $2\times 2$ matrices over the integers, 
${\bf M}_{2,2}({\bZ})$.  
Let $M$ be the abelian group, ${\bf M}_{2,1}({\bZ})$, of $2\times 1$ matrices under addition.  Then $(1)$ and $(2)$ correspond to the distributive law for matrix multiplication, 
$(3)$ is the associative law, and $(4)$ is multiplication on the left by the $2\times 2$ identity matrix.  Thus, ${\bf M}_{2,1}({\bZ})$ is an ${\bf M}_{2,2}({\bZ})$-module.
If $x\in {\bf M}_{2,1}({\bZ})$ then, obviously, $x+x=x$ implies that $x=\theta_{21}$ 
where $\theta_{21}$ is the zero matrix in ${\bf M}_{2,1}({\bZ})$.  To see this, just add $-x$ to both sides of $x+x=x$.  This fact is true in any module for the same reason.  
In particular, in any module if $z\in M$ and $0\in R$ is the identity of $(R, +)$, then
$0z = (0+0)z= 0z + 0z$ and, taking $x=0z$ in the identity $x+x=x$, $0z=\theta$, the zero in $(M,+)$.  This fact is obvious in our ${\bf M}_{2,2}({\bZ})$-module ${\bf M}_{2,1}({\bZ}).$
Likewise for any module, if $\alpha \in R$ and $\theta$ is the zero in $(M,+)$, then
$\alpha\theta = \alpha(\theta + \theta)= \alpha\theta + \alpha\theta$ implies that $\alpha\theta=\theta$.
In the modules that we will be interested in (i.e., vector spaces~\ref{def:vectoralgebra}) it is true that
for $\alpha \in R$ and $z\in M$, $\alpha z = \theta$ implies that either $\alpha=0$ or  $z=\theta$. 
This assertion is not true in our  ${\bf M}_{2,2}({\bZ})$-module ${\bf M}_{2,1}({\bZ}):$
\begin{equation*}
\alpha z =
\left(
\begin{array}{cc}
1&0\\
1&0
\end{array}
\right)
\left(
\begin{array}{c}
0\\
1
\end{array}
\right)
=
\left(
\begin{array}{c}
0\\
0
\end{array}
\right)
=\theta_{21}.
\end{equation*}
\hspace*{1 in}
\end{rem}
\index{module!free}
\begin{rem}[\bfseries Free modules]
\label{rem:fremod}
Of special interest to us are certain $R$ modules,  $R^n$ (see \ref{rem:basicsets}), where $R$ is a ring with identity $1_R$.  The abelian groups of these modules consist 
of $n$-tuples (n-vectors)  of elements in $R$ where addition is  component-wise. 
The multiplication of n-tuples $(x_1, \ldots, x_n)$ by elements $\alpha$ of $R$  is  defined component wise:
 $\alpha(x_1, \ldots, x_n)= (\alpha x_1, \ldots, \alpha x_n)$.  Such modules are called {\em free modules of rank} $n$ over $R$.  For a careful discussion see Wikipedia  ``Free module.''
\end{rem}
%END REMARK
%DEFINITION
\index{vector space and algebra!definition}
\begin{defn}[\bfseries Vector space and algebra]
\label{def:vectoralgebra}
If an abelian group $(M, + )$ is an $F$-module where $F$ is a field (\ref{rem:ringfield}), then we say $(M, + )$ (or, simply, $M$) is a {\em vector space over} $F$ (or $M$ is an $F$ {\em vector space}).
Suppose $(M, + , \cdot )$ is a ring where $(M, + )$ is a vector space over $F$. Then $(M, + , \cdot )$ is an {\em algebra over} $F$ (or $M$ is an $F$ {\em algebra}) if the following scalar rule holds:
\begin{description}
\item[scalar rule] for all $\alpha \in F$, $a, b \in M$ we have  $\alpha(a\cdot b)=(\alpha a)\cdot b = a \cdot (\alpha b)$.\\
\end{description}
\end{defn}
%END DEFINITION
%REMARK
\begin{remark}[\bfseries Vector spaces versus modules]
\label{rem:vecspamod}
We will use certain  modules over special Euclidean domains (e.g., integers and polynomials)
as  a tool to understand properties of finite dimensional vector spaces.  One basic difference between finite dimensional vector spaces and  general modules is that a proper subspace of such a vector space always has lower dimension (rank) than the vector  space itself -- not so in general for modules.  
As an example, consider the integers  $\bZ.$ 
The ordered pairs $\bZ\times \bZ$ is   a $\bZ$ module with the usual operations  on ordered pairs (free module of rank $2$ over $\bZ$).  
The natural ``module basis'' is $E=\{(0,1), (1,0)\}$ so this module has ``rank $2$''.
Take $E'=\{(0,1), (2,0)\}$.  
The span of $E'$ is a proper submodule of $\bZ\times \bZ$ over
the integers $\bZ$ since the first  component of every element in the span is even, but the ``module rank'' of this proper submodule is still $2$.  If we had used the field of rational numbers $\bQ$ instead and regarded $E'$ as a set in the vector space 
$\bQ\times \bQ$ over $\bQ$ then the span of $E'$ is the entire vector space, not a proper subspace as in the case of $\bZ$.
We are not defining our terms here, but the basic idea should be clear.
\end{remark}
%REMARK
\begin{rem}[\bfseries Complex matrix algebra]
\label{rem:vectoralgebra}
Let ${\bC}$ be the field of complex numbers and let $M$ be ${\bf M}_{2,2}({\bC})$, the 
additive abelian group of $2\times 2$ matrices with complex entries.  
Conditions (1) to (4) of \ref{def:module} are familiar properties of multiplying matrices by scalars (complex numbers).
Thus, $M$ is a complex vector space or, alternatively, $M$ is a vector space over the field of complex numbers, ${\bC}$.  If we regard $M$ as the ring, ${\bf M}_{2,2}({\bC})$, of 
$2\times 2$ complex matrices using the standard multiplication of matrices, then
it follows from the definitions of matrix multiplication and multiplication by scalars that the scalar rule of \ref{def:vectoralgebra} holds, and ${\bf M}_{2,2}({\bC})$ is an algebra
over $\bC$.
\end{rem}
%END REMARK 
\index{rings $\bK$ used in this book}
\begin{remark}[\bfseries Rings $\bK$ of interest to us]
\label{rem:specialringsfields}
Fields $\bF$ of interest will be $\bQ$ (rational numbers), $\bR$ (real numbers), 
$\bC$ (complex numbers) and the fields of rational functions (ratios of polynomials) 
$\bQ(z)$, $\bR(z)$ and $\bC(z)$.  
Thus,  $\bF\in\{\bQ, \bR, \bC, \bQ(z), \bR(z),\bC(z)\}$. 
Let $\bK\in \{\bZ, \bF[x], \bF\}$ where $\bZ$ denotes the integers and $\bF[x]$ the  polynomials over $\bF$.  
All of these rings are of characteristic zero (Definition~\ref{def:characteristic}).
Note that $\bK$ is a Euclidean domain with valuation absolute value (integers), degree (polynomials), or identically $1$ for all nonzero elements (field).
The fields $ \bQ(z), \bR(z),\bC(z), \bQ$  are the {\em quotient fields} for the Euclidean domains
$\bQ[z], \bR[z], \bC[z], \bZ,$  respectively.  We will not have much interest in the polynomials
$\bF[x]$ or rational functions $\bF(x)$ where $\bF = \bQ(z), \bR(z),\bC(z).$

\end{remark}
\index{gcd and lcm!divisors and multiples}
\begin{remark}[\bfseries Greatest common divisor, least common multiple]
\label{rem:gcdlcmmeetjoin}
Suppose $d$ is a common divisor of $a, b\neq 0$ in $\bK$.
If for all $c$, $c\Mid d$ whenever $c\Mid a$ and $c\Mid b$ then $d$ is a {\em greatest common divisor}  of $a$ and $b$  ($d$ is a $\gcd(a,b)$).  
If $d$ is a  $\gcd (a,b)$ then $td$ is a  $\gcd(ua, vb)$ for any units $t, u, v\in \bK$ 
(i.e., the $\gcd$ is defined ``up to units''). 
If $a\neq 0$ we adopt the convention that $a$ is a $\gcd(a,0)$.  
An element $c\in \bK$ is a {\em common multiple} of $a$ and $b$ if $a|c$ and $b|c$.    
If for all $x\in \bK$,   $a|x$ and $b|x$ implies $c|x$ then $c$ is a 
{\em least common multiple} of $a$ and $b$ (c is an \rm{lcm}(a,b)).
The $\mathrm{lcm}(a,b)$ is determined up to units.
We also write $\rm{lcm}(a,b)=a \lor b$ (``join'' of $a$ and $b$) 
and $\rm{gcd}(a,b)=a \land b$ (``meet'' of $a$ and $b$).  

If $\bK=\bF$ then all nonzero elements  are units  so $d$ is a $\gcd(a,b)$ for any nonzero $d,a,b$.  Likewise,  $d$ is an $\mathrm{lcm}(a,b)$ for any nonzero $d,a,b$.   
For $\bK=\bZ$, the units are $\{+1, -1\}$, and for $\bK=\bF[x]$, the units are the nonzero constant polynomials.
Thus, we focus on the cases $\bK=\bZ\; {\rm or} \; \bF[x]$:   

Suppose $a=p_1^{e_1} \cdots p_m^{e_m}$ and $b=p_1^{f_1} \cdots p_m^{f_m}$ are
prime factorizations of $a$ and $b$.  
Then  $ p_1^{\min (e_1,f_1)} \cdots p_m^{\min (e_m,f_m)}$ is a $\rm{gcd}(a,b)$,
and  $p_1^{\max (e_1,f_1)} \cdots p_m^{\max (e_m,f_m)}$ is an $\rm{lcm}(a,b)$.
Note that $ab=\gcd(a,b)\mathrm{lcm}(a,b)=(a \land b)(a \lor b).$
If $\bK=\bZ$ then let $d>0$ be the largest positive divisor  of $a$ and $b$.  
The notation $d=\gcd(a,b)$ is sometimes used to indicate that this is the canonical choice up to units for a $\gcd$ of the integers $a$ and $b$.   For $\bK=\bF[x]$ the canonical choice up to units is often taken to be the monic polynomial $d$ (coefficient  of highest power one).
\end{remark}

%END ALGEBRA
\section*{Sets, lists, multisets and functions}
\label{sec:setsandfunctions}
We consider collections of objects where order and repetition play different roles.
The concept of a {\em function} and the terminology for specifying various types of functions are discussed.
 %Begin Remark 
\index{notation!lists, multisets, functions}
\index{multisets}
\begin{rem}[\bfseries Notation for sets, lists and multisets]
\label{rem:notation}
The empty set is denoted by $\emptyset$. Sets are specified by braces: $A=\{1\}$, $B=\{1,2\}$.  They are unordered, so
$B=\{1,2\}=\{2,1\}$.  Sets $C=D$ if  $x\in C\; \mathrm{implies}\; x\in D$ 
(equivalently, $C\subseteq D$) and
$x\in D\; \mathrm{implies}\; x\in C$.  If you write $C=\{1,1,2\}$ and $D=\{1,2\}$ then, by the definition of set equality, $C=D$.  A {\em list}, {\em vector}  or {\em sequence} (specified by parentheses) is ordered:
$C' = (1,1,2)$ is not the same as $(1,2,1)$ or $(1,2)$.  Two lists (vectors, sequences), 
$(x_1, x_2, \ldots , x_n) = (y_1, y_2, \ldots , y_m)$, are equal if and only if $n=m$ and
$x_i = y_i$ for $i=1, \ldots , n$.  
A list such as $(x_1, x_2, \ldots , x_n)$ is also written $x_1, x_2, \ldots , x_n$, without the parentheses.
There are occasions where we discuss collections of objects where, like sets,
order doesn't matter but, like lists, repetitions do matter.  
These objects are called {\em multisets}.  
A  multiset can be specified by giving the elements of the multiset with each repeated a certain number of times.  
For example, $\{1,1,2,2,2,3,3,3\}$ is a multiset with $1$ twice, $2$ three, and $3$ three times.  In this case, $\{1,1,2,2,2,3,3,3\}\neq \{1,2,3\}$ but
$\{1,1,2,2,2,3,3,3\} = \{1,2,1,3,2,3,3,2\}.$  
We say $2$ is a member of $\{1,1,2,2,2,3,3,3\}$ with repetition $3$.
The size of the multiset $\{1,1,2,2,2,3,3,3\}$ is $2+3+3=8$ (sum of the distinct repetition numbers).
The use of braces to define multisets is like the use of braces to define sets.  You must make clear in any discussion whether you are discussing sets or multisets.
The union of two multisets combines their elements and their multiplicities:
$\{1,1,2,2,2,3\} \cup \{1,2,2,3,3\} = \{1,1,1, 2,2,2,2,2,3,3,3\}$.

If $S$ is a finite set, then $|S|$ is the number of elements in $S$.
Thus, $|\{1,1,2\}|=|\{1,2\}|=2.$  $\bP(S)$ is the set of all subsets of $S$, 
and $\bP_k(S)$ is all subsets of  $S$ of size $k$.  
If $|S|=n$ then $\left|\bP(S)\right| = 2^n$    and   
$\left|\bP_k(S)\right| = \binom{n}{k}$
%\left(\begin{array}{c} n\\k\end{array}\right)$ 
(binomial coefficient).
We use {\em underline} notation, $\underline{n}= \{1, 2, \ldots, n\}.$
We sometimes leave off the underline when the meaning is clear: $\bP_k(\underline{n}) = \bP_k(n).$
\end{rem}
\index{partition of set}
\begin{defn}[\bfseries Partition of a set]
\label{def:setpartition}
A partition of a set $S$ is a collection, ${\mathcal B}(S)$, of nonempty subsets, $X$, of $S$ such that each element of $S$ is contained in {\em exactly one}  set
$X\in{\mathcal B}(S)$. 
The sets $X\in{\mathcal B}(S)$ are called the {\em blocks} of the partition
 ${\mathcal B}(S)$ .  
A set $D\subseteq S$ consisting of exactly one   element from each block is called a  {\em system of distinct representatives} (or ``SDR'') for the partition.
\end{defn}

\begin{remark}[\bfseries Partition examples]
\label{rem:partitionexamples}
${\mathcal B}(S) = \{\{a, c\}, \{b,d,h\}, \{e,f,g\}\}$ is a partition of the set 
$S=\{a,b,c,d,e,f,g,h\}$. The set $\{c,d,g\}$ is an SDR for this partition.
Note that $\{\{a, c\}, \{a,c\}, \{b,d,h\}, \{e,f,g\}\}$ is a partition of $S$ and is the same as
$\{\{a, c\}, \{b,d,h\}, \{e,f,g\}\}$.  (Recall that repeated elements in a description of a set count as just one element.)
\end{remark}
\index{equivalence relations}
\begin{rem}[\bfseries Equivalence relations]
\label{rem:equivrel}
Consider the partition of~\ref{rem:partitionexamples}.  
For each pair $(x, y)\in S$, write $x\sim y$ if $x$ ``is in the same block as'' $y$ and 
$x\nsim y$ if $x$ ``is not in the same block as'' $y$.  For all $x, y, z \in S$ we have
\begin{equation}
{\rm (1)}\;x\sim x\;\;\;{\rm(2)}\;x\sim y\implies y\sim x \;\;\;(3)\;
x\sim y\;{\rm and}\;y\sim z \implies x\sim z.
\end{equation}
Condition (1) is called  {\em reflexive}, (2) is called {\em symmetric}, and (3) is called 
{\em transitive}.  Any relation defined for all $(x,y)\in S$ which can be either true or
false can be written as $x\sim y$ if true  and $x\nsim y$ if false.  
If such a relation is reflexive, symmetric and transitive then it is called an
{\em equivalence relation}.   Every equivalence relation can be thought of as
a partition of $S$.   
As an example, take the  ``is in the same block as'' equivalence relation.  Suppose we just knew how to check if two things were in the same block but didn't know the blocks.  We could reconstruct the set of blocks (i.e., the partition) by taking each $x\in S$ and constructing the set
$E_x=\{y\Mid x\sim y\}$.  
\index{equivalence class}
This block is called the ``equivalence class'' of $x$. 
The  partition could be reconstructed as the set of equivalence classes:
\begin{equation}
\label{eq:setequivclasses}
{\mathcal B}(S) = \{E_x\Mid x\in S\}=\{E_a, E_b, E_c,  E_d, E_e, E_f, E_g, E_h\}.
\end{equation}
In this list we have duplicate blocks (e.g., $E_a=E_c$).  But duplicates count as the same 
element in set notation: $\{1,1,2\}=\{1,2\}$ (\ref{rem:notation}). You should carry out the construction and proof of \ref{eq:setequivclasses} for the general equivalence class.  
You will need to use Definition~\ref{def:setpartition}.
Wikipedia has a good article.
\end{rem}
\index{function!definition}
\begin{defn}[\bfseries \bf Functions] 
\label{def:function}
Let $A$ and $B$ be sets.   
A function $f$ from $A$ to $B$ is a rule that assigns to each element 
$x\in A$ a unique element $y \in B$.  We write $y=f(x)$.
Two functions $f$ and $g$ from $A$ to $B$ are equal if $f(x)=g(x)$ for all $x\in A$.
\end{defn}

This definition is informal as it uses ``rule,'' ``assign'' and ``unique'' intuitively, but that is good enough for us as we shall give many examples.  
Given a function $f$ from $A$ to $B$, we can define a set $F\subseteq A\times B$ by \index{function!graph}
\begin{equation}
\label{eq:graphfunction}
F=\{(x,f(x))\Mid x\in A\}
\end{equation}
We call $F$ the {\em graph} of $f$, denoted by ${\rm Graph}(f)$. 
 A subset $F\subseteq A\times B$ is the graph  of a function from $A$ to $B$ if and only if it satisfies 
 the following two conditions:
\begin{equation}
{\rm G1}:\;\;(x,y)\in F\;\;\;  {\rm and}\;\;\; (x,y')\in F \implies y=y' 
\end{equation}
\begin{equation}
{\rm G2}:\;\;\{x\Mid (x,y)\in F\} = A.
\end{equation}
Condition G$1$ makes the idea of ``unique'' more precise, and G$2$ specifies what is meant by
``assigns to each.''  Two functions, $f$ and $g$, are equal if and only if their graphs are equal as sets:
${\rm Graph}(f) = {\rm Graph}(g).$
\index{function!domain, range}
The set $A$ is called the {\em domain} of $f$ (written $A={\rm domain}(f)$), and $B$ is called the {\em range} of $f$ (written $B={\rm range}(f)$). 
The notation $f:A \rightarrow B$ is used to denote that $f$ is a function with domain $A$
and range $B$.
\index{function! image}
For $S\subseteq A$, define $f(S)$ (image of $S$ under $f$) by
$f(S)\equiv \{f(x)\Mid x \in S\}$. 
In particular, $f(A)$ is called the {\em image} of $f$ (written $f(A)={\rm image}(f)$). 
The set of {\em all} functions with domain $A$ and range $B$ can be written  $\{f\Mid f:A \rightarrow B\}$ or simply as $B^A$.  If $A$ and $B$ are finite then $\left| A^B\right|$ is  $|A|^{|B|}$.
\index{function!indicator or characteristic}
The {\em characteristic or indicator function} of a set $S\subseteq A$, 
$\mathcal{X}_S:A\rightarrow \{0,1\}$, is defined by
\begin{equation}
\label{eq:charfunc}
\mathcal{X}_S(x) = 1 \,\,{\rm if\;and\;only\;if\,\,} x\in S.
\end{equation} 
\index{function!restriction, composition}
The {\em restriction} $f_S$ of  $f:A\rightarrow B$ to a subset $S\subseteq A$ is defined by
\begin{equation}
\label{eq:restriction}
f_S: S \rightarrow B\;\;{\rm where}\;\;f_S(x)=f(x)\;\;{\rm for\;\;all\;\;} x\in S.
\end{equation}

If $f:A\rightarrow B$ and $g:B \rightarrow C$ then the {\em composition} of $g$ and $f$, denoted
by $gf: A\rightarrow C$, is defined by 
\begin{equation}
\label{eq:composition}
gf(x) = g(f(x))\;\; {\rm for}\;\; x \in A.
\end{equation}
Note that composition of functions is associative: $h(gf)=(hg)f$ if 
$f:A\rightarrow B$, $g:B \rightarrow C$ and $h:C\rightarrow D$.
In some discussions, the product of functions, also denoted $gf$,  is defined by $gf(x)=g(x)f(x)$.
Another notation for composition of functions is $g\circ f$.
Thus, $g\circ f (x)\equiv g(f(x))$.

In the following six examples (\ref{eq:fnctypes}), the sets are specified by listing vertically their elements.  Arrows collectively define the rule. In the first example, $x=3$ is in $A = \{1,2,3,4\}$, and $f(x) = a$ is defined by the arrow 
from $3$ to $a\,$.

\begin{minipage}{\textwidth}
\begin{equation}
\label{eq:fnctypes}
{\bf Examples\;\,of\;\,Functions}
\end{equation}
\begin{center}
\includegraphics{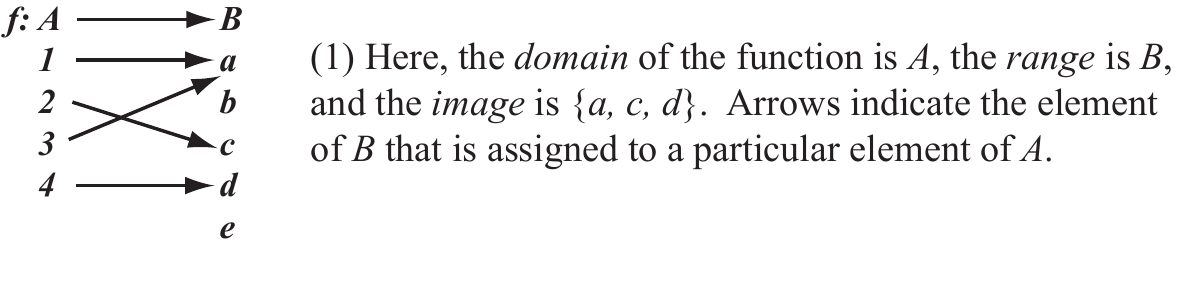}
\end{center}
\end{minipage}
%(2)
\begin{center}
\includegraphics{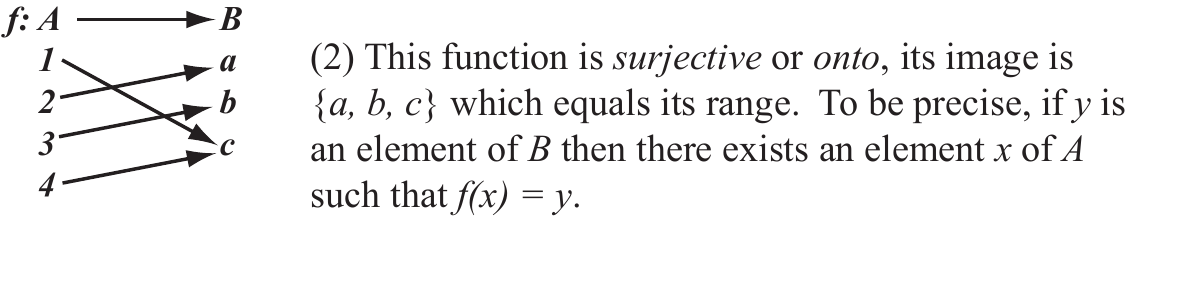}
\end{center}
%(3)
\begin{center}
\includegraphics{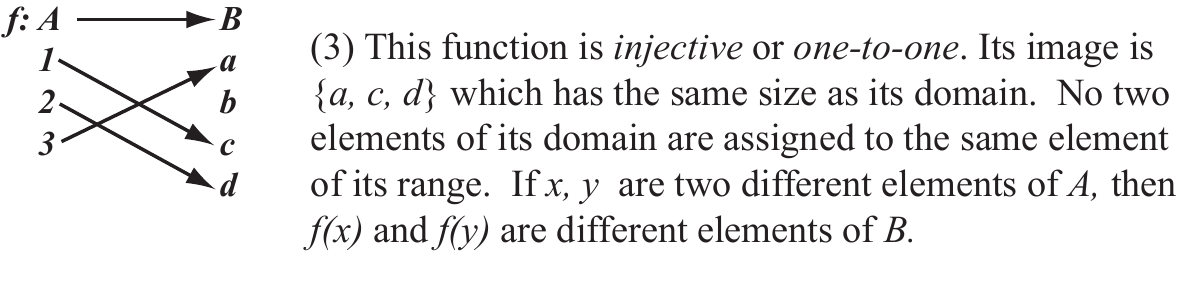}
\end{center}
%(4)
\begin{center}
\includegraphics{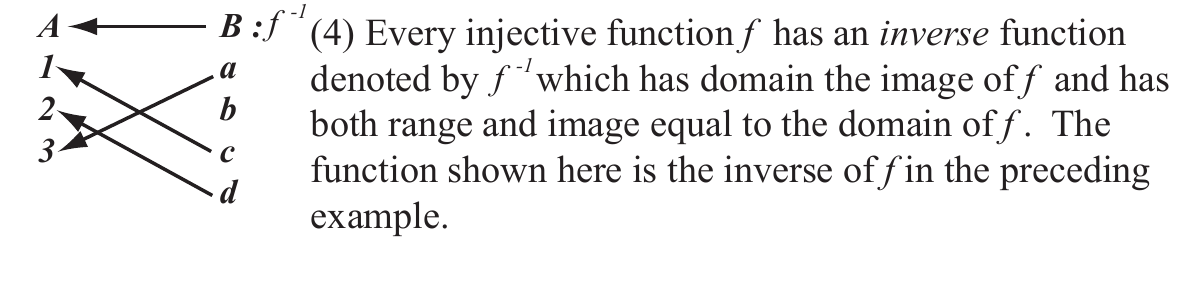}
\end{center}
%(5)
\begin{center}
\includegraphics{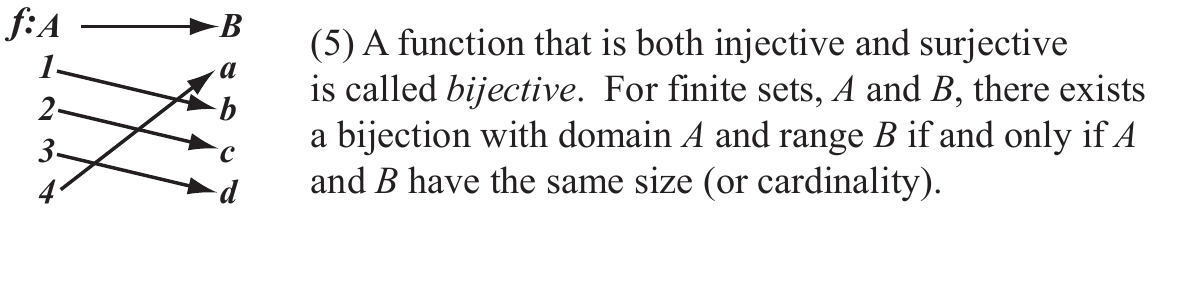}
\end{center}
%(6)
\begin{center}
\includegraphics{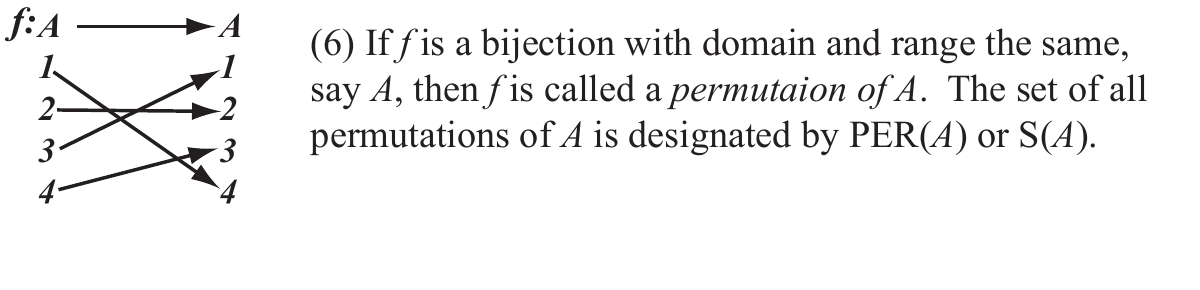}
\end{center}
\label{fig:p1ex1f}
%\end{figure}

There are many ways to describe a function.  Any such description must specify the domain, the range, and the rule for assigning some range element to each domain element.\
From the discussion following Definition~\ref{def:function}, you could specify the first function in the above examples (\ref{eq:fnctypes})  using set notation:  the domain is the set $\{1,2,3,4\}$, the range to be the set $\{a, b, c, d, e\}$, and the function $f$ is the set: 
${\rm Graph}(f)=\{(1,a)\,(2,c)\,(3,a)\,(4,d)\}$. 
\index{function!two line description} 
Alternatively, you could describe the same function by giving the range as $\{a, b, c, d, e\}$ and using {\em two line notation}
\begin{equation}
\label{eq:twolinedis}
f=
\left(
\begin {array}{cccc}
1 & 2 & 3 & 4 \\
a & c & a & d 
\end{array}
 \right).
\end{equation}
If we assume the domain, in order, is $1\,2\,3\,4$, then~\ref{eq:twolinedis} can be abbreviated  to 
{\em one line}: $a\,c\,a\,d$.\\ 
%\end{equation}
%Note {cccc} specifies centered columns.  We could use {llll}, {rrrr} or mixed {lrrl}.
%To omit the right parenthesis, write \right 
\index{function!coimage}
\begin{defn}[\bfseries Coimage partition]
\label{def:coimage}
Let $f:A\rightarrow B$ be a function with domain $A$ and range $B$.  
Let ${\rm image}(f) = \{f(x)\Mid x\in A\}$ (\ref{def:function}).  The inverse image of an
element $y\in B$ is the set $f^{-1}(y)\equiv \{x\Mid f(x)=y\}.$
The {\em coimage} of  $f$ is the set of subsets of $A:$ 
\begin{equation}
{\rm coimage}(f) = \{f^{-1}(y)\Mid y\in {\rm image}(f)\}.
\end{equation}
The ${\rm coimage}(f)$ is a partition of $A$ (\ref{def:setpartition}) called the
{\em coimage partition of $A$ induced by }$f$.
\end{defn}

For the function $f$ of \ref{eq:twolinedis}, we have ${\rm image}(f) = \{a,c,d\}.$
Thus, the coimage of $f$ is 
\begin{equation}
\label{eq:examplecoimage}
{\rm coimage}(f) = \{f^{-1}(a), f^{-1}(c), f^{-1}(d)\}=\{\{1,3\},\{2\},\{4\}\}.
\end{equation}
\hspace*{1 in}
%SETS OF FUNCTIONS
\index{functions!strictly increasing SNC}
\index{functions!weakly increasing WNC}
\index{functions!injective SNC}
\index{functions!permutations PER}
\begin{defn}[\bfseries Sets of functions] 
\label{def:setsfunc}
Let ${\underline n}=\{1, 2, \ldots, n\}$ and let
${\underline p}^{\underline n}$ be all functions with domain  
${\underline n}$, range ${\underline p}$.  Define
\[ 
{\rm SNC}(n,p) = \{f\Mid f\in {\underline p}^{\underline n}, i<j \implies f(i)<f(j)\}
\;\;{\bf (strictly\;increasing)}
\]
\[ 
{\rm WNC}(n,p) = \{f\Mid f\in {\underline p}^{\underline n}, i<j \implies f(i)\leq f(j)
\;\;{\bf (weakly\;increasing)}
\]
\[ 
{\rm INJ}(n,p) = \{f\Mid f\in {\underline p}^{\underline n}, i\neq j \implies f(i)\neq f(j)\}
\;\;{\bf (injective)}
\]
\[
{\rm PER}(n) = {\rm INJ}(n,n) 
\;\;{\bf (permutations\;of\;}{\underline n}).
\] 
From combinatorics, $|{\rm INJ}(n,p)|=(p)_n=p(p-1)\cdots (p-n+1)$, |PER(n)| = n!,
\[
|{\rm SNC}(n,p)|= \left(\begin{array}{c} p\\n\end{array}\right)
\;\;{\rm and}\;\;
|{\rm WNC}(n,p)|= \left(\begin{array}{c} p+n-1\\n\end{array}\right).    
\]

More generally, if $X\subseteq {\underline n}$ and $Y\subseteq {\underline p}$, then
${\rm SNC}(X,Y)$ denotes the strictly increasing functions from $X$ to $Y$.  We define
${\rm WNC}(X,Y)$ and ${\rm INJ}(X,Y)$ similarly.
%\label{def:fncsets}
\end{defn}
Sometimes  ``increasing'' is used instead of  ``strictly increasing'' or
``nondecreasing''   instead of  ``weakly increasing'' for the functions of \ref{def:setsfunc}.
 
 The next identity expresses the set ${\rm INJ}(n,p)$ as a composition of functions 
 in  ${\rm SNC}(n,p)$ and ${\rm PER}(n)$. and is illustrated in the table (\ref{eq:injsncpermp}) 
 that follows for $n=3$, $p=4$ ($f$, $g$ and the table entries, $fg$, are in one line notation):
\begin{equation}
\label{eq:injsncper}
{\rm INJ}(n,p) = \{fg\Mid f\in {\rm SNC}(n,p),\; g\in {\rm PER}(n) \}. 
\end{equation}

%Table for INJ(3,4)
\noindent%must be here to prevent indenting the whole minipage
\begin{minipage}{\textwidth}
\begin{equation}
\label{eq:injsncpermp}%this label can be found even if in minipage
{\bf Table\; for\;\;}{\rm INJ}(3,4) = \{fg\Mid  f\in {\rm SNC}(3,4),\;g\in {\rm PER}(3)\}
\end{equation}
\begin{center}	
\includegraphics{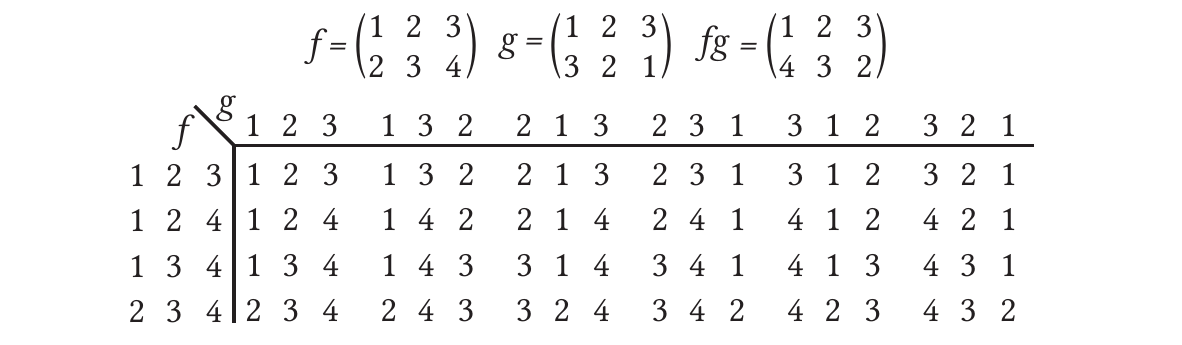}
\end{center}
\end{minipage}
\\%
\\%
\begin{remark}[\bfseries Example of an ${\rm INJ}(n,p)$ composition] 
Consider $h\in INJ(3,5)$ where
$h=\left(\begin{array}{ccc}1&2&3\\4&1&3\end{array}\right).$
Following~\ref{eq:injsncper} and~\ref{eq:injsncpermp}, we want to write $h$ as the composition $fg$
where $f\in {\rm SNC}(3,5)$ and $g\in {\rm PER}(3).$ 
Obviously, ${\rm image}(h) = {\rm image}(f) = \{1,3,4\}$ (see discussion following~\ref{def:function}),
and thus $f$ is uniquely defined:  
$f=\left(\begin{array}{ccc}1&2&3\\1&3&4\end{array}\right).$
To write $h=fg$ we only need the permutation $g=f^{-1}h$ where 
$f^{-1}=\left(\begin{array}{ccc}1&3&4\\1&2&3\end{array}\right).$
Thus, $g=\left(\begin{array}{ccc}1&2&3\\3&1&2\end{array}\right).$
\end{remark}

\section*{Permutations}
\index{permutations!cycle notation}
\begin{equation}
\label{eq:twolinecycle}
{\bf Two\;\,line,\;one\;\,line\;\,and\;\,cycle \;\,notation}
\end{equation}
%Subsection {Two line, one line and cycle notation}
%Use minipage if necessary.
\nobreak
There are three standard notations we shall use for writing permutations.
The permutation $f$ specified in the sixth arrow diagram of \ref{eq:fnctypes} can be written in two line notation as
\[
f=
\left(
\begin {array}{cccc}
1 & 2 & 3 & 4 \\
4 & 2 & 1 & 3 
\end{array}
\right)\,.
\]
The same permutation could be written as   

\[
f=
\left(
\begin {array}{cccc}
4 & 3 & 2 & 1\\
3 & 1 & 2 & 4 
\end{array}
\right)\,
\]
since the same rule of assignment is specified.
In the first of the above examples, the second line is $(4,2,1,3)$.  Since we know this is the second line of a two line notation for that permutation, we know that $D=R=\{4,2,3,1\} = \{1,2,3,4\}$.  In some discussions, the permutations are referred to by only the second line in the two line notation.  This shorthand representation of permutations is called the {\em one line} notation.  When using one line notation, the order of the elements in the missing first line (read left to right) must be specified in advance.  If that is done, then the missing first line can be used to construct the two line notation and thus the rule of assignment.  This seems trivially obvious, but sometimes the elements being permuted have varied and complex natural linear orders. In our example $f$, above, there are just two natural orders on the set $D$, increasing and decreasing.  If the agreed order is increasing then $f=(4,2,1,3)$ is the correct one line notation.  If the agreed order is decreasing then $f=(3,1,2,4)$ is the correct one line notation.

Our third notation for permutations, {\em cycle} notation, is more interesting. 
In cycle notation, the permutation $f$ of the previous paragraph is written
\begin{equation}
\label{eq:lengthcycle}
f = (1,4,3)(2)
\end{equation}
The ``cycle'' $(1,4,3)$ is read as ``1 is assigned to 4, 4 is assigned to 3, and 3 is assigned to 1.''  This cycle has three {\em elements}, 1, 4, and 3, and thus it has {\em length 3},  the number of elements in the cycle.
The cycle (2) has length one and is read as ``2 is assigned to 2''.  The usual convention with cycle notation is to leave out cycles of length one.  Thus, we would write
\[ f= (1,4,3)\]
and we would specify that the set D being permuted is $D=\{1,2,3,4\}$.
This latter information is required so we can reconstruct the missing cycles of length one.
%Product, inverse and identity
\index{permutations!product, inverse}
\begin{equation}
\label{eq:prodinv}
{\bf Product,\;\,inverse\;\;and\;\;identity}
\end{equation}
%Subsection {Product, inverse and identity}
\nobreak
If $f$ and $g$ are permutations, then the product, $h=fg$, of $f$ and $g$ is their {\em composition}:
$h(x) = (fg)(x) = f(g(x))$.  
For example, take $f=(1,4,3)$ and $g=(1,5,6,2)$ as permutations of 
${\underline 6} = \{1,2, \ldots, 6\}$.  Take $x=5$ and compute $g(5)=6$ and $f(6)=6$.
Thus, $(fg)(5)=6$. Continuing in this manner we get
\begin{equation}
\label{eq:compose}
h=fg =
\left(
\begin {array}{cccccc}
1 & 2 & 3 & 4 & 5 & 6\\
5 & 4 & 1 & 3 & 6 & 2 
\end{array}
\right)
\;\;{\rm or}\;\;
fg = (1,5,6,2,4,3)\,.
\end{equation}
In computing  $fg$, we mixed cycle notation and two line notation.
You should compose $f=(1,4,3)$ and $g=(1,5,6,2)$ to get $h=fg = (1,5,6,2,4,3)$, working entirely in cycle notation.  
Then put both $f$ and $g$ into two line notation and compose them to get the first identity in \ref{eq:compose}.

The identity permutation, $e$,  on $D$ is defined by $e(x)=x$ for all $x\in D$.
For $D={\underline 6}$ we have 
\begin{equation}
\label{eq:identityperm}
e=
\left(
\begin {array}{cccccc}
1 & 2 & 3 & 4 & 5 & 6\\
1 & 2 & 3 & 4 & 5 & 6
\end{array}
\right)
\;\;{\rm or}\;\;
e = (1)(2)\cdots (6)\,.
\end{equation}
For obvious reasons, we ignore the convention of leaving out cycles of length one when writing the identity permutation in cycle notation.

Every permutation, $h$, has an inverse, $h^{-1}$, for which $hh^{-1}=h^{-1}h=e$. 
For example, if we take $h=fg$ of \ref{eq:compose}, then
\begin{equation}
\label{eq:inverseh}
h^{-1} =
\left(
\begin {array}{cccccc}
5 & 4 & 1 & 3 & 6 & 2\\ 
1 & 2 & 3 & 4 & 5 & 6
\end{array}
\right)
\;\;{\rm or}\;\;
h^{-1} = (3,4,2,6,5,1)\,.
\end{equation}
Note that the two line representation of $h=fg$ in \ref{eq:compose} was just ``flipped over''
to get the two line representation of $h^{-1}$ in \ref{eq:inverseh}. 
It follows that $(h^{-1})^{-1} = h$ (which is true in general).
The cycle representation of  $h=fg$ in \ref{eq:compose} was written in reverse order to get
the cycle representation of $h^{-1}$ in \ref{eq:inverseh}. 
\begin{minipage}{\textwidth} 
\begin{equation}
\label{eq:propcyc}
{\bf Properties\;\,of\;\,cycles\;\,and\;\,transpositions}
\end{equation}
\index{permutations!cycles, transpositions}
Two cycles are {\em disjoint} if they  have  no entries in common (i.e., are disjoint as sets).  Thus, the permutation $(2517)(346)$ of $\underline {7}$ is written as the product of two disjoint cycles.  We leave out the commas in cycle notation, writing $(2517)$ rather than $(2,5,1,7)$, when the meaning is clear.
\end{minipage}
When a permutation is the product of disjoint cycles $f=c_1c_2\ldots c_t$, these cycles can be reordered in any manner (e.g., $f=(12)(34)(56) = (34)(12)(56)=(56)(12)(34)$, etc.) without changing $f$.  Also, the order of the entries in a cycle can be shifted around cyclically without changing the permutation:
\[
(2517)(346) = (5172)(346) = (1725)(346) = (7251)(463) = (2517)(634).
\]

A cycle of length two, like $(27)$, is called a {\em transposition}.
A cycle such as $c=(1423)$ can be written as a product of transpositions in a number of ways:
\begin{equation}
\label{eq:difftranspositions}
c=(13)(12)(14) =(13)(12)(14)(43)(42)(14)(23)(14)(24).  
\end{equation}
The arguments to the  permutations are on the right. Thus, the function $c$ evaluated at the integer $4$ is
$c{\bf[4]}=(13)(12)(14){\bf[4]}=(13)(12){\bf[1]}=(13){\bf[2]}=2.$
Note that the number of transpositions in the first ``transposition product'' representation of 
$c$ above (\ref{eq:difftranspositions}) is 3 and in the second representation is 9.  Although a given cycle $c$ can be written as a product of transpositions in different ways, say 
$t_1t_2\cdots t_p$ and  $s_1s_2\cdots s_q$, the number of transpositions, $p$ and $q$, are always both even or both odd (they have the same {\em parity} or, equivalently, 
$p\equiv q \pmod{2}$).  
In the example \ref{eq:difftranspositions}, we have $p=3$ and $q=9$, both odd.  We will discuss this further below.

Given a cycle $c=(x_1x_2\cdots x_k)$ of length $|c|=k$, $c$ can always be written as the product of $k-1$ transpositions. For example,
\begin{equation}
\label{eq:minprodtrans}
c=(x_1x_2\cdots x_k) = (x_1,x_k)(x_1,x_{k-1})\cdots (x_1,x_3)(x_1,x_2).
\end{equation}
Obviously, $c$ cannot be written as a product of less than $k-1$ transpositions.
\index{permutations!index defined}
\begin{defn}[\bfseries \bf Index of permutation] 
 \label{def:indexperm}
Let $f=c_1c_2\cdots c_p$ be a permutation written as a product of disjoint cycles 
$c_t$, where $|c_t| = k_t, t=1,2, \ldots, p$.
We define $I(f)$, the {\em index} of $f$, by
\begin{equation}
({\bf Index})\;\;I(f) = \sum_{t=1}^p (k_t - 1).
\end{equation}
\end{defn}
It is easy to see that the index, $I(f)$, is the smallest number of transpositions in any transposition product representation of $f$.  As another example, take $f=c_1c_2=(1234)(567)$ to be a permutation
on {\underline 9}.  The index, $I(f)=5$.  Let's check what happens to the index when we multiply $f$ by a transposition $\tau = (a,b)$, depending on the choice of $a$ and $b$.  

If $\tau = (a,b) = (89)$, then $\tau f =(89)(1234)(567)$ and $I(\tau f) = I(f)+1$.  
In this case, neither $a$ nor $b$ is in either cycle $c_1$ or $c_2$. 
%a and b not in any cycle
If $f=c_1c_2\cdots c_k$ where $c_i$ are the disjoint cycles ($|c_i|>1$) and if 
 $a$ and $b$ are not in any of the $c_i$, then 
$\tau f = (a\,b)f$ satisfies $I(\tau f) = I(f)+1$:
\begin{equation}
\label{eq:anotinbnotin}
{\bf a\,not\,in,\,b\,not\,in:\;}\tau f = (a\,b)c_1c_2\cdots c_k\,.
\end{equation}
$$I(\tau f) = I(f)+1$$

Let $\tau=(a\,b) = (59)$ with $a=5$, and let $f = (1234)(567)$ be a permutation on {\underline 9}.  In this case, $a=5$ is in a cycle (the cycle $(567)$),  but $b$ is not in any cycle.  Since the cycles commute, let's put the cycle containing $a$ first, $f = (567)(1234)$ (this is just a notational convenience).
We compute $\tau f = (59)(567)(1234) = (5679)(1234)$. 
Applying Definition~ \ref{def:indexperm}, we get $I(\tau f) = I(f)+1$.

%a is in a cycle but b is not
The general situation is to take 
$f=c_1c_2\cdots c_k$ ($|c_i|>1$), and take $\tau = (a,b)$ where  $a$ is in a cycle but $b$ is not.  We assume, without loss of generality, that $a$ is in 
$c_1=(a\,\,x_1\,\,x_2\,\ldots x_r)$.  We compute,
$(a\,b)(a\,\,x_1\,\ldots x_r)\;=(a\,\,x_1\,\ldots x_r\,b)$. Thus,
$I(\tau f)=r+1 + I(c_2\cdots c_k)$, $I(f) = r +  I(c_2\cdots c_k)$ and  $I(\tau f) = I(f)+1\,.$ 
To summarize:
\begin{equation}
\label{eq:ainbnot}
{\bf a\,in,\,b\,not\,in:\;} (a\,b)(a\,\,x_1\,\ldots x_r)\;=(a\,\,x_1\,\ldots x_r\,b)
\end{equation}
$$I(\tau f) = I(f)+1$$

%a and b in same cycle
Next we consider the case where 
$f=c_1c_2\cdots c_k$ ($|c_i|>1$), and  $\tau = (a,b)$ where $a$ is in a cycle and $b$ is in the same cycle.  
We assume, without loss of generality, that $a$ and $b$ are in 
$c_1=(a\,\,x_1\,\ldots x_r\,b\,y_1\ldots y_s )$.  We compute,
$(a\,b)(a\,\,x_1\,\ldots x_r\,b\,y_1\ldots y_s )\;=(a\,\,x_1\,\ldots x_r)\,(b\,y_1\ldots y_s )$. 
Thus 
$I(\tau f)=r+s + I(c_2\cdots c_k)$, $I(f) = r +s+1+  I(c_2\cdots c_k)$ and 
$I(\tau f) = I(f)-1\,.$
To summarize,
\begin{equation}
\label{eq:ainbinsame}
{\bf a,b\,same:\;}
(a\,b)(a\,\,x_1\,\ldots x_r\,b\,y_1\ldots y_s )\;=(a\,\,x_1\,\ldots x_r)\,(b\,y_1\ldots y_s )  
\end{equation}
$$I(\tau f) = I(f)-1$$

%a and b in different cycles
Finally, we consider the reverse of equation~\ref{eq:ainbinsame} where $a$ and $b$ are in different cycles. We compute,
$(a\,b)(a\,\,x_1\,\ldots x_r)\,(b\,y_1\ldots y_s )\;=(a\,\,x_1\,\ldots x_r\,b\,y_1\ldots y_s )$.
Thus, 
$I(\tau f)=r+s+1 + I(c_2\cdots c_k)$, $I(f) = r + s  + I(c_2\cdots c_k)$ and
$I(\tau f) = I(f)+1\,.$  To summarize,
\begin{equation}
\label{eq:ainbindiff}
{\bf a,b\,diff:\;}
(a\,b)(a\,\,x_1\,\ldots x_r)\,(b\,y_1\ldots y_s )\;=(a\,\,x_1\,\ldots x_r\,b\,y_1\ldots y_s ) 
\end{equation}
$$(\tau f) = I(f)+1\,$$

\begin{defn}[\bfseries Parity]
We say that $m$ and $n$ in ${\bZ}=\{0, \pm 1, \pm 2, \ldots \}$ have the same {\em parity} if $m-n$ is even.  
Equivalently,  we can write $m \equiv  n \pmod 2$ or $(-1)^m = (-1)^n$.
\end{defn}
\index{permutations!parity!index, transpositions}
Recall the index function, $I(f)$, of Definition~ \ref{def:indexperm}. 
\begin{lem}[\bfseries Parity of index and transposition count]
\label{lem:indtransparity}
Let $f$ be a permutation.  Suppose that 
$f=\tau_1\tau_2 \cdots \tau_q$ is any representation of $f$ as a product of transpositions $\tau_i$. 
Then the parity of $q$ is the same as the parity of the index, $I(f)$, of $f$.
\end{lem}
\begin{proof}
Note that 
\begin{equation}
\label{eq:invastrans}
e=(\tau_q\cdots \tau_2\tau_1)(\tau_1\tau_2\cdots \tau_q) = \tau_q\cdots \tau_2\tau_1f.
\end{equation}
Equations~\ref{eq:anotinbnotin} through \ref{eq:ainbindiff} state that multiplying an arbitrary permutation
by a transposition $\tau$ either increases or decreases that permutation's index by one. 
Thus, $I(\tau_1f)=I(f) \pm 1$, $\tau_2(\tau_1f) = I(f) \pm 1 \pm 1$, etc.
Applying this observation to~\ref{eq:invastrans} inductively gives
\[
0 = I(e) = I(f) + p -n
\]
where $p$ is the number of times a transposition in the sequence $(\tau_q,\ldots ,\tau_2,\tau_1)$ 
results in a ``+1'' and $n$ the number of times a transposition results in a ``-1''.
Thus, $I(f)=n-p$. But, $n-p\equiv n+p\pmod 2$.  Thus, $I(f)\equiv n+p \pmod 2$.  
But $n+p = q$, the number of transpositions.  This completes the proof that $I(f)\equiv q \pmod 2$ or, equivalently,
$I(f)$ and $q$ have the same parity.
\end{proof}

\begin{defn}[\bfseries \bf Sign of a permutation]  Let $f$ be a permutation.
\label{def:sign} 
\index{permutations!sign of} 
The {\em sign} of $f$ is defined as ${\rm sgn}(f)=(-1)^{I(f)}$.
\end{defn}

From Lemma~\ref{lem:indtransparity}, an alternative definition is  ${\rm sgn}(f)=(-1)^q$  where $q$ is the number of transpositions in any representation of $f$ as a product of transpositions: $f=\tau_1\tau_2\cdots \tau_q$.
As an example, consider the permutation $c=(1423)$ of~\ref{eq:difftranspositions}.
\begin{equation}
\label{eq:ninetranspositions}
c=(1423) =(13)(12)(14)(43)(42)(14)(23)(14)(24).
\end{equation}
The index, $I(c)$ is 3 so ${\rm sgn}(c) = (-1)^3 = -1$.  
The number of transpositions in the  product of transpostions in~\ref{eq:ninetranspositions} is 9. 
Thus, ${\rm sgn}(c) = (-1)^9 = -1.$
\index{permutations!inversions defined}
\begin{defn}[\bfseries \bf Inversions of a permutation]  Let $f$ be a permutation.
\label{def:inversions}  
Let $S$ be a set, $|S|=n$, for which there is an agreed upon ordering 
$s_1<s_2<\cdots < s_n$ of the elements.   Let $f$ be a permutation of $S$.
An {\em inversion} of $f$ with respect to this ordering is a pair of elements, 
$(x,y)\in S\times S$
where $x<y$ but $f(x)>f(y)$.  Let ${\rm Inv}(f)$ denote the set of all inversions of $f$ with respect to the specified ordering.  
An inversion of the form $(s_i, s_{i+1})$, $1\leq i < n$, is called an {\em adjacent} inversion ($s_i$ and $s_{i+1}$ are ``adjacent'' or ``next to each other'' in the ordering on $S$). 
\end{defn}
\index{permutations!parity!index, inversions}
For notational convenience, we take $S={\underline n}$ with the usual order on integers 
($1<2<\cdots <n$).  Figure~\ref{eq:inversiongrid} gives an example of a permutation $f$ which is given in two line notation and also in one line notation (at the base of the {\em inversion grid} used to plot the set ${\rm Inv}(f)$).

\begin{minipage}{\textwidth}
\begin{equation}
\label{eq:inversiongrid}
{\bf Figure:\;Inversion\;grid}
\end{equation}
\begin{center}
\includegraphics{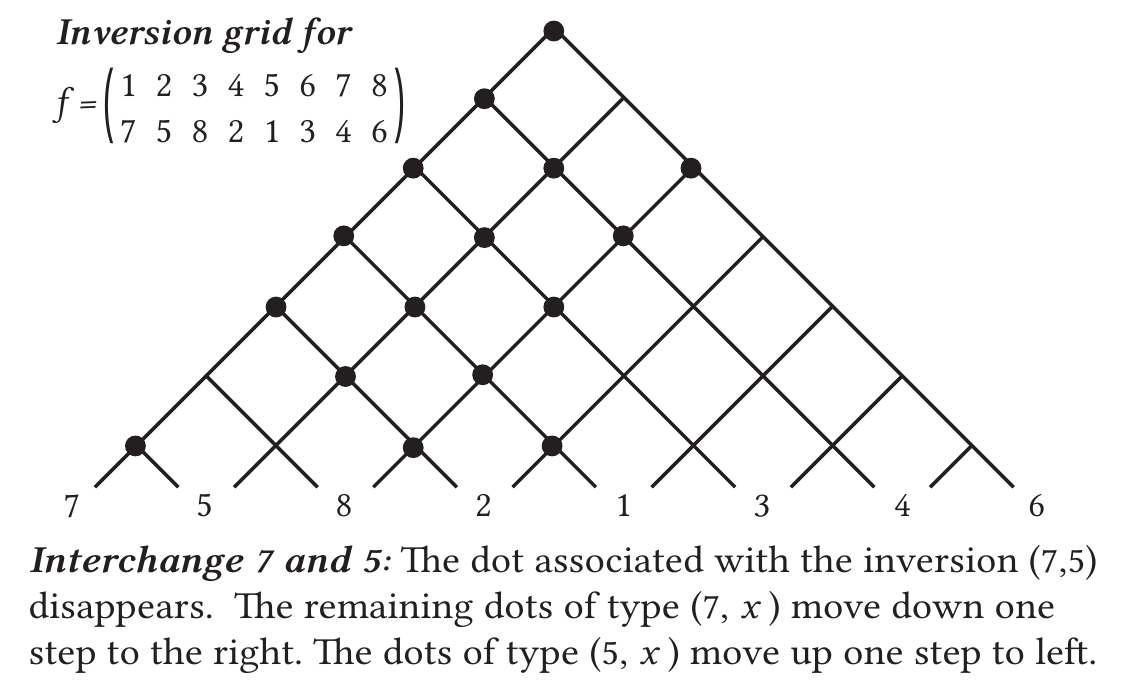}
\end{center}
\end{minipage}

Each intersection of two lines in the inversion grid (~\ref{eq:inversiongrid}) corresponds to a pair of integers in the one line notation for $f$. The intersection points corresponding to inversions are marked with solid black dots.  
There are 16 solid black dots, so $|{\rm Inv}(f)| =16$.  If we write $f$ in disjoint cycle form,
we get $f=(17425)(386)$ so the index $I(f)=4+2=6$.  
Thus, the parity of $I(f)$ (the index)  and $|{\rm Inv}(f)|$ (the inversion number) are the same in this case.  It turns out that they are always the same parity.

To understand why $I(f)\equiv |{\rm Inv}(f)| \pmod 2$, we take a close look at 
Figure~\ref{eq:inversiongrid}.  
Look at the adjacent inversion, $(1, 2)$, corresponding $f(1)=7>f(2)=5$.  Imagine what will happen when you switch 7 and 5 in the one line notation.  This switch corresponds to multiplying $f\tau_1=(17425)(386)(12)=(15)(274)(386)$ or $f\tau_1=5\,7\,8\,2\,1\,3\,4\,6$ in one line notation. Note, ${\rm Inv}(f\tau_1) = {\rm Inv}(f) -\{(1,2)\}$ as explained in 
Figure~\ref{eq:inversiongrid}.  The inversion $(1,2)$ is removed from ${\rm Inv}(f)$, all others remain.
The new function $f\tau_1$ has an adjacent inversion, $(3,4)$, corresponding to 8 and 2 in one line notation.  Remove it by multiplying by an adjacent transposition, $\tau_2 = (34)$ in this case.  In this way, you remove one adjacent inversion 
(from the original set ${\rm Inv}(f)$) by one transposition multiplication each time until you have $f\tau_1\tau_2\cdots \tau_q = e$ where $q=|{\rm Inv}(f)|$.   But we know from 
Lemma~\ref{lem:indtransparity} that
$q\equiv I(f) \pmod 2$.  Thus, $I(f)\equiv |{\rm Inv}(f)| \pmod 2$. To summarize,
\begin{equation}
\label{eq:invtransindex}
(-1)^{I(f)} = (-1)^{|{\rm Inv}(f)|} = (-1)^q = {\rm sgn}(f). 
\end{equation}

%EXERCISES
\section*{Exercises}

\begin{ex}
Which permutation on {\underline n} has the most inversions?  How many?
\end{ex}

\begin{ex}
Show that ${\rm sgn}(fg)={\rm sgn}(f){\rm sgn}(g)$ where $f$ and $g$ are permutations on a finite set $D$.  Hint: Use the identity ${\rm sgn}(f) = (-1)^q$ where $q$ is the number of transpositions in any representation of $f$
as a product of transpositions .  See equation ~\ref{eq:invtransindex}.
\end{ex}

\begin{ex}
In our discussion of Figure~\ref{eq:inversiongrid}, we proved that the permutation shown there was the product of $q$ adjacent transpositions,  $q = |{\rm Inv}(f)|$.  We assumed that each time we eliminated an inversion, if another inversion remained, we could choose an {\em adjacent} inversion to be eliminated.  Prove that this choice is always possible.
\end{ex}

\chapter{Matrices and Vector Spaces}
%SECTION
\section*{Review}
You should be familiar with Section~\ref{sec:setsandfunctions}.  We summarize a few key definitions:
If $K$ and $J$ are sets, then we use $f: K\rightarrow J$ to indicate that $f$ is a function with domain 
$K$ and range $J$. 
The notation, $\{f\Mid f: K\rightarrow J\}$ stands for the set of all functions with domain $K$ and range $J$ (also stated, ``set of all $f$ from $K$ to $J$'').
Alternatively, we use exponential notation, $J^K$, to denote the set of all functions from $K$ to $J$.
For a finite set $S$, we use $|S|$ to denote the cardinality (number of elements) of $S$.  If $S$ and $K$ are sets then $S\times K =\{(s,t)\Mid s\in S, t\in K\}$.
We note that $|J^K|$ = $|J|^{|K|}$ if $K$ and $J$ are finite.

\begin{rem}[\bfseries Underline notation]
\label{rem:underline}
We use the notation ${\underline n} = \{i\Mid 1\leq i \leq n\}$. Thus,
${\underline n} \times {\underline m} = \{(i,j)\Mid 1\leq i \leq n, 1\leq j \leq m\}$.
\end{rem}
\index{matrix!defined}
\begin{definition}[\bfseries Matrix]
\label{def:mat}
Let $m, n$ be positive integers.  An $m$ by $n$ matrix with entries in a set $S$ is a function 
$f: \underline{m}\times \underline{n} \rightarrow S$.
The set of all such $f$ is denoted by $\mathbf{M}_{m,n}(S)$.
\end{definition}

In some applications,  a matrix is defined as a function 
$f:\underline{\mu}\times \underline{\nu}\rightarrow S$ where 
$\underline{\mu}$ and $\underline{\nu}$ are linearly ordered sets. 
We won't need that generality in what follows.
The subject of ``matrix theory''  involves the use of matrices and operations defined on them as
{\em conceptual data structures} for understanding various ideas  in algebra and combinatorics.
For such purposes, matrices have a long history in mathematics and are surprisingly useful.

Next we recall and specialize Definitions~\ref{def:module} and \ref{def:vectoralgebra}:

\index{vector space and algebra!axioms}
\begin{defn}[\bfseries \bf Summary of vector space and algebra axioms] Let $\bF$ be a field and let
$(M, + )$ be an abelian group.  
Assume there is an operation,  $\bF\times M\rightarrow M$, which takes $(r,x)$ to $rx$ (juxtaposition of $r$ and $x$).  
To show that $(M, + )$ is a {\em vector space over} $\bF$, we show the following four things
hold for every $r,s \in \bF$ and $x,y \in M$:
$${\bf(1)}\;r(x+ y) = rx+  ry\;\;{\bf(2)}\\;(r+s)x = rx+  sx\;\;{\bf(3)}\\;(rs)x=r(sx)\;\;{\bf(4)}\\;1_\bF\,x=x$$
where $1_\bF$ is the multiplicative identity in $\bF$.  If $(M, + , \cdot )$ is a ring for which 
$(M, + )$ is a vector space over $\bF$, then $(M, + , \cdot )$ is an {\em algebra over} $\bF$ if the following scalar rule holds:
\begin{description}
\item[${\bf(5)}$ scalar rule] for all $\alpha \in \bF$, $x, y \in M$ we have  
$\alpha(xy)=(\alpha x)y = x(\alpha y)$.\\
\end{description}
%\ref{def:vspaceaxioms} 
\label{def:vspaceaxioms}
\end{defn}

The standard example of a vector space is as follows:

Let $J=\bR$, the field of real numbers (real number field).  
If $f$ and $g$ are in $\bR^K$, we define $h=f+g$ by
$h(x)=f(x)+g(x)$
for all $x \in K$.  
For $\alpha \in \bR$ and $f \in \bR^K$, we define the product, $h=\alpha f$, of the number 
$\alpha$ and 
the function $f$ by 
$h(x) = \alpha f(x)$ for all $x\in K$.  
Alternatively stated, $(\alpha f)(x) = \alpha f(x)$ for $x$ in $K$
($\alpha \beta$ denotes  $\alpha$ times $\beta$ in $\bR$).

The set $M=\bR^K$, together with these standard rules for adding and multiplying various things, satisfies the following four conditions (see \ref{def:module}) that make $M$ a module over $\bR$:
For all $r,s \in \bR$ and $f,g \in M$:
\[
(1)\;r(f+ g) = rf+  rg\;\;(2)\;(r+s)f = rf+  sf\;\;(3)\;(rs)f=r(sf)\;\;(4)\;1\,f=f.
\]
Since $\bR$ is a field, $M$ is a  {\em vector  space} by Definition~\ref{def:vectoralgebra}.

Next, define the {\em pointwise} product of $f, g\in \bR^K$ by  $(fg)(x)=f(x)g(x)$.  To show that the 
vector space $M=\bR^K$ with this rule of multiplicaton is an algebra, we need to verify  the scalar rule of
Definition~\ref{def:vectoralgebra}:
\begin{description}
\item[scalar rule] for all $\alpha \in \bR$, $f, g \in \bR^K$ we have  
$\alpha(fg)=(\alpha f)g = f(\alpha g)$.
\end{description}
This scalar rule follows trivially from the rules for multiplying functions and scalars.
Thus, $\bR^K$, under the standard rules for function products and multiplication by scalars, {\em is an  algebra over} $\bR$, the field of real numbers.  

A vector space $(M, + )$  over $\bF$ is a {\em real} vector space if
$\bF = \bR$ and a complex vector space if $\bF = \bC$. Examples are $\bR^K$ and $\bC^K$.  As noted above, these vector spaces are algebras under
pointwise multiplication of functions.

%Begin Example 1.1.2
\begin{examp} {\bf Examples} 
\label{example:colvec}
%\ref{example:colvec}
Let $V$ be the vector space of column vectors ($n\times 1$ matrices) with real entries:
%------------------------
\[
\left[
\begin{array}{c}
a_1\\
a_2\\
\vdots\\
a_n
\end{array}
\right]
+
\left[
\begin{array}{c}
b_1\\
b_2\\
\vdots\\
b_n
\end{array}
\right]
%-------------------------
=
%-------------------------
\left[
\begin{array}{c}
a_1+b_1\\
a_2+b_2\\
\vdots\\
a_n+b_n
\end{array}
\right]
\;\;\;\;{\rm and}\;\;\;\;\mu
\left[
\begin{array}{c}
a_1\\
a_2\\
\vdots\\
a_n
\end{array}
\right]
%---------------------------
=
%----------------------
\left[
\begin{array}{c}
\mu a_1\\
\mu a_2\\
\vdots\\
\mu a _n
\end{array}
\right]
\]
%------------------------
Or let $V$ be the vector space, $\mathbf{M}_{m,n}(\bR)\,$, of $m\times n$ matrices over $\bR$ where
%------------------------
\[
\left[
\begin{array}{cccc}
a_{1,1} & a_{1,2} & \hdots & a_{1,n} \\
a_{2,1} & a_{2,2} & \hdots & a_{2,n} \\
\vdots & \vdots & \vdots & \vdots\\
a_{m,1} & a_{m,2} & \hdots & a_{m,n} 
\end{array}
\right]
+
\left[
\begin{array}{cccc}
b_{1,1} & b_{1,2} & \hdots & b_{1,n} \\
b_{2,1} & b_{2,2} & \hdots & b_{2,n} \\
\vdots & \vdots & \vdots & \vdots\\
b_{m,1} & b_{m,2} & \hdots & b_{m,n} 
\end{array}
\right]
%--------------------------
=
%------------------------
\]										
\\
\[
\left[
\begin{array}{cccc}
a_{1,1}+b_{1,1} & a_{1,2}+b_{1,2}& \hdots & a_{1,n}+b_{1,n} \\
a_{2,1}+b_{2,1} & a_{2,2}+b_{2,2} & \hdots & a_{2,n}+b_{2,n} \\
\vdots & \vdots & \vdots & \vdots\\
a_{m,1}+b_{m,1} & a_{m,2}+b_{m,2} & \hdots & a_{m,n} +b_{m,m}
\end{array}
\right]
\]
and
%-------------------------
\[
\mu
\left[
\begin{array}{cccc}
a_{1,1} & a_{1,2} & \hdots & a_{1,n} \\
a_{2,1} & a_{2,2} & \hdots & a_{2,n} \\
\vdots & \vdots & \vdots & \vdots\\
a_{m,1} & a_{m,2} & \hdots & a_{m,n} 
\end{array}
\right]
=
\left[
\begin{array}{cccc}
\mu a_{1,1} &  \mu a_{1,2} & \hdots & \mu a_{1,n} \\
\mu a_{2,1} & \mu a_{2,2} & \hdots & \mu a_{2,n} \\
\vdots & \vdots & \vdots & \vdots\\
\mu a_{m,1} & \mu a_{m,2} & \hdots & \mu a_{m,n} 
\end{array}
\right].
\]
%----------------------------

\end{examp}
%End Example1.1.2

Note that the two examples in~\ref{example:colvec} are of the standard form 
$V=\bR^K$.  In the first case, $K= \underline{n} = \{1,2,\cdots,n\}$ and $f(t)=a_t$ for $t\in K$.  
In the second example, $K = \underline{n}\times\underline{m}$, a cartesian product of two sets, and
$f(i,j) = a_{i,j}$ (usually, $a_{i,j}$ is written $a_{ij}$ without the comma).\\    

%Begin review
\index{linear algebra!review}
\begin{review}[\bfseries Linear algebra concepts]
\label{rev:vecspace}
We review some concepts from a standard first course in linear algebra. Let $V$ be a vector space over $\bR$ and let $S$ be a nonempty subset of $V$.  
The {\em span} of  $S$, denoted ${\rm Span}(S)$ or ${\rm Sp(S)}$, is the set of all finite linear combinations of elements in $S$.  
That is, ${\rm Sp}(S)$ is the set of all vectors of the form $\sum_{x\in S} c_x x$ where $c_x \in \bR$ and $|\{x\Mid c_x\neq 0\}|<\infty$ ({\em finite support} condition).
Note that ${\rm Sp}(S)$ is a subspace of $V.$  
If ${\rm Sp}(S) = V$ then $S$ is a {\em spanning set} for $V$.  
$S$ is  {\em linearly independent} over 
$\bR$ if $\sum_{x\in S} c_x x \in {\rm Sp}(S)$ and $\sum_{x\in S} c_x x= \theta$ (the zero vector) then $c_x = 0$ for all $x$.  
If $S$ is linearly independent and spanning, then $S$ is a {\em basis} for $V$.   
If $S$ is a basis and finite ($|S|<\infty$), then $V$ is {\em finite dimensional}.  
The cardinality, $|S|$, of this basis is called the {\em dimension} of $V$ (any two bases for $V$ have the same cardinality).  
The zero vector is never a member of a basis.
The ${\rm Span}(x_1, x_2, \ldots , x_n)$ of a {\em sequence} of vectors is the set 
$\{\sum_{i=1}^n c_ix_i\Mid c_i\in \bR\}$.
A {\em sequence} of vectors, $(x_1, x_2, \ldots , x_n)$, is linearly independent if
$$c_1x_1 + c_2x_2 + \cdots + c_nx_n = 0$$ 
implies that $c_1=c_2=\cdots =c_n=0$.  Otherwise, the sequence is linearly dependent.
Thus, the sequence of  nonzero vectors $(x, x)$ is linearly dependent 
due to the repeated vector $x$.  The set $\{x,x\} = \{x\}$ is linearly independent.
However, ${\rm Span}(x,x) = {\rm Span}(\{x\}).$   
\end{review}
Most vector spaces we study will arise as subspaces of finite dimensional vector spaces already familiar to us. 
The vector space $V=\bR^n$ of n-tuples of real numbers will be a frequently used example
(\ref{def:vspaceaxioms}).
The elements of $V$ are usually written as $(\alpha_1, \alpha_2, \ldots , \alpha_n)$ were each $\alpha_i$ is a real number.  This sequence of real numbers can also be regarded as a function $f$ from the set of indices ${\underline n} = \{1,2,\ldots,n\}$ to $\bR$ where $f(i)=\alpha_i$.  
In this sense, $V=\bR^{\underline n}$ (using exponential notation for the set of functions).
$V$ can also be denoted by  $\times^n\bR$. 
%\REMARK
\begin{rem}[\bfseries Delta notation]
\label{rem:deltanotation}
\index{delta!examples} 
The $n$ vectors $e_i = (0, \ldots, 0, 1, 0, \ldots 0)$, where the single $1$ occurs in position $i$, are a basis (the ``standard basis'') for $V=\bR^n$.  In the function notation, 
$\bR^{\underline n}$, we can simply say that $e_i(j) = \delta (i=j)$ where $\delta({\bf Statement})$ is $1$ if {\bf Statement} is true and  $0$ if {\bf Statement} is false (\ref{rem:basicsets}).  In this case, 
{\bf Statement} is ``$i=j$''. One also sees $\delta_{ij}$ for the function $e_i(j)$.
\end{rem}

%Begin Exercises
\index{vector space!subspace examples}
\begin{defn}[\bfseries \bf Subspace]
Let $V$ be a vector space over $\bR$, the real numbers, and $H$ a nonempty 
subset of $V$. 
If for any $x, y \in H$ and $\alpha \in \bR$, $x+y \in H$ and $\alpha x\in H$, then H is a
{\em subspace} of $V$.  Similarly, we define a subspace of a vector space $V$ over the complex numbers or any field $\bF$.
\label{sub:subspace}
Note that a subspace $H$ of $V$ satisfies (or ``inherits'') all of the conditions of 
Definition~\ref{def:vspaceaxioms}
and thus $H$ is itself a vector space over $\bR$.
\end{defn}

\section*{Exercises: subspaces} 
\label{exe:subspaces}
%\ref{exe:subspaces}
%Subspace
.
\begin{ex}
Let $K_0$ be a subset of $K$.  Let $U$ be all $f$ in $\bR^K$, the real valued functions with domain $K$, such that $f(t)= 0$ for $t\in K_0$.  
Let $W$ be all $f$ in $\bR^K$ such that for $i,j\in K_0$, $f(i)=f(j)$.  Show that $U$ and $W$ are subspaces of the vector space $V=\bR^K$.
\label{exe:1.3.1}
%\ref{exe:1.3.1}
\end{ex}

\begin{ex}
In each of the following problems determine whether or not the indicated  subset is a subspace of the given vector space.   The vector space of $n$-tuples of real numbers is denoted by $\bR^n$ or $\times^n \bR$ (Cartesian product or $n$-tuples).
\label{exe:1.3.2}
%\ref{exe:1.3.2}
\begin{enumerate}
\item  $V=\bR^2$; $\;\;H = \{(x,y)\Mid y\geq 0\}$
\item $V=\bR^3$; $\;\;H = \{(x,y,z)\Mid z= 0\}$
\item $V=\bR^2$; $\;\;H = \{(x,y)\Mid x=y\}$
\item $V=C(K)$ where $C(K)$ is all continuous real valued functions on the interval 
$K=\{x\Mid 0<x<2\}$; $\;\;H= \{f\Mid f\in V, f(1)\in \bQ\}$ where $\bQ$ is the rational numbers.
\item $V=C(K)$; $\;\;H$ is all $f\in C(K)$ such that $3\frac{df}{dx}=2f$.
\item $V=C(K)$; $\;\;H$ is all $f\in C(K)$ such that there exists real numbers $\alpha$ and $\beta$ (depending on $f$) such that $\alpha\frac{df}{dx}=\beta f$.
\item $V=\bR^{N_0}$ where $N_0=\{0,1,2,3,\ldots\}$ is the set of nonnegative integers; 
$H=\{f\Mid f\in V,\; |\{t\Mid f(t)\neq 0\}|<\infty\}$ is the set of all functions in $V$ with ``finite support.''
\end{enumerate}
\end{ex}

\begin{ex}
Let $\mathbf{M}_{n,n}(\bR)\,$ be the real vector space of  $n\times n$ matrices 
(Example~\ref{example:colvec}).  
Let $H_1$ be the subspace of {\em symmetric} matrices 
$B=(b_{ij}), 1\leq i,j\leq n$, where $b_{ij}=b_{ji}$ for all $i$ and $j$.
Let $H_2$ be the subspace of {\em lower triangular} matrices 
$B=(b_{ij}), 1\leq i,j\leq n$, where $b_{ij}=0$ for all $i<j$.  
Show that for any two subspaces, $H_1$ and $H_2$, of a vector space $V$, $H_1\cap H_2$ 
is a subspace.  What is $H_1 \cap H_2$ for this example?  What is the smallest subspace
containing $H_1\cup H_2$ in this example?
\label{exe:1.3.3}
%\ref{exe:1.3.3}
\end{ex}
\index{vector space!spanning sets, dimension}
\section*{Exercises: spanning sets and dimension}
\label{sec:spandim}
\begin{ex}
\label{old:25.2}
Show that the matrices
\[
 \left( \begin{array}{cc}
  	1  & 1\\
	0  & 1  
   \end{array}\right)
   \;,\;
\left( \begin{array}{cc}
       -1  & 1\\
	0  & 1  
   \end{array}\right)
      \;,\;
\left( \begin{array}{cc}
  	0  & 1\\
	0  & 0  
   \end{array}\right)
\]
do not span the vector space of all $2\times 2$ matrices over $\bR$.\\
\end{ex}

\begin{ex}
\label{old:25.3}
Let $V$ be a vector space and $x_1, x_2, \ldots, x_n$ be a sequence (ordered list) of vectors in $V$.   
If 
\[
x_1\neq 0\,, x_2\notin {\rm Sp}(\{x_1\})\,, x_3\notin {\rm Sp}(\{x_1, x_2\})\,,
\ldots\,, x_n \notin {\rm Sp}(\{x_1, x_2,\ldots x_{n-1}\})\,,
\]
show that the vectors $x_1, x_2, \ldots, x_n$ are linearly independent.\\
\end{ex}

\begin{ex}
\label{old:25.4}
Let $V$ be a vector space and $x_1, x_2, \ldots, x_n$ be  linearly independent 
vectors in $V$.   
Let
\[
y= \alpha_1x_1 + \alpha_2x_2 + \cdots + \alpha_nx_n\,.
\]
What condition on the scalars $\alpha_i$ will guarantee that for each $i$,  the vectors $x_1, x_2, \ldots, x_{i-1}, y, x_{i+1}, \ldots x_n$ are linearly independent?\\
\end{ex}

\begin{ex}
\label{old:25.5}
Show that the vectors $s_t=\sum_1^te_i, t=1, \ldots, n$, are a basis for $\bR^{\underline n}$ where $e_i(j)=\delta (i=j)$, $i=1, \ldots, n$ (see the Remark~\ref{rem:deltanotation} ).\\
\end{ex}

\begin{ex}
\label{old:25.6}
Recall the basics of matrix multiplication. The matrices
\[
\sigma_x=
 \left( \begin{array}{cc}
  	0  & 1\\
	1  & 0  
   \end{array}\right)
   \;,\;
\sigma_y=
\left( \begin{array}{cc}
        0  & -i\\
	 i  & 0  
   \end{array}\right)
      \;,\;
\sigma_z=
\left( \begin{array}{cc}
  	1  & 0\\
	0  & -1  
   \end{array}\right)
\]
are called the {\em Pauli spin matrices}.
Show that these three matrices plus the identity matrix $I_2$~(\ref{def:matmult}) form a basis for the vector space of $2\times 2$ matrices over the complex numbers, $\bC$.
Show also that 
\[
\sigma_x \sigma_y = - \sigma_y \sigma_x\;\;\;\;{\rm and}\;\;\;\;
\sigma_x \sigma_z = - \sigma_z \sigma_x\;\;\;\;{\rm and}\;\;\;\;
\sigma_y \sigma_z = - \sigma_z \sigma_y\;.\\
\]
\end{ex}

\begin{ex}
\label{old 25.9}
Let $A_1, A_2, \ldots, A_k$ be a sequence of  $m\times n$ matrices.  Let $X \neq \theta_{n\times 1}$ be an $n\times 1$ matrix ($\theta_{n\times 1}$ is  the $n\times 1$ zero matrix).
Show that if $A_1X = A_2X = \cdots A_kX =  \theta_{m\times 1}$ then the matrices 
$A_1, A_2, \ldots, A_k$ do not form a basis for the vector space ${\bf M}_{m,n}(\bR)$ of all $m\times n$ matrices.\\
\end{ex}

\begin{ex}
\label{old 25.13}
Find a basis for the vector space of $n\times n$ matrices, ${\bf M}_{n,n}(\bR)$, that consists only of matrices 
$A$ that satisfy $A^2=A$ (these are called {\em idempotent} matrices).  
Hint:  For the case $n=2$, here is such a basis:
\[
\left( \begin{array}{cc}
  	1  & 0\\
	0  & 0  
   \end{array}\right)
 \;\;
 \left( \begin{array}{cc}
  	0  & 0\\
	0  & 1  
   \end{array}\right)
 \;\;
 \left( \begin{array}{cc}
  	0  & 0\\
	1  & 1  
   \end{array}\right)
 \;\; 
  \left( \begin{array}{cc}
  	0 & 1\\
	0 & 1  
   \end{array}\right)\;.
\]
\end{ex}
\hspace*{1in}\\
\begin{ex}
\label{old 25.14}
Show that if $A_iA_j = A_jA_i$ for all $1\leq i< j \leq k$, then the $n\times n$ matrices 
$A_1, A_2, \ldots, A_k$ are not a basis for ${\bf M}_{n,n}(\bR)$ ($n>1$). \\  
\end{ex}

\begin{ex}
\label{old 25.16}
({\bf Trace of matrix}) Show that if $A_i = B_iC_i - C_iB_i$ for  $1\leq i \leq k$, where $A_i$, $B_i$ and $C_i$ are 
$n\times n$ matrices,  then the 
$A_1, A_2, \ldots, A_k$ do not form a basis for the vector space ${\bf M}_{n,n}(\bR)$ ($n>1$).
\index{matrix!trace, {\rm Tr}}
Hint: 
Recall that the trace of a matrix $A$, ${\rm Tr}(A)$,  is the sum of the diagonal entries of $A$ (i.e., $A(1,1)+A(2,2)+\cdots A(n,n)$).  It is easy to show that 
${\rm Tr}(A+B)={\rm Tr}(A) + {\rm Tr}(B)$ and ${\rm Tr}(AB)={\rm Tr}(BA)$. 
This latter fact, ${\rm Tr}(AB)={\rm Tr}(BA)$, implies that if $A$ and
$B$ are similar matrices, $A=SBS^{-1}$, then ${\rm Tr}(A)={\rm Tr}(B)$
(a fact not needed for this exercise).
The proof  is trivial: ${\rm Tr}(S(BS^{-1}))={\rm Tr}((BS^{-1})S)={\rm Tr}(B)$.\\
\end{ex}

\begin{ex}
Is it possible to span the vector space ${\bf M}_{n,n}(\bR)$ ($n>1$). using the powers of a single matrix:
$I_n, A, A^2, \ldots, A^t, \ldots $?\\
\end{ex}

\begin{ex}
%\label{old 25.17}
Show that any $n\times n$ matrix with real coefficients satisfies a polynomial equation 
$f(A)=0$, where $f(x)$ is a nonzero polynomial with real coefficients.
Hint: Can the matrices $I_n, A, A^2, \ldots, A^p$ be linearly independent for all $p$?\\
\end{ex}

\vfil
\section*{Matrices -- basic stuff}
We first discuss some notational issues regarding matrices.\\
%\begin{minipage}{\textwidth}
%\begin{equation}
%\label{eq:notationalissues}
%{\bf Notational\,\,issues}\\
%\end{equation} 
%\end{minipage} 
%\hspace*{1 in}\\

\index{matrix!index-to-entry function}
\begin{remark}[\bfseries Index-to-entry function]
\label{rem:indtoent}
The matrix $A$ in~\ref{eq:indextoentry} is shown in the standard general form for an $m\times n$ matrix.
Assume the entries of $A$ are from some set $S$.
From \ref{def:mat}, $A$ is the  function $(i,j)\mapsto a_{ij}$ whose domain is 
$\underline{m}\times \underline{n}$ and whose range is $S$.
As we shall see, matrices can be interpreted as functions in other ways.  
We refer to the basic definition (\ref{def:mat}) as the index-to-entry representation of $A$. 
The standard rectangular presentation is as follows:
\begin{equation}
\label{eq:indextoentry}
A=
\left[
\begin{array}{cccc}
a_{11} & a_{12} & \hdots & a_{1n} \\
a_{21} & a_{22} & \hdots & a_{2n} \\
\vdots & \vdots & \vdots & \vdots\\
a_{m1} & a_{m2} & \hdots & a_{mn} 
\end{array}
\right]
\end{equation}
%Call this function, $A: \underline{m}\times \underline{n} \rightarrow S$, 
%an {\em index-to-entry} function.\\
\end{remark}

In~\ref{eq:twointerpretations} we see two representations of the same  
index-to-entry function, $A$.
The first representation is the standard two-line description of a function.
The domain, $\underline{2}\times \underline{2}$, is listed in lexicographic order as the first line;
the values of the function are the second line.
In the second representation, the domain values 
\index{matrix!two line interpretation}
are not shown explicitly but are inferred by the standard rule for indexing the elements of a matrix, 
$(i,j)\rightarrow A(i,j)$, 
where $i$ is the row index and $j$ the column index. 
Thus, $A(1,1)=2$, $A(1,2)=2$, $A(2,1)=3$, $A(2,2)=4$.  
\begin{equation}
\label{eq:twointerpretations}
A=
\left(
\begin{array}{cccc}
(1,1)&(1,2)&(2,1)&(2,2)\\
2&2&3&4
\end{array}
\right)
\hspace{0.5cm}
{\rm or}
\hspace{0.5cm}
A=
\left(
\begin{array}{cc}   
2&2\\
3&4
\end{array}
\right)
\end{equation}
The second representation of the {\em index-to-entry} function $A$ in~\ref{eq:twointerpretations} is the one most used in matrix theory.

\begin{rem}[\bfseries Matrices and function composition]
\label{rem:matricescomposition}
Matrices as index-to-entry functions can be composed (\ref{eq:composition}) with other functions. 
\index{matrix!function composition}
Here is $A$, first in two line and then in standard matrix form: 
\end{rem}
%MATRIX TWO LINE
\begin{equation}
\label{eq:matrixa}
A=
\left(
\begin{array}{cccc}
(1,1)&(1,2)&(2,1)&(2,2)\\
2&2&3&4
\end{array}
\right)=
\left(
\begin{array}{cc}   
2&2\\
3&4
\end{array}
\right),
\end{equation}
%SIGMA TWO LINE
Next is a permutation $\sigma$ of $\underline{2}\times \underline{2}$:
\begin{equation}
\label{eq:permutationsigma}
\sigma=
\left(
\begin{array}{cccc}
(1,1)&(1,2)&(2,1)&(2,2)\\
(1,2)&(2,1)&(2,2)&(1,1)
\end{array}
\right).
\end{equation}
%A SIGMA
Next compose $A$ with $\sigma$ (work with two line forms):
\begin{equation}
\label{eq:acircsigma}
A\circ \sigma \equiv A\sigma=
\left(
\begin{array}{cccc}
(1,1)&(1,2)&(2,1)&(2,2)\\
2&3&4&2
\end{array}
\right)
=
\left(
\begin{array}{cc}   
2&3\\
4&2
\end{array}
\right).
\end{equation}
%TAU
Compose $\sigma$ with a second permutation $\tau$
\begin{equation}
\label{eq:permutationtau}
\tau=
\left(
\begin{array}{cccc}
(1,1)&(1,2)&(2,1)&(2,2)\\
(1,1)&(2,1)&(1,2)&(2,2)
\end{array}
\right)
\end{equation}
to obtain $\sigma\tau$ in two line notation:
\begin{equation}
\label{eq:permutationsigmatau}
\sigma\tau=
\left(
\begin{array}{cccc}
(1,1)&(1,2)&(2,1)&(2,2)\\
(1,2)&(2,2)&(2,1)&(1,1)
\end{array}
\right).
\end{equation}
Finally, compose $A$ with $\sigma\tau$:
\begin{equation}
\label{eq:acircsigmatau}
A(\sigma\tau)=(A\sigma)\tau=
\left(
\begin{array}{cccc}
(1,1)&(1,2)&(2,1)&(2,2)\\
2&4&3&2
\end{array}
\right) =
\left(
\begin{array}{cc}   
2&4\\
3&2
\end{array}
\right).
\end{equation}
Converting $\sigma$ (\ref{eq:permutationsigma}) to cycle notation, we get 
$\sigma = \left((1,1), (1,2), (2,1), (2,2)\right)$ which is a cycle of length four.
In cycle notation, $\tau$ (\ref{eq:permutationtau}) is $\tau = \left((1,2), (2,1)\right)$
which is a transposition.  In fact, $A\tau$ is called the {\em transpose} (\ref{def:transpose}) of $A$.
The permutation $\sigma\tau= ((1,1),(1,2),((2,2))$ is a three cycle.\\
Note that the matrix $A$ (\ref{eq:twointerpretations}) can be composed with functions that are not permutations:
\begin{equation}
\label{eq:notaperm}
f=\left(
\begin{array}{cc}
(1,1)&(1,2)\\
(1,2)&(2,1)
\end{array}
\right)\;\;
Af=\left(
\begin{array}{cc}
(1,1)&(1,2)\\
2&3
\end{array}
\right)\;=\;
\left(
\begin{array}{cc}
2&3
\end{array}
\right).
\end{equation}
In ~\ref{eq:notaperm}, the function $f$ transforms a $2\times 2$ matrix 
$A$ into a $1\times 2$ matrix $Af$.

\begin{rem}[\bfseries Function terminology applied to matrices]
\label{rem:functiontermsmatrices}
Consider the matrix $A$ of \ref{rem:matricescomposition}:
\begin{equation}
\label{eq:considermatrixa}
A=
\left(
\begin{array}{cccc}
(1,1)&(1,2)&(2,1)&(2,2)\\
2&2&3&4
\end{array}
\right)=
\left(
\begin{array}{cc}   
2&2\\
3&4
\end{array}
\right)
\end{equation}
Recall the terminology for functions associated with Definition~\ref{def:function} through
\index{matrix!coimage, image}
Definition~\ref{def:coimage} (coimage).
The ${\rm image}(A)=\{2,3,4\}$.  The ${\rm domain}(A)= \{(1,1), (1,2), (2,1),(2,2)\}$.
The range of $A$ could be any set containing the image (e.g., $\underline{4}$).
The ${\rm coimage}(A) = \{\{(1,1), (1,2)\}, \{(2,1)\}, \{(2,2)\}\}$.
\end{rem}

\begin{defn}[\bfseries Basic matrix notational conventions]
\label{def:matrixbasics}
Let $A=(a_{ij})$ and $A'=(a'_{ij})$ be two $m\times n$ matrices with entries in a set $S$.
We use the notation $A(i,j)\equiv a_{ij}$.
Two matrices are equal, $A=A'$, if $A(i,j)=A'(i,j)$ for all $(i,j)\in \underline{m}\times \underline{n}$.
$A_{(i)} = (a_{i1} \ldots a_{in})$, $1\leq i\leq m$,  designates row $i$ of $A$.  $A_{(i)}$ is a $1\times n$ matrix called a 
{\em row vector} of $A$.   $A^{(j)}$, $1\leq j\leq n$,  designates column $j$ of $A$.  
$A^{(j)}$ is an $m\times 1$ matrix called a  {\em column vector} of $A$.
\end{defn}

\hspace*{1 in}\\
\begin{minipage}{\textwidth}
The range $S$ of the index-to-entry function of a matrix can be quite general.
Figure~\ref{eq:NCmatrices} shows two matrices, $N$ and $C$, which have matrices as entries.\\
\index{matrix!general range}
\begin{equation}
\label{eq:NCmatrices}
{\bf Figure: Matrices\;\;with\;\;matrices\;\;as\;\;entries:}
\end{equation}
\begin{center}
\includegraphics{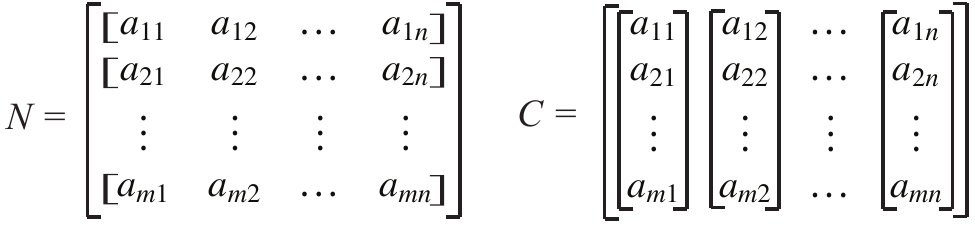}
\end{center}
\end{minipage}

The matrix $N$ (\ref{eq:NCmatrices}) is an $m\times 1$ matrix with each entry a 
$1\times n$ row vector of a matrix $A$~(\ref{def:matrixbasics}). 
Thus, $N(i,1)$, $1\leq i\leq m$, is the $1\times n$ {\em matrix} consisting of row $i$ of $A$ which, 
in terms of $A$, is designated $A_{(i)}$.

The matrix $C$ (\ref{eq:NCmatrices}) is a $1\times n$ matrix with each entry an $m\times 1$ matrix called a {\em column vector} of $A$.
Thus, $C(1,j)$, $1\leq j\leq n$, is the $m\times 1$ {\em matrix} consisting of column $j$ of $A$ which, in terms of $A$, is designated $A^{(j)}$.

\begin{rem}[\bfseries Equality of row and column vectors]
If $A$ is an $m\times n$ matrix, a row vector $A_{(i)}$ can never be equal to a column vector $A^{(j)}$ 
(unless $m=n=1$).  Two matrices can be equal only if they have the same number of rows and the same number of columns.
Sometimes you will see a statement that ``row $i$ equals column $j$.''   Such a statement might be made, for example,
if $m=n$ and the sequence of numbers in $A_{(i)}$ equals the sequence of numbers in $A^{(j)}$.  In this case the row vector equals the column vector as  sequences of numbers, not as matrices.  
\end{rem}
%Endline needed for minipage
%Subsection {Multiplication of matrices}
\begin{minipage}{\textwidth}
\begin{equation}
{\bf Multiplication\;\; of\;\; matrices}
\label{eq:interpret}
\end{equation}
The relationship between linear transformations and matrices is the primary 
(but not only) motivation for the following definition of matrix multiplication:
\end{minipage}
\index{matrix!multiplication definition}
\begin{defn}[\bfseries Matrix multiplication]
\label{def:matmult} 
Let $A\in \mathbf{M}_{m,p}(\bK)$ and $B\in \mathbf{M}_{p,n}(\bK)$ be matrices (\ref{rem:specialringsfields} for $\bK$ of interest here).
Let $A(i,j)=a_{ij}$, $(i,j)\in {\underline m}\times {\underline p}$, and $B(i,j) = b_{ij}$,
$(i,j)\in {\underline p}\times {\underline n}$.
The product $D=AB$ is defined by 
\begin{equation}
d_{ij}=D(i,j) = \sum_{k=1}^pA(i,k)B(k,j)=\sum_{k=1}^p a_{ik}b_{kj},\;\;(i,j)\in {\underline m}\times {\underline n}.
\end{equation}
The $q\times q$ matrix, $I_q$ defined by $I_q(i,j) = 1$ if $i=j$ and $I_q(i,j) = 0$ if $i\neq j$
is called the $q \times q$ {\em identity matrix}.  For an $m\times n$ matrix $M$, 
\begin{equation}
\label{eq:identity}
I_mM = MI_n = M
\end{equation}
\end{defn}
From Definitions~\ref{def:matrixbasics} and \ref{def:matmult}  we have
\begin{equation}
\label{eq:rowcolentry}
(M_{(i)})^{(j)} = (M^{(j)})_{(i)}\equiv M_{(i)}^{(j)} = (M(i,j)).
\end{equation}

Another way to write the sum in Definition~\ref{def:matmult} is using the 
{\em summation convention}
\begin{equation}
\label{eq:sumcon}
D(i,j)=A(i,k)B(k,j)
\end{equation}
where the two consecutive $k$ indices imply the summation.

Using Definition~\ref{def:matmult}, it easy to show that if $A$ is $n_0\times n_1$, $B$ and $C$ are $n_1\times n_2$ and $D$  is $n_2\times n_3$, then the {\em distributive laws} hold:
\begin{equation}
A(B+C) = AB + AC
\;\;{\rm and}\;\;
(B+C)D = BD + CD.
\end{equation}
\[
\framebox{\textbf{\emph{The associative law for matrix multiplication }}}
\]

\index{matrix!associativity of multiplication}
An important property of matrix multiplication is that it is {\em associative}.
If $A$ is an $n_0\times n_1$ matrix, $B$ an $n_1\times n_2$ matrix, and $C$ an 
$n_2\times n_3$ matrix, then $(AB)C = A(BC)$.  
An aficionado of the summation convention would give a short proof:
%\begin{equation}
$$((AB)C)(i,j) = (A(i,t_1)B(t_1,t_2))C(t_2,j)$$ 
$$= A(i,t_1)(B(t_1,t_2)C(t_2,j)) = (A(BC))(i,j).$$
%\end{equation}
This summation-convention proof uses the commutative, distributive, and associative laws.   Here is the longer proof of the associative law for matrix multiplication using explicit summation notation:

\begin{equation}
\label{eq:associative}
((AB)C)(i,j) = \sum_{t_2=1}^{n_2} \left( \sum_{t_1=1}^{n_1}  A(i,t_1)B(t_1,t_2)\right) C(t_2,j)=
\end{equation}
$$
\sum_{t_2=1}^{n_2} \sum_{t_1=1}^{n_1}  A(i,t_1)B(t_1,t_2)C(t_2,j) =
$$
$$
 \sum_{t_1=1}^{n_1} \sum_{t_2=1}^{n_2} A(i,t_1)B(t_1,t_2)C(t_2,j) =
$$
$$
 \sum_{t_1=1}^{n_1}A(i,t_1)\left( \sum_{t_2=1}^{n_2} B(t_1,t_2)C(t_2,j)\right) = (A(BC))(i,j).
$$

The associativity of matrix multiplication is a powerful combinatorial tool.
Suppose we are to compute the product $C=A_1A_2A_3A_4$ of four matrices assuming, of course, that the product is defined.  For example, suppose 
$A_1$ is $n\times n_1$, $A_2$ is $n_1\times n_2$, $A_3$ is $n_2\times n_3$, 
$A_4$ is $n_3\times m$.  We can express $C(i,j)$, an entry in the $n\times m$ matrix $C$, as
\[
C(i,j)=\sum_{t_1,t_2,t_3} A_1(i,t_1)A_2(t_1,t_2)A_3(t_2,t_3)A_4(t_3,j)
\]
where the sum is over all $(t_1,t_2,t_3) \in {\underline n_1}\times {\underline n_2}\times {\underline n_3}$ in any order. 

Let $A$, $L$ and $R$ be $n\times n$ matrices. Suppose that $LA=I_n$ and $AR=I_n$. 
Then  $R=(LA)R = L(AR)=L$, and, hence, $R=L$  (see \ref{rem:monoid}). 
The matrix $B=R=L$ is called the {\em inverse} of $A$ if it exists.  
\index{matrix!inverse}
We use the notation,
$B=A^{-1}$ for the inverse of $A$: 
\begin{equation}
\label{eq:inverse}
A^{-1}A = AA^{-1} = I_n.
\end{equation}
\index{matrix!invertible, unit, nonsingular}
If a matrix $A$ has an inverse, we say that $A$ is {\em nonsingular} or {\em invertible} or a {\em unit} in the ring  $R={\bf M}_{n,n}({\bK})$ (\ref{rem:specialrings}).

Again, assume $D=AB$ where $A$ is  an $m\times p$ matrix and $B$ a $p\times n$ matrix. 
Note that Definition~\ref{def:matmult} also implies that the $1\times n$ row matrix $D_{(i)}$ and the $m\times 1$ column matrix $D^{(j)}$ satisfy
\index{matrix!$(AB)_{(i)} = A_{(i)}B$ and!$(AB)^{(j)}= AB^{(j)}$}
\begin{equation}
\label{eq:colrowmatrix}
D_{(i)} = (AB)_{(i)} = A_{(i)}B \;\;\;\;{\rm and}\;\;\;\;D^{(j)} = (AB)^{(j)}= AB^{(j)}.
\end{equation}
Explicitly,  for the row version we have
\begin{equation}
\label{eq:lincomrows}
D_{(i)} = A_{(i)}B = \sum_{t=1}^p A(i,t)B_{(t)}.
\end{equation}
The column version of \ref{eq:lincomrows} is
\begin{equation}
\label{eq:lincomcols}
D^{(j)} = AB^{(j)} = \sum_{t=1}^p A^{(t)}B(t,j).
\end{equation}
Equation~\ref{eq:lincomrows} states that row $i$ of the product $AB$ is a linear combination of the rows of $B$ with coefficients from row $i$ of $A$.
Equation~\ref{eq:lincomcols} states that column $j$ of the product $AB$ is a linear combination of the columns of $A$ with coefficients from column $j$ of $B$.
Here is an example:\\
\begin{minipage}{\textwidth}
\begin{equation}
\label{eq:matrixtimesvector}
{\bf Figure: Matrix\;times\;a\;vector}
\end{equation}
\begin{center}
\includegraphics{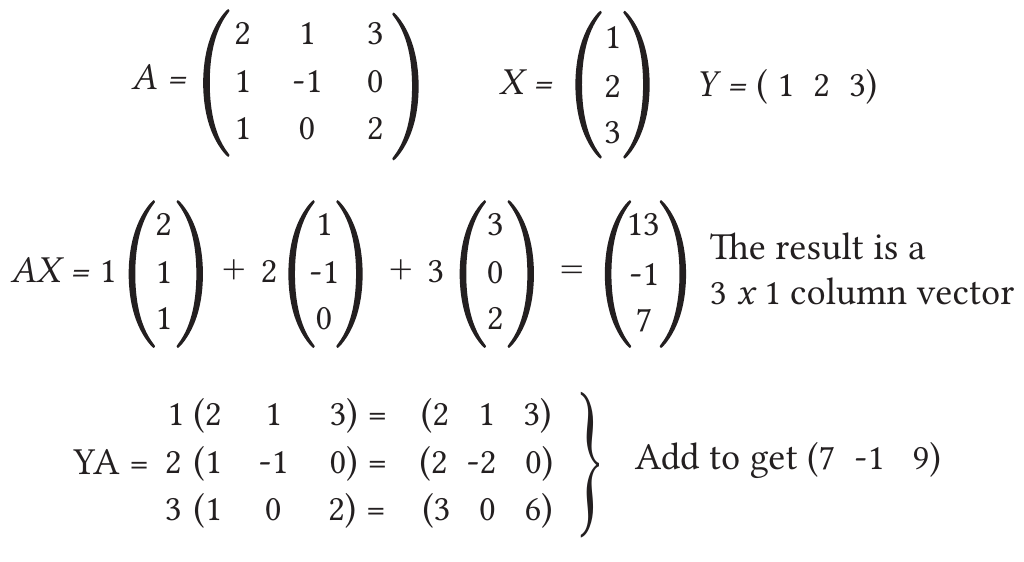}
\end{center}
\end{minipage}
\index{matrix!submatrix notation}
We need some notation for {\em submatrices} of a matrix.
Figure~\ref{eq:submatrixfncs} gives some examples ($p=5$, $n=3$) of what is needed.  Note that $f$ and $g$ are functions with domain $\underline 3$ and range $\underline 5$ (i.e., elements of ${\underline 5}^{\underline 3}$).  The function $c$ is a permutation of $\underline 3$ and $fc$ denotes the composition of $f$ and $c$ 
(i.e., $fc\in {\underline 5}^{\underline 3}$).

\begin{minipage}{\textwidth}
\begin{equation}
\label{eq:submatrixfncs}
{\bf Figure:\,Submatrix\; notation:\;}  f,\,g\in {\underline 5}^{\underline 3}\;{\bf and\;}
c = (1\,3).
\end{equation}
\begin{center}
\includegraphics{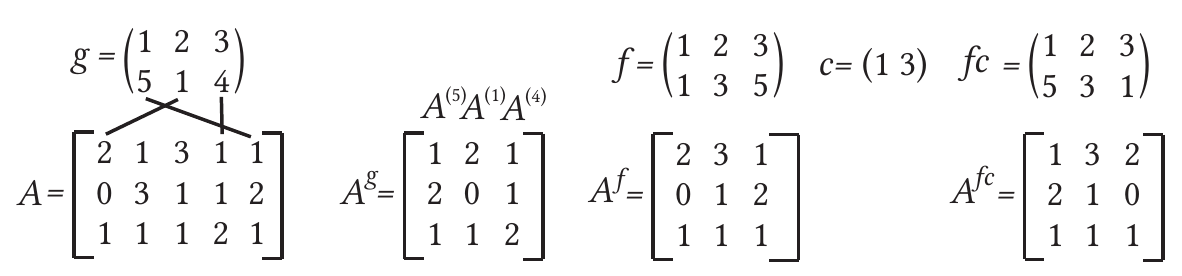}
\end{center}
\end{minipage}
\hspace*{1 in}

What we call ``submatrices'' is an extension the usual usage. Here is the formal definition:\\
%Definition submatrices 2
\index{matrix!submatrix notation}
\begin{defn}[\bfseries Submatrix notation]
\label{def:submatrices2}
Let $X$ be an $m\times n$ matrix and let $f\in {\underline m}^{\underline r}$ and
$g\in {\underline n}^{\underline s}$ be functions, $r,s >0$.
We define the $r\times s$ matrix 
\begin{equation}
\label{eq:submatfncs}
Y\equiv X^g_f \equiv X[f\Mid g] \equiv X[f(1), \ldots , f(r)\Mid g(1), \ldots, g(s)]
\end{equation} 
by $Y(p,q) = X(f(p), g(q))$.  
Suppose ${\alpha} = \{a_1, \ldots , a_r\}\subseteq {\underline m}$ and 
${\beta}= \{b_1, \ldots , b_s\}\subseteq {\underline n}$ are subsets of size $r$ and  $s$ where the $a_1<\cdots<a_r$ and $b_1<\cdots <b_s$ are in increasing order. 
Define
\begin{equation}
\label{eq:submatsets}
Y\equiv X[{\alpha}\Mid {\beta}] \equiv X[a_1, \ldots ,a_r \Mid b_1,\ldots, b_s] 
\end{equation}
by $Y(i,j)=X(a_i, b_j)$ for 
all $(i,j)\in {\underline r}\times {\underline s}$.
If ${\alpha}'={\underline m} - {\alpha}$ and ${\beta}'={\underline n} - {\beta}$ are the 
ordered complements of ${\alpha}$ and ${\beta}$, then define 
\begin{equation}
\label{eq:setsinout}
X({\alpha} \Mid {\beta}) = X[{\alpha}' \Mid {\beta}'],\;X({\alpha} \Mid {\beta}] = X[{\alpha}' \Mid {\beta}],\;
X[{\alpha} \Mid {\beta}) = X[{\alpha} \Mid {\beta}'].
\end{equation}
\end{defn}
%End definition submatrices 2
Note that in Definition~\ref{def:submatrices2} we have 
\begin{equation}
\label{eq:updowncommute}
X^g_f = (X^g)_f=(X_f)^g
\end{equation}
%\pagebreak
\begin{rem}[\bfseries Example of submatrix notation]
\label{rem:exampmatnot}
Let $m=2$ and $n=3$ with 
\begin{equation}
\label{eq:indtoent1}
X=
\left(
\begin{array}{ccc}
 1 & 2&3 \\
 4 & 5& 6
\end{array} 
\right).
\end{equation}
Let $r=3$ and $s=4$ with $f=(1\, 2\, 1)$ and $g=(2\, 3\, 2\, 1)$.  Then
\begin{equation}
\label{eq:indtoent2}
Y=X[f\Mid g]=
\left(
\begin{array}{cccc}
  2 & 3 & 2 & 1 \\
  5 & 6 & 5 & 4 \\
  2 & 3 & 2 & 1 
\end{array}
\right).
\end{equation}
\end{rem}
The index-to-entry function of the submatrix~\ref{eq:indtoent2}
refers to the domain values of $f$ and $g$. 
If $Y=X[f\Mid g]$ then $Y(2,4)=4$ but $X(2,1)=4$.\\

\begin{rem}[\bfseries Submatrices as sets or functions]
\label{rem:setsfncsrelate}
Note that \ref{eq:submatsets}  is a special case of~\ref{eq:submatfncs}.
Let ${\alpha}$ and ${\beta}$ be as in \ref{eq:submatsets}, and define $f_{\alpha} \in {\rm SNC}(r,m)$ 
and $f_{\beta} \in {\rm SNC}(s,n)$ (\ref{def:setsfunc}) by  
 ${\rm image}(f_{\alpha})={\alpha}$ and ${\rm image}(f_{\beta})={\beta}$.  
 Then 
 $$X[{\alpha}\Mid {\beta}] = X[a_1, \ldots ,a_r \Mid b_1,\ldots, b_s] = X[f_{\alpha}\Mid f_{\beta}].$$
 \end{rem}
Using Definition~\ref{def:submatrices2}, we can generalize identity~\ref{eq:colrowmatrix}.
Assume now that $D=AB$ where $A$ is  an $a\times p$ matrix and $B$ a $p\times b$ matrix. 
Let $g\in {\underline a}^{\underline m}$ and $h\in  {\underline b}^{\underline n}$ be functions.
We can think of $g$ as a ''row selection'' function so that $A_g$ is an $m\times p$ matrix and
$h$ as a ``column selection'' function so that $B^h$ is a $p\times n$ matrix.
Thus, the product $A_gB^h$ is a $m\times n$ matrix.  Then, we have
\begin{equation}
\label{eq:gencolrowmatrix}
D[g\Mid h] \equiv D_g^h \equiv (AB)_g^h = A_gB^h.
\end{equation}

%START OF CHAPTER 2
\chapter{Determinants}
We now define the determinant of an $n\times n$ matrix. In the discussion of determinants that follows, assume the  matrices have entries in the rings, $\bK$,  described in Remark~\ref{rem:specialringsfields}, all of which are Euclidean domains.   If you are interested in more generality, review  the discussions of Section~1, specifically \ref{def:ring}, \ref{rem:specialrings}, \ref{def:characteristic}, \ref{rem:unitassoc}, and do a web search for ``rings determinants.''
Recall the definition of the sign, ${\rm sgn}(f)$, of a permutation (Definition~\ref{def:sign}) and the discussion that follows that definition, including identity~\ref{eq:invtransindex}
%DETERMINANT DEFINITION
\begin{defn}[\bfseries \bf Determinant]
\index{determinant!definition}
\label{def:determinant1}
Let $A$ be an $n\times n$ matrix with entries $A(i,j)$.  The determinant, $\det(A)$, 
is defined by
\[
\det (A) = \sum_f {\rm sgn}(f)\prod_{i=1}^n A(i, f(i))
\]
where the sum is over all permutations of the set ${\underline n}=\{1,2,\ldots,n\}.$\\
\end{defn}%
The terms of the product, $\prod_{i=1}^n A(i, f(i))$, commute.  Thus, the product can be taken
in any order over the set, ${\rm Graph}(f)$ (\ref{eq:graphfunction}):
\begin{equation}
\label{eq:detovergraph}
\det(A) = \sum_f {\rm sgn}(f)\prod_{(i,j)\in {\rm Graph}(f)}A(i, j).
\end{equation}
In particular, note that 
\begin{equation}
\label{eq:graphunordered}
{\rm Graph}(f) = \{(i,f(i))\Mid i\in {\underline n}\} = \{(f^{-1}(i),i)\Mid i\in {\underline n}\}.
\end{equation}
Thus, we have
\begin{equation}
\label{eq:detunordered}
\det(A) = \sum_f {\rm sgn}(f)\prod_{{\rm Graph}(f)}A(i, j)=\sum_f {\rm sgn}(f)\prod_{i=1}^n A(f^{-1}(i), i). 
\end{equation}

Summing over all $f$ is the same as summing over all $f^{-1}$ and ${\rm sgn}(f) = {\rm sgn}(f^{-1})$.
Therefore, the second sum in \ref{eq:detunordered} can be written
\[
\sum_{f^{-1}} {\rm sgn}(f^{-1})\prod_{i=1}^n A(f^{-1}(i), i)= \sum_{f} {\rm sgn}(f)\prod_{i=1}^n A(f(i), i).
\]
Thus, we have the important identity
%EQUATION ROW AND COLUMN DETERMINANTS
\begin{equation}
\label{eq:detrowcol}
\index{determinant!row, column forms}
\det(A) = \sum_f {\rm sgn}(f)\prod_{i=1}^n A(i, f(i)) = \sum_{f} {\rm sgn}(f)\prod_{i=1}^n A(f(i), i).
\end{equation}
The first sum in \ref{eq:detrowcol} is called the row form of the determinant and the second is called the
column form.  In the first sum, the domain of  $f$ is the set of row indices and the range is the set of column indices.  In the second, the domain is the set of column indices and the range is the set of row indices.
%DEFINITION TRANSPOSE OF MATRIX
\begin{definition}[\bfseries Transpose of a matrix]
\label{def:transpose}
\index{matrix!transpose}
Let $A$ be an $n\times n$ matrix with entries $A(i,j)$.  
The {\em transpose} of $A$ is the matrix $A^T$ defined by $A^T(i,j) = A(j,i)$.\\
\end{definition}
%REMARK
\begin{remark}[\bfseries Transpose  basics]
\label{rem:transposeprod}
The transpose of 
$
A=
\left(
\begin{array}{cc}
a_{11}  &  a_{12}  \\
a_{21}  &  a_{22}   
\end{array}
\right)
$
is
$
A^T=
\left(
\begin{array}{cc}

a_{11}  &  a_{21}  \\
a_{12}  &  a_{22}   
\end{array}
\right).
$
Note that the transpose of a product is the product of the transposes in reverse order: $(AB)^T = B^TA^T$.
Recall remark~\ref{rem:indtoent} concerning the index-to-entry function and
note that $A^T(2,1)=a_{12}=A(1,2)$.  
Suppose we take $A=(a_{ij})$ to be a $4\times 4$ matrix.  
Let $X=A(2\Mid 3)$ be the  $3\times 3$ submatrix
\[
X=
\left(
\begin{array}{ccc}
a_{11}&a_{12}&a_{14}\\
a_{31}&a_{32}&a_{34}\\
a_{41}&a_{42}&a_{44}
\end{array}
\right)
\;\;\mathrm{and}\;\;
X^T=
\left(
\begin{array}{ccc}
a_{11}&a_{31}&a_{41}\\
a_{12}&a_{32}&a_{42}\\
a_{14}&a_{34}&a_{44}
\end{array}
\right).
\]
Recall remark~\ref{rem:indtoent} concerning the index-to-entry function and
note that $X(1,3)=a_{14}$ and $X^T(1,3)=a_{41}$. 
The index-to-entry function has domain $\underline{3}\times \underline{3}$
for these submatrices.
Thus, $X(i,j) = X^T(j,i)$, $1\leq i,j \leq 3$ as required by definition~\ref{def:transpose}.
Starting with $A$, we have $(A(2\Mid 3))^T = A^T(3\Mid 2)$.
You can first take a submatrix of $A$ and then transpose that or first transpose $A$ 
and then take the appropriate submatrix.
Using the submatrix notation of~\ref{def:submatrices2}, the rule is
\begin{equation}
\label{eq:trasubmatcom}
\index{matrix!transpose and submatrices}
(X[f\Mid g])^T = X^T[g\Mid f]
\end{equation}
As an example, consider 
\[
X=
\left(
\begin{array}{ccc}
 1 & 2&3 \\
 4 & 5& 6
\end{array} 
\right)
\;\;\mathrm{and}\;\;
X^T=
\left(
\begin{array}{cc}
 1 & 4 \\
 2 & 5 \\
 3 &6
\end{array} 
\right).
\]
Let $r=3$ and $s=4$ with $f=(1\, 2\, 1)$ and $g=(2\, 3\, 2\, 1)$.  Then
\begin{equation}
\label{eq:indtoent3}
(X[f\Mid g])^T=
\left(
\begin{array}{cccc}
  2 & 3 & 2 & 1 \\
  5 & 6 & 5 & 4 \\
  2 & 3 & 2 & 1 
\end{array}
\right)^T
= 
\left(
\begin{array}{ccc}
2&5&2\\
3&6&3\\
2&5&2\\
1&4&1
\end{array}
\right)
= X^T[g\Mid f].
\end{equation}
\hspace*{1in}\\
\end{remark}
%END REMARK
%THEOREM DETERMINANT OF TRANSPOSE
\begin{theorem}[\bfseries Determinant of transpose]
\label{thm:detoftranpose}
\index{determinant!of transpose}
Let $A$ be an $n\times n$ matrix with entries $A(i,j)$ and let $A^T$ be its transpose.
Then
\[ \det(A) = \det(A^T). \]
\begin{proof}
We use \ref{eq:detrowcol}.
\[ \det(A) = \sum_f {\rm sgn}(f)\prod_{i=1}^n A(i, f(i)) = \sum_{f} {\rm sgn}(f)\prod_{i=1}^n A^T(f(i), i) = \det(A^T).\]
\end{proof}
\end{theorem}

%\begin{equation}
%\label{eq:propdet}
%{\bf Elementary\;\,properties\;\,of\;\,determinants}
%\end{equation}
\section* {Elementary properties of determinants}

We now derive some ``elementary'' properties of the determinant -- properties that follow directly from the definition (\ref{def:determinant1}) using routine (but not necessarily short) computations.  It is common to use Greek letters for permutations so we switch to that convention.  

First, note that if $\varphi$ and $\gamma$ are permutations of {\underline n} then
\begin{equation}
\label{eq:differentorder}
\prod_{i=1}^n A(i, \varphi(i)) = \prod_{i=1}^n A(\gamma(i), \varphi(\gamma(i)))=\prod_{(i,j)\in {\rm Graph}(\varphi)}A(i, j).
\end{equation}

Let $A$ be an $n\times n$ matrix and let $\gamma$ be a permutation on ${\underline n}$.  As a sequence of columns, we write $A=(A^{(1)}, \dots, A^{(i)}, \ldots, A^{(n)})$.  We define
\[
A^\gamma = (A^{\gamma(1)}, \dots, A^{\gamma(i)}, \ldots, A^{\gamma(n)}),
\]
This notation is a special case of Definition~\ref{def:submatrices2} ($r=s=m=n$, $g=\gamma$, $f$ the identity).
From~\ref{eq:detrowcol},
\[
\det(A^\gamma)=\sum_\varphi {\rm sgn}(\varphi)\prod_{i=1}^n A^\gamma(\varphi(i), i)
\]
where 
\[
\prod_{i=1}^n A^\gamma(\varphi(i), i) = \prod_{i=1}^n A(\varphi(i), \gamma(i))= 
\prod_{i=1}^n A(\varphi\gamma^{-1}(i), i).
\]
Thus,   
\begin{equation}
\label{eq:offbyperm}
 \det(A^\gamma)=\sum_\varphi {\rm sgn}(\varphi)\prod_{i=1}^n A(\varphi \gamma^{-1}(i), i).
\end{equation}
From~\ref{eq:offbyperm}, we get the very important symmetry property of the determinant function
under permutation of columns (by a similar argument, rows) which states that 
$\det(A^\gamma)={\rm sgn}(\gamma) \det(A)\;$ 
(for rows, $\det(A_\gamma)={\rm sgn}(\gamma) \det(A)$):
\begin{equation}
\label{eq:detsymmetry}
\index{determinant!symmetry properties}
\det(A^\gamma)=  
{\rm sgn}(\gamma)\sum_\varphi  {\rm sgn}(\varphi \gamma^{-1})\prod_{i=1}^n A(\varphi \gamma^{-1}(i), i))
={\rm sgn}(\gamma)\det(A).
\end{equation}

The next definition is fundamental to the study of determinants.
%DEFINITION MULTILINEAR
\begin{defn}[\bfseries \bf Multilinear function]  
\label{def:multilinear}
\index{function!multilinear definition}
\index{function!linear}
\index{linear function}
Let $V_1, \ldots, V_n$ be  vector spaces over a field 
$\bF$ and let $W = \times_1^n V_i = \{(x_1, \ldots, x_n)\Mid x_i \in V_i, i=1, \ldots, n\}$ be the
direct (Cartesian) product of these $V_i$. 
A function $\Phi$ from $W$ to $\bF$ is {\em multilinear} if it is linear separately in each variable:  
For  $c, d\in{\bF}$ and for $t=1, \ldots, n$ ,
\[ 
\Phi(x_1,\ldots, (cx_t + dy_t), \ldots x_n) =
c\Phi(x_1,\ldots, x_t, \ldots x_n) + d\Phi(x_1,\ldots, y_t, \ldots x_n).
\]
If $n=1$ then $\Phi$ is a {\em linear function} from $V_1$ to $\bF$.
\end{defn}
\[
\framebox{\textbf{\emph{The determinant is a multilinear function. }}}
\]
\index{determinant!multilinear function}
If $A$ is an $n\times n$ matrix over $\bF$, we can regard $A$ as an ordered sequence,  
$(x_1, \ldots, x_n)$, of vectors in $\bF^n$ where either $x_i = A^{(i)}, 1=1,\ldots n\,,$ are the columns of $A$ or $x_i=A_{(i)}, 1=1,\ldots n$, are the rows of $A$.  In either case, rows or columns, $\det(A) = \det(x_1, \ldots, x_n)$ is a multilinear function from 
$W = \times_1^n V_i$ to $\bF$.

To verify multilinearity (row version), let 
$B_{(t)} = (B(t,1), \ldots, B(t,n))$ for a fixed $t$.
Replace row $t$ of $A$ with $B_{(t)}$, to get $\hat{A}$:
$$\hat{A} =   (A_{(1)}, \ldots , B_{(t)},  \ldots A_{(n)})$$
Replace row $t$ of $A$ with $cA_{(t)} + dB_{(t)}$, to get $\tilde{A}$:
$$\tilde{A} = (A_{(1)}, \ldots , (cA_{(t)} + dB_{(t)}),  \ldots A_{(n)}).$$
Using the definition of the determinant we compute
$$\det (\tilde{A})=\sum_\varphi {\rm sgn}(\varphi)A(1,\varphi(1))\cdots 
(cA(t, \varphi(t))+dB(t,\varphi(t))\cdots A(n,\varphi(n))=$$
$$c\sum_\varphi {\rm sgn}(\varphi)A(1,\varphi(1))\cdots A(t, \varphi(t))\cdots A(n,\varphi(n))+$$
$$d\sum_\varphi {\rm sgn}(\varphi)A(1,\varphi(1))\cdots B(t,\varphi(t))\cdots A(n,\varphi(n)).$$
Thus we have
\begin{equation}
\label{eq:multidetrow}
\det (\tilde{A})= c\det{A} + d\det(\hat{A})
\end{equation}
which verifies that of the determinant is a multilinear function (\ref{def:multilinear}).

\begin{defn} [\bf Alternating multilinear]
\label{def:alternating}
\index{function!alternating multilinear}
A multilinear function $\Phi$ from $\times^n V$ to $\bF$ is {\em alternating} if for any transposition 
$\tau=(s\;t)$ on ${\underline  n}$,
\[
\Phi(x_1, \ldots, x_i, \ldots x_n) = - \Phi(x_{\tau(1)}, \ldots, x_{\tau(i)}, \ldots x_{\tau(n)}).
\]
Or, equivalently,
\[
\Phi(x_1, \ldots , x_s, \ldots, x_t, \ldots, x_n) = -\Phi(x_1, \ldots , x_t, \ldots, x_s, \ldots, x_n).
\]
\end{defn}
In particular, note that if $\Phi$ is alternating and $x_s=x_t=x$ then the identities
of~\ref{def:alternating} become
\begin{equation}
%alternating property
\label{eq:altmult}
\Phi(x_1, \ldots , x, \ldots, x, \ldots, x_n) = -\Phi(x_1, \ldots , x, \ldots, x, \ldots, x_n).
\end{equation}
Thus, $\Phi(x_1, \ldots , x, \ldots, x, \ldots, x_n)=0$. 
We use the fact that $\bF$ (or $\bK$) is of characteristic $0$ 
(Definition~\ref{def:characteristic} and Remark~\ref{rem:specialringsfields}).
In particular, if any pair of vectors, $x_s$ and $x_t$, are linearly dependent then
\begin{equation}
\Phi(x_1, \ldots , x_s, \ldots, x_t, \ldots, x_n) = 0
\end{equation}
For suppose $x_t$ and $x_s$ are nonzero and $x_t=cx_s$ where $0\neq c\in \bF$ (or $\bK$). 
Then
\[
\Phi(x_1, \ldots , x_s, \ldots, x_t, \ldots, x_n) = c\Phi(x_1, \ldots , x_s, \ldots, x_s, \ldots, x_n)=0.
\]

\[
\framebox{\textbf{\emph{The determinant is alternating multilinear}}}
\]
\index{determinant!alternating multilinear}
Identity \ref{eq:multidetrow} shows that $\det(x_1, \ldots, x_n)$ is a multilinear function of its 
rows or columns.  Identity~\ref{eq:detsymmetry} implies that if $A$ is an $n\times n$ matrix,
$\gamma$ is a permutation on ${\underline n}$ and (column form)
$A^\gamma = (A^{\gamma(1)}, \dots, A^{\gamma(i)}, \ldots, A^{\gamma(n)})$ or (row form)
$A_\gamma = (A_{\gamma(1)}, \dots, A_{\gamma(i)}, \ldots, A_{\gamma(n)})$  then
\begin{equation}
\label{eq:sgndet}
\det(A^\gamma) = {\rm sgn}(\gamma) \det (A)\;\,{\rm and\;\,}
\det(A_\gamma) = {\rm sgn}(\gamma) \det (A).
\end{equation}
Thus, if we take $\gamma=\tau = (s\;t)$ we get (using ${\rm sgn}(\tau) = -1$)
\begin{equation}
\label{eq:detisalt}
\det(A^{\tau(1)}, \dots, A^{\tau(n)}) = -\det(A^{(1)}, \dots, A^{(n)}).
\end{equation}
This identity shows that the determinant is an alternating (multilinear) function of its columns
(and, similarly, its rows).  Thus, $\det(A)=0$ if any two rows or columns are the same (or are 
linearly dependent).\\
\begin{minipage}{\textwidth}
\begin{equation}
\label{eq:gendirectsums}
{\rm\bf General\;\; direct\;\;sums}
\end{equation}
We start with the standard definition. See \ref{def:submatrices2} for related notation.\\
\end{minipage}
%DIRECT SUM - SIMPLE
\begin{defn}[\bfseries  Direct sum of matrices]
\label{def:dirsummatrices}
\index{matrix!direct sum}
An $n\times n$ matrix $A=B\oplus C$ is called the direct sum of  a $r \times r$ matrix,  
$B$, and  $s\times s$ matrix, $C$,  if  $0<r<n$, $n=r+s$ and
\[
A[{\underline r}\Mid {\underline r}] = B\;\;\;\;
A({\underline r}\Mid {\underline r}) = C\;\; \;\;
A[{\underline r}\Mid {\underline r}) = \Theta_{r,s}\;\;\;\;
A({\underline r}\Mid {\underline r}] = \Theta_{s,r}
\]
where $\Theta_{r,s}$ is the $r\times s$ zero matrix.
\end{defn}

The determinant of a direct sum is the product of the determinants of the summands: 
\begin{equation}
\label{eq:dirsummatrices}
\det(B\oplus C)=\det(A) = 
\det(A[{\underline r}\Mid {\underline r}])\det(A({\underline r}\Mid {\underline r})= \det(B)\det(C).
\end{equation}
This identity can be proved directly from the definition of the determinant.

As an example of~\ref{eq:dirsummatrices}, let $A$ be a $4\times 4$ matrix of integers ($A\in {\bf M}_{4,4}(\bZ)):$
\begin{equation}
\label{eq:dirsumsimple}
A=
\left(
\begin{array}{cccc}
1 & 2 & 0 & 0\\
2 & 3 & 0 & 0 \\
0 & 0 & 3 & 4\\
0 & 0 & 4 & 1
\end{array}
\right)
=
\left(\begin{array}{cc}1&2\\2&3 \end{array}\right) \oplus 
\left(\begin{array}{cc}3&4\\4&1 \end{array}\right).
\end{equation}
Then $\det(A)=\det(B)\det(C) = (-1)(-13) = 13.$\\

\begin{defn}[\bfseries General direct sum]
\label{def:gendirsum}
\index{matrix!direct sum - general}
Let $A$ be an $n\times n$ matrix and let $X,Y\in \bP_r(n)$ where $0<r<n$ (\ref{rem:notation}). 
We say that $A$ is a {\em general direct sum} relative to $X$ and $Y$ 
of  an $r \times r$ matrix $B$ and  $s\times s$ matrix $C$  if  $n=r+s$ and
\begin{equation*}
A[X\Mid Y]=B,\;\;
A(X\Mid Y)=C,\;\;
A[X\Mid Y)=\Theta_{r,s}\;\;{\rm and}\;\;
A(X\Mid Y] =\Theta_{s,r}.
\end{equation*} 
We write 
\[
A = B \oplus_X^Y C.
\]
\end{defn}  

As an example of~\ref{def:gendirsum}, let $A$ be a $4\times 4$ matrix of integers ($A\in {\bf M}_{4,4}(\bZ)):$
\begin{equation}
\label{eq:dirsumgeneral}
A=
\left(
\begin{array}{cccc}
1 & 0 & 2 & 0\\
0 & 3 & 0 & 4 \\
2 & 0 & 3 & 0\\
0 & 4 & 0 & 1
\end{array}
\right)
=
\left(\begin{array}{cc}1&2\\2&3 \end{array}\right) \oplus_X^Y 
\left(\begin{array}{cc}3&4\\4&1 \end{array}\right)\;\;\;X=Y=\{1,3\}.
\end{equation}
Note that the direct sum of Definition~\ref{def:dirsummatrices} is  a special case of
Definition~\ref{def:gendirsum} (take $X=Y=\underline{r}$). 
The matrix of example \ref{eq:dirsumgeneral} can be transformed by row and column interchanges to that of \ref{eq:dirsumsimple}, thus the determinants differ by only the sign.
Direct computation gives $\det(A)=13$ for $A$ in \ref{eq:dirsumgeneral} .
An example of the transformation process by row and column interchanges is given in 
Figure~\ref{eq:laplacerowcol}.  We take a different  approach in order to develop  precise combinatorial and analytic tools for future use.\\
\begin{minipage}{\textwidth}
\begin{equation}
\label{eq:laplacerowcol}%this label can be found even if in minipage
{\bf Figure:\,Reduce\; submatrix}\; A[X\Mid Y]\; {\bf to\;initial\;position} 
\end{equation}
\begin{center}	
\includegraphics{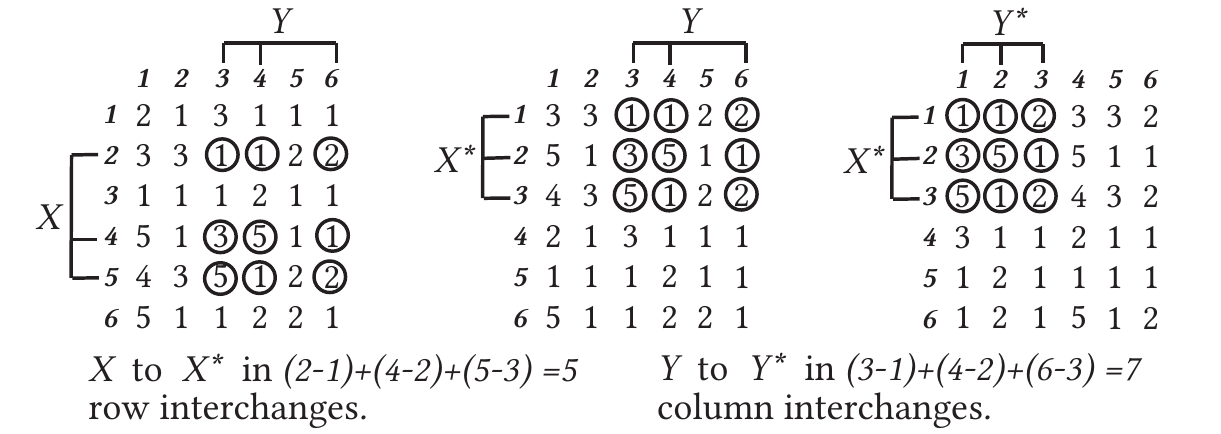}
\end{center}
\end{minipage}
\hspace*{1 in}\\
\begin{defn}[\bfseries The set $S_X^Y$]
\label{def:SXY}
\index{determinant!general direct sum!begin discussion}
Let $S={\rm PER(n)}$ be the permutatons of $\underline{n}$, and let $X, Y\in \bP_k(n)$ be subsets of ${\underline n}$ of size $k$.
Let 
\begin{equation*}
S_X^Y = \{\sigma\Mid \sigma \in S\; {\rm and}\; \sigma(X)=Y\}.
\end{equation*}
Let $\gamma\in S_X^Y$ and let $X'$ and $Y'$ denote the complements of $X$ and $Y$.  Suppose the restrictions (\ref{eq:restriction}) $\gamma_X$ and   
$\gamma_{X'}$ of  $\gamma$ are strictly increasing:
\[
\gamma_X\in {\rm SNC}(X,Y)\:\; {\rm and}\;\;
\gamma_{X'}\in {\rm SNC}(X',Y'). 
\]
Then $\gamma$ is called the {\em canonical} element of $S_X^Y$.
Note that $\gamma_X$ and $\gamma_{X'}$  are unique since $|X|=|Y|$ and $|X'|=|Y'|$.
\end{defn}
%REMARK
\begin{remark}[\bfseries Example of $S_X^Y$]
\label{rem:exampleblock}
Take $S={\underline 6}$ and $X, Y\in \bP_3(6)$ where $X=\{2,4,5\}$ and $Y=\{3,4,6\}$.
The set $S_X^Y\,$  consists of all permutations 
$\sigma$ such that $\sigma(\{2,4,5\})=\{3,4,6\}$ (i.e., $\sigma(X)=Y$). 
This implies (from the definition of a permutation) that
$\sigma(\{1,3,6\})=\{1,2,5\}$ (i.e., $\sigma(X')=Y'$).
In other words,
\begin{equation}
\label{eq:graphsigma}
{\rm Graph}(\sigma)\subseteq (X\times Y) \cup (X'\times Y').
\end{equation}
\end{remark}
%\hspace{1 in}
%END REMARK
%REMARK
\begin{remark}[\bfseries Canonical element, $\gamma$]
We continue Remark~\ref{rem:exampleblock}.
Designate the elements of $X$ in order as $(x_1, x_2, x_3) = (2,4,5)$ and 
$Y$ in order as $(y_1,y_2,y_3)=(3,4,6)$.  
Similarly, $X'$ in order is $(x'_1, x'_2, x'_3) = (1,3,6)$ and 
$Y'$ in order is $(y'_1, y'_2, y'_3) = (1,2,5)$.
In two line notation, let
\begin{equation}
\label{eq:sampleperm}
\gamma =
\left(
\begin{array}{cccccc}
2 & 4  & 5 & 1 & 3 & 6\\
3 & 4  & 6 & 1 & 2 & 5
\end{array}
\right)
=
\left(
\begin{array}{cccccc}
1 & 2  & 3 & 4 & 5 & 6\\
1 & 3  & 2 & 4 & 6 & 5
\end{array}
\right).
\end{equation}
The permutation $\gamma$ is the canonical element of $S_X^Y$ (\ref{def:SXY}).
\end{remark}
%CANONICAL EXAMPLE
%REMARK
\begin{remark}[\bfseries Restrictions]
The $\gamma$ of \ref{eq:sampleperm} has the following restrictions to $X$ and $X'$:
\begin{equation}
\label{eq:xyrestrict}
\gamma_X =
\left(
\begin{array}{ccc}
2 & 4  & 5 \\
3 & 4  & 6    
\end{array}
\right)
\;\;{\rm and}\;\;
\gamma_{X'}=
\left(
\begin{array}{ccc}
1 & 3 & 6\\
1 & 2 & 5    
\end{array}
\right).
\end{equation}
Note from \ref{eq:xyrestrict} that $\gamma_X\in {\rm SNC}(X,Y)$ and
$\gamma_{X'}\in {\rm SNC}(X',Y')$ as required by the definition of  
$\gamma$ in \ref{def:SXY}.
\end{remark}
%END REMARK
%REMARK TYPICAL ELEMENT
\begin{remark}[\bfseries  Typical element]
A typical element $\sigma\in S_X^Y$ has restrictions to $X$ and $X'$ that are injections:
\begin{equation}
\label{eq:xyrestrict2}
\sigma_X =
\left(
\begin{array}{ccc}
2 & 4  & 5 \\
6& 3  & 4    
\end{array}
\right)
\;\;{\rm and}\;\;
\sigma_{X'}=
\left(
\begin{array}{ccc}
1 & 3 & 6\\
5 & 2 & 1    
\end{array}
\right).
\end{equation}
Using the fact that $(y_1,y_2,y_3)=(3,4,6)$ and 
$(y'_1, y'_2, y'_3) = (1,2,5)$, the second lines of $\sigma_X$ and $\sigma_{X'}$ (\ref{eq:xyrestrict2}) 
can  be specified as permutations of $(y_1,y_2,y_3)$ and $(y'_1, y'_2, y'_3)$ : 
\begin{equation}
\label{eq:nuspecified}
(6,3,4) = (y_3,y_1,y_2)=(y_{\nu(1)}, y_{\nu_1(2)}, y_{\nu(3)})\;\;
\;\; \nu=\left(\begin{array}{ccc} 1&2&3\\ 3&1&2\end{array}\right)
\end{equation}
and
\begin{equation}
\label{eq:muspecified}
(5,2,1) = (y'_3,y'_2,y'_1)=(y'_{\mu(1)}, y'_{\mu(2)}, y'_{\mu(3)})\;\;\;\;\mu=\left(\begin{array}{ccc} 1&2&3\\ 3&2&1\end{array}\right).
\end{equation}
\end{remark}
%END REMARK
%\hspace{1 in}\\
%LEMMA
\begin{lemma}[\bfseries Description of $S_X^Y$]
\label{lem:descriptsxy}
Let $X$ be ordered $x_1 <\cdots< x_k$ (i.e.,  {\em the sequence $x_1, \ldots, x_k$ is ordered as integers})  and $X'$ ordered $x'_1< \cdots < x'_{n-k}$.
Similarly, let $Y$ be ordered $y_1< \cdots < y_k$ and $Y'$ ordered $y'_1 < \cdots < y'_{n-k}$.
Let $\gamma$ be the canonical representative of $S_X^Y$.
\[
\gamma =
\left(
\begin{array}{cccccc}
x_1&\ldots&x_k &x'_1&\ldots&x'_{n-k}\\
y_1&\ldots&y_k &y'_1&\ldots&y'_{n-k}
 \end{array}
\right).
\]
For $\sigma \in S_X^Y$, define $\nu\in {\rm PER}(k)$ by 
$(\sigma(x_1), \ldots, \sigma(x_k))=(y_{\nu(1)}, \ldots , y_{\nu(k)}).$  
Likewise, define $\mu\in {\rm PER}(n-k)$ by
$(\sigma(x'_1), \ldots, \sigma(x'_{n-k}))=(y'_{\mu(1)}, \ldots , y'_{\mu(n-k)}).$
Then $S_X^Y$ is the set of all permutations of the form
\begin{equation}
\label{eq:generateblock}
\gamma^{\nu\mu} =
\left(
\begin{array}{cccccc}
x_1&\ldots&x_k &x'_1&\ldots&x'_{n-k}\\
y_{\nu(1)}&\ldots&y_{\nu(k)} &y'_{\mu(1)}&\ldots&y'_{\mu(n-k)} \end{array}
\right)
\end{equation}
for $\nu\in {\rm PER}(k)$ and $\mu\in {\rm PER}(n-k).$
%PROOF
\begin{proof}
%The definition of $\gamma^{\nu\mu}$ is a routine generalization of the discussion of identities 
%\ref{eq:sampleperm},
%\ref{eq:xyrestrict},
%\ref{eq:xyrestrict2},
%\ref{eq:nuspecified} and
%\ref{eq:muspecified}.
The fact that the set 
\[\{\gamma^{\nu\mu}\Mid  \nu\in {\rm PER}(k)\; {\rm and}\; \mu\in {\rm PER}(n-k)\} = S_X^Y\]
follows from the requirement that $\sigma \in S_X^Y$ if and only if the restrictions,
$\sigma_X$ and $\sigma_{X'}$, are injections. 
Thus, $(\sigma_X(x_1),\ldots, \sigma_X(x_k))=(y_{\nu(1)}, \dots, y_{\nu(k)})$ uniquely defines $\nu$
and $(\sigma_X'(x'_1),\ldots, \sigma_X'(x'_{n-k})=(y_{\mu(1)}, \dots, y_{\mu(n-k)})$ uniquely defines 
$\mu$.
\end{proof}
\end{lemma}
%END LEMMA
%LEMMA
\begin{lemma}[\bfseries Signs of $S_X^Y$ elements]
\label{lem:signxyelement}
We use the terminology of Lemma~\ref{lem:descriptsxy}.
Let $X$ be ordered as integers $x_1, \ldots, x_k$ and $X'$ ordered $x'_1, \ldots, x'_{n-k}$.
Similarly, let $Y$ be ordered $y_1, \ldots, y_k$ and $Y'$ ordered $y'_1, \ldots, y'_{n-k}$.
Let $\gamma$ be the canonical representative.

\[
\gamma=
\left(
\begin{array}{cccccc}
x_1&\ldots&x_k &x'_1&\ldots&x'_{n-k}\\
y_1&\ldots&y_k &y'_1&\ldots&y'_{n-k}
 \end{array}
\right),
\]
and let 
\[
\gamma^{\nu\mu}=
\left(
\begin{array}{cccccc}
x_1&\ldots&x_k &x'_1&\ldots&x'_{n-k}\\
y_{\nu(1)}&\ldots&y_{\nu(k)} & y'_{\mu(1)}&\ldots&y'_{\mu(n-k)} \end{array}
\right).
\]
Then 
\begin{equation}
\label{eq:signxyelement}
{\rm sgn}(\gamma^{\nu\mu})={\rm sgn}(\gamma){\rm sgn}(\nu){\rm sgn}(\mu).
\end{equation}
%PROOF
\begin{proof}
Note that the second line of $\gamma^{\nu\mu}$ can be converted to the second line of $\gamma$ by first transposition sorting
$\nu$ to transform $(y_{\nu(1)}, \ldots, y_{\nu(k)})$ to $(y_1, \ldots, y_k)$ and then transposition sorting
$(y'_{\mu(1)}, \ldots, y'_{\mu(n-k)})$ to $(y'_1, \ldots, y'_{n-k})$.
\end{proof}
\end{lemma}

%NEW LEMMA
\begin{lemma}[\bfseries Sign of $\gamma_y$]
\label{lem:signfactor}
We use the terminology of Lemma~\ref{lem:signxyelement}.
Let $Y$ be ordered as integers $y_1, \ldots, y_k$ and $Y'$ ordered $y'_1, \ldots, y'_{n-k}$ ($1\leq k <n$).
Let $\gamma_y$ be the permutation
%\[
\[
\gamma_y =
\left(
\begin{array}{cccccc}
1   &\ldots &k     &k+1&\ldots&n\\
y_1&\ldots&y_k &y'_1&\ldots&y'_{n-k}
 \end{array}
\right).
\]
Then ${\rm sgn}(\gamma_y) = (-1)^{\sum_{i=1}^k (y_i-i)}=(-1)^{\sum_{i=1}^k y_i}(-1)^{k(k+1)/2}$.
\begin{proof}
The proof is by induction on $k$.  Suppose $k=1$. 
Then 
\[
\gamma_y =
\left(
\begin{array}{cccc}
1    &2 &\ldots&n\\
y_1&y'_1&\ldots&y'_{n-1}
 \end{array}
\right).
\]
Since $y_1'< \cdots < y_{n-1}'$, we must do exactly $y_1-1$ transpositions of $y_1$ with the $y'_i$ to get to the
sequence $1, 2, \ldots, n$.  Thus, for $k=1$, ${\rm sgn}(\gamma_y) = (-1)^{y_1-1}$ which proves the lemma
for $k=1$.  Next, assume the case $k-1$ and consider
\[
\gamma_y =
\left(
\begin{array}{cccccc}
1   &\ldots &k     &k+1&\ldots&n\\
y_1&\ldots&y_k &y'_1&\ldots&y'_{n-k}
 \end{array}
\right).
\]
First, insert $y_k$ into its proper  position in $y_1'< \cdots < y_{n-k}'$.  This can be done in
$y_k - 1 - (k-1) = y_k-k$ transpositions since $y_1<\cdots<y_{k-1}<y_k$ and thus $y_{k}$ doesn't have to be transposed with 
the $k-1$ numbers $y_1, \ldots, y_{k-1}$. 
Thus, we have ${\rm sgn}(\gamma_y) = (-1)^{(y_k - k)} {\rm sgn}(\gamma_y')$ where
\[
\gamma_y' =
\left(
\begin{array}{cccccc}
1   &\ldots &k-1      &k    &\ldots&n\\
y_1&\ldots&y_{k-1} &y'_1&\ldots&y'_{n-(k-1)}
 \end{array}
\right).
\]
Applying the induction hypothesis (case $k-1$) to $\gamma_y'$ gives 
${\rm sgn}(\gamma_y) = (-1)^{\sum_{i=1}^k (y_i-i)}$ which 
was to be shown.  The term $(-1)^{k(k+1)/2}$ in the statement of the lemma comes from writing
$\sum_{i=1}^k (y_i-i) = \sum_{i=1}^k y_i - \sum_{i=1}^k i$ and using  the fact that 
$\sum_{i=1}^k i = k(k+1)/2$.
\end{proof}
\end{lemma}
%END LEMMA
%\hspace{1 in}\\
%EXAMPLE COMPUTING SGN GAMMA_Y
\begin{remark}[\bfseries Example of computing ${\rm sgn}(\gamma_y)$]
As an example of Lemma~\ref{lem:signfactor},
take $k=3$ and $n=6$.  Take
$(y_1, y_2, y_3) = (3, 4, 6)$ and  $(y'_1, y'_2, y'_3) = (1, 2, 5)$. 
Thus, 
\[
\gamma_y =
\left(
\begin{array}{cccccc}
1  & 2  & 3 & 4 & 5 & 6 \\
3  & 4  & 6 & 1 & 2 & 5     
\end{array}
\right).
\]
From Lemma~\ref{lem:signfactor}, 
${\rm sgn}(\gamma_y) = (-1)^{\sum_{i=1}^k (y_i-i)} = (-1)^{(3-1)+(4-2)+(6-3)}=-1.$
In cycle form, $\gamma_y = (1,3,6,5,2,4)$.
\end{remark}
%END EXAMPLE

The next lemma computes ${\rm sgn}(\gamma)$ and ${\rm sgn}(\gamma^{\nu\mu}).$   We use the terminology of  Lemmas~\ref{lem:signxyelement} and~\ref{lem:signfactor}.
\vfil
%LEMMA
\begin{lemma}[\bfseries Signs of  $\gamma$ and $\gamma^{\nu\mu}$]
\label{lem:signcanon}
Let $\gamma$ and $\gamma^{\nu\mu}$ be as in \ref{lem:signxyelement}.
Let $X$ be ordered as integers $x_1, \ldots, x_k$, $1\leq k< n$, and $X'$ be ordered $x'_1, \ldots, x'_{n-k}$.
Similarly, let $Y$ be ordered $y_1, \ldots, y_k$ and $Y'$ be ordered $y'_1, \ldots, y'_{n-k}$.
For simplicity, let ${\sum_{i=1}^k x_i}\equiv \sum X$ and ${\sum_{i=1}^k y_i}=\sum Y.$
Then,
\begin{equation*}
{\rm sgn}(\gamma) = (-1)^{\sum X}(-1)^{\sum Y}\;\;\;{\rm and}\;\;\;
{\rm sgn}(\gamma^{\nu\mu}) = (-1)^{\sum X}(-1)^{\sum Y}{\rm sgn}(\nu){\rm sgn}(\mu).
\end{equation*}
%PROOF
\begin{proof}
%Write $\gamma = \gamma_y \,\gamma_x$ where
Recall
\[
\gamma=
\left(
\begin{array}{cccccc}
x_1&\ldots&x_k &x'_1&\ldots&x'_{n-k}\\
y_1&\ldots&y_k &y'_1&\ldots&y'_{n-k}
 \end{array}
\right).
\]
Let
\begin{equation}
\label{eq:gammay}
\gamma_y =
\left(
\begin{array}{cccccc}
1   &\ldots &k     &k+1&\ldots&n\\
y_1&\ldots&y_k &y'_1&\ldots&y'_{n-k}
 \end{array}
\right)
\end{equation}
and  let
\begin{equation}
\label{eq:gammax}
\gamma_x =
\left(
\begin{array}{cccccc}
1   &\ldots &k     &k+1&\ldots&n\\
x_1&\ldots&x_k &x'_1&\ldots&x'_{n-k}
 \end{array}
\right).
\end{equation}
Note that $\gamma\,\gamma_x = \gamma_y$.  
Applying Lemma~\ref{lem:signfactor} to $\gamma_x$ and $\gamma_y$ we get
\[
{\rm sgn}(\gamma) \,(-1)^{\sum_{i=1}^k  x_i}(-1)^{k(k+1)/2}=(-1)^{\sum_{i=1}^k y_i}(-1)^{k(k+1)/2}
\]

and thus
${\rm sgn}(\gamma) = (-1)^{\sum_{i=1}^k x_i}(-1)^{\sum_{i=1}^k y_i} = (-1)^{\sum X}(-1)^{\sum Y}.$
From \ref{eq:signxyelement} we obtain
\begin{equation}
\label{eq:signgammanumu}
{\rm sgn}(\gamma^{\nu\mu}) = (-1)^{\sum X}(-1)^{\sum Y}{\rm sgn}(\nu){\rm sgn}(\mu).
\end{equation}
\end{proof}
%END PROOF
\end{lemma}
%END LEMMA
%BEGIN THEOREM GENERAL DIRECT SUMS
\begin{theorem}[\bfseries Determinants of general direct sums]
\label{thm:gendirsums}
\index{determinant!general direct sum!theorem}
Let $A$ be an $n\times n$ matrix and let $X, Y\in \bP_k(n)$ be  fixed subsets of $\{1, \ldots ,n\}$ 
of size $k$, $1\leq k<n$.
Let $\sum X = \sum_{x\in X} x$ and let $X'$ be the complement of $X$ in $\{1, \ldots ,n\}$ 
(similarly for $Y, Y'$).  
%Then, summing over just 
%$S_X^Y = \{\sigma\Mid \sigma \in S\; {\rm and}\; \sigma(X)=Y\}$ (\ref{def:SXY}), we have
%ROW FORM
If  $A=B \oplus_X^Y C$ then $B=A[X\Mid Y]$ and  $C=A(X\Mid Y)$ and
\begin{equation}
\det(A) =  (-1)^{\sum X} (-1)^{\sum Y} \det A[X\Mid Y] \det A(X\Mid Y).
\end{equation}
%COLUMN FORM
%PROOF
\begin{proof}
For any $n\times n$ matrix $A$, we show that the {\em restricted determinant sum} 
\begin{equation}
\label{eq:restrictedsum}
\Delta(X,Y, A)=\sum_{\sigma \in S_X^Y} {\rm sgn}(\sigma)\prod_{i=1}^n  A(i,\sigma(i))
\end{equation}
(sum over $S_X^Y$ only) satisfies
\begin{equation}
\label{eq:deltaequalsresdetsum}
\Delta(X,Y, A)=
(-1)^{\sum X} (-1)^{\sum Y} \det(A[X\Mid Y]) \det(A(X\Mid Y)).
\end{equation}
If $A=B \oplus_X^Y C$ then the restricted determinant sum gives $\det(A)$ and thus proves the result.
Use the characterization of $S_X^Y$ given in Lemma~\ref{lem:descriptsxy}:
\[
S_X^Y = \{\gamma^{\nu\mu}\Mid \nu\in {\rm PER}(k),\; \mu\in {\rm PER}(n-k)\}
\]
where
\begin{equation}
\gamma^{\nu\mu} =  \left(
\begin{array}{cccccc}
x_1&\ldots&x_k &x'_1&\ldots&x'_{n-k}\\
y_{\nu(1)}&\ldots&y_{\nu(k)} &y'_{\mu(1)}&\ldots&y'_{\mu(n-k)} 
\end{array}
\right).
\end{equation}
The restricted determinant sum (~\ref{eq:restrictedsum}) becomes
\begin{equation}
\label{eq:restrictedsumbecomes}
 \Delta(X,Y, A)=\sum_{\nu} \sum_{\mu}
{\rm sgn}(\gamma^{\nu\mu}) \prod_{i=1}^k A(x_i,y_{\nu(i)})\prod_{i=1}^{n-k} A(x'_i,y'_{\mu(i)})
\end{equation}
where  $\nu\in {\rm PER}(k)$ and $\mu\in {\rm PER}(n-k)$.
From  \ref{lem:signcanon} we obtain
\begin{equation}
\label{gammanumu2}
{\rm sgn}(\gamma^{\nu\mu}) = (-1)^{\sum X}(-1)^{\sum Y}{\rm sgn}(\nu){\rm sgn}(\mu).
\end{equation}
Thus, \ref{eq:restrictedsumbecomes} becomes $\Delta(X,Y, A)=$
\[
\sum_{\nu} \sum_{\mu}
(-1)^{\sum X}(-1)^{\sum Y}
{\rm sgn}(\nu){\rm sgn}(\mu) \prod_{i=1}^k A(x_i,y_{\nu(i)})\prod_{i=1}^{n-k} A(x'_i,y'_{\mu(i)})=
\]
\[
(-1)^{\sum X}(-1)^{\sum Y}\left(\sum_\nu {\rm sgn}(\nu)\prod_{i=1}^k A(x_i,y_{\nu(i)})\right)
\left(\sum_{\mu}{\rm sgn}(\mu)\prod_{i=1}^{n-k} A(x'_i,y'_{\mu(i)})\right) =
\]
\[
(-1)^{\sum X}(-1)^{\sum Y}\det (A[X\Mid Y]) \det(A(X\Mid Y)).
\]
If $A=B \oplus_X^Y C$ then $\Delta(X,Y, A)= \det(A)$ which completes the proof.
\end{proof}
\end{theorem}

\begin{corollary}[\bfseries Restricted determinant sums]
\label{cor:resdetsums}
\index{determinant!restricted determinant sum}
We use the terminology of Theorem~\ref{thm:gendirsums}.
Let $A$ be an $n\times n$ matrix.   Let
\[
\Delta(X,Y, A)=\sum_{\sigma \in S_X^Y} {\rm sgn}(\sigma)\prod_{i=1}^n  A(i,\sigma(i))
\]
denote the restricted determinant sum (\ref{eq:restrictedsum}).
Then
\[
\Delta(X,Y, A)= (-1)^{\sum X}(-1)^{\sum Y}\det (A[X\Mid Y]) \det(A(X\Mid Y)).
\]
\begin{proof}
The proof is developed in the process of proving Theorem~\ref{thm:gendirsums}.
\end{proof}
\end{corollary}

%SECTION LAPLACE
\section*{Laplace expansion theorem}
\index{Laplace expansion!general case}
We derive the general Laplace expansion theorem.  Our proof is valid for matrices with entries in a commutative ring (e.g., the Euclidean domains, $\bK$, \ref{rem:specialringsfields}) and is based on Corollary~\ref{cor:resdetsums}. We use Definitions~\ref{def:setpartition}, \ref{def:setsfunc} and \ref{def:submatrices2}.
We also use the  notation discussed in  Remark~\ref{rem:notation}.\\  
%BEGIN DEF
\begin{defn}[\bfseries Laplace partition and canonical SDR]
\label{def:laplacepartition}
\index{Laplace expansion!partition!canonical SDR}
Let $S={\rm PER(n)}$ be the permutatons of $\underline{n}$, and let $X, Y\in \bP_k(n)$ be subsets of ${\underline n}$ of size $k$.
Let 
\begin{equation}
\label{eq:laplaceblock}
S_X^Y = \{\sigma\Mid \sigma \in S\; {\rm and}\; \sigma(X)=Y\}.
\end{equation}
For a fixed $X$, the collection of sets
\begin{equation}
\label{eq:laplacepartition}
\bL^X=\{S_X^Y\Mid Y\in \bP_k(n)\}
\end{equation}
is the {\em Laplace partition} of $S={\rm PER(n)}$ associated with $X$.
Let $X'$ and $Y'$ denote the complements of $X$ and $Y$ in ${\underline n}.$
For $\sigma \in S_X^Y$, let $\sigma_X$ and $\sigma_{X'}$ be the restrictions 
(\ref{eq:restriction}) of  $\sigma$ to $X$ and $X'$.  For $X$ fixed, the set
%CANONICAL DEF
\begin{equation}
\label{eq:laplacesdr}
\index{Laplace expansion!partition!canonical SDR}
D^X = \{\gamma \Mid \gamma_X\in {\rm SNC}(X,Y)\,,\; 
\gamma_{X'}\in {\rm SNC}(X',Y')\,, Y\in \bP_k(n)\} 
\end{equation}
is the {\em canonical system of distinct representatives (SDR)} for the Laplace partition of $S$ associated with $X$.
\end{defn}
%END DEF
Note that $\gamma_X$ and $\gamma_{X'}$ in \ref{eq:laplacesdr}  are unique since
$|X|=|Y|$ and $|X'|=|Y'|$.  
The number of blocks in the  partition $\varmathbb{L}^X$ is 
$|\varmathbb{L}^X|=\left( \begin{array}{c}n\\k\end{array} \right).$\\

%THEOREM LAPLACE EXPANSION
\begin{theorem}[\bfseries Laplace expansion theorem]
\label{thm:laplaceexpansion}
\index{Laplace expansion!theorem}
Let $A$ be an $n\times n$ matrix and let $X\in \bP_k(n)$ be a fixed subset of $\{1, \ldots ,n\}$ 
of size $k$, $1\leq k<n$.
Let $\sum X = \sum_{x\in X} x$.  Then the fixed-rows form of the Laplace expansion is
%ROW FORM
\begin{equation}
\label{eq:rowlaplace}
\det(A) = (-1)^{\sum X} \sum_{Y\in \bP_k(n)} (-1)^{\sum Y} \det(A[X\Mid Y]) \det(A(X\Mid Y))
\end{equation}
%COLUMN FORM
and the fixed-columns form is
\begin{equation}
\label{eq:collaplace}
\det(A) = (-1)^{\sum X} \sum_{Y\in \bP_k(n)} (-1)^{\sum Y} \det(A[Y\Mid X]) \det(A(Y\Mid X)).
\end{equation}
%PROOF
\begin{proof}
By definition, 
\begin{equation}
\det(A) = \sum_{\sigma\in {\rm PER}(n)} \prod_{i=1}^n {\rm sgn}(\sigma) A(i,\sigma(i)).
\end{equation}
Let $\varmathbb{L}^X=\{S_X^Y\Mid Y\in \bP_k(n)\}$ be the Laplace partition of ${\rm PER}(n)$ (\ref{def:laplacepartition}).
Then
\begin{equation}
\label{eq:blockdeterminant}
\det(A) = \sum_{Y\in \bP_k(n)} \sum_{\sigma\in S_X^Y} {\rm sgn}(\sigma) \prod_{i=1}^n A(i,\sigma(i)).
\end{equation}
By definition~\ref{eq:restrictedsum}, the inner sum of~\ref{eq:blockdeterminant} is the
restricted determinant sum $\Delta(X,Y,A).$  Thus,
\begin{equation}
\label{eq:sumdeltaxya}
\det(A) = \sum_{Y\in \bP_k(n)} \Delta(X,Y,A)
\end{equation}
By Corollary~\ref{cor:resdetsums}
\[
\Delta(X,Y, A)= (-1)^{\sum X}(-1)^{\sum Y}\det (A[X\Mid Y]) \det(A(X\Mid Y)).
\]
%and thus we have
%\begin{equation}
%\det(A) = (-1)^{\sum X} \sum_{Y\in \bP_k(n)} (-1)^{\sum Y} \det(A[X\Mid Y]) \det(A[X'\Mid Y']).
%\end{equation}
The column form~\ref{eq:collaplace} follows by replacing $A$ by its transpose $A^T$  (\ref{def:transpose}).
This completes the proof.
\end{proof}
\end{theorem}
%END LAPLACE THEOREM

Let $X'$ be the complement of $X$ in $\{1, \ldots ,n\}$.  
Using the notation (\ref{eq:setsinout}) $A[X'\Mid Y'] = A(X\Mid Y)$ we can write 
the fixed-row Laplace expansion (\ref{thm:laplaceexpansion}) as
\begin{equation}
\label{eq:laplacealtnotation}
\det(A) = (-1)^{\sum X} \sum_{Y\in \bP_k(n)} (-1)^{\sum Y} \det(A[X\Mid Y]) A[X'\Mid Y']
\end{equation}
and the fixed-column form as 
\begin{equation}
\label{eq:laplacealtnotationcol}
\det(A) = (-1)^{\sum X} \sum_{Y\in \bP_k(n)} (-1)^{\sum Y} \det(A[Y\Mid X]) A[Y'\Mid X'].
\end{equation}

%Bijection between sets and strictly increasing functions
The set of subsets of {\underline n} of size k, $\bP_k(n)$, corresponds bijectively to the set, ${\rm SNC}(k,n)$,
of strictly increasing functions from {\underline k} to {\underline n}. 
The natural bijection is for $Y\in \bP_k(n)$, $Y = \{y_1, \ldots, y_k\}$ with $y_1< \cdots < y_k$, to correspond 
to 
$
\left(
\begin{array}{ccc}
1   &\ldots &k\\
y_1&\ldots&y_k
 \end{array}
\right).
$
Thus, we can rewrite the Laplace expansion theorem 
in terms of functions.
We use notation like that of \ref{eq:setsinout}.
\index{Laplace expansion!function version}
Let $g \in {\rm SNC}(k,n)$ be fixed, let $A=(a_{ij})$ be an $n\times n$ matrix and  let 
$\sum g = \sum_{i=1}^k g(i)$.
Then the fixed-row Laplace expansion expressed in terms of functions is 
%Laplace expansion in terms of functions
\begin{equation}
\label{eq:laplacealtnotation2row}
\det(A) = (-1)^{\sum g} \sum_{f \in {\rm SNC}(k,n)} (-1)^{\sum f} \det(A[g\Mid f]) \det(A(g\Mid f))
\end{equation}
and the fixed-column Laplace expansion expressed in terms of functions is 
%Laplace expansion in terms of functions
\begin{equation}
\label{eq:laplacealtnotation2col}
\det(A) = (-1)^{\sum g} \sum_{f \in {\rm SNC}(k,n)} (-1)^{\sum f} \det(A[f\Mid g] \det(A(f\Mid g)).
\end{equation}

%REMARK
\begin{remark}[\bfseries Example of Laplace expansion]
\index{Laplace expansion!examples}
Let $A$ be a $4\times 4$ matrix of integers ($A\in {\bf M}_{4,4}(\bZ))$
\[
A=
\left(
\begin{array}{cccc}
1&  2 & 3 & 4\\
2&  3 & 4 & 1 \\
3&  4&  1 & 2\\
4 & 1 & 2 & 3
\end{array}
\right).
\]
Use the fixed-row Laplace expansion ~(\ref{eq:laplacealtnotation2row}) 
with the fixed function $g=(1,3)$ (in one line notation).
Take the variable functions, $f$, lexicographically in one-line notation:
$(1,2), (1,3), (1,4), (2,3), (2,4), (3,4).$  Then $\det(A) \equiv |A| =$
\[
-\left|\begin{array}{cc}1&2\\3&4 \end{array}\right| \left|\begin{array}{cc}4&1\\2&3 \end{array}\right|
+\left|\begin{array}{cc}1&3\\3&1 \end{array}\right| \left|\begin{array}{cc}3&1\\1&3 \end{array}\right|
-\left|\begin{array}{cc}1&4\\3&2 \end{array}\right| \left|\begin{array}{cc}3&4\\1&2 \end{array}\right|
\]
\[
-\left|\begin{array}{cc}2&3\\4&1 \end{array}\right| \left|\begin{array}{cc}2&1\\4&3 \end{array}\right|
+\left|\begin{array}{cc}2&4\\4&2 \end{array}\right| \left|\begin{array}{cc}2&4\\4&2 \end{array}\right|
-\left|\begin{array}{cc}3&4\\1&2 \end{array}\right| \left|\begin{array}{cc}2&3\\4&1 \end{array}\right|
= 160.
\]
\end{remark}
%END REMARK
%SIMPLE LAPLACE

The next corollary is the version of the Laplace expansion theorem that is most often stated and proved in elementary courses:
\begin{corollary}[\bfseries Simple Laplace expansion]
\label{cor:laplacesimple}
\index{Laplace expansion!simple standard}
Let $A$ be an $n\times n$ matrix, $n>1$, and let $1\leq i \leq n$.
Then
\begin{equation}
\label{eq:laplacesimple}
\det(A) = (-1)^i \sum_{j=1}^n (-1)^j a_{i j} \det(A(i\Mid j)) 
\end{equation}
\begin{equation}
\label{eq:laplacesimplecol}
\det(A) = (-1)^i \sum_{j=1}^n (-1)^j a_{j i} \det(A(j\Mid i)).
\end{equation}

\begin{proof}  We use~\ref{eq:laplacealtnotation}:
\[ \det(A) = (-1)^{\sum X} \sum_{Y\in \bP_k(n)} (-1)^{\sum Y} \det(A[X\Mid Y]) \det(A(X\Mid Y)). \]
In this case, 
\[ \sum X = i,\;\;\sum Y = j,\;\;
A[X\Mid Y]= A[i\Mid j] = a_{i j},\;\;
{\rm and\;\;}
 \det(A[X\Mid Y]) = a_{i j}.
 \] 
 Substituting these values and noting that $A(X\Mid Y) = A(i\Mid j)$ proves~\ref{eq:laplacesimple}.
 Identity~\ref{eq:laplacesimplecol} follows by using~\ref{eq:laplacealtnotationcol}
 instead of~\ref{eq:laplacealtnotation}.
\end{proof}
\end{corollary}
%END SIMPLE LAPLACE
Note that for
 $A=(a_{ij})$ then $A(i,j)$ or $A[i\Mid j]$ can be used
in place of $a_{ij}$ in~\ref{eq:laplacesimple}.
%SIMPLE LAPLACE EXTENDED
Recall the ``delta'' notation: $\delta({\bf statement}) = 0$ if ${\bf statement}$ is false,  $1$ if {\bf statement} is true.

\begin{corollary}[\bfseries Simple Laplace extended]
\label{cor:laplaceinverse}
\index{Laplace expansion!simple extended}
Let $A$ be an $n\times n$ matrix, $n>1$, and let $1\leq i \leq n$. Then
\begin{equation}
\label{eq:laplaceinverse}
\delta(k = i)\det(A) = (-1)^k \sum_{j=1}^n (-1)^j a_{i j} \det(A(k\Mid j)).
\end{equation}
or, alternatively, 
\begin{equation} 
\label{eq:cofactorinverse0}
AB_A(i,k) = \sum_{j=1}^n A(i,j)B_A(j,k) = \delta(i=k)\det(A)
\end{equation}
where $B_A=(b_{jk})$ and $b_{jk} = (-1)^{(k+j)} \det(A(k\Mid j)).$
\begin{proof}
The case where $k=i$ becomes
\begin{equation}
\label{eq:kequali}
\det(A) = (-1)^i \sum_{j=1}^n (-1)^j a_{i j} \det(A(i\Mid j))
\end{equation}
which is a simple Laplace expansion (\ref{eq:laplacesimple}).

Consider $k\neq i$.   Take the matrix $A$ with rows
\[
A=(A_{(1)}, \ldots,  A_{(k-1)}, A_{(k)}, A_{(k+1)}, \ldots, A_{(n)})
\]
and replace row $k$, $A_{(k)}$, with row $i$, $A_{(i)}$, to obtain a matrix $A'$:
\[
A'=(A_{(1)}, \ldots,  A_{(k-1)}, A_{(i)}, A_{(k+1)}, \ldots, A_{(n)}).
\] 
  
The matrix $A'$ thus has two identical rows and hence $\det(A')=0$.

Apply equation~\ref{eq:laplacesimple} to $A'$ to get
\begin{equation}
\label{eq:fakerow}
 \det(A') = (-1)^k \sum_{j=1}^n (-1)^j A'(k,j) \det(A'(k\Mid j)).
\end{equation}
By definition of $A'$, $A'(k,j)=A(i,j)=a_{ij}$.  Also by definition of $A'$, the  $(n-1)\times (n-1)$ matrix 
$A'(k\Mid j)=A(k\Mid j).$
Thus by \ref{eq:fakerow} we have
\begin{equation}
\label{eq:knotequali}
0=\det(A') = (-1)^k \sum_{j=1}^n (-1)^j a_{i j} \det(A(k\Mid j))\;\;\;{\rm where}\;\;\;k\neq i.
\end{equation}
Combining \ref{eq:kequali} and \ref{eq:knotequali} gives \ref{eq:laplaceinverse}.

 Rewrite equation~\ref{eq:laplaceinverse} as follows:
\begin{equation}
\label{eq:laplacecofactor}
\sum_{j=1}^n a_{i j} [(-1)^{k+j} \det(A(k\Mid j))]= \delta(i = k)\det(A).
\end{equation}
If we define a matrix $B_A=(b_{jk})$ by $b_{jk} = (-1)^{(k+j)} \det(A(k\Mid j))$ then we obtain
\begin{equation} 
\label{eq:cofactorinverse}
AB_A(i,k) = \sum_{j=1}^n A(i,j)B_A(j,k) = \delta(i=k)\det(A).
\end{equation}
This completes the proof.
\end{proof}
\end{corollary}

\begin{definition}[\bfseries Signed cofactor matrix]
\label{def:signcofactor}
Let $A$ be an $n\times n$ matrix.  For $1\leq i, j \leq n$,  
define $c_{ij} = (-1)^{i+j} \det A(i\Mid j).$ We call $c_{ij}$ the 
{\em signed cofactor} of $A(i,j)$.  
The $n\times n$ matrix $C_A=(c_{ij})$ is the {\em signed cofactor matrix} of $A.$
The transpose, $B_A=C_A^T$, of $C_A$ is sometimes called the {\em adjugate} of $A$ and
written $\mathrm{adj}(A)$.
\end{definition}

The next corollary restates~\ref{eq:cofactorinverse}. We use $\det(XY)=\det(X)\det(Y)$ (\ref{cor:detproduct}).

\begin{corollary}
\label{cor:transposesignedcofactor}
\index{Laplace expansion!near inverse form}
Let $C_A=(c_{ij})$ be the signed cofactor matrix of $A$.
Let $B_A = C^T$ be the transpose of  $\,C_A$ (i.e., $B_A=\mathrm{adj}(A)$) .  Then
\begin{equation}
\label{eq:nearinverse}
AB_A = B_AA = (\det A) I_n
\end{equation}
where $I_n$ is the $n \times n$ identity matrix. If $\det(A)$ is a unit in $\bK$ then
$A$ is a unit in ${\bf M}_{n,n}(\bK)$ (\ref{rem:specialringsfields})
and $A^{-1}=B_A(\det(A))^{-1}$. 
Thus, $A$ is unit in ${\bf M}_{n,n}(\bK)$ if and only if $\det(A)$ is a unit in $\bK$.
\end{corollary}
\begin{proof}
If $\det(A)\neq 0$ then $A^{-1}=B_A(\det(A))^{-1}$ follows from~\ref{eq:cofactorinverse}.  
The converse follows from the fact that if $AA^{-1} = I_n$ then $\det(A)\det(A^{-1}) = 1$ 
so $(\det(A))^{-1} = \det(A^{-1}).$
If $\det(A)\neq 0$ then $AB_A=B_AA$ follows from the fact that $A^{-1}A=AA^{-1}$.
The statement $AB_A=(\det A) I_n$ is exactly the same as 
equation~\ref{eq:cofactorinverse} and does not require $A$ to be nonsingular.
In general, commutivity, $AB_A=B_AA$, follows from \ref{cor:laplaceinverse} by replacing $A$ by $A^T$ in equation~\ref{eq:cofactorinverse}:  
\begin{equation}
\label{eq:laplaceinversetra1}
\delta(k = i)\det(A^T) = \sum_{j=1}^n (-1)^{k+j} a^T_{i j} \det(A^T(k\Mid j)).
\end{equation} 
$\det(A^T(k\Mid j))=\det((A(j\Mid k))^T) = \det(A(j\Mid k))$ (\ref{rem:transposeprod})
and \[(-1)^{k+j}\det(A(j\Mid k))= B_A(k,j)\]  thus using \ref{eq:laplaceinversetra1} 
\begin{equation}
\label{eq:laplaceinversetra2}
\delta(k = i)\det(A) = \sum_{j=1}^n B_A(k,j) a_{j i} = B_AA.
\end{equation}
Thus, $AB_A=B_AA$ in all cases.\\
\end{proof}

\begin{remark}[\bfseries Example of signed cofactor and adjugate matrices]
\label{rem:sgncofactexample}
Let $A$ be a $4\times 4$ matrix of integers ($A\in {\bf M}_{4,4}(\bZ))$
with $A(i,j) = \delta(i\leq j)$.  The matrices $A$, $C_A$, and $B_A=\mathrm{adj}(A)$ of
Corollary~\ref{cor:transposesignedcofactor} are as follows:
\[
A=
\left(
\begin{array}{cccc}
1&  1 & 1 & 1\\
0&  1 & 1 & 1 \\
0&  0&  1 & 1\\
0 & 0 & 0 & 1
\end{array}
\right)\;\;
C_A=
\left(
\begin{array}{rrrr}
+1&  0 & 0 & 0\\
-1&  +1 & 0 & 0 \\
0&  -1&  +1 & 0\\
0 & 0 & -1 & +1
\end{array}
\right)\;\;
B_A=
\left(
\begin{array}{rrrr}
+1&  -1 & 0 & 0\\
0&  +1 & -1 & 0 \\
0&  0&  +1 & -1\\
0 & 0 & 0 & +1
\end{array}
\right).
\]
\end{remark}
%END REMARK
If $A$ is an $n\times n$ matrix over a field $\bF$ and $A^{-1}$ is the inverse of $A$,
then $\det(A\,A^{-1})= \det(A)\det(A^{-1}) = \det(I_n)=1.$
If $X$ is $n\times 1$, we write $X(t) = X(t,1)$, $t=1, \ldots, n$.
Putting together Definition~\ref{def:signcofactor} and 
Corollary~\ref{cor:transposesignedcofactor} we get the following corollary.
%COROLLARY CRAMER
\begin{corollary}[\bfseries Cramer's rule]
\label{cor:cramerrule}
\index{Laplace expansion!Cramer's rule}
Let $AX=Y$ where $A$ is $n\times n$, $X$ is $n\times 1$ and $Y$ is $n\times 1$
(entries in a field $\bF$).
Designate $A=(A^{(1)}, \ldots , A^{(i)}, \ldots, A^{(n)})$ as a sequence of columns.
Define ${\hat A}=(A^{(1)}, \ldots, A^{(i-1)}, Y, A^{(i+1)}, \ldots, A^{(n)})$ to be the matrix
$A$ with column $A^{(i)}$ replaced by $Y$. Then
\[X(i) = \det({\hat A})/\det(A).\]
%PROOF
\begin{proof}
$AX=Y$ implies $X=A^{-1}Y$ and $X(i)=\sum_{j=1}^n A^{-1}(i,j)Y(j).$
From
Definition~\ref{def:signcofactor} and Corollary~\ref{cor:transposesignedcofactor} 
we get
%EQUATION  
\begin{equation}
\label{eq:inversecofactorspecific}
A^{-1}(i,j) = \left[ (-1)^{i+j}\frac{\det(A(j\Mid i))}{\det(A)}\right].
\end{equation}
%END EQUATION
%EQUATION
Thus
\begin{equation}
\label{eq:usedeflemma1}
X(i)=\sum_{j=1}^n A^{-1}(i,j)Y(j) = \sum_{j=1}^n \left[ (-1)^{i+j}\frac{\det(A(j\Mid i))}{\det(A)}\right]Y(j) 
\end{equation}
%END EQUATION
and
%EQUATION
\begin{equation}
\label{eq:usedeflemma2}
X(i)=\frac{1}{\det(A)}\sum_{j=1}^n (-1)^{i+j}\det(A(j\Mid i))Y(j). 
\end{equation}
%END EQUATION
Observe that the matrix $A(j\Mid i) = {\hat A}(j\Mid i)$ since only columns $i$ differ
between $A$ and ${\hat A}$.  Note also that $Y(j) = {\hat A}(j,i)$.
Thus, the sum of \ref{eq:usedeflemma2} becomes
\[
\sum_{j=1}^n (-1)^{i+j}\det(A(j\Mid i))Y(j) = 
(-1)^i \sum_{j=1}^n (-1)^{j}{\hat A}(j,i)\det({\hat A}(j\Mid i)) = \det({\hat A})
\]
using  \ref{cor:laplacesimple} (column form, \ref{eq:laplacesimplecol}).
This completes the proof.
%EQUATION
\end{proof}
\end{corollary}
%END COROLLARY
\begin{remark}[\bfseries Example of Cramer's rule]
\label{rem:cramer}
The two equations
\[
\begin{array}{ccccc}
x_1  & +  &  x_2 & = & 3 \\
x_1  & -  &  x_2 & =  & 1
\end{array}  
\]
can be expressed by the equation $AX=Y$ where 
\[
A=\left(
\begin{array}{rr}
1&1\\
1&-1
\end{array}
\right)
\;\;\;\;
X=\left(
\begin{array}{c}
x_1\\
x_2
\end{array}
\right)
\;\;\;\;
Y=\left(
\begin{array}{c}
3\\
1
\end{array}
\right).
\]
Applying ~\ref{cor:cramerrule} twice, to $X(1)=x_1$ and $X(2) = x_2$, and noting
that $\det(A)=-2$ gives
\[
x_1= \frac{\det\left(\begin{array}{rr}3&1\\1&-1\end{array}\right)}{-2} = 2
\;\;\;\;{\rm and}\;\;\;\;
x_2= \frac{\det\left(\begin{array}{rr}1&3\\1& 1\end{array}\right)}{-2} = 1.
\]
\end{remark}

%SECTION 3
\section*{Cauchy-Binet theorem}
\index{Cauchy-Binet!discussion}
We need some notational conventions for describing product - sum interchanges.   
Consider the following example:
\[
(x_{11}+x_{12})(x_{21}+x_{22})=x_{11}x_{21}+x_{11}x_{22}+x_{12}x_{21}+x_{12}x_{22}.
\]
Look at the second integers in each pair of subscripts:
\[
x_{1{\underline 1}}x_{2{\underline 1}}+x_{1{\underline 1}}x_{2{\underline 2 }}+
x_{1{\underline 2}}x_{2{\underline 1}}+x_{1{\underline 2}}x_{2{\underline 2 }}.
\]
The pairs of underlined integers are, in order:
\[
11, 12, 21, 22.
\]
These pairs represent  (in one line notation) all of the functions in ${\underline 2}^{\underline 2}$.
Thus, we can write
\begin{equation}
\prod_{i=1}^2\left(\sum_{k=1}^2x_{ik}\right) = 
\sum_{f\in {\underline 2}^{\underline 2}} \prod_{i=1}^2 x_{i\,f(i)}.
\end{equation}
The general form of this identity is
\begin{equation}
\label{eq:changesumprod}
\prod_{i=1}^n\left(\sum_{k=1}^p x_{ik}\right) = 
\sum_{f\in {\underline p}^{\underline n}} \prod_{i=1}^n x_{i\,f(i)}.
\end{equation}
This product-sum-interchange identity \ref{eq:changesumprod} is important to what follows.

We now prove the Cauchy-Binet theorem using \ref{eq:sgndet}, \ref{def:submatrices2},
\ref{eq:changesumprod} and \ref{eq:injsncper}.  
The proof  is valid for matrices with entries in a commutative ring (e.g., a Euclidean domain).
\begin{thm}[\bf Cauchy-Binet] Let $A$ be an $n\times p$ and $B$ a $p\times n$ matrix.
\label{thm:cauchybinet}
\index{Cauchy-Binet!theorem}
Then 
\[
\det(AB)= \sum_{f\in {\rm SNC}(n,p)} \det(A^f)\det(B_f)
\]
where ${\rm SNC}(n,p)$ denotes the strictly increasing functions from ${\underline n}$ to 
${\underline p}$ (see~\ref{def:setsfunc}); $A^f$ denotes the submatrix of $A$ 
with columns selected by $f$;
and $B_f$ denotes the submatrix of $B$ with rows selected by $f$ (see \ref{def:submatrices2}).
\end{thm}
\begin{proof}
%id 1
\[
\det(AB)= \sum_{\gamma \in {\rm PER}(n)} {\rm sgn}(\gamma)\prod_{i=1}^n(AB)(i,\gamma(i))\;\;
({\bf definition\;\det})
\]
%id 2
\[
\det(AB)= \sum_{\gamma\in {\rm PER}(n)} {\rm sgn}(\gamma)\prod_{i=1}^n
\sum_{k=1}^p A(i,k)B(k,\gamma(i))
\;\;({\bf definition\;}AB)
\]
%id 3
\[
\det(AB)= \sum_{\gamma\in {\rm PER}(n)} {\rm sgn}(\gamma)\sum_{h\in {\underline p}^{\underline n}}
\prod_{i=1}^n A(i,h(i))B(h(i), \gamma(i))
\;\;({\bf identity\;}\ref{eq:changesumprod})
\]
%id 4
\[
\det(AB)= \sum_{h\in {\underline p}^{\underline n}}
\prod_{i=1}^n A(i,h(i))\sum_{\gamma\in {\rm PER}(n)} {\rm sgn}(\gamma)
\prod_{i=1}^n B(h(i), \gamma(i))
\;\;({\bf algebra\;rules})
\]
%id 5
\begin{equation}
\label{ngreaterp}
\det(AB)= \sum_{h\in {\underline p}^{\underline n}}
\prod_{i=1}^n A(i,h(i)) \det(B_h)
\;\;({\bf definition}\;B_h\; \ref{def:submatrices2})
\end{equation}
Now observe that $\det(B_h)=0$ if $h\not\in {\rm INJ}(n,p)$ by \ref{eq:detisalt}.  Thus,
%id 6
\[
\det(AB)= \sum_{h\in {\rm INJ}(n,p)}
\prod_{i=1}^n A(i,h(i)) \det(B_h)
\;\;({\bf determinant\;property\;} \ref{eq:altmult}).
\]
Thus,
%id 7
\[
\det(AB)=\sum_{f\in {\rm SNC}(n,p)}
\sum_{\gamma\in {\rm PER}(n)}\left(\prod_{i=1}^n A(i,f\gamma(i))\right) \det(B_{f\gamma})
\;\;({\bf set\;identity\;} \ref{eq:injsncper})
\]
Note that  $A(i,f\gamma(i))= A^f(i,\gamma(i))$.
The matrix $B_{f\gamma}$ of the previous equation is the same as $(B_f)_\gamma$.
Identity~\ref{eq:sgndet} implies that 
$\det((B_f)_\gamma)={\rm sgn}(\gamma) \det (B_f)$.  Thus we obtain
\[
\det(AB)=\sum_{f\in {\rm SNC}(n,p)}
\left(\sum_{\gamma\in {\rm PER}(n)}{\rm sgn}(\gamma)\prod_{i=1}^n A^f(i,\gamma(i))\right) \det(B_f).
\;\;({\bf algebra\;rules})
\]
Finally, applying the definition of the determinant to $A^f$ we get
\[
\det(AB)= \sum_{f\in {\rm SNC}(n,p)} \det(A^f)\det(B_f).
\]
\end{proof}
\begin{rem}[\bfseries Zero determinant of product]
\label{rem:whenzero}
Equation~\ref{ngreaterp} of the preceding theorem is a sum of the form
\[
\det(AB)= \sum_{h\in {\underline p}^{\underline n}}
C(h) \det(B_h)
\]
where $C(h)$ depends on $A$ and $h$.
If $n>p$ then every $h\in {\underline p}^{\underline n}$ has $h(s)=h(t)$ for some pair of
values $s<t$ (i.e., the set ${\rm INJ}(n,p)$ is empty).  For every such $h$,  $\det(B_h)=0$ and 
hence $\det(AB)=0$ if $n>p$.
\end{rem}

\begin{cor}[\bf Determinant of product]
\label{cor:detproduct}
\index{Cauchy-Binet!determinant of product}
If $A$  and $B$ are $n\times n$ matrices, then $\det(AB)=\det(A)\det(B)$.
\begin{proof}
Apply the Cauchy-Binet theorem (\ref{thm:cauchybinet}) with $p=n$.  In that case,
$\det(AB)= \sum_{f\in {\rm SNC}(n,n)} \det(A^f)\det(B_f)$.  ${\rm SNC}(n,n)$ has only one element, the identity function $f(i)=i$ for all $i\in {\underline n}$.  Thus, $A^f=A$ and
$B_f= B$.  
\end{proof}
\end{cor}

Corollary~\ref{cor:gencauchybinet}  is a useful restatement of Theorem~\ref{thm:cauchybinet}.\\
\begin{cor}[\bfseries General Cauchy-Binet]
\label{cor:gencauchybinet}
\index{Cauchy-Binet!general form}
Let $A$ be an $a\times p$ matrix and $B$ a $p\times b$ matrix.  
Let  $g\in {\underline a}^{\underline n}$ and $h\in {\underline b}^{\underline n}$.
Let $C=AB$.
Then
\begin{equation}
\label{eq:gencauchybinet}
\det(C[g\Mid h]) =  \sum_{f\in {\rm SNC}(n,p)} \det(A[g\Mid f])\det(B[f\Mid h])
\end{equation}
\end{cor}
\begin{proof}
We apply Theorem~\ref{thm:cauchybinet} to  $A_g$ ($n\times p$) 
and $B^h$ ($p\times n$):
\[
\det(A_gB^h)= \sum_{f\in {\rm SNC}(n,p)} \det((A_g)^f)\det((B^h)_f).
\]
This becomes
\[
\det((AB)[g\Mid h])= \sum_{f\in {\rm SNC}(n,p)} \det(A[g\Mid f])\det(B[f\Mid h]).
\]
Substituting $C=AB$ gives \ref{eq:gencauchybinet}.\\
\end{proof}
\begin{remark}[\bfseries General Cauchy-Binet]
\label{rem:gencauchybinet}
If $n>p$ then ${\rm SNC}(n,p) = \emptyset$ (empty set).  Thus, 
$n\leq p$ is the more interesting case of \ref{eq:gencauchybinet}. 
Likewise, if $n>a$ or $n>b$ then the right hand side of \ref{eq:gencauchybinet} is
zero.  
Thus,  we are most interested in the case $n\leq \min(a,b,p)$.  Even if $n\leq \min(a,b,p)$ then we still need $g$ and $h$ to be injections to make the corresponding determinants nonzero.  
These observations lead to the following version:
\end{remark}
\begin{corollary}[\bfseries Extended Cauchy-Binet]
\label{cor:extcauchybinet}
Let $A$ be an $a\times p$ matrix and $B$ a $p\times b$ matrix.  
Let  $g\in {\underline a}^{\underline n}$ and $h\in {\underline b}^{\underline n}$.
Assume $n\leq \min(a,b,p)$ and $g$ and $h$ are injective (\ref{def:setsfunc}).
Let $C=AB$.
Then 
\begin{equation}
\label{eq:extcauchybinet}
\det(C[g\Mid h]) =  \sum_{f\in {\rm SNC}(n,p)} \det(A[g\Mid f])\det(B[f\Mid h])
\end{equation}
\begin{equation}
\label{eq:extcauchybinetfnc}
\det(C_g^h) =  \sum_{f\in {\rm SNC}(n,p)} \det(A_g^f)\det(B_f^h)\;\;\;\;
\left(C_g^h=(AB)_g^h=A_gB^h\right).
\end{equation}
\begin{proof}
These statements are a special case of \ref{cor:gencauchybinet}.\\
\end{proof}
\end{corollary}
\begin{rem}[\bfseries $g$ and $h$ strictly increasing]
In Corollary \ref{cor:extcauchybinet} the additional assumption that  $g$ and $h$ are  strictly increasing is often made:
$g \in {\rm SNC}(n,a)$ and $h \in {\rm SNC}(n,b)$ (\ref{def:setsfunc}).  This assumption implies the standard "set" version of the Cauchy-Binet theorem.  We identify a function, $g \in {\rm SNC}(n,a)$, with the set $G={\rm image(g)}.$
\end{rem}
%COROLLARY
\begin{cor}[\bfseries Cauchy-Binet set version]
\label{cor:cauchybinetset}
\index{Cauchy-Binet!set version}
Let $A$ be an $a\times p$ matrix and $B$ a $p\times b$ matrix. Let $C=AB$.
Assume $n\leq \min(a,b,p)$.  Let $G\in \bP_n({\underline a})$ and $H\in \bP_n({\underline b})$ be subsets of
size $n$.
Then
\begin{equation}
\label{eq:cauchybinetset}
\det(C[G\Mid H]) =  \sum_{F\in \bP_n({\underline p})} \det(A[G\Mid F])\det(B[F\Mid H]).
\end{equation}
Alternatively, let $g\in {\rm SNC}(n,a)$ and  $h\in {\rm SNC}(n,b)$ be defined by  ${\rm image}(g)=G$ and ${\rm image}(h)=H$.
Then
\begin{equation}
\label{eq:caubinexatwoplus}
\det(A_gB^h)= \sum_{f\in {\rm SNC}(n,p)} \det(A_g^f)\det(B^h_f)
\end{equation}
where $A_g$ is $n\times p$, $B^h$ is $p\times n$ while $A_g^f$ and $B^h_f$ are both $n\times n$.
\begin{proof}
This statement is a special case of \ref{cor:extcauchybinet}.
We use \ref{rem:setsfncsrelate}.
We also use equation~\ref{eq:gencolrowmatrix} to select rows and columns of matrices:
\begin{equation}
\label{eq:gencolrowmatrixrep}
C[g\Mid h] \equiv C_g^h \equiv (AB)_g^h = A_gB^h.\\
\end{equation}
\end{proof}
\end{cor}
\begin{remark}[\bfseries Discussion of theorem~\ref{thm:cauchybinet}]
\index{Cauchy-Binet!examples}
Take $A$ to be a  $2\times 6$ matrix (i.e. $n=2$ and $p=6$), and $B$ to be
a $6\times 2$ matrix as follows:
\begin{equation}
\label{eq:exacaubinthm}
A=
  \begin{blockarray}{ccccccc}
   & \scriptstyle1& \scriptstyle2& \scriptstyle3& \scriptstyle4& \scriptstyle5& \scriptstyle6\\
    \begin{block}{c[cccccc]}
       \scriptstyle1& 2&  5&  4&  -2&  2 & 1 \\ 
       \scriptstyle2& 2&  6&  5& 0& -1  & 0 \\
    \end{block}
  \end{blockarray}
  \;\;\;  \;\;\;
  B=
   \begin{blockarray}{ccc}
   & \scriptstyle1& \scriptstyle2\\
    \begin{block}{c[cc]}
     \scriptstyle1& 2&5\\ 
      \scriptstyle2&1&1\\
      \scriptstyle3&2&3\\
      \scriptstyle4&2&1\\ 
      \scriptstyle5&3&1\\
      \scriptstyle6&2&2\\
    \end{block}
  \end{blockarray}.
\end{equation}
Note that $C=AB$ is a $2\times 2$ matrix and thus $\det(C)$ is defined.
However, $\det(A)$ and $\det(B)$ are not defined.  
We have (by \ref{thm:cauchybinet})
\begin{equation}
\label{eq:caubinexaone}
\det(AB)= \sum_{f\in {\rm SNC}(2,6)} \det(A^f)\det(B_f)
\end{equation}
where ${\rm SNC}(2,6)$ denotes the strictly increasing functions from ${\underline 2}$ to 
${\underline 6}$ (see~\ref{def:setsfunc}); $A^f$ denotes the submatrix of $A$ 
with columns selected by $f$,
and $B_f$ denotes the submatrix of $B$ with rows selected by $f$.
For example, take $f=(1,3)$ in one-line notation.  
Then 
\[
A^f=A^{(1,3)}=
\left[
\begin{array}{cc}
2&4\\
2&5
\end{array}
\right]\;\;\;\;{\rm and};\;\;\;
B_f=B_{(1,3)}=
\left[
\begin{array}{cc}
2&5\\
2&3
\end{array}
\right].
\]
Thus, $\det(A^f)\det(B_f)=2\cdot(-4)= -8$ is one of $15$ terms in
the sum of~\ref{eq:caubinexaone}.
\end{remark}

\begin{remark}[\bfseries Discussion of \ref{cor:cauchybinetset}]
\label{rem:discaubinset}
Take $A$ to be a $4\times 6$ ($a\times p$) and $B$ to be a $6\times 3$
($p\times b$) matrix as follows:
\begin{equation}
A=
  \begin{blockarray}{ccccccc}
   &\scriptstyle1&\scriptstyle2&\scriptstyle3&\scriptstyle4&\scriptstyle5&\scriptstyle6\\
    \begin{block}{c[cccccc]}
      \scriptstyle1& 2&  5&  4&  -2&  2 & 1 \\
      \scriptstyle2& 0&  1&  1&   0& -1 & 0 \\
      \scriptstyle3& 2&  6&  5& 0& -1  & 0 \\
      \scriptstyle4& 2&  1&  1& -1& -1  & 0 \\        
    \end{block}
  \end{blockarray}
  \;\;\; \;\;\; \;\;\;
  B=
   \begin{blockarray}{cccc}
   &\scriptstyle1&\scriptstyle2&\scriptstyle3\\
    \begin{block}{c[ccc]}
     \scriptstyle1& 2&5&5\\ 
     \scriptstyle2&1&1&1\\
     \scriptstyle3&2&3&5 \\
     \scriptstyle4&2&1&4\\ 
     \scriptstyle5&3&1&1\\
     \scriptstyle6&2&2&5 \\
    \end{block}
  \end{blockarray}.
\end{equation}
Note that $C=AB$ is a $4\times 3$ matrix and thus $\det(C)$ is not defined.
From \ref{eq:cauchybinetset}
\begin{equation}
\label{eq:cauchybinetsetexa1}
\det(C[G\Mid H]) =  \sum_{F\in \bP_n({\underline p})} \det(A[G\Mid F])\det(B[F\Mid H]).
\end{equation}
We choose $n\leq \min(a,b,p)=\min(4,3,6)=3$  to be $n=2$.  Choose $G \in \bP_2(\underline{4})$
to be $G=\{1,3\}$ and $H \in \bP_2(\underline{3})$ to be $H=\{1,2\}$.
\begin{equation*}
\label{eq:cauchybinetsetexa2}
\det(C[\{1,3\}\Mid \{1,2\}]) =  
\sum_{F\in \bP_2({\underline 6})} \det(A[\{1,3\}\Mid F])\det(B[F\Mid \{1,2\}]).
\end{equation*}
We can rewrite this equation using \ref{eq:gencolrowmatrixrep}.  
\begin{equation}
\label{eq:caubinexatwo}
\det(A_gB^h)= \sum_{f\in {\rm SNC}(2,6)} \det(A_g^f)\det(B^h_f)
\end{equation}
where $g=(1,3)$ and $h=(1,2)$ in one-line notation.
Note that $A_g$ is the matrix $A$ and $B^h$ is the matrix $B$ of \ref{eq:exacaubinthm}.
From \ref{eq:gencolrowmatrix} we have $A_gB^h= (AB)_g^h=C_g^h$.
Thus,  the matrices of Corollary~\ref{cor:cauchybinetset} are  ``containers'' for 
many instances where Theorem~\ref{thm:cauchybinet} can be applied.\\
\end{remark}
%DEFINITION
\begin{minipage}{\textwidth}
\begin{equation}
\index{matrix!rank}
\index{Cauchy-Binet!rank of matrix}
\label{eq:secranmat}
{\bf Rank\;\;of\;\;a\;\;matrix}
\end{equation}
\begin{definition}[\bfseries Rank of a matrix]
\label{def:ranmat} %labdef Rank of a matrix
Recall the notation for sets of subsets, $\bP_n({\underline a})$ (\ref{rem:notation}).
Let $C\in {\bf M}_{a,b}(\bK)$ be an $a\times b$ matrix.  
The {\em rank} $\rho(C)$ is the size of the largest nonzero sub-determinant of $C$:
\[
\rho(C)=
\max\{n\Mid n\in \bN_0,\;G\in \bP_n({\underline a}),\;H\in \bP_n({\underline b})\;\;
\det(C[G\Mid H])) \neq 0\}.
\]
If $C=\Theta_{a,b}$ is the zero matrix, then $\rho(C)=0$.
\end{definition}
\end{minipage}
\hspace* {1in}\\
\begin{remark}[\bfseries Alternative definitions of rank]
\label{rem:ranmat}
The notion of  ``rank''  for modules was discussed \ref{rem:modulerank}.  
Definition~\ref{def:ranmat} defines rank for a matrix.  You will recall from your linear algebra courses that the rank of a matrix $C\in {\bf M}_{a,b}(\bK)$, $K$ a field, is the same as the dimension of the row space of $C$ which is the same as the dimension of the column space of $C$.  
This dimension is equal to the maximum number of linearly independent rows or columns of $C$.
The following is a technically another corollary of theorem~\ref{thm:cauchybinet}.
\end{remark}
%COROLLARY
\begin{cor}[\bfseries Rank of a product]
\label{cor:ranpro}
\index{Cauchy-Binet!rank of product}
\index{matrix!rank of product}
Let $A\in {\bf M}_{a,p}(\bK)$ and $B\in {\bf M}_{p,b}(\bK)$.
Let $C=AB$.  Then  the rank $\rho(C)$ satisfies 
\begin{equation}
\label{eq:ranlesmin}
\rho(C)=\rho(AB)\leq \min\{\rho(A),\rho(B)\}.
\end{equation}
If $B$ is nonsingular then $\rho(AB)=\rho(A)$ and if $A$ is nonsingular
$\rho(AB)=\rho(B)$.
\begin{proof}
Let $r=\rho(AB)$. 
From the definition of rank, there exists  $G\in \bP_r({\underline a})$ and 
$H\in \bP_r({\underline b})$  such that  $\det(AB[G\Mid H])\neq 0.$
From \ref{eq:cauchybinetset}, 
\begin{equation}
\label{eq:caubinsetnot}
\det(AB[G\Mid H]) =  \sum_{F\in \bP_r({\underline p})} \det(A[G\Mid F])\det(B[F\Mid H])\neq 0.
\end{equation}
To be able to choose the subsets $F$, $G$, and $H$, we have $r\leq \min(a,b,p)$. 
Suppose, without loss of generality, that $\rho(A)<\rho(AB)=r$.
Then we have $\det(A[G\Mid F])=0$ for every term in the sum and hence $\det(AB[G\Mid H])=0$, contrary to assumption.
Thus, $\rho(A)\geq \rho(AB)$.  Similarly, $\rho(B)\geq \rho(AB).$
Let $C=AB$.   
We have shown that $\rho(C)\leq \rho(B)$ whether or not $A$ is nonsingular.
If $A$ is nonsingular, let $B=A^{-1}C$ and apply \ref{eq:ranlesmin} again to get
$\rho(B)\leq \rho(C)$.  Thus, $\rho(B)=\rho(C)$ if $A$ is nonsingular.
The argument to show $\rho(A)=\rho(C)$ if $B$ is nonsingular is the same.\\
\end{proof}
\end{cor}
%REMARK
\begin{remark}[\bfseries Function notation for \ref{cor:ranpro} proof]
\label{rem:fncnotpro}
The set notation for submatrices used in the proof of \ref{cor:ranpro}
is standard in the literature.  
The equivalent ``function notation'' is, however,
more expressive of what is going on.
Let $A\in {\bf M}_{a,p}(\bK)$ and $B\in {\bf M}_{p,b}(\bK)$ and
let $r=\rho(AB)$. 
Let $g\in {\rm SNC}(r,a)$ and  $h\in {\rm SNC}(r,b)$ be strictly 
increasing functions such that   $\det((AB)_g^h)\neq 0.$   
Note that the $r\times r$ matrix $(AB)_g^h=A_gB^h$
where $A_g$ is $r\times p$ and  $B^h$ is  $p\times r$.
Analogous to \ref{eq:caubinexatwoplus} and \ref{eq:caubinexatwo} we can write
\begin{equation}
\label{eq:fncnotpro}
\det((AB)_g^h)=\det(A_gB^h)= \sum_{f\in {\rm SNC}(r,p)} \det(A_g^f)\det(B^h_f)
\end{equation}
and use this identity instead of  \ref{eq:caubinsetnot} in the proof of Corollary~\ref{cor:ranpro}.

\end{remark}
%EXERCISES CAUCHY BINET LAPLACE
\section*{Exercises: Cauchy Binet and Laplace}
\label{sec:cauchylaplace}
\begin{definition}[\bfseries Greatest common divisor]
\label{def:gcd}
Let $S\subset \bZ$ be a finite set of integers containing at least one nonzero integer.  The set of greatest common divisors of $S$ is $\{-d, d\}$ where $d$ is the largest positive integer that divides all of the integers in $S$. 
We call $d$ {\em the} greatest common divisor, $d={\rm gcd}(S)$.  
See~\ref{rem:gcdlcmmeetjoin}.  \\
\end{definition}
%EXERCISE DETERMINANTAL DIVISORS
\begin{exercise}
\label{ex:detdivcauchybinet}
Let $A\in {\bf M}_{a,p}(\bZ)$ and $B\in {\bf M}_{p,b}(\bZ).$  Let $C=AB$.
Assume $n\leq \min(a,b,p)$.
Suppose that $\det(C[G\Mid H]) \neq 0$ for some 
$G\in \bP_n(\,{\underline a}\,),\;H\in \bP_n(\,{\underline b}\,).$
Let
\[{\mathcal C}_n = {\rm gcd}(\{\det(C[G\Mid H])\Mid G\in \bP_n(\,{\underline a}\,),\;H\in \bP_n(\,{\underline b}\,)\})
\] 
\[{\mathcal A}_n = {\rm gcd}(\{\det(A[G\Mid F])\Mid G\in \bP_n(\,{\underline a}\,),\;F\in \bP_n(\,{\underline p}\,)\})
\]
\[{\mathcal B_n} = {\rm gcd}(\{\det(B[F\Mid H])\Mid F\in \bP_n(\,{\underline p}\,),\;H\in \bP_n(\,{\underline b}\,)\})
\]
where ${\rm gcd}$ denotes the greatest common divisor.
Prove that  $\det(A[G\Mid F])\neq 0$ for some $G$ and $F$,
$\det(B[F\Mid H])\neq 0$ for some $F$ and $H$ and
${\mathcal A}_n$ and ${\mathcal B}_n$  
both divide ${\mathcal C}_n$ (recall \ref{cor:cauchybinetset}).\\
\end{exercise}
%EXERCISE COFACTORS AND UNITS
\begin{exercise}
Repeat example~\ref{rem:sgncofactexample} with
\[
A=
\left(
\begin{array}{rrr}
1&  1 & 1 \\
-1&  -2 & -2  \\
1&  2&  1 
\end{array}
\right)
\]
by finding $C_A$, the signed cofactor matrix, and $B_A=C_A^T$, the transpose of $C_A$.
Verify the identity of \ref{cor:transposesignedcofactor}.  Is this matrix, $A$,
an invertible element (i.e., a unit \ref{rem:unitassoc} ) in the ring ${\bf M}_{3,3}(\bZ)$?\\
\end{exercise}

%EXERCISE CRAMER
\begin{exercise}
Use Cramer's rule (\ref{cor:cramerrule}) to find $X$ in the equation $AX=Y$ where
\[
A=\left(
\begin{array}{rrr}
1&0&1\\
-1&1&0\\
0&-1&2
\end{array}
\right)
\;\;\;\;
X=\left(
\begin{array}{c}
x_1\\
x_2\\
x_3
\end{array}
\right)
\;\;\;\;
Y=\left(
\begin{array}{c}
1\\
0\\
1
\end{array}
\right).
\]
\end{exercise}

%CHAPTER 3   
\chapter{Hermite/echelon forms}
\index{Hermite form!left-unit equivalence}
\section*{Row equivalence}
\label{sec:rowequ}
A system of distinct representatives (SDR, \ref{def:setpartition})  for a partition or equivalence relation (\ref{rem:equivrel}) is sometimes called a set of ``canonical forms.''  As an example, we will discuss  equivalence relations on
${\bf M}_{n,m}(\bZ)$, the set of all $n\times m$ matrices over the integers $\bZ$.
%DEFINITION

Review~\ref{rem:specialringsfields} for a list of the rings of primary interest to us and the class of rings that we designate by $\bK$.
%END DEFINITION
\index{matrix!unit, nonsingular, invertible}
Recall that a matrix $Q$ that has a multiplicative inverse, $Q^{-1}$, is called a ``unit'' or ``nonsingular''  or ``invertible.''  These terms are used interchangeably.
Thus, $Q$ is a unit in $\bf M_{n,n}(\bK)$ if and only if $\det(Q)$ is a unit in $\bK$
 (\ref{cor:transposesignedcofactor}).
 The units of a ring with identity form a group called the {\em group of units} (\ref{rem:specialrings}):
 
 %Definition
 \begin{definition}[\bfseries Group of units of $ {\bf M}_{n,n}(\bK)$]
 \index{matrix!group of units!${\rm GL}(n,\bK)$}
 \index{general linear group!${\rm GL}(n,\bK)$}
 \label{def:grouni}
 The group of units of the ring $ {\bf M}_{n,n}(\bK)$ is denoted by ${\rm GL}(n,\bK)$
and is called the {\em general linear group} of ${\bf M}_{n,n}(\bK).$
\end{definition}
%DEFINITION
\begin{definition}[\bfseries Left-unit equivalence]
\label{def:leftuniteq}
\index{matrix!left-unit equivalence}
Define two matrices $A$ and $B$ in ${\bf M}_{n,m}(\bK)$ to be {\em left-unit equivalent} 
if there exists a unit $Q\in {\bf M}_{n,n}(\bK)$ such that $A=QB.$  It is easy to show that left-unit equivalence is an 
equivalence relation on the set ${\bf M}_{n,m}(\bK).$
\end{definition}
%END DEFINITION
 We may think of multiplying on the left by an element $Q$ of 
  ${\rm GL}(n,\bK)$ as a function $f_Q: \bf M_{n,m}(\bK) \rightarrow \bf M_{n,m}(\bK)$
  where $f_Q(A)=QA$.  The function $f_Q$ is bijective with the  range equal to the domain and is thus a permutation of $\bf M_{n,m}(\bK).$
 From \ref{eq:ranlesmin} we have that the function $f_Q$ preserves rank:
 $\rho(f_Q(A)) = \rho (QA) = \rho(A).$   
 We will also discuss right-unit equivalence: $A$ and $B$ in ${\bf M}_{n,m}(\bK)$ are {\em right-unit equivalent}  if there exists a unit $Q\in {\rm GL}(m,\bK))$ 
 such that $A=BQ.$  
Our results for left-unit equivalence transform in a trivial way to right-unit equivalence, including the preservation of rank.
%\hspace*{1 in}\\
 
\begin{lemma}[\bfseries Preserving column relations]
\label{lem:lincolpreserved}
\index{matrix!left-unit equivalence!preserve column relations} 
Let $A, B\in {\bf M}_{n,m}(\bK)$ be left-unit equivalent as in Definition~\ref{def:leftuniteq}. 
Thus, $A=QB$, $Q\in {\rm GL}(n,\bK)).$
Then
$\sum_{j=1}^n c_j  A^{(j)} = \Theta_{n,1}$ (zero matrix) if and only if  $\;\sum_{j=1}^n c_j  B^{(j)} = \Theta_{n,1}$.
\begin{proof}
 \begin{equation}
 \label{eq:preservelincol}
 \left(\sum_{j=1}^n c_j  A^{(j)}\right)=
 \sum_{j=1}^n c_j  (QB)^{(j)} = \sum_{j=1}^n c_j  Q(B^{(j)}) =
 Q\left(\sum_{j=1}^n c_j  B^{(j)}\right).
 \end{equation}
 Thus, $\left(\sum_{j=1}^n c_j  A^{(j)}\right)=Q\left(\sum_{j=1}^n c_j  B^{(j)}\right)$ and
 hence  $\left(\sum_{j=1}^n c_j  B^{(j)}\right)= \Theta_{n,1}$  
 implies that $\sum_{j=1}^n c_j  A^{(j)}=\Theta_{n,1}.$
 Since $Q$ is invertible, 
 $Q^{-1}\left(\sum_{j=1}^n c_j  A^{(j)}\right)=\left(\sum_{j=1}^n c_j  B^{(j)}\right).$
 Hence  $\left(\sum_{j=1}^n c_j  A^{(j)}\right)= \Theta_{n,1}$ 
 implies that $\sum_{j=1}^n c_j  B^{(j)}=\Theta_{n,1}.$\\
 \end{proof}
 \end{lemma}
 
 %REMARK PRESERVE COL RELATIONS
 \begin{remark}[\bfseries Preserving column relations]
 \label{rem:preservecolrels}
Take
\[
 A=
\left(
\begin{array}{cccccr}
1&0&2&0&2&2 \\
 0&1&1&0&1&1\\
 0&0&0&1&1&-1  
\end{array}
\right)
\;\;{\rm and}\;\;
 B=
\left(
\begin{array}{cccccr}
1&1&3&0&3&3 \\
 0&1&1&1&2&0\\
 0&1&1&2&3&-1  
 \end{array}
\right).
\]
Note that $QA=B$ for the unit matrix
 \[
Q=
\left(
\begin{array}{ccc}
1 & 1 & 0  \\
0 & 1 & 1  \\
0& 1 &  2 
 \end{array}
\right)
\;\;{\rm where}\;\;
Q^{-1}=
\left(
\begin{array}{ccc}
1 & -2 & +1  \\
0 & +2 & -1  \\
0& -1 &  +1 
 \end{array}
\right).
\]
Thus, $A^{(3)} + A^{(4)} = A^{(5)}$ implies that $B^{(3)} + B^{(4)} = B^{(5)}$
and conversely. 
Likewise, $A^{(3)} - A^{(4)} = A^{(6)}$ implies that $B^{(3)} -B^{(4)} = B^{(6)}$ and
conversely.
In other words,  linear relations among columns is 
a left-unit {\em equivalence class invariant.}  In particular, this invariance implies that if $A$ and $B$ are in the same 
left-unit equivalence class then columns $A^{j_1}, \ldots, A^{j_r}$ are linearly independent
if and only if columns $B^{j_1}, \ldots, B^{j_r}$ are linearly independent.  
Importantly, columns $A^{(2)}, A^{(3)}, A^{(4)}$ are {\em obviously} independent so
columns $B^{(2)}, B^{(3)}, B^{(4)}$ are independent.
In particular, the fact that linear relations among columns is 
a left-unit equivalence class invariant shows that matrix $D$ is not 
left-unit equivalent to $A$ and $B$: 
 \[
D=
\left(
\begin{array}{cccccc}
1&0&2&0&2&2 \\
 0&1&1&0&1&1\\
 0&0&0&1&1& 1  
\end{array}
\right)
\]
because $D^{(3)} - D^{(4)} \neq D^{(6)}.$\\
\end{remark}
%END REMARK
%Left-unit equivalence preserves linear relations on columns ( \ref{rem:preservecolrels}).  Next we show that it preserves the span (\ref{rev:vecspace}) of the rows (as vectors over $\bK$).\\
Let $A\in {\bf M}_{n,m}(\bK)$ and let  $(A_{(1)}, \ldots , A_{(n)})$ be the rows of $A$.
We use the notation ${\rm Span}((A_{(1)}, \ldots , A_{(n)})$ to denote all linear combinations of the rows of $A$: 
\[
%{\rm Span}((A_{(1)}, \ldots , A_{(n)}) =
\left\{\sum_{i=1}^n a_iA_{(i)}\Mid a_i\in \bK, i=1, \ldots, n \right\}.
\]
\begin{lemma}
\label{lem:preservespan}
\index{matrix!left-unit equivalence!preserve span rows}
Let $A, B\in {\bf M}_{n,m}(\bK)$ be left-unit equivalent as in Definition~\ref{def:leftuniteq}. 
Thus, $A=QB$, $Q\in {\bf M}_{n,n}(\bK)$ a unit.  
Let $(A_{(1)}, \ldots , A_{(n)})$ and $(B_{(1)}, \ldots , B_{(n)})$ be the sequences of
row vectors of $A$ and $B$.
Then \[{\rm Span}((A_{(1)}, \ldots , A_{(n)}) = {\rm Span}((B_{(1)}, \ldots , B_{(n)}).\]
\begin{proof}
Let $\sum_{i=1}^n a_iA_{(i)}$ be a linear combination of the rows of $A$.
We show that there is a linear combination of the rows of $B$ such that
\[ \sum_{i=1}^n b_iB_{(i)} = \sum_{i=1}^n a_iA_{(i)}.\]
In matrix terms, $(b_1, \ldots , b_n) B = (a_1, \ldots , a_n) A.$ 
This latter identity can be written
\[
(b_1, \ldots , b_n) (QA) = (a_1, \ldots , a_n) A 
\]
which can be solved by taking
\[
(b_1, \ldots , b_n) = (a_1, \ldots , a_n)Q^{-1}.
\]
Thus, \[{\rm Span}((A_{(1)}, \ldots , A_{(n)}) \subseteq {\rm Span}((B_{(1)}, \ldots , B_{(n)}).\]
The reverse inclusion follows from $B=Q^{-1}A$ and the same argument.
\end{proof}
\end{lemma}

We now define elementary row and column operations on a matrix.   
Let $\bK$ be as in (\ref{rem:specialringsfields}).
%DEFINITION
\begin{definition}[\bfseries Elementary row and column operations]
\label{def:elrowcolops} %labdef Elementary row and column operations
\index{matrix!row/column operations}
Let $A\in{\bf M}_{n,m}(\bK)$ be an $n\times m$ matrix.  
Define three types of functions from ${\bf M}_{n,m}(\bK)$ to ${\bf M}_{n,m}(\bK)$ called
 {\em elementary row operations}:
\begin{description}
\item[(Type I)] ${\hat R}_{[i][j]}(A)$ interchanges row  $A_{(i)}$ with row $A_{(j)}.$
\item[(Type II)] ${\hat R}_{[i]+c[j]}(A)$  replaces row $A_{(i)}$ with $A_{(i)} + cA_{(j)}$,  $c\in \bK.$
\item[(Type III)]  ${\hat R}_{u[i]}(A)$ replaces row $A_{(i)}$ with $uA_{(i)}$, $u$ a unit in $\bK$.
\end{description}
Let ${\hat C}_{[i][j]}(A)$, ${\hat C}_{[i]+c[j]}(A)$, ${\hat C}_{u[i]}(A)$ be the corresponding elementary column operations.
\end{definition}
%END DEFINITION
%REMARK
\begin{remark}[\bfseries Elementary row operations as matrices]
\label{rem:rowopstomatrix}
\index{matrix!row/column matrices}
Using \ref{eq:colrowmatrix} (first identity), we know that for any $n\times n$ matrix, $Q$,
and $A\in{\bf M}_{n,m}(\bK)$, the  rows 
\[((QA)_{(1)}, \ldots, (QA)_{(n)}) = (Q_{(1)}A, \ldots , Q_{(n)}A).\] 
In particular, define an $n\times n$ matrix by $Q = {\hat R}_{[i]+c[j]}(I)$ where $I$ is the $n\times n$ identity  matrix.
Let $I_{(t)}$ denote row $t$ of the identity matrix.  Then, in terms of rows, 
\[Q=(I_{(1)}, \ldots , I_{(i-1)},\left[I_{(i)} + cI_{(j)}\right], I_{(i+1)}, \ldots I_{(n)})\]
\[QA=(I_{(1)}A, \ldots , I_{(i-1)}A,\left[I_{(i)} + cI_{(j)}\right]A, I_{(i+1)}A, \ldots I_{(n)}A)\]
\[QA=(A_{(1)}, \ldots , A_{(i-1)},\left[A_{(i)} + cA_{(j)}\right], A_{(i+1)}, \ldots A_{(n)})\]
Thus, $QA = {\hat R}_{[i]+c[j]}(A)$ so left multiplication of $A$ by $Q$ is the same as 
applying the elementary row operation, ${\hat R}_{[i]+c[j]}$ to $A$.  
Instead of $Q = {\hat R}_{[i]+c[j]}(I)$, we use the notation 
$R_{[i]+c[j]}= {\hat R}_{[i]+c[j]}(I)$ (remove the hat).  
Thus, $R_{[i]+c[j]}\in{\bf M}_{n, n}(\bK)$ is a nonsingular matrix such that
$R_{[i]+c[j]}A = {\hat R}_{[i]+c[j]}(A)$.\\
 \end{remark}
%END REMARK

\begin{lemma}
\label{lem:elementmatrices}
Let $A\in{\bf M}_{n,m}(\bK)$.  
For each of the elementary row (or column) operations, ${\hat R}$
(or $\hat C$), there is an invertible matrix  $R\in{\bf M}_{n, n}(\bK)$ (or $C\in{\bf M}_{m,m}(\bK)$) that when left-multiplied (or right-multiplied) with $A$ results in the same matrix as 
${\hat R}(A)$ 
(or ${\hat C}(A)$). In each case, $R = {\hat R}(I_n)$ (or $C={\hat C}(I_m)$).
\begin{proof}
The argument in each case is similar to that given for 
${\hat R}_{[i]+c[j]}$ in~\ref{rem:rowopstomatrix}.
\end{proof}
\end{lemma}

An $n\times n$ elementary row matrix $R_{[i]+c[j]}$ acts by left multiplication on any 
$n\times m$ matrix $A\in{\bf M}_{n,m}(\bK)$ (for any $m\geq 1$).  
The corresponding elementary row operation $\hat{R}_{[i]+c[j]}$ is defined on 
any matrix $A\in{\bf M}_{n,m}(\bK)$ where $m\geq 1$  and $n \geq \max(i,j)$.
This difference in natural domains between the functions $\hat{R}$ and the matrices $R$ 
(by left multiplication) needs to be kept in mind in some discussions.

%REMARK
\begin{remark}[\bfseries Identities for elementary row matrices]
\label{rem:rowcolmatrices}
\index{matrix!row/column identities}
Check the following for $n=3$ - the case for $n\times n$ matrices is the same idea.
\begin{description}
\item[(Type I)]  $R_{[i][j]}^{-1} = R_{[i][j]},\;\;R_{[i][j]}^T =R_{[i][j]},\;\; \det(R_{[i][j]})=-1$
\item[(Type II)] $R_{[i]+c[j]}^{-1} = R_{[i]-c[j]},\;\;R_{[i]+c[j]}^T = R_{[j]+c[i]},\;\;\det(R_{[i]+c[j]})=1$
\item[(Type III)] $R_{u[i]}^{-1} = R_{u^{-1}[i]},\;\;R_{u[i]}^T = R_{u[i]},\;\;\det(R_{u[i]}) = u$
\end{description}
For example,
$
R_{[2]+[3]}=\left(
\begin{array}{ccc}
1&0&0\\
0&1&1\\
0&0&1
\end{array}\right)
$
and
$
R_{[2]+[3]}^{-1} = R_{[2]-[3]}=\left(
\begin{array}{ccc}
1&0&0\\
0&1&-1\\
0&0&1
\end{array}\right).
$ 
\end{remark}
%\hspace*{1 in}\\
%END REMARK
\begin{remark}[\bfseries Euclidean algorithm and greatest common divisors]
\label{rem:euclideanalg}
\index{gcd and lcm!Euclidean algorithm}
We recall the Euclidean algorithm for computing $r_k=\gcd(r_0,r_1)$ where $r_0$ and 
$r_1$ are nonzero elements of $\bZ$.  
The same algorithm works for  any Euclidean domain, in particular for $\bK.$
The algorithm is usually described by a layout  representing successive divisions.
The layout for $\bK=\bF$ is trivial: $r_0 = q_1r_1$ where $q_1=r_0r_1^{-1}$.
Here is the general pattern:
\[
\begin{array}{ccc}
r_0&=&q_1r_1+r_2\\
r_1&=&q_2r_2+r_3\\
r_2&=&q_3r_3+r_4\\
\vdots & \vdots & \vdots\\
r_{k-3}&=&q_{k-2}r_{k-2}+r_{k-1}\\
r_{k-2}&=&q_{k-1}r_{k-1}+r_k\\
r_{k-1}&=&q_kr_k + 0
\end{array}
\]
If $r_2\neq 0$, the remainders, $r_2, \ldots , r_k$, have strictly decreasing valuations
and thus must terminate with zero -- in this case, $r_{k+1}=0$.
The last nonzero remainder, $r_k$ in this case, is the $\gcd(r_0, r_1)$. 
In fact, the set of all divisors of $r_k$ satisfies:
 $
 \{x: x\Mid r_k\} = \{x: x\Mid r_0\}\cap \{x: x\Mid r_1\}.
 $
 This fact is easily seen (or proved by induction) from the layout above.
 \end{remark}

\begin{remark}[\bfseries Greatest common divisor as linear combination]
\label{rem:gcdaslincomb}
\index{gcd and lcm!$gcd(a,b)=sa+tb$}
Referring to Remark~\ref{rem:euclideanalg}, the second to the last identity in the successive division layout, 
 $r_{k-2}=q_{k-1}r_{k-1}+r_k$
 can be solved for $r_k$ to get $r_k = r_{k-2}-q_{k-1}r_{k-1}=s_{k-2}r_{k-2}+t_{k-1}r_{k-1}$
 (this defines $s_{k-2}$ and $t_{k-1}$).
 Using $r_{k-3}-q_{k-2}r_{k-2}=r_{k-1}$ to eliminate $r_{k-1}$ gives
 $r_k=s_{k-3}r_{k-3}+t_{k-2}r_{k-2}.$  
 Repeating this process (or using induction) gives $r_k=s_0r_0+t_1r_1.$
 The standard theorem from basic algebra is that if $a$ and $b$ are nonzero elements
 of a Euclidean domain $\bK$ and $d=\gcd(a,b)$ then there exists $s, t\in \bK$ such that
 \[d=sa+tb.\]  
This theorem is easily proved without using the Euclidean algorithm by using the fact that the Euclidean domain $\bK$ is also a principle  ideal domain (\ref{rem:pid}).
\end{remark}

\begin{remark}[\bfseries Matrix versions of Euclidean algorithm]
\label{rem:matrixeuclideanalg}
\index{gcd and lcm!matrix version!Euclidean algorithm}
The sequence of remainders displayed in Remark~\ref{rem:euclideanalg} can be represented by a sequence of matrix multiplications as follows:
\begin{equation}
\label{eq:matrixprod4gcd}
\left(\begin{array}{c}r_{t}\\r_{t+1}\end{array}\right) =
\left(\begin{array}{cc} 0&1\\1&-q_t\end{array}\right)
\left(\begin{array}{c}r_{t-1}\\r_t\end{array}\right) 
\;\;t=1, \ldots, k\;\;(r_{k+1}=0).
\end{equation}
The matrix of \ref{eq:matrixprod4gcd} is a product type I and II elementary row matrices
(\ref{def:elrowcolops}):
\begin{equation}
\label{eq:prodtwoelrowmatrices}
R_{[1][2]}R_{[1]-q_t[2]} =\hat{R}_{[1][2]}\hat{R}_{[1]-q_t[2]}I_2 =
\left(\begin{array}{cc} 0&1\\1&-q_t\end{array}\right).
\end{equation}
Let 
\begin{equation}
\label{eq:productelrowsforgcd}
Q_2=\prod_{t=1}^k R_{[1][2]}R_{[1]-q_t[2]}=
\left(\begin{array}{cc} r_{11}&r_{12}\\r_{21}&r_{22}\end{array}\right).
\end{equation}
From the Euclidean algorithm, \ref{rem:euclideanalg}, we have
\begin{equation}
\label{eq:onematrixtogcd}
\left(\begin{array}{c}d\\0 \end{array}\right) =
\left(\begin{array}{cc} r_{11}&r_{12}\\ r_{21}& r_{22}\end{array}\right)
\left(\begin{array}{c}a\\ b\end{array}\right).
\end{equation}

An alternative point of view follows from~\ref{rem:gcdaslincomb}. If $a$ and $b$ are nonzero elements
 of a Euclidean domain $\bK$ and $d=\gcd(a,b)$ then there exists $s, t\in \bK$ such that
 $sa+tb=d$ and hence
\begin{equation}
\label{eq:onematrixprod4gcd}
\left(\begin{array}{c}d\\0 \end{array}\right) =
\left(\begin{array}{cc} s&t\\ \frac{-b}{d}& \frac{+a}{d} \end{array}\right)
\left(\begin{array}{c}a\\ b\end{array}\right).
\end{equation}
Note that $\det \left(\begin{array}{cc} s&t\\ \frac{-b}{d}& \frac{+a}{d} \end{array}\right) = 1$
so this matrix is  a unit in ${\bf M}_{2,2}(\bK).$\\
\end{remark}

If we take 
\[
\left(\begin{array}{c}a\\ b\end{array}\right) = \left(\begin{array}{c}18\\ 12\end{array}\right)\;\;
s=+1,\;\;t=-1,\;\;d=6,\;\;
\left(\begin{array}{cc} s&t\\ \frac{-b}{d}& \frac{+a}{d} \end{array}\right)=
\left(\begin{array}{cc} +1&-1\\-2& +3 \end{array}\right)
\]
then
\begin{equation}
\label{eq:matrixq2}
\left(\begin{array}{cc} +1&-1\\-2& +3 \end{array}\right)
 \left(\begin{array}{c}18\\ 12\end{array}\right) = 
\left(\begin{array}{c}6\\0\end{array}\right)=
\left(\begin{array}{c}d\\ 0\end{array}\right).
\end{equation}
%REMARK
\begin{remark}[\bfseries Examples of matrix versions]
\label{rem:examplesmatrixver}
\index{gcd and lcm!matrix version!examples}
Take $\left(\begin{array}{c}a\\ b\end{array}\right) = \left(\begin{array}{c}18\\ 12\end{array}\right).\;\;$
The Euclidean algorithm has two steps: 
$18=1\cdot 12 + 6$ ($q_1=1$) and $12=2\cdot 6$ ($q_2 = 2$).
Thus, \ref{eq:prodtwoelrowmatrices} and \ref{eq:productelrowsforgcd} become
\begin{equation}
\label{eq:matrixeuclid2by2}
Q_2=
\left(\begin{array}{cc} 0&1\\1& -2 \end{array}\right)
\left(\begin{array}{cc} 0&1\\1& -1 \end{array}\right) =
\left(\begin{array}{cc} +1&-1\\-2& +3 \end{array}\right)=
\left(\begin{array}{cc} r_{11}&r_{12}\\ r_{21}& r_{22}\end{array}\right).
\end{equation}
Thus, \ref{eq:onematrixtogcd} becomes
\begin{equation}
\label{eq:onematrixtogcdex}
\left(\begin{array}{c}6\\0 \end{array}\right) =
\left(\begin{array}{cc} +1&-1\\-2& +3 \end{array}\right)
\left(\begin{array}{c}18\\ 12\end{array}\right).
\end{equation}
To see how these $2\times 2$ matrices are used in general, recall the notation of \ref{def:gendirsum}
and take $Q_4=Q_2\oplus_X^X I_2$ to be the
general direct sum corresponding to $X=Y=\{2,4\}.$ Let $A$ be a $4\times 6$ matrix as shown:
\begin{equation}
\label{eq:matrixq4}
Q_4=
\left(
\begin{array}{cccc}
1&  0 & 0 & 0\\
0&  \underline{+1} & 0 & \underline{-1} \\
0&  0&  1 & 0\\
0 & \underline{-2} & 0 & \underline{+3}
\end{array}
\right)\;\;
A=
\left(
\begin{array}{cccccc}
+1&  -1 & 0 & 0&3&4\\
+2&  0&  \underline{18} & -1&2&6\\
0&  +1 & 2 & 2&5&3 \\
-5 & -4 & \underline{12} & 5&2 &8
\end{array}
\right).
\end{equation}
Using \ref{thm:gendirsums}, we see that $Q_4$ is a unit matrix in ${\bf M}_{4,4}(\bZ)$: 
\[
\det(Q_4)=(-1)^{2\sum X}\det(Q_4[X\Mid X])\det(Q_4(X\Mid X))=\det(Q_2)\det(I_2)=+1.
\]
Thus, $Q_4$ is a unit in ${\bf M}_{4,4}(\bK).$ 
Since $Q_2$ is a product of elementary row operations (or matrices) so is $Q_4$.
In fact, only type I and II matrices are needed.
Consider $Q_4A$.
\begin{equation}
\label{eq:q4timesa}
Q_4A =
\left(
\begin{array}{cccccc}
+1&  -1 & 0 & 0&3&4\\
7&  4&  \underline{6} & -6&0&-2\\
0&  +1 & 2 & 2&5&3 \\
-19 & -12 & \underline{0} & 17&2 &12
\end{array}
\right).
\end{equation}
\end{remark}
Note how $Q_4$ transforms the underlined entries in $A$
shown in \ref{eq:matrixq4} to those shown in \ref{eq:q4timesa} and compare these
transformations with \ref{eq:onematrixtogcdex}.
Note also, that setting $A'=Q_4A$, we have $A'(X\Mid \underline{m}] = A(X\Mid \underline{m}]$
where $X=\{2,4\}$ and $m=6$. 
This discussion leads to the following lemma.
%LEMMA
\begin{lemma}
\label{lem:qntimesa}
\index{gcd and lcm!left-unit equivalence}
Let $A\in {\bf M}_{n,m}(\bK),$  let $(A(i_1, j),  \ldots, A(i_k,j))$, 
$i_1<\cdots < i_k$, be 
$k$ nonzero entries in column $j$ of $A$, and let $s\in \{i_1, \ldots, i_k\}$ be specified.  
There exists a unit $Q\in {\bf M}_{n,n}(\bK)$ such that $QA=A'$ satisfies
\begin{equation}
\label{eq:colsamplereduction1}
(A'(i_1, j), \ldots, A'(s,j), \ldots , A'(i_k,j)) = (0, \ldots ,d, \ldots, 0)
\end{equation}
\begin{equation}
\label{eq:colsamplereduction2}
d=\gcd(A(i_1, j), A(i_2, j), \ldots, A(i_k,j)) 
\end{equation} 
\begin{equation}
\label{eq:colsamplereduction3}
A'(i_1, \ldots i_k \Mid \underline{m}] = A(i_1, \ldots i_k \Mid \underline{m}]\;\;
({\rm notation}\;\;\ref{def:submatrices2}).
\end{equation}
Furthermore, $Q$ can be chosen to be a product of type I and II row operations.
\begin{proof}
The lemma is a restatement of ideas discussed  in 
Remarks~\ref{rem:matrixeuclideanalg} and \ref{rem:examplesmatrixver}.
It suffices to consider $s=i_1$ since repositioning $d$ can be done by one elementary 
(type II) row operation.
The proof is by induction on $k$. The case $k=2$ is discussed in~\ref{rem:matrixeuclideanalg}.
Assume there is a product of elementary row operations $\tilde{Q}$ such that 
$\tilde{Q}A=\tilde{A}$ satisfies
\[
(\tilde{A}(i_2, j), \ldots,  \tilde{A}(i_k,j)) = (\tilde{d}, 0, \ldots , 0),
\]
\[
\tilde{d} = \gcd(A(i_2, j), \ldots,  A(i_k,j)), {\;\;\rm and}
\] 
\[
\tilde{A}(i_2, \ldots i_k \Mid \underline{m}] = A(i_2, \ldots i_k \Mid \underline{m}].
\]
Thus, $(\tilde{A}(i_1,j), \tilde{A}(i_2,j)) = (A(i_1,j), \tilde{d})$ is the case $k=2$. 
Let $Q_n=Q_2\oplus_X^X I_{n-2}$ where $X=\{i_1,i_2\}.$  
$Q=Q_n\tilde{Q}$ is the required unit $Q$ such that $QA=A'$ has the properties 
 \ref{eq:colsamplereduction1}, \ref{eq:colsamplereduction2}, 
 and \ref{eq:colsamplereduction3} (with $s=i_1$).
\end{proof}
\end{lemma}
%SECTION ROW ECHELON FORMS
\section*{Hermite form, canonical forms and uniqueness}
\label{sec:rowechherfor}
Let $A\in{\bf M}_{n,m}(\bK)$ (\ref{rem:specialringsfields}).
We are interested in characterizing certain ``nicely structured'' matrices $H$ which are left-unit equivalent to $A$ (\ref{def:leftuniteq}).
In particular, we study those with the  general structure shown in Figure~\ref{eq:hermiteform}.
Such matrices are called {\em row echelon forms} or {\em Hermite forms}.\\ 
\begin{minipage}{\textwidth}
\index{Hermite form!general structure}
\begin{equation}
\label{eq:hermiteform}
{\bf Figure:\;Hermite\;\,form}
\end{equation}
\begin{center}
\includegraphics{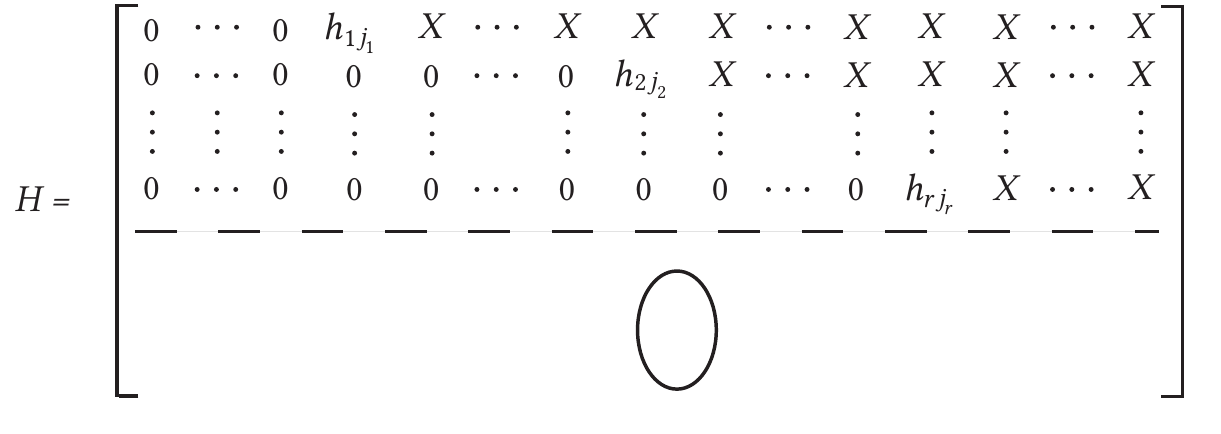}
\end{center}
\end{minipage}
%Definition Hermite form
\begin{definition}[\bfseries Row Hermite or row echelon form]
\label{def:hermiteform}
\index{Hermite form!definition}
A matrix $H\in{\bf M}_{n,m}(\bK)$ is in {\em row Hermite form} if it is the zero matrix, 
$\Theta_{nm}$, or it is nonzero and looks like the matrix in
Figure~\ref{eq:hermiteform}. 
Specifically, for a nonzero $H$ the following hold:
\begin{description}  
\item[$1$] For some $1\leq r\leq n$ the first $r$ rows are nonzero; the rest are zero. 
\item[$2$] In each  nonzero row $i$ the first  or {\em primary} nonzero  row entry is $h_{ij_i}$.
\item[$3$] The {\em primary column indices} are $j_1<j_2\cdots <j_r$.
\end{description}
The number $r$ of nonzero rows is the rank $\rho(H)$ of $H$ which is also the dimension of the vector space spanned by the rows or columns of $H$ over the quotient field of $\bK$.
Note that $\det H[1,\ldots r\Mid j_1,  \ldots,  j_r] \neq 0$, and any $k\times k$
sub-determinant of $H$ with $k>r$ has determinant zero.  
Thus $r$ is the rank of $H$ in the sense of~\ref{def:ranmat} and also 
the rank of any matrix $A$ left-unit equivalent to $H$ (\ref{eq:ranlesmin}).\\
\end{definition}

In most discussions it will be clear if we are talking about ``row'' Hermite form or the alternative
``column'' Hermite form.
We will prove that for $n, m\geq 1$, any matrix in ${\bf M}_{n,m}(\bK)$ is left-unit equivalent to a matrix in Hermite form.  
The general proof for matrices in ${\bf M}_{n,m}(\bK)$ is by induction on $m$, having established the case $m=1$ for $n\geq 1$.\\

\begin{theorem}[\bfseries Hermite form]
\label{thm:hermiteformgen}
\index{Hermite form!existence proof}
Let $A\in{\bf M}_{n,m}(\bK)$ (\ref{rem:specialringsfields}).  
There exists a unit $Q\in{\bf M}_{n,n}(\bK)$ such that $QA=H$ where
$H$ is a Hermite (row echelon) form (\ref{def:hermiteform}).
$Q$ is a product of type I and II elementary row matrices.
\begin{proof}
Suppose $m=1$ and $n\geq 1$.  If $A=\Theta_{n,1}$ then it is a Hermite form by definition.
Suppose $A\neq \Theta_{n,1}$  contains $k$ nonzero entries. Apply~\ref{lem:qntimesa} with $j=1$, $s=i_1$ where 
$(A(i_1, 1),  \ldots, A(i_k,1))$ contains all of the nonzero entries in column $1$.
Thus we obtain a unit $Q'$ such that $A'=Q'A$ has    
$(A'(i_1, 1), \ldots, \ldots , A'(i_k,1)) = (d, \ldots, 0).$ 
Apply $\hat{R}_{[1][i_1]}I_n$ (if $1<i_1$) to obtain the unit matrix $Q\in{\bf M}_{n,n}(\bK)$ such
that $QA=H$ is a nonzero $n\times 1$ matrix in Hermite form.
$Q$ is a product of type  I and II elementary row matrices.

By induction on $m$, assume there is unit 
matrix $P$, a product of type I and II elementary row matrices, such that $PB$ is a Hermite form for any matrix $B\in{\bf M}_{n,m-1}(\bK)$,  $n\geq 1$. 
Now suppose $A\in{\bf M}_{n,m}(\bK)$ where $m>1$.  
Either ${\rm\bf(1)}\,{A}^{(1)}=\Theta_{n1}$ or 
${\rm\bf(2)}\,$ there is unit matrix $\tilde{Q}$, 
a product of type I and II elementary row matrices, 
such that $\tilde{Q}A=\tilde{A}$ has $\tilde{A}^{(1)} = (d, 0, \ldots, 0)$, $d\neq 0$.

In case ${\rm\bf (1)}$, by the induction hypothesis, there is a unit $Q\in {\bf M}_{n,n}(\bK)$ such that $QA[1, \ldots, n\Mid 2, \ldots, m]$ is a Hermite form and thus
$Q{A}$ is a Hermite form.  
$Q$ is a product of type I and II elementary row matrices.

In case ${\rm\bf (2)}$, there is a unit $\tilde{P}\in {\bf M}_{n-1,n-1}(\bK)$ such 
that $\tilde{P}\tilde{A}[2, \ldots, n\Mid 2, \ldots, m]$ is a Hermite form.
Thus, $Q=((1)\oplus \tilde{P})\tilde{Q}$ is such that $QA$ is a Hermite form.
By the induction hypothesis, $\tilde{P}$ (hence $(1)\oplus \tilde{P}$) and $\tilde{Q}$ are products of type I and II elementary row 
operations  and hence so is $Q$.\\
\end{proof} 
\end{theorem}

\begin{remark}[\bfseries Hermite form of a unit matrix]
\label{rem:hermiteofunit}
\index{Hermite form!of unit matrix}
Suppose in \ref{thm:hermiteformgen} $m=n$ and the matrix $A\in{\bf M}_{n,n}(\bK)$ is a unit.
Thus any Hermite form, $H=Q'A$, must be a unit and upper triangular: $H(i,j)=0$ if $i>j$.  
Thus all diagonal elements of $H$ are units.
By using additional type I elementary row  operations there is a unit $P$ such that $PH=H'$ is a diagonal matrix with units along the diagonal.  
By using additional type III elementary row operations, there is a unit $P'$ such that $P'H'= I_n$.
Thus, for any unit matrix $A$ there is a product of elementary row operations of type I and II that reduces $A$ to a diagonal matrix and a product of type I, II, and III elementary row operations that reduces $A$ to the identity.  In particular, $A\in{\bf M}_{n,n}(\bK)$ is a unit if and only if it is a product of elementary row operations (or matrices).  The statement "row operations" can be 
replaced by "column operations."
\end{remark}
%COROLLARY
\begin{corollary}
\label{cor:hermiteofunit}
\index{Hermite form!of unit matrix}
$A\in{\bf M}_{n,n}(\bK)$ is a unit if and only if it is a product of elementary row operations (or matrices). 
\begin{proof}
See the discussion of  \ref{rem:hermiteofunit}.  Note that, in general, type I, II, and III
elementary row operations are required.
\end{proof}
\end{corollary}

\begin{minipage}{\textwidth}
\begin{equation}
{\bf Uniqueness\;results\; for\; Hermite\; forms}
\end{equation}
\index{Hermite form!uniqueness results}
We now discuss  additional structural conditions on Hermite forms that
make them unique (or {\em canonical}).
Note that if $A\in{\bf M}_{n,m}(\bK)$  and   
$\tilde{P},\tilde{Q}\in{\bf M}_{n,n}(\bK)$ are units such that 
$\tilde{P}A=K$ and $\tilde{Q}A=H$  where
$K$ and $H$  are Hermite forms (\ref{def:hermiteform}). 
Then,  $PH=K$ where $P=\tilde{P}\tilde{Q}^{-1}.$
\end  {minipage}

%Lemma
\begin{lemma}[\bfseries Left-unit equivalence of Hermite forms]
\label{lem:unitequiv}
Suppose the Hermite forms $H, K\in{\bf M}_{n,m}(\bK)$  
are left-unit equivalent: $PH=K$.
Then the primary column  indices of $H$ and $K$ are the same.
\begin{proof}
We use Lemma~\ref{lem:lincolpreserved} which states that linear relations among column vectors are preserved under left-unit equivalence. 
Assume first that $\bK$ is  field.
Let $(j_1, \ldots, j_r)$ and $(i_1, \ldots, i_s)$ be the primary column indices for
$H$ and $K$ respectively. 
If   $(j_1, \ldots, j_r)$ does not equal $(i_1, \ldots, i_r)$, then, without loss of generality, let $t$ be the first index such that $j_t<i_t$ (possibly $t=1$.)
Then columns $H^{(j_1)}, \dots, H^{(j_t)}$ are linearly independent but columns
$K^{(j_1)}, \dots, K^{(j_t)}$ are not linearly independent.  
Thus,
$(j_1, \ldots, j_r) = (i_1, \ldots, i_r)$.
If $\bK$ is not a field, apply the same argument to the quotient field of $\bK$,
noting that $PH=K$ is a valid identity in the quotient field.
\end{proof}
\end{lemma}
%End lemma
Recall Definition~\ref{def:setpartition}and the discussion following it.
%Definition
\begin{definition}[\bfseries Canonical SDR for associates]
\label{def:unitsdr}  %labdef Canonical SDR for associates
\index{Hermite form!canonical SDR associates}
For $s,t\in \bK$, define  an equivalence relation  (~\ref{rem:equivrel})  
by $s\sim t$ if $us=t$ for some unit $u\in\bK$.  If $s\sim t$ then $s$ and $t$ are
{\em associates} in $\bK$.  Otherwise, $s$ and $t$ are {\em nonassociates}.
If $\bK=\bZ$ define the SDR for associates to be the set $\{n\Mid n\geq 0\}.$  
If $\bK=\bF[x]$ the SDR for associates is the zero polynomial and all monic polynomials  
(i.e., $a_kx^k + \cdots + a_0$ with $a_k=1$, $1x^0 = 1$). 
If $\bK=\bF$ is a field the  SDR for associates is $\{0,1\}$.
These SDRs are also called {\em complete systems of nonassociates.} 
\end{definition}
%Definition
\begin{definition}[\bfseries Canonical SDR for residues]
\label{def:residuesdr}
\index{Hermite form!canonical SDR residues}
Given $0\neq m\in \bK$, define an equivalence relation on $\bK$ by $s\sim t$
if $m\Mid (s-t)$ ($m$ divides $s-t$).  If $s\sim t$ then $s$ and $t$ are equivalent modulo $m$.
If $\bK=\bZ$ the  {\em (canonical)} SDR for residues modulo $m\;$ is $\{0, 1, \ldots, m-1\}.$   
If  $\bK=\bF[x]$ the SDR for residues modulo $m\;$ is $\{0\} \cup \{p(x)\Mid \deg(p(x))<\deg(m(x))\}.$
If  $\bK=\bF$ the SDR for residues modulo $m\;$ is $\{0\}$. 
These SDRs are also called {\em complete systems of residues modulo $m$}.\\
\end{definition}
 %Definition
 \begin{definition}[\bfseries Hermite canonical form -- row version]
 \label{def:hermitecanonical}
 \index{Hermite form!canonical}
Let $H\in{\bf M}_{n,m}(\bK)$ be a row Hermite form (\ref{def:hermiteform}).
Suppose the primary row entries,
$h_{tj_t}$, $1\leq t\leq r$,  are elements of the 
canonical SDR for associates for $\bK$ (\ref{def:unitsdr}) and the $h_{ij_t},\;i<t,\;$
are elements of the canonical SDR for residues modulo $h_{tj_t}$, $1\leq t\leq r$.
Then the Hermite form $H$ is called a Hermite {\em canonical} form.
 \end{definition}
% \hspace*{ 1 in}\\
%Remark
\begin{remark}[\bfseries Computing Hermite canonical form]
\label{rem:normalizedhermite}
\index{Hermite form!computing}
Let $Q\in{\bf M}_{n,n}(\bK)$ be a unit and suppose that $QA=H$ where
$H$ is a Hermite form (\ref{def:hermiteform}).
By using  elementary row operations of the form $\hat{R}_{u[t]}$, $u$ a unit,
we can transform the primary row entries,
$h_{tj_t}$, $1\leq t\leq r$, such that they are elements of the 
SDR for associates for $\bK$ (\ref{def:unitsdr}).
Next, by using elementary row operations of the form $\hat{R}_{[i]-q_{ij_t}[t]}$, 
$1\leq i<t$ ($2\leq t\leq r$), on the resulting $H$, we can arrange that $h_{i,j_t}$ 
 are in the SDR for residues modulo $h_{tj_t}$ (\ref{def:residuesdr}).
The  SDR for associates phase is done first, then the SDR for residues  phase. 
The residues are computed left to right. Figure~\ref{eq:hermitecanonical} shows
this computation where the residues already computed are indicated by $h'.$ 
The last residues,  corresponding to dividing by $h_{rj_r}$  (by applying $\hat{R}_{[i]-q_{ij_r}[r]}$ as needed), are yet to be computed. 
\end{remark}
\begin{minipage}{\textwidth}
\begin{equation}
\label{eq:hermitecanonical}
{\bf Figure:\;Compute\;residue\;phase} 
\end{equation} 
\begin{equation*}
\left[
\begin{array}{ccccccccccc}
0&\cdots &h_{1j_1} &\ast\ast\ast   &h'_{1j_2} &\ast\ast\ast & h'_{1j_{r-1}}&\ast\ast\ast&{h}_{1j_r} &\ast\ast\ast&\ast \\
0&\cdots &0&\cdots  &h_{2j_2} &\ast\ast\ast & h'_{2j_{r-1}}&\ast\ast\ast &{h}_{2j_r}&\ast\ast\ast&\ast \\
\vdots & \vdots & \vdots &\vdots & \vdots & \vdots & \vdots& \vdots& \vdots&\vdots&\vdots\\
0&\cdots &0&\cdots  &0& \cdots & h_{(r-1)j_{r-1}} & \ast\ast\ast &{h}_{(r-1)j_r}&\ast\ast\ast&\ast\\
0&\cdots &0&\cdots  &0& \cdots &0 & \cdots & h_{rj_r}&\ast\ast\ast&h_{rm}  
\end{array}
\right]   
\end{equation*}  
\end{minipage} 
%Examples Hermite 
\hspace*{1 in}\\
\begin{remark}[\bfseries Examples of Hermite canonical forms]
\label{rem:examplescanonical}
\index{Hermite form!examples}
We give three examples.
In these examples, we omit the initial zero columns  and the terminal zero rows.
The initial, nonzero, elements of the rows (the $h_{ij_i}$) are referred to as
the primary row entries (sometimes called ``pivots''). For the definitions of the SDR for associates and the SDR for residues see~\ref{def:unitsdr}, \ref{def:residuesdr}.

In our first example, \ref{eq:fieldnormalform}, let $\bK=\bF$ be a field.  In the Hermite canonical form for a field, 
the primary row entries are all $1$.
The elements above the these pivots are all zero since the SDR for residues is $\{0\}$ in a field.
\begin{equation}
\label{eq:fieldnormalform}
\left[
\begin{array}{ccccccccccc}
1 &\ast\ast\ast  &0 & \ast\ast\ast& 0&\ast\ast\ast&0 &\ast\ast\ast&\ast \\
0 &\cdots  &1 & \ast\ast\ast & 0&\ast\ast\ast &0&\ast\ast\ast&\ast \\
\vdots &\vdots & \vdots & \vdots & \vdots& \vdots& \vdots&\vdots&\vdots\\
0&\cdots  &0& \cdots & 1 & \ast\ast\ast &0&\ast\ast\ast&\ast\\
0&\cdots  &0& \cdots &0 & \cdots & 1&\ast\ast\ast&\ast  
\end{array}
\right]
\;\;\;\;\bK=\bF
\end{equation}

In the second example, \ref{eq:intnormalform}, let $\bK=\bZ$ be the integers.  
In the Hermite canonical form for $\bZ$, the pivots are all nonzero 
and belong to the  SDR for associates: $\{n\Mid n\geq 0\}$.  
The elements above the pivots $h_{tj_t}$  are all in the  SDR 
for residues modulo $h_{tj_t}$: $\{0, 1, \ldots, h_{tj_t}-1\}.$

\begin{equation}
\label{eq:intnormalform}
\left[
\begin{array}{ccccccccccc}
3 &\ast\ast\ast  &5 & \ast\ast\ast& 1&\ast\ast\ast&6 &\ast\ast\ast&\ast \\
0 &\cdots  &6 & \ast\ast\ast & 3&\ast\ast\ast &0&\ast\ast\ast&\ast \\
\vdots &\vdots & \vdots & \vdots & \vdots& \vdots& \vdots&\vdots&\vdots\\
0&\cdots  &0& \cdots & 5 & \ast\ast\ast &8&\ast\ast\ast&\ast\\
0&\cdots  &0&\cdots &0 & \cdots & 9&\ast\ast\ast&\ast  
\end{array}
\right]
\;\;\;\;\bK=\bZ
\end{equation}

In the third example, \ref{eq:polynormalform}, let $\bK=\bF[x]$ be the polynomials
over the field $\bF$.  
In the Hermite canonical form for $\bF[x]$, the pivots are all nonzero monic polynomials.
The elements above the pivots $h_{tj_t}$  are  in the SDR 
for residues modulo $h_{tj_t}$: $\{0\} \cup \{p(x)\Mid \deg(p(x))<\deg(h_{tj_t})\}.$

\begin{equation}
\label{eq:polynormalform}
\left[
\begin{array}{ccccccc}
x^3+3 &\ast\ast\ast  &x^5 & \ast\ast\ast& 2x+1&\ast\ast\ast&x^2+6 \\
0 &\cdots  &0 & \cdots & x^2-1&\ast\ast\ast &-2x^2\\
\vdots &\vdots & \vdots & \vdots & \vdots& \vdots& \vdots\\
0&\cdots  &0& \cdots & 0 & \cdots &x^3-8
\end{array}
\right]
\;\;\bK=\bF[x]
\end{equation}
\end{remark}
\hspace*{1 in}

Recall Lemma~\ref{lem:unitequiv} which showed that if $H$ and $K$ are
Hermite forms and $PH=K$ then the primary
column indices of $H$ and $K$ are the same.
If $H$ and $K$ are Hermite {\em canonical} forms then the result is much stronger.
In what follows, $P\in{\bf M}_{n,n}(\bK)$ is a unit  and 
$H, K\in{\bf M}_{n,m}(\bK)$ are Hermite canonical forms.  The primary column
indices of $H$ and $K$ are $j_1<j_2<\cdots < j_r$ where $r$ is the number of non zero rows in $H$ and $K$.  If $\bK$ is a field, $r$ is the rank (row rank equals column rank) 
of $H$ and $K$.
%Uniqueness example
\begin{remark}[\bfseries Uniqueness of Hermite canonical form with $r=1$]
\label{rem:uniquerequal1}
\index{Hermite form!uniqueness results}
Assume $PH=K$ where $H$ and $K$ are $n\times m$ Hermite canonical forms and
$P$ is a unit in ${\bf M}_{n,n}(\bK)$  (which is equivalent to $\det(P)$ a unit in $\bK$).
If $K=\Theta_{nm}$ then $PH=K=\Theta_{nm}$ implies $H=\Theta_{nm}$
for any unit $P\in{\bf M}_{n,n}(\bK)$. 
Thus, we consider the case where $K\neq \Theta$ and take $r=1$,  $P=(p_{ij})$ and $PH=K$ where
%Equation
\begin{equation}
\label{eq:signgammamhermite}
K=
\left[
\begin{array}{cccccccc}
0&\cdots &0&k_{1{j_1}} &\ast &\cdots & \ast &k_{1m}\\
0&\cdots &0&0 &0 &\cdots & 0&0\\
\vdots&\cdots &\vdots&\vdots &\vdots &\cdots &\vdots&\vdots\\
0&\cdots &0&0 &0 &\cdots & 0&0
\end{array}
\right].
\end{equation}
Note that 
%Equation
\begin{equation}
\label{eq:sdrforunitsneeded}
(PH)^{(j_1)} = PH^{(j_1)} = 
\left[
\begin{array}{c}
p_{11}h_{1j_1}\\
p_{21}h_{1j_1}\\
\vdots\\
p_{n1}h_{1j_1}
\end{array}
\right]
=
\left[
\begin{array}{c}
k_{1j_1}\\
0\\
\vdots\\
0
\end{array}
\right].
\end{equation}
Identity \ref{eq:sdrforunitsneeded} implies two important facts:

($1$)~Since $p_{s1}h_{1j_1}=0$, $s=2, \ldots, n$, and $h_{1j_1}\neq 0$, we have
$p_{s1}=0$, $s=2, \ldots, n$.
Thus, ($1$)  implies that $\det(P) =p_{11}\det(P(1\Mid 1)).$  Since $\det(P)$ is a unit in $\bK$, both $p_{11}$  and $\det(P(1\Mid 1))$ are units (\ref{rem:unitassoc}).

($2$) Since $p_{11}h_{1j_1} = k_{1j_1}$ and $p_{11}$ is a unit, the fact that 
$h_{1j_1}$ and $k_{1j_1}$ belong to the same SDR for associates for $\bK$ implies that 
$p_{11}=1$.

Hence,  $P$ has the following structure:
%Equation
\begin{equation}
P=\left[
\begin{array}{cccc} 
1&p_{12}&\cdots&p_{1n}\\
0&\cdot&\cdots&\cdot\\
\vdots&\vdots&P(1|1)&\vdots\\
0&\cdot&\cdots&\cdot
\end{array}
\right]=
\left[
\begin{array}{cc}
I_1&P[1|1)\\
\Theta_{n-1,1}&P(1|1)
\end{array}
\right].
\end{equation}
where $P(1|1) \in {\bf M}_{n-1,n-1}(\bK)$ is a unit and $P[1|1)\in {\bf M}_{1,n-1}(\bK)$ is arbitrary.
%Equation
Thus, since $H$ is a rank one Hermite canonical form,
\begin{equation}
PH=
\left[
\begin{array}{cccccccc}
0&\cdots &0&h_{1{j_1}} &\ast &\cdots & \ast &h_{1m}\\
0&\cdots &0&0 &0 &\cdots & 0&0\\
\vdots&\cdots &\vdots&\vdots &\vdots &\cdots &\vdots&\vdots\\
0&\cdots &0&0 &0 &\cdots & 0&0
\end{array}
\right] = K.
\end{equation}
implies that $h_{1t}= k_{1t}$ for $t=j_1, \ldots, m$ and thus  $H=K$.
\end{remark}

\begin{remark}[\bfseries Uniqueness of Hermite canonical form $2\times m$ case]
\label{rem:uniquerequal2}
\index{Hermite form!uniqueness results}
Assume $PH=K$ where $P\in{\bf M}_{2,2}(\bK)$ is a unit  and 
$H, K\in{\bf M}_{2,m}(\bK)$ are $2\times m$ Hermite canonical forms.
Thus, in this example, $n=r=2$, 
$P=\left[\begin{array}{cc} p_{11}&p_{12}\\p_{21}&p_{22}\end{array}\right]$ 
and $PH=K$ where
\begin{equation}
\label{eq:twobyemhermite}
H=
\left[
\begin{array}{ccccccccccc}
0&\cdots &0&h_{1{j_1}} &\ast &\cdots& \ast&h_{1{j_2}}&\ast  &\cdots & \ast\\
0&\cdots &0&0&0 &\cdots& 0& h_{2{j_2}}&\ast &\cdots& \ast
\end{array}
\right].
\end{equation}
Note that the  first primary column $K^{(j_1)}=(PH)^{(j_1)}=P(H^{(j_1)}).$ 
Thus, 
\begin{equation}
\label{eq:2x2firstprimcol}
\left[
\begin{array}{c}
k_{1{j_1}}\\
0
\end{array}   
\right]
=
\left[
\begin{array}{c}
p_{11}h_{1{j_1}}\\  
p_{21}h_{1{j_1}}
\end{array}
\right]
\end{equation}
which implies (since $h_{1{j_1}}\neq 0$) that $p_{21} = 0$.
Since $\det(P)=p_{11}p_{22}$ is a unit, $p_{11}$ and $p_{22}$ are units in $\bK$.
Since $k_{1{j_1}}$ and $h_{1{j_1}}$ are assumed to be from the same
SDR for associates of $K$ (\ref{def:unitsdr}), we have $p_{11}=1$  and
$k_{1{j_1}}=h_{1{j_1}}.$
At this point, $P=\left[\begin{array}{cc} 1&p_{12}\\0&p_{22}\end{array}\right]$
and
\begin{equation}
\label{eq:topline}
PH=
\left[
\begin{array}{ccccccccccc}
0&\cdots &0&k_{1{j_1}} &\ast &\cdots& \ast&h_{1{j_2}}+p_{12}h_{2{j_2}}&\ast  &\cdots & \ast\\
0&\cdots &0&0&0 &\cdots& 0&p_{22}h_{2{j_2}}&\ast &\cdots& \ast
\end{array}
\right].
\end{equation}
From $PH=K$ we get that $p_{22}h_{2{j_2}}=k_{2{j_2}}$ and since $h_{2{j_2}}$and $k_{2{j_2}}$
belong the same SDR for associates, the unit $p_{22}=1$ and $h_{2{j_2}}=k_{2{j_2}}$.
Thus, we have  $h_{1{j_2}}+p_{12}k_{2{j_2}}=k_{1{j_2}}.$
But, $k_{1{j_2}}$ is in the SDR for residues modulo $k_{2{j_2}}$ (\ref{def:residuesdr}).  
Thus, $p_{12}=0$, $P=I_2$ and $H=K$.
\end{remark}

Remarks \ref{rem:uniquerequal1} and \ref{rem:uniquerequal2} illustrate all of the ideas needed for the general proof. 
We use the standard submatrix notation~\ref{def:submatrices2}.
%THEOREM UNIQUENESS OF HERMITE 
\begin{theorem}[\bfseries Uniqueness of Hermite canonical form]
\label{thm:uniquehermformgen}
\index{Hermite form!uniqueness results}
If $PH=K$ where $P\in{\bf M}_{n,n}(\bK)$ is a unit  and 
$H, K\in{\bf M}_{n,m}(\bK)$ are  Hermite canonical forms
then $H=K$ and $P$ is of the form
\begin{equation}
\label{eq:uniquehermiteformofp}
P=\left[
\begin{array}{cc}
I_r&P[\underline{r}|\underline{r})\\
\Theta_{n-r,r}&P(\underline{r}|\underline{r})
\end{array}
\right]
\end{equation}
where $r$ is the rank of $H$ and $K$, 
$P(\underline{r}|\underline{r}) \in {\bf M}_{n-r,n-r}(\bK)$ is a unit 
and $P[\underline{r}|\underline{r})\in {\bf M}_{r,n-r}(\bK)$.
\begin{proof}
The proof is by induction on $r$ where $r=0$ is trivial.  
The case $r=1$ was proved in Remark~\ref{rem:uniquerequal1}. 
Assume $r>1$ and the theorem is true for the case $r-1$.  
We know the primary row indices, $j_1 < j_2<\cdots <j_r$, are the same for $H$ and $K.$
Thus,
\begin{equation}
\label{eq:sdrforunitsneeded2}
(PH)^{(j_1)} = PH^{(j_1)} = 
\left[
\begin{array}{c}
p_{11}h_{1j_1}\\
p_{21}h_{1j_1}\\
\vdots\\
p_{n1}h_{1j_1}
\end{array}
\right]
=
\left[
\begin{array}{c}
k_{1j_1}\\
0\\
\vdots\\
0
\end{array}
\right].
\end{equation}
Equation~\ref{eq:sdrforunitsneeded2} implies (see \ref{eq:2x2firstprimcol} for the idea)
that 

($1$)~Since $p_{s1}h_{1j_1}=0$, $s=2, \ldots, n$, and $h_{1j_1}\neq 0$, we have
$p_{s1}=0$, $s>1.$ Thus $p_{11}\det(P(1\Mid1))=\det(P)$.  
This implies that $p_{11}$ and $\det(P(1\Mid1))$ are units (\ref{rem:eud}).
We also have that 

($2$)~$p_{11}h_{1j_1}=k_{1j_1}$, and since $h_{1j_1}$ and $k_{1j_1}$ belong to the same SDR for associates (\ref{def:unitsdr}), $p_{11}=1$.

Thus, ($1$) and ($2$) imply that $\det(P) =\det(P(1\Mid 1))$  (\ref{eq:setsinout}). 
Hence, $P(1|1)$ is a unit in ${\bf M}_{n-1,n-1}(\bK)$ (\ref{cor:transposesignedcofactor}). 
$P$ has the following structure:
%Equation
\begin{equation}
P=\left[
\begin{array}{cccc} 
1&p_{12}&\cdots&p_{1n}\\
0&\cdot&\cdots&\cdot\\
\vdots&\vdots&P(1|1)&\vdots\\
0&\cdot&\cdots&\cdot
\end{array}
\right]=
\left[
\begin{array}{cc}
I_1&P[1|1)\\
\Theta_{n-1,1}&P(1|1)
\end{array}
\right].
\end{equation}
where $P(1|1) \in {\bf M}_{n-1,n-1}(\bK)$ is a unit and $P[1|1)\in {\bf M}_{1,n-1}(\bK)$. 
Thus we have
\begin{equation}
\label{eq:callcaserminus1}
PH=
\left[
\begin{array}{cc}
H[1|1,\ldots, j_1]&(PH)[1|1,\ldots, j_1)\\
\Theta_{n-1,j_1}&P(1|1)H(1|1,\ldots, j_1)
\end{array}
\right] = K.
\end{equation}

Equation \ref{eq:callcaserminus1} implies that $H[1|1,\ldots, j_1] = [0 \cdots 0\;k_{j_1}],$
$(PH)[1|1,\ldots, j_1) = K[1|1,\ldots, j_1)$ and
\begin{equation}
\label{eq:phequalkinduction}
P(1|1)H(1|1,\ldots, j_1)=K(1|1,\ldots, j_1).
\end{equation}
By the induction hypothesis, \ref{eq:phequalkinduction} gives
$H(1|1,\ldots, j_1)=K(1|1,\ldots, j_1)$ since $P(1|1) \in {\bf M}_{n-1,n-1}(\bK)$ is a unit
and  $H(1|1,\ldots, j_1)$ and $K(1|1,\ldots, j_1)$ are Hermite canonical forms in
${\bf M}_{n-1,m-j_1}(\bK).$ 

We claim that $(PH)[1|1,\ldots, j_1) = K[1|1,\ldots, j_1)$ implies that
$p_{12}=\cdots = p_{1r} = 0$.  
Otherwise, let $t$ be the first integer such that $p_{1t}\neq 0$, $2\leq t\leq r$.
This implies (see  \ref{eq:topline} for basic idea)
\[k_{1t}=PH(1,t) = p_{11}h_{1t} + p_{1t}k_{tj_t}= h_{1t} + p_{1t}k_{tj_t}\;\; ({\rm as}\;p_{11}=1)\]
which contradicts the fact that $k_{1t}$ and $h_{1t}$ are in the canonical 
SDR for residues modulo $k_{tj_t}=h_{tj_t}$.    

Finally, note that  \ref{eq:phequalkinduction} implies that by the induction
hypothesis, $P(1|1)$ has the structure of \ref{eq:uniquehermiteformofp} with
$n$ replaced by $n-1$, $r$ replaced by $r-1$ and $P$ replaced by $P(1|1)$.
This completes the proof.
\end{proof}
\end{theorem}
%END THEOREM

\section*{Stabilizers of ${\rm GL}(n,\bK)$; column Hermite forms}
%${\rm GL}(n,\bK)$ acting on ${\bf M}_{n,n}
%$\varmathbb{K}$
Recall that the group of units of the ring  ${\bf M}_{n,n}(\bK)$ is denoted 
by ${\rm GL}(n,\bK)$ and is called the {\em general linear group} 
of ${\bf M}_{n,n}(\bK)$ (\ref{def:grouni}).
%DEFINITION
\begin{definition}[\bfseries Stability subgroups for left unit multiplication]
\label{def:stasublefmul}
\index{group!stabilizers ${\rm GL}(n,\bK)$}
The subgroup $\{Q\Mid Q\in {\rm GL}(n,\bK), QA=A\}$ is called the {\em stabilizer} or
{\em stability} subgroup of ${\rm GL}(n,\bK)$ at $A\in {\bf M}_{n,m}(\bK)$.
We denote this subgroup by ${\rm GL}_A(n, \bK).$
\end{definition}
%REMARK
\begin{remark}[\bfseries Conjugate stability subgroups]
\label{rem:constasub}
\index{group!conjugate stabilizers ${\rm GL}(n,\bK)$! on ${\bf M}_{n,m}(\bK)$}
Note that if $QA=B$ then ${\rm GL}_A(n, \bK) =Q {\rm GL}_B(n, \bK)Q^{-1}$.
To prove this identity, note that $X\in {\rm GL}_B(n, \bK)$ if and only if $XB=B$ if and only if
$XQA = QA$ if and only if $Q^{-1}XQA=A$ if and only if $Q^{-1}XQ\in {\rm GL}_A(n, \bK).$
Thus, if $A$ and $B$ are left unit equivalent ($QA=B$), their stability subgroups are 
conjugate: ${\rm GL}_A(n, \bK) =Q {\rm GL}_B(n, \bK)Q^{-1}$.
Alternatively stated, if $A \sim B$ under left unit equivalence than the stability subgroups of $A$ and $B$ are conjugate.
\end{remark}
%REMARK
\begin{remark}[\bfseries Stability subgroup of Hermite canonical form]
\label{rem:stasubgroherfor}
\index{Hermite form!stabilizer ${\rm GL}_H(n, \bK)$}
Let $H\in {\bf M}_{n,m}(\bK)$ be a Hermite canonical form. As a consequence of  
Theorem~\ref{thm:uniquehermformgen} we know that if $PH=H$ then $P$ is of the form 
\ref{eq:uniquehermiteformofp}.  Conversly, if $P$ is of that form then $PH=H$.  Thus, we have characterized ${\rm GL}_H(n, \bK).$  Block multiplication is as follows: $PQ=$
\begin{equation*}
\left[
\begin{array}{cc}
I_r&P[\underline{r}|\underline{r})\\
\Theta_{n-r,r}&P(\underline{r}|\underline{r})
\end{array}
\right]
\left[
\begin{array}{cc}
I_r&Q[\underline{r}|\underline{r})\\
\Theta_{n-r,r}&Q(\underline{r}|\underline{r})
\end{array}
\right]=
\left[
\begin{array}{cc}
I_r&Q[\underline{r}|\underline{r})+P[\underline{r}|\underline{r})Q(\underline{r}|\underline{r})\\
\Theta_{n-r,r}&P(\underline{r}|\underline{r})Q(\underline{r}|\underline{r})
\end{array}
\right].
\end{equation*}
The inverse of $P$ is constructed as follows:
\begin{equation*}
P^{-1}=
\left[
\begin{array}{cc}
I_r&-P[\underline{r}|\underline{r})P^{-1}(\underline{r}|\underline{r})\\
\Theta_{n-r,r}&P^{-1}(\underline{r}|\underline{r})
\end{array}
\right].
\end{equation*}
\end{remark}
We summarize with the following corollary.
%COROLLARY
\begin{corollary}[\bfseries Characterization of stabilizer ${\rm GL}_A(n, \bK)$]
\label{cor:stasublefmul}
Let $A\in {\bf M}_{n,m}(\bK)$ and let $QA = H$  where $H$ is the Hermite
canonical form of $A$.  Assume the rank of $A$ (and hence $H$) is $r$.  The 
stabilizer ${\rm GL}_A(n, \bK) =Q {\rm GL}_H(n, \bK)Q^{-1}$ where ${\rm GL}_H(n, \bK)$ is the 
set of all matrices $P$ of the form
\[
P=
\left[
\begin{array}{cc}
I_r&P[\underline{r}|\underline{r})\\
\Theta_{n-r,r}&P(\underline{r}|\underline{r})
\end{array}
\right]
\]
where $P(\underline{r}|\underline{r}) \in {\bf M}_{n-r,n-r}(\bK)$ is an arbitrary unit 
and $P[\underline{r}|\underline{r})\in {\bf M}_{r,n-r}(\bK)$ is an arbitrary matrix.
\begin{proof}
The proof follows from the discussion in~\ref{rem:constasub} and~\ref{rem:stasubgroherfor}. \\
\end{proof}
\end{corollary}

\begin{minipage}{\textwidth}
\begin{equation}
\label{eq:subseccolher}
{\bf Right\;unit\;equivalence\;and\;column\;Hermite\;form}
\end{equation}
All of the results concerning left unit equivalence and row Hermite forms have direct analogs 
for right unit equivalence and column Hermite forms.  We discuss the key results here.
\end{minipage}

%REMARK
\begin{remark}[\bfseries Column Hermite form]
\label{rem:colherfor}
\index{Hermite form!canonical!column version}
Suppose $\tilde{A}\in {\bf M}_{n,m}(\bK)$ and  $\tilde{Q}\in {\rm GL}_{n}(\bK)$. 
Suppose
$\tilde{Q}\tilde{A}=\tilde{H}\in {\bf M}_{n,m}(\bK)$  is a {\em row} Hermite 
form of $\tilde{A}$.
Taking transposes, $\tilde{A}^T=A$, $\tilde{Q}^T=Q$,  $\tilde{H}^T=H$, we have
${A}{Q}={H}$ where  ${Q}\in {\rm GL}_{n}(\bK)$.  The matrix ${H}\in {\bf M}_{m,n}(\bK)$  is a {\em column} Hermite form of ${A}\in {\bf M}_{m,n}(\bK)$.
The structure  of ${H}$ is shown in \ref{eq:colherfor} (columns and rows of zeros in bold type can be repeated):

\begin{equation}
\label{eq:colherfor}
{H}=
\left[
\begin{array}{ccccc}
   {\bf 0}         &  {\bf 0}     &\cdots&    {\bf 0}    &     {\bf 0}\\%endline
h_{j_11}&    0    &\cdots&    0    &    {\bf 0} \\
\begin{array}{c}\ast \\\ast \end{array}   &\vdots &\cdots&\vdots&\vdots\\%endline
\ast &0 &\cdots&0 &{\bf 0}\\%endline
h_{j_21}   &h_{j_22}&\cdots&0 &{\bf 0}\\%endline
\begin{array}{c}\ast \\\ast \end{array}   &\begin{array}{c}\ast \\\ast \end{array}   
&\cdots&\vdots&\vdots\\%endline
\begin{array}{c}\ast \\\ast \end{array}   &\begin{array}{c}\ast \\\ast \end{array}
&\cdots&\vdots&\vdots\\%endline
h_{j_r1}   &h_{j_r2}&\cdots&h_{j_rr}&{\bf 0}\\
\begin{array}{c}\ast \\\ast \end{array}   &\begin{array}{c}\ast \\\ast \end{array} 
 &\cdots&\begin{array}{c}\ast \\\ast \end{array} &\vdots\\%endline
h_{m1}   &h_{m2}&\cdots&h_{mr}&{\bf 0}%last line
\end{array}
\right].
\end{equation}
\end{remark}

Referring to \ref{eq:colherfor} , we rephrase definition ~\ref{def:hermiteform}:
\begin{definition}[\bfseries  Column Hermite  or column echelon form]
\label{def:colhermiteform}
A matrix $H\in{\bf M}_{m,n}(\bK)$ is in {\em column Hermite form} if it is the zero matrix, 
$\Theta_{mn}$, or it is nonzero and looks like the matrix in
Figure~\ref{eq:colherfor}. 
Specifically, for a nonzero $H$ the following hold:
\begin{description}  
\item[$1$] For some $1\leq r\leq n$ the first $r$ columns are nonzero; the rest are zero. 
\item[$2$] In each  nonzero column $i$ the first  nonzero or {\em primary} column entry is $h_{j_ii}.$
\item[$3$] The {\em primary row indices} satisfy $j_1<j_2\cdots <j_r$.
\end{description}
The number $r$ of nonzero columns is the rank $\rho(H)$ of $H$ which is also the dimension of the vector space spanned by the colums or rows of $H$ over the quotient field of $\bK$.
Note that $\det H[j_1,  \ldots,  j_r\Mid 1,\ldots, r] \neq 0$, and any $k\times k$
sub-determinant of $H$ with $k>r$ has determinant zero.  
Thus $r$ is the rank of $H$ in the sense of~\ref{def:ranmat} and also 
the rank of any matrix $A$ right-unit equivalent to $H$ (\ref{eq:ranlesmin}).
\end{definition}

\begin{definition}[\bfseries Hermite canonical form -- column version]
 \label{def:colhermitecanonical}
 \index{Hermite form!column canonical}
Let $H\in{\bf M}_{m,n}(\bK)$ be a column Hermite form -- column version (\ref{def:colhermiteform}).
Suppose the primary column entries,
$h_{j_tt}$, $1\leq t\leq r$,  are elements of the 
canonical SDR for associates for $\bK$ (\ref{def:unitsdr}) and the $h_{j_t i},\;i<t,\;$
are elements of the canonical SDR for residues modulo $h_{j_tt}$, $1\leq t\leq r$.
Then the Hermite form $H$ is called a Hermite {\em canonical} form
-- column version.
 \end{definition}
 %REMARK
\begin{remark}[\bfseries Row column canonical form uniqueness]
\label{rem:rowcolcanforuni}
If we reduce $A \in{\bf M}_{m,n}(\bK)$ to Hermite canonical form -- column version, $AQ=H$ 
(\ref{def:hermitecanonical}), then the primary column entries $h_{j_i,i}$ are uniquely determined.
If we now reduce $H$ to Hermite canonical form -- row version, $PH$, 
( \ref{def:colhermitecanonical}) we get a matrix of the following form:
$PAQ=$
\begin{equation}
\label{eq:kayrowcolnorfor}
\left[
\begin{array}{cccccccccccc}
d_1  &0 & \cdots& 0&0 &\cdots&0 \\
0  &d_2 & \cdots & 0 &0&\dots&0 \\
\vdots & \vdots & \ddots & \vdots& \vdots&\vdots&\vdots\\
0 &0& \cdots & d_{r-1} &0&\cdots&0\\
0&0& \cdots &0 &d_r&\cdots&0 \\
{\bf 0} &{\bf 0}& \cdots &{\bf 0} & {\bf 0}&\cdots&{\bf 0}  
\end{array}
\right] = 
\left[
\begin{array}{ccc}
D_r&\vline&\Theta_{r,n-r}\\ \hline
\Theta_{m-r,r}&\vline&\Theta_{m-r,n-r}
\end{array}
\right]
\end{equation}
where $r=\rho(A)$, the rank of $A$. 
The matrix $D_r$ is uniquely determined with its entries in the canonical SDR for associates in $\bK$ (\ref{def:unitsdr}).  By further use of elementary row and column operations, we can put the diagonal entries of $D_r$ in any order.
Diagonalization of  matrices will be discussed in section~\ref{sec:diacanforsmifor}.
\end{remark}
\begin{remark}[\bfseries Row/column canonical form when $\bK=\bF$]
\label{rem:rowcancolcan}
Let $A\in {\bf M}_{m,n}(\bK)$ where $\bK=\bF$ is a field have rank $r$. From~\ref{eq:fieldnormalform} we know that the row Hermite canonical form has the  
structure shown in~\ref{eq:fierowcolnorfor0} (where boldface rows and columns of zeros can be repeated):
%Equation
\begin{equation}
\label{eq:fierowcolnorfor0}
\left[
\begin{array}{cccccccccccc}
{\bf 0}&1 &\ast\ast\ast  &0 & \ast\ast\ast& 0&\ast\ast\ast&0 &\ast\ast\ast&\ast \\
{\bf 0}&0 &\cdots  &1 & \ast\ast\ast & 0&\ast\ast\ast &0&\ast\ast\ast&\ast \\
\vdots&\vdots &\vdots & \vdots & \vdots & \vdots& \vdots& \vdots&\vdots&\vdots\\
{\bf 0}&0&\cdots  &0& \cdots & 1 & \ast\ast\ast &0&\ast\ast\ast&\ast\\
{\bf 0}&0&\cdots  &0& \cdots &0 & \cdots & 1&\ast\ast\ast&\ast \\
{\bf 0}&{\bf 0}&\cdots  &{\bf 0}& \cdots &{\bf 0} & \cdots & {\bf 0}&\cdots&{\bf 0}  
\end{array}
\right]
\end{equation}
%Equation
Using elementary column operations (\ref{def:elrowcolops}) of type I and II
(i.e., ${\hat C}_{[i][j]}$  and ${\hat C}_{[i]+c[j]}$) we can reduce \ref{eq:fierowcolnorfor0}
to \ref{eq:fierowcolnorfor1}:
\begin{equation}
\label{eq:fierowcolnorfor1}
\left[
\begin{array}{cccccccccccc}
1  &0 & \cdots& 0&0 &\cdots&0 \\
0  &1 & \cdots & 0 &0&\dots&0 \\
\vdots & \vdots & \ddots & \vdots& \vdots&\vdots&\vdots\\
0 &0& \cdots & 1 &0&\cdots&0\\
0&0& \cdots &0 & 1&\cdots&0 \\
{\bf 0} &{\bf 0}& \cdots &{\bf 0} & {\bf 0}&\cdots&{\bf 0}  
\end{array}
\right] = 
\left[
\begin{array}{ccc}
I_r&\vline&\Theta_{r,n-r}\\ \hline
\Theta_{m-r,r}&\vline&\Theta_{m-r,n-r}
\end{array}
\right]
\end{equation}
\end{remark}
\section*{Finite dimensional vector spaces and Hermite forms}
In this section the ring $\bK$ is a field $\bF$.  
\index{Hermite form!linear transformations}
%We first consider matrices $A\in {\bf M}_{n,m}(\bK)$ where $\bK=\bF$ is a field.
Suppose $V$ is a vector space over  $\bF$ with finite dimension $\dim (V)$.  
If $U\subset V$ is a proper subspace of $V$ then $\dim(U)<\dim(V)$; 
a proper subspace of a vector space has dimension {\em strictly} smaller than that space.
This strict decrease of dimension (rank) with  proper  inclusion of subspaces (submodules) 
is not generally true for modules over rings  (see \ref{rem:vecspamod}).      

Let $V$ and $W$ be vector spaces and let $\bL(V,W)$ denote the linear transformations from $V$ to $W$.  
$\bL(V,W)$ is also designated by ${\rm Hom}(V,W)$ (vector space homomorphisms).
%Definition
\begin{definition}[\bfseries Matrix of a pair  of bases]
\label{def:matpaibas}
\index{matrix!of base pair}
\index{linear transformation!matrix given bases}
Let ${\bf v}=(v_1, \ldots v_q)$ and ${\bf w}=(w_1, \ldots, w_r)$ be ordered bases for $V$ and $W$
respectively.
Suppose for  each $j$, $1\leq j\leq q$, $T(v_j)=\sum_{i=1}^r a_{ij}w_i$. 
The matrix $A=(a_{ij})$ is called the matrix of $T$ with respect to the  base pair $({\bf v},{\bf w})$.
We write $[T]_{\bf v}^{\bf w}$ for $A$.
\end{definition}

For example, let $V = \bR^2$ and $W=\bR^3$.   
Let ${\bf w}= \{w_1, w_2, w_3\}$ and ${\bf v}=\{v_1, v_2\}$.  
Define $T$ by $T(v_1) = 2w_1 + 3w_2 - w_3$ and $T(v_2) = w_1 + 5w_2 + w_3$.
Then
\[[T]_{\bf v}^{\bf w} = \left[\begin{array}{cc} 2&1\\3&5\\-1&1\end{array}\right]\]
is the matrix of $T$ with respect to the  base pair $({\bf v},{\bf w})$.
%Theorem
\begin{theorem}[\bfseries Composition of $T$ and $S$  as matrix multiplication]
\label{thm:commatmul}
\index{matrix!composition of transformations}
\index{linear transformations!composition!$[TS]_{\bf u}^{\bf w} = [T]_{\bf v}^{\bf w}[S]_{\bf u}^{\bf v}$}
Let $S\in \bL(U,V)$ and $T\in \bL(V,W)$.  Let ${\bf u}=(u_1, \ldots, u_p)$, 
${\bf v}=(v_1, \ldots, v_q)$ and ${\bf w}=(w_1, \ldots, w_r)$ be bases for $U,\, V,\, W$.
Then 
\begin{equation*}
[TS]_{\bf u}^{\bf w} = [T]_{\bf v}^{\bf w}[S]_{\bf u}^{\bf v}.
\end{equation*}
\begin{proof}
For $j=1, \ldots , p$,  $(TS)(u_j) =  T(S(u_j))$.   Let $A=[T]_{\bf v}^{\bf w}$ and let
$B= [S]_{\bf u}^{\bf v}$.   Let $C=AB$.   Thus,
\[ 
(TS)(u_j) =  T(S(u_j)) = T\left(\sum_{i=1}^q B(i,j)v_i\right).
\]
By linearity of $T$ we obtain
\[
T\left(\sum_{i=1}^q B(i,j)v_i\right) = \sum_{i=1}^q B(i,j)T(v_i)=
 \sum_{i=1}^q B(i,j)\left(\sum_{t=1}^rA(t,i)w_t\right)=
 \]
\[
\sum_{i=1}^q\sum_{t=1}^r B(i,j)A(t,i)w_t= 
\sum_{t=1}^r\left(\sum_{i=1}^qA(t,i)B(i,j)\right)w_t=
\sum_{t=1}^rC(t,j)w_t. 
\]
Thus  $(TS)(u_j)= \sum_{t=1}^rC(t,j)w_t$ for $j=1, \ldots , p$ or
$[TS]_{\bf u}^{\bf w} = [T]_{\bf v}^{\bf w}[S]_{\bf u}^{\bf v}$.\\
\end{proof}
\end{theorem}

Theorem~\ref{thm:commatmul} has an interesting special case when 
$T\in\bL(V,V)$.

%Corollary
\begin{corollary}[\bfseries Change of basis for $T\in\bL(V,V)$]
\label{cor:chabaslintra}
\index{matrix!change of basis}
\index{linear transformations!change of bases}
Let ${\bf v}=(v_1, \ldots v_q)$ and ${\bf w}=(w_1, \ldots, w_q)$ be ordered bases for 
the $q$-dimensional vector spaces $V$ and $W$ over the field $\bF$ and
let $T\in \bL(V,V)$. Then
\begin{equation}
\label{eq:chabaslintra}
[I]_{\bf v}^{\bf w}[T]_{\bf v}^{\bf v}[I]_{\bf w}^{\bf v}=[T]_{\bf w}^{\bf w}.
\end{equation}
\begin{proof}
Apply Theorem~\ref{thm:commatmul}.\\
\end{proof}
\end{corollary}

%Remark
\begin{remark}[\bfseries Similarity of matrices and change of bases]
\label{rem:simmatchabas}
\index{matrix!base change - similarity}
\index{linear transformations!identity re base pair!$[I]_{\bf v}^{\bf w}I_q[I]_{\bf w}^{\bf v}=I_q$}
If $T$ is the identity transformation, $T(x)=x$ for all $x\in V$, and
$I_q\in {\bf M}_{n,n}(\bF)$  is the identity matrix, then \ref{eq:chabaslintra}
 becomes $[I]_{\bf v}^{\bf w}I_q[I]_{\bf w}^{\bf v}=I_q.$
 Thus, if $S=[I]_{\bf w}^{\bf v}$ then $S$ is nonsingular (invertible, unit) and 
$S^{-1}=[I]_{\bf v}^{\bf w}$.  
Given any basis  ${\bf v}=(v_1, \ldots v_q)$ for $V$ and any matrix $A\in {\bf M}_{n,n}(\bF)$,  $A$ and $\bf v$ define $T\in \bL(V,V)$ by $A=[T]_{\bf v}^{\bf v}.$
All $S\in {\rm GL}_n(\bF)$ can be interpreted as $[I]_{\bf w}^{\bf v}$ for selected bases (see discussion of this fact in remark~\ref {rem:equbasnonmat}).
Thus, in matrix terms, \ref{eq:chabaslintra} becomes $S^{-1}AS=B$ where
$B=[T]_{\bf w}^{\bf w}$ with respect to the basis ${\bf w}$ which is defined by 
$S^{-1}=[I]_{\bf v}^{\bf w}$.  
When $A$ and $B$ are related by $S^{-1}AS=B$ for nonsingular $S$, they
are called {\em similar} matrices.
 \end{remark}

%Remark
%\hspace*{1 in}\\
\begin{remark}[\bfseries Equivalence of bases and nonsingular matrices]
\label{rem:equbasnonmat}
Let $\mathcal{B}=\{\mathbf{u}\Mid {\bf u}=(u_1, \dots, u_q) {\rm \;a\;basis\;for\;V}\}$
be the set of all ordered bases for the $q$-dimensional vector space $V$ over $\bF$.
Let ${\rm GL}_n(\bF)$   be all $n\times n$ nonsingular matrices over $\bF$ (i.e., the general linear
group).
For a fixed basis ${\bf v}=(v_1, \dots, v_q)$, the correspondence ${\bf u}\mapsto [I]_{\bf u}^{\bf v}$
(alternatively, we could work with ${\bf u}\mapsto [I]_{\bf v}^{\bf u}$) is a bijection between $\mathcal{B}$ and ${\rm GL}_n(\bF)$. 
In other words, given a fixed basis ${\bf v}$,  we have that $P\in {\rm GL}_n(\bF)$  if and only if there is a basis ${\bf u}$ such that
$P=[I]_{\bf u}^{\bf v}$.

Suppose $V=\bR^3$ and ${\bf v}=(v_1, v_2, v_3)=((1,0,0), (0,1,0), (0,0,1))$.
Let 
\[
P= \left[\begin{array}{ccc} 1&0&1\\-1&1&0\\0&-1&2 \end{array}\right] \equiv [I]_{\bf u}^{\bf v}
\]
define a basis ${\bf u}=(u_1, u_2, u_3)$ where
$u_1=1v_1 + (-1)v_2 + 0v_3 = (1, -1, 0)$ and thus the transpose 
$(1, -1, 0)^T=P^{(1)}$.
Likewise, columns $P^{(2)}$ and $P^{(3)}$ define $u_2$ and $u_3$.
It follows that    
\[
[I]_{\bf v}^{\bf u} = P^{-1}= \left[\begin{array}{ccc} 2/3&-1/3&-1/3\\2/3&2/3&-1/3\\1/3&1/3&1/3 \end{array}\right] 
\]
Because of this equivalence between matrices and linear transformations, most  concepts in linear  algebra have a ``matrix version'' and a ``linear transformation'' (or ``operator'') version.  
Going back and forth between these points of view can greatly simplify proofs.
\end{remark}

\begin{definition}[\bfseries Image and kernel of a linear transformation]
\label{def:imakerlintra}
\index{linear transformation!image, kernel}
Let $T\in \bL(V,W)$ be a linear transformation.  
The  set $\{T(x)\Mid x\in V\}$ is a subspace of $W$ called the {\em image} of $T$.   
It is denoted by ${\rm image}(T)$ or ${\rm Im}(T)$.   
The dimension, $\dim({\rm Im}(T))$, is called the {\em rank}, $\rho(T)$, of $T$ 
(see \ref{def:ranmat} and \ref{rem:ranmat}).
The set $N(T)=\{x\Mid x\in V, T(x)=0_W\}$ where $0_W$ is the zero vector  of $W$ is called
the {\em kernel} or {\em null space} of $T$.  
The $\dim(N(T))=\eta(T)$ is called the {\em nullity} of $T$.
\end{definition}

\begin{remark}[\bfseries Representing $T(x)=y$ in matrix form:  $Ax=y$]
\label{rem:matfor}
Let $T\in \bL(V,W)$, $\dim(V)=n$, $\dim(W)=m$.   
Let ${\bf v}=(v_1, \ldots, v_n)$ be a basis for $V$ and
${\bf w}=(w_1, \ldots, w_m)$ a basis for $W$.   
If $x\in V$ write $x=\sum_{i=1}^n x_i v_i.$
Let $T(x)=T\left(\sum_{i=1}^n x_i v_i\right)=\sum_{j=1}^n x_jT(v_j)$.
Defiine $(a_{ij})$ by  $T(v_j) = \sum_{i=1}^m a_{ij}w_i$. 
Let $[x]^{\bf v}$ denote the column vector $[x_1, \ldots , x_n]^T$. 
Thus 
$A[x_1, \ldots , x_n]^T=$
\begin{equation}
\label{eq:matforbascha}
\left[
\begin{array}{cccc}
a_{11} & a_{12} & \hdots & a_{1n} \\
\vdots & \vdots & \vdots & \vdots\\
a_{m1} & a_{m2} & \hdots & a_{mn} 
\end{array}
\right]
\left[\begin{array}{c}x_1\\ \vdots\\ x_n\end{array}\right] =
[T]_{\bf v}^{\bf w}[x]^{\bf v}=[T(x)]^{\bf w}.
\end{equation}
Note that multiplication of \ref{eq:matforbascha} on the left by a nonsingular $m\times m$ matrix $P$ results in expressing $T(x)$ in a different basis:
\[
PA[x_1, \ldots , x_n]^T= 
[I]_{\bf w}^{\bf u}[T(x)]^{\bf w}=[T(x)]^{\bf u}
\]
since for any ${\bf w}$ and matrix $P$ there is a unique  basis ${\bf u}=(u_1, \ldots u_m)$ such that 
$P=[I]_{\bf w}^{\bf u}\in {\rm GL}_{m}(\bF)$ (see \ref{rem:equbasnonmat}).
As an example, suppose $[T]_{\bf v}^{\bf w}=A$ as follows:
\begin{equation}
\label{eq:examat}
A=
  \begin{blockarray}{cccccc}
    \stackrel{n_1}{\downarrow}&\stackrel{n_2}{\downarrow}&&\stackrel{n_3}{\downarrow}&&\\
    \begin{block}{[cccccc]}
     2&  5&  4&  -2&  2 & 1  \\ 
     0&  1&  1&   0& -1 & 0  \\
      2&    6&  5& 0& -1 & 0  \\
    \end{block}
  \end{blockarray}.
\end{equation}
By a sequence of elementary row operations 
\begin{equation}
\label{eq:seqelerowope}
P=R_{[1]+[3]}R_{(1/2)[3]}R_{[1]-(5/2)[2]}R_{(1/2)[1]}R_{[3]-[2]}R_{[3]-[1]}
\end{equation}
we can construct a nonsingular $P$ such that
$PA=H$ is in Hermite canonical form:

\[
PA=H=
  \begin{blockarray}{cccccc}
   \stackrel{n_1}{\downarrow}&\stackrel{n_2}{\downarrow}&&\stackrel{n_3}{\downarrow}&&\\
    \begin{block}{[cccccc]}
     1&  0&  \frac{-1}{2}&  0&  \frac{5}{2} & 0  \\ 
     0&  1&  1&   0& -1 & 0  \\
      0&    0&  0& 1& -1 & \frac{-1}{2}  \\
    \end{block}
  \end{blockarray}
\]
where $n_1=1$, $n_2=2$ and $n_3=4$ are the primary column indices (\ref{def:hermiteform}) and
\[
P=[I]_{\bf w}^{\bf u}=
\left[
\begin{array}{ccc}
0&-3&\frac{1}{2}\\
0&1&0\\
\frac{-1}{2}&\frac{-1}{2}&\frac{1}{2}
\end{array}
\right].
\]
Thus, the matrix $H= [T]_{\bf v}^{\bf u}$ is the matrix of $T\in \bL(V,W)$ with
respect to the bases ${\bf v}=(v_1, \ldots, v_n)$  for $V$ and
${\bf u}=(u_1, \ldots, u_m)$  for $W$.
To  obtain ${\bf u}$ explicitly in terms of $\bf w$, we need to compute
\[
P^{-1}=[I]_{\bf u}^{\bf w}=
\left[
\begin{array}{ccc}
2&5&-2\\
0&1&0\\
2&\-6&0
\end{array}
\right]
\]
and thus obtain: $u_1=2w_1 + 2w_3$, $u_2=5w_1+w_2+6w_3$ and $u_3=-2w_1$.\\
\end{remark}

%Remark
\begin{remark}[\bfseries Solving equations, rank and nullity]
\label{rem:solequranul}
\index{linear equations!Hermite form and solutions}
\index{linear equations!rank and nullity}
We discuss the equation $T(x)=a$ where 
$T\in \bL(V,W)$.  Let  ${\bf v}=(v_1, \ldots, v_n)$ be a basis for $V$ and
${\bf w}=(w_1, \ldots, w_m)$ a basis for $W$. 
Using the notation of the previous discussion (\ref{rem:matfor}) we solve the 
equivalent matrix equation
\[A[x_1, \ldots , x_n]^T=[a_1, \ldots, a_m]^T\] where $A\in {\bf M}_{m, n}(\bF)$
and $a=[a_1, \ldots, a_m]^T$ (see \ref{eq:matforbascha}):
\begin{equation}
\label{eq:matforbaschaexa}
\left[
\begin{array}{cccc}
a_{11} & a_{12} & \hdots & a_{1n} \\
\vdots & \vdots & \vdots & \vdots\\
a_{m1} & a_{m2} & \hdots & a_{mn} 
\end{array}
\right]
\left[\begin{array}{c}x_1\\ \vdots\\ x_n\end{array}\right] =
\left[\begin{array}{c}a_1\\ \vdots\\ a_m\end{array}\right].
\end{equation} 
The strategy is to multiply both sides of  equation~\ref{eq:matforbaschaexa} 
by a matrix $P\in{\rm GL}_m(\bF)$ to reduce $A$ to its Hermite canonical form  $H$.
Thus, $PA[x_1, \ldots , x_n]^T=P[a_1, \ldots, a_m]^T$ becomes
$H[x_1, \ldots , x_n]^T=[a'_1, \ldots, a'_m]^T$ where $PA=H$ and
$P[a_1, \ldots, a_m]^T = [a'_1, \ldots, a'_m]^T$.
Left unit equivalence preserves linear relations among 
columns (\ref{lem:lincolpreserved} ) therefore
 \[PA[x_1, \ldots , x_n]^T=P[a_1, \ldots, a_m]^T\] 
 if and only if
\[H[x_1, \ldots , x_n]^T=[a'_1, \ldots, a'_m]^T.\]
To help with the computation it is customary to form the {\em augmented}
matrix of the system: $[A\Mid a]$. For
\begin{equation}
\label{eq:exaofequtosol}
A=
\left[
\begin{array}{cccccc}
 2&  5&  4&  -2&  2 & 1  \\ 
 0&  1&  1&   0& -1 & 0  \\
 2&    6&  5& 0& -1 & 0  
\end{array}
\right]
\;\;{\rm and}\;\;a=
\left[
\begin{array}{c}
5\\
-2\\
-3
\end{array}
\right]
\end{equation}
%Augmented matrix
the augmented matrix is
\begin{equation}
\label{eq:augmat}
[A|a]=
  \begin{blockarray}{ccccccc}
    \stackrel{n_1}{\downarrow}&\stackrel{n_2}{\downarrow}&&\stackrel{n_3}{\downarrow}&&&\\
    \begin{block}{[cccccc|c]}
     2&  5&  4&  -2&  2 & 1 &   5 \\ 
     0&  1&  1&   0& -1 & 0 &  -2 \\
      2&    6&  5& 0& -1 & 0 &  -3 \\
    \end{block}
  \end{blockarray}.
\end{equation}

By a sequence of elementary row operations (\ref{eq:seqelerowope}) we 
construct a nonsingular $P$ such that
\begin{equation}
\label{eq:rannulexa}
P[A|a]=[PA\Mid Pa]=[H|h]=
  \begin{blockarray}{ccccccc}
   \stackrel{n_1}{\downarrow}&\stackrel{n_2}{\downarrow}&&\stackrel{n_3}{\downarrow}&&&\\
    \begin{block}{[cccccc|c]}
     1&  0&  \frac{-1}{2}&  0&  \frac{5}{2} & 0 &   \frac{9}{2} \\ 
     0&  1&  1&   0& -1 & 0 &  -2 \\
      0&    0&  0& 1& -1 & \frac{-1}{2} &  -3 \\
    \end{block}
    \begin{block}{ccccccc}
    x_1&x_2&z_1&x_3&z_2&z_3&\\
    \end{block}
  \end{blockarray}.
\end{equation}
where $H$ is the Hermite canonical form of $A$ and
$n_1=1$, $n_2=2$ and $n_3=4$ are the primary column indices.
We can easily solve the equation
\[H[x_1, x_2, x_3, x_4, x_5, x_6]^T=[\frac{9}{2}, -2, -3]^T\]
by taking 
\[[x_1, x_2, x_3, x_4, x_5, x_6]=[\frac{9}{2}, -2,0,-3,0,0].\]
This same solution solves the original equation
\[A[x_1, x_2, x_3, x_4, x_5, x_6]^T=[5, -2, -3]^T.\]
where $A$ and $a$ are specified in~\ref{eq:exaofequtosol}.
Thinking of $A$ as a linear function from $\bF^6$ to $\bF^3$,
the image of $A$ is all of $\bF^3$ and the rank of  $A$ is $3$,
the dimension of the image.  
Consider all vectors of the form (see~\ref{eq:rannulexa} bottom line right)
\[z=[x_1, x_2, z_1, x_3, z_2, z_3].\]
Choose $z_1, z_2, z_3$ arbitrarily.  Then choose $x_1=(\frac{-1}{2})z_1- (\frac{5}{2})z_2)$, $x_2=-z_1+z_2$ and 
$x_3=z_2+z_3/2$ to satisfy $Az=0$ (for arbitrary $[z_1, z_2, z_3].$)
Thus, the null space, $\{z\Mid Az=0\}$ has dimension $3$.
The rank plus nullity, $\rho(A)+\eta(A) = 3 + 3 = 6 = \dim(\bF^6).$
These ideas extend easily to the general case and show constructively why the rank plus nullity of a linear transformation on a vector space $V$ equals $\dim(V)$.
\end{remark}

%\section*{Diagonal forms}

\section*{Diagonal canonical forms -- Smith form}
\label{sec:diacanforsmifor}
\index{matrix!diagonal form!Smith form}
\index{Smith form, weak Smith form!equivalence of matrices}
\index{diagonal forms!Smith form!weak Smith form}
For the rings $\bK$ see \ref{rem:specialringsfields}.
Start with a matrix $A\in {\bf M}_{m,n}(\bK)$ and reduce it to the column 
canonical form $AQ=H$ shown in~\ref{rem:colherfor} where 
$Q\in {\rm GL}_n(\bK).$
Next reduce the matrix $H$ to row canonical form with $P\in {\rm GL}_m(\bK)$ to get
$PAQ=$
\begin{equation}
\label{eq:kayrowcolnorfor1}
\left[
\begin{array}{cccccccccccc}
d_1  &0 & \cdots& 0&0 &\cdots&0 \\
0  &d_2 & \cdots & 0 &0&\dots&0 \\
\vdots & \vdots & \ddots & \vdots& \vdots&\vdots&\vdots\\
0 &0& \cdots & d_{r-1} &0&\cdots&0\\
0&0& \cdots &0 &d_r&\cdots&0 \\
{\bf 0} &{\bf 0}& \cdots &{\bf 0} & {\bf 0}&\cdots&{\bf 0}  
\end{array}
\right] = 
\left[
\begin{array}{ccc}
D_r&\vline&\Theta_{r,n-r}\\ \hline
\Theta_{m-r,r}&\vline&\Theta_{m-r,n-r}
\end{array}
\right]
\end{equation}
where $r=\rho(A)$.
We will show that by further row and column operations the matrix $PAQ$
of~\ref{eq:kayrowcolnorfor1} can be reduced to a  diagonal matrix  in which
$d_1|d_2|\cdots |d_r$ (i.e., the $d_i$ form a divisibility chain).

\begin{definition}[\bfseries Smith form, weak Smith form]
\label{def:smifor}
\index{Smith form, weak Smith form!definitions}
A matrix $D\in {\bf M}_{m,n}(\bK)$ of the form shown in~\ref{eq:kayrowcolnorfor1}
is a {\em Smith form} if $r=1$ or if $r>1$ and the diagonal elements $d_1, d_2, \ldots d_r$  form a divisibility chain: $d_1|d_2|\cdots |d_r$.  
If $d_1, d_2, \ldots d_r$ satisfy the weaker condition that 
$d_1|d_i$ for $i=1,\ldots, r$ (i.e., $d_1$ divides all of the rest of the $d_i$) then we call $D$ a {\em weak Smith form}.\\ 
\end{definition}

\begin{definition}[\bfseries Equivalence of matrices]
\label{def:equmat}
\index{equivalence of matrices!definition}
Matrices $A, B\in {\bf M}_{m,n}(\bK)$ are {\em left-right equivalent}  if there exists
$Q\in {\rm GL}_n(\bK)$ and $P\in {\rm GL}_m(\bK)$ such that $PAQ=B$.
Note that left-right equivalence is an equivalence relation (\ref{rem:equivrel}).
Often, we refer to left-right equivalent matrices  as just ``equivalent matrices.''
\end{definition}

We will show that every $A\in {\bf M}_{m,n}(\bK)$ is  equivalent (left-right) to a Smith form.
If $\bK$ is a field then~\ref{eq:kayrowcolnorfor1} is a Smith form since all $d_i$ are nonzero and hence units of the field.
Thus, the case of interest will be when $\bK$ is not a field.
We first discuss the case  $\rho(A)=2$.

\begin{remark}[\bfseries The case $\rho(A)=2$] 
\label{rem:thecastwo} 
\index{Smith form!rank 2 case}
Let $A\in {\bf M}_{m,n}(\bK)$ have rank $\rho(A)=2$.
The matrix $PAQ$ of~\ref{eq:kayrowcolnorfor1} becomes
\[
D=
\left[
\begin{array}{ccc}
D_2&\vline&\Theta_{2,n-2}\\ \hline
\Theta_{m-2,2}&\vline&\Theta_{m-2,n-2}
\end{array}
\right] =
D_2 \oplus \Theta_{m-2,n-2}
\]  

where $D_2=\left[\begin{array}{cc} d_1&0\\ 0&d_2\end{array}\right]$. 
Since $\rho(D_2)=2$,  $d_1$ and $d_2$ are nonzero.
 %in the canonical SDR for associates (\ref{def:unitsdr}). 
Suppose $\bK$ is a Euclidean domain but not a field (e.g., $\bZ$ or $\bF(x)$, see~\ref{rem:specialringsfields}).
Let ${\rm gcd}(d_1,d_2)=\delta$ and let ${\rm lcm}(d_1,d_2)=\lambda$ be the 
greatest common divisor and least common multiple of $d_1$ and $d_2$.   From
basic algebra, we have $\delta\lambda = d_1d_2$ (\ref{rem:gcdlcmmeetjoin}).
Choose $s, t\in \bK$ such that $sd_1+td_2=\delta$.
Let $Q_2=\left[\begin{array}{cc} s&-d_2/\delta\\ t&d_1/\delta \end{array}\right]$ and 
$P'_2=\left[\begin{array}{cc} 1&1\\ 0&1 \end{array}\right]. $ 
We have
\[
P'_2D_2Q_2 = 
\left[\begin{array}{cc} 1&1\\ 0&1 \end{array}\right]
\left[\begin{array}{cc} d_1&0\\ 0&d_2\end{array}\right]
\left[\begin{array}{cc} s&-d_2/\delta\\ t&d_1/\delta \end{array}\right]
\]
where 
\[
\left[\begin{array}{cc} 1&1\\ 0&1 \end{array}\right]
\left[\begin{array}{cc} d_1&0\\ 0&d_2\end{array}\right] =
\left[\begin{array}{cc} d_1&d_2\\ 0&d_2\end{array}\right]
\]
and 
\begin{equation}
\label{eq:finalstep}
\left[\begin{array}{cc} d_1&d_2\\ 0&d_2\end{array}\right]
\left[\begin{array}{cc} s&-d_2/\delta\\ t&d_1/\delta \end{array}\right] =
\left[\begin{array}{cc} \delta&0\\ td_2&\lambda \end{array}\right]
\end{equation}
Finally, noting that $\delta \Mid td_2$ (in fact, $\delta \Mid d_2$),  applying the elementary row matrix $R_{[2] - c[1]}$ to \ref{eq:finalstep} where
$c=td_2/\delta$ results in 
\begin{equation}
\label{eq:2by2smith} 
P_2D_2Q_2 =
\left[\begin{array}{cc} \delta&0\\ 0&\lambda \end{array}\right]=
 \left[\begin{array}{cc} d_1\land d_2&0\\ 0&d_1\lor d_2 \end{array}\right]=
\hat{D}_2
\end{equation}
where $P_2=R_{[2] - c[1]}P'_2$ and 
$Q_2=\left[\begin{array}{cc} s&-d_2/\delta\\ t&d_1/\delta \end{array}\right].$
We use the notation $d_1\land d_2 = \gcd(d_1,d_2)$ and $d_1\lor d_2 = \rm{lcm}(d_1,d_2)$
and the fact that $d_1d_2/\gcd(d_1,d_2) = \rm{lcm}(d_1,d_2)$ (\ref{rem:gcdlcmmeetjoin}).
Note that the diagonal matrix $\hat{D}_2$ has rank $2$, and 
$\hat{D}_2(1,1)$  divides $\hat{D}_2(2,2)$.  $\hat{D}_2$ is a Smith form for $D_2$.\\
\end{remark}
%Lemma
\begin{lemma}[\bfseries Weak Smith form]
\label{lem:weasmifor}
\index{Smith form, weak Smith form!proof of weak case}
Let $A\in {\bf M}_{m,n}(\bK)$ with $\rho(A)=r>0$.  
There exists $Q\in {\rm GL}_n(\bK)$ and  $P\in {\rm GL}_m(\bK)$ such 
$PAQ=\left[
\begin{array}{ccc}
D_r&\vline&\Theta_{r,n-r}\\ \hline
\Theta_{m-r,r}&\vline&\Theta_{m-r,n-r}
\end{array}
\right]$ 
where $D_r=\rm{diag}(d_1, \ldots, d_r)$
and $d_1|d_i$ for $i=1\ldots r$ (i.e., weak Smith form). 
%Proof 
\begin{proof}
The proof is by induction on $r$.  The case $r=1$ is trivial.  
The case $r=2$ was shown in remark~\ref{rem:thecastwo}. 
Let $r>2$ and assume the lemma is true for $r-1$. 
We can choose $P$ and $Q$ such that 
\begin{equation}
\label{eq:weasmifor}
PAQ=
\left[
\begin{array}{ccc}
D_r&\vline&\Theta_{r,n-r}\\ \hline
\Theta_{m-r,r}&\vline&\Theta_{m-r,n-r}
\end{array}
\right] 
\end{equation}
as in~\ref{eq:kayrowcolnorfor1} where
$D_r=\rm{diag}(d_1, \ldots d_r).$  
By the induction hypothesis, we can further apply left-right multiplications by nonsingular
matrices (or, equivalently, row and column operations)
such that, using the same notation of~\ref{eq:weasmifor}, $d_1|d_i$ for $i=1,\ldots,r-1$.
Next, by applying the result for $r=2$ (\ref{rem:thecastwo} ), we can construct 
$\hat{P}=P_2\oplus_{\{1,r\}}^{\{1,r\}} I_{m-2}$ and 
$\hat{Q}=Q_2\oplus_{\{1,r\}}^{\{1,r\}} I_{n-2}$ (notation~\ref{def:gendirsum}) 
such that 
\[
\hat{P}D_r \hat{Q} = \hat{D}_r = \rm{diag}(d_1\land d_r, d_2, \ldots d_{r-1}, d_1\lor d_r)
\]
where  now, $\hat{D}_r(1,1)|\hat{D}_r(i,i)$ for $i=1\ldots r$.
This completes the proof.\\
\end{proof}
\end{lemma}
%Theorem
\begin{theorem}[\bfseries Smith form]
\label{thm:smifor}
\index{Smith form, weak Smith form!proof of Smith form case}
Let $A\in {\bf M}_{m,n}(\bK)$, $\rho(A)=r>0$.  
There exists $Q\in {\rm GL}_n(\bK)$ and  $P\in {\rm GL}_m(\bK)$ such 
$PAQ=\left[
\begin{array}{ccc}
D_r&\vline&\Theta_{r,n-r}\\ \hline
\Theta_{m-r,r}&\vline&\Theta_{m-r,n-r}
\end{array}
\right]$ 
where $D_r=\rm{diag}(d_1, \ldots, d_r)$ and $d_1|d_2|\cdots |d_r$.  
Thus, every $A\in {\bf M}_{m,n}(\bK)$ is equivalent to a Smith form.
Alternatively, any $A\in {\bf M}_{m,n}(\bK)$ can be transformed into a Smith form
by elementary row and column operations (\ref{cor:hermiteofunit}).
\begin{proof}
By lemma~\ref{lem:weasmifor},  there exists $\tilde{P}$ and $\tilde{Q}$
such that 
\[
\tilde{P}A\tilde{Q}=\left[
\begin{array}{ccc}
\tilde{D}_r&\vline&\Theta_{r,n-r}\\ \hline
\Theta_{m-r,r}&\vline&\Theta_{m-r,n-r}
\end{array}
\right]
\;\;\;{\rm (weak\;Smith\;form)}
\]
where $\tilde{D}_r=\rm{diag}(\tilde{d}_1,\tilde{d}_2, \ldots, \tilde{d}_r)$ and
$\tilde{d}_1|\tilde{d}_i$ for $i=1,\dots, r$.

Thus, we need only show that there are $\hat{P}\in {\rm GL}_{r}(\bK)$
and   $\hat{Q}\in {\rm GL}_{r}(\bK)$ such that  $\hat{P}\tilde{D}_r\hat{Q}$
is in Smith form.
The proof is by induction. 
The case $r=1$ is trivial.  
The case $r=2$ was shown in remark~\ref{rem:thecastwo}.
Note, in particular, equation~\ref{eq:2by2smith}.  
Let  $r>2$ and assume the theorem is true for $r-1$.
The induction hypothesis applied to 
$D'_{r-1}=\tilde{D}_r(1|1)=\rm{diag}(\tilde{d}_2, \ldots, \tilde{d}_r)$
implies there exists $P'\in {\rm GL}_{r-1}(\bK)$ and $Q'\in {\rm GL}_{r-1}(\bK)$
such that $P'D'_{r-1}Q'=\rm{diag}(d_2, d_3, \ldots, d_r)$ where 
$d_2|d_3|\cdots |d_r$.   Note that since $\tilde{d}_1$ divides all entries of 
$D'_{r-1}=\rm{diag}(\tilde{d}_2, \ldots, \tilde{d}_r)$,
$\tilde{d}_1$ also divides all entries of $P'D'_{r-1}Q'$ (easily seen for multiplication by elementary row and column operations) and hence divides $d_2$.  
Thus, setting $\tilde{d}_1 \equiv d_1$ we have $d_1|d_2|\cdots |d_r$. 
Taking $\hat{P}=(1)\oplus P'$ and $\hat{Q}=(1)\oplus Q'$ gives 
$\hat{P}\tilde{D}_r\hat{Q} = \rm{diag}(d_1, d_2, \ldots, d_r)$ which is a Smith form.\\
\end{proof}
\end{theorem}

\begin{corollary}[\bfseries Pairwise relatively prime diagonal entries]
\label{cor:relpridiaent}
\index{Smith form!pairwise relatively prime!diagonal entries}
Let 
\[
\tilde{D}=\rm{diag}(\tilde{d}_1,\tilde{d}_2, \ldots, \tilde{d}_n)\in {\bf M}_{n,n}(\bK),
\;\;n>1,
\]
and suppose that $\gcd(\tilde{d}_i,\tilde{d}_j) = 1$ for $1\leq i<j\leq n$.
Then
\[
D=(1,\ldots,1,(\tilde{d}_1\tilde{d}_2 \cdots \tilde{d}_n))
\]
is a Smith form for  $\tilde{D}$.
\begin{proof}
We use the notation $a\land b= \gcd(a,b)$ and $a\lor b = \rm{lcm}(a,b)$ 
(\ref{rem:gcdlcmmeetjoin}).
The case $n=2$ follows from \ref{eq:2by2smith}: 
$\tilde{D}=\rm{diag}(\tilde{d}_1,\tilde{d}_2)$ is 
equivalent to 
\[
D=
\left[
\begin{array}{cc} 
\tilde{d}_1\land \tilde{d}_2&0\\ 
0&\tilde{d}_1\lor \tilde{d}_2
\end{array}
\right]
=
\left[
\begin{array}{cc} 
1&0\\ 
0&\tilde{d}_1 \tilde{d}_2
\end{array}
\right].
\]
Assume the case $n-1$.   Then  
$\tilde{D}=\rm{diag}(\tilde{d}_1,\tilde{d}_2, \ldots, \tilde{d}_n)$ is equivalent to
$D'=\rm{diag}(1,\ldots,1,(\tilde{d}_1\tilde{d}_2 \cdots \tilde{d}_{n-1}),\tilde{d}_n).$
Using the case $2$ again on the last two entries we get the result.
\end{proof}
\end{corollary}
%CHAPTER
\chapter{Similarity and equivalence}
\label{ch:simequ}
% ~\ref{thm:smifor} applied to special cases. \\

In this chapter we focus on the Euclidean domains (\ref{def:euclidean}) $\bK\in\{\bZ,\bF[x]\}$, $\bF$ a field as specified in remark~\ref{rem:specialringsfields}.
%Section
\section*{Determinantal divisors and related invariants}
%DEFINTION k determinantal divisor
\begin{definition}[\bfseries $k^\text{th}$ order determinantal divisor]
\label{def:kthdetdiv} %labdef  $k^\text{th}$ order determinantal divisor
\index{determinantal divisors!definition}
Let $A\in {\bf M}_{m,n}(\bK)$ where $\bK$ is a Euclidean domain 
($\bK\in\{\bZ,\bF[x]\}$, $\bF$ a field as specified in remark~\ref{rem:specialringsfields}). % Rings of interest to us
Let  $1\leq k\leq \min(m,n)$.
Let $f_k$ denote a greatest common divisor of all $k\times k$ subdeterimants of $A$:
\begin{equation}
\label{eq:detdiv}
f_k=\gcd \{\det(A_\omega^\gamma)\Mid \omega\in {\rm SNC}(k,m), \gamma\in {\rm SNC}(k,n)\}
\end{equation}
where ${\rm SNC}(k,m)$ denotes the strictly increasing functions from 
$\underline{k}\mapsto \underline{m}$ and $A_\omega^\gamma$ is the submatrix of $A$ with
rows selected by $\omega$ and columns by $\gamma$ (\ref{def:submatrices2}).  
We call $f_k$ a $k^\text{th}$ {\em order determinantal divisor} of $A$;  it is determined up to units in $\bK$. We define $f_0=1$ (the multiplicative identity in $\bK$).
\end{definition}
%DEFINITION
\begin{definition}[\bfseries Determinantal divisor sequences]
\label{def:detdivseq} %labdef Determinantal divisor sequences
\index{determinantal divisors!sequences of}
Let $A\in {\bf M}_{m,n}(\bK)$ with $\bK$  as in remark~\ref{rem:specialringsfields}. 
Let $f_k$, $0\leq k \leq \min(m,n)$, be as in definition~\ref{def:kthdetdiv}. 
Define 
\[
(f_0^A, f_1^A, \ldots , f_{\min(m,n)}^A) 
\]
to be a {\em sequence of determinantal divisors} of $A$.  
Let $\rho(A)$ (\ref{def:ranmat}) be the rank of $A$.   
Define
\[
(f_0^A, f_1^A, \ldots, f_{\rho(A)}^A)
\]
to be a {\em maximal sequence of nonzero determinantal divisors} of $A$.
\end{definition}
%REMARK -EXAMPLE
\begin{remark}[\bfseries Example of determinantal divisor sequences]
\index{determinantal divisors!examples of sequences}
\label{rem:exadetdivseq} %labrem Example of determinantal divisor sequences
To save space, we sometimes write $\det X = \left| X \right|$.
Let $A\in {\bf M}_{3,4}(\bZ)$.
\begin{equation}
\label{eq:exadetdivseq1} %labeq Example of determinantal divisor sequences
A=
\left[
\begin{array}{cccccc}
 0&  4&  6&  2\\ 
 8&  2&  10& 8\\
 2&  0&  4& 4  
\end{array}
\right]
\end{equation}
By definition, $f^A_0=1$.  Obviously, $f^A_1=2$ or $-2$.  
Let's choose $f^A_1=2$, using the canonical SDR for associates for $\bZ$ (\ref{def:unitsdr}).
Since all $2\times 2$ subdeteriminants  of $A$ have even entries, 
all such subdeterminants are divisible by $4$.
The determinant of $A[2,3\Mid 1,2]$ is $-4$.   Thus, $f^A_2= 4$ , again choosing
from the canonical SDR.  
Clearly, $\rho(A)=3$ since $\det A[1,2,3\Mid 1,2,3]= -72\neq 0$.
There are four possible $3\times 3$ subdeterminants:
$\left| A^\gamma\right|$ where 
\[\gamma = (1,2,3),\,(1,2,4),\,(1,3,4)\,\,\textrm{or}\,\,(2,3,4).\]
Note that columns $A^{(2)} + A^{(4)} = A^{(3)}$.  
This implies that $| A^\gamma |= 0$ if $\gamma = (2,3,4)$.
Check that $|A^{(1,2,3)}|=-72$, $|A^{(1,2,4)}|=-72$, $|A^{(1,3,4)}|=-72$ and, thus, $f^A_3 = 72$.  We have,
\[
(f^A_0, f^A_1, f^A_2, f^A_3)=(1, 2, 4, 72)
\]
is both a sequence of determinantal divisors and a maximal sequence of nonzero determinantal divisors of $A$.
\end{remark}

Recall the definition of left-right equivalence (or just equivalence) of  
matrices~\ref{def:equmat}.  
By corollary~\ref{cor:ranpro} we know that if $A$ and $B$ are equivalent,
$A=PBQ$, then their ranks are equal: $\rho(A)=\rho(B)$.
%Lemma
\begin{lemma}[\bfseries Determinantal divisors of equivalent matrices]
\label{lem:detdivequmat}
\index{determinantal divisors!of equivalent matrices}
Let $A, B\in {\bf M}_{m,n}(\bK)$   
and let $Q\in {\rm GL}_n(\bK)$, $P\in {\rm GL}_m(\bK)$ be such that 
$A=PBQ$.   Let $r=\rho(A)=\rho(B)$.  
Then the determinantal divisor sequences satisfy
\begin{equation}
\label{eq:detdivequmat1}
(f_0^A, f_1^A, \ldots , f_{r}^A) = (f_0^B, u_1f_1^B, \ldots , u_rf_{r}^B)
\end{equation}
where $u_i$, $1\leq i\leq r$, are units in $\bK$.   
If the determinantal divisors $f^A_i$ and $f^B_i$, $1\leq i\leq r$,
come from the same SDR for associates in $\bK$ (\ref{def:unitsdr}) then 
\begin{equation}
\label{eq:detdivequmat2}
(f_0^A, f_1^A, \ldots , f_{r}^A) = (f_0^B, f_1^B, \ldots , f_{r}^B).
\end{equation}
%Proof
\begin{proof}
If $A=\Theta_{m,n}$ then \ref{eq:detdivequmat1} and  \ref{eq:detdivequmat2} 
are trivial: $(1)=(1).$
Assume $A\neq\Theta_{m,n}.$
We use Cauchy-Binet, corollary~\ref{cor:cauchybinetset} (equation~\ref{eq:caubinexatwoplus}).
Let $1\leq k\leq r$ and choose $g\in \rm{SNC}(k,m)$ and $h\in \rm{SNC}(k,n).$ 
From \ref{eq:caubinexatwoplus} with $X=PB\in {\bf M}_{m,n}(\bK)$ we have
\begin{equation}
\det(X_g^h)=\det((PB)_g^h)=\det(P_gB^h)= 
\sum_{f\in {\rm SNC}(k,m)} \det(P_g^f)\det(B^h_f)
\end{equation}
where $X_g^h\in {\bf M}_{k,k}(\bK)$, $P_g\in {\bf M}_{k,m}(\bK)$,  
$ B^h\in {\bf M}_{m,k}(\bK)$ and $P_g^f,\, B^h_f\in {\bf M}_{k,k}(\bK)$.
Note that the $k^{\rm{th}}$ ($1\leq k\leq r$) determinantal divisor $f^B_k\neq 0$ divides  
$\det(B^h_f)$ for all $f\in {\rm SNC}(k,m)$.  
Hence  $f^B_k$ divides $\det(X_g^h)$  for all 
$g\in \rm{SNC}(k,m)$ and $h\in \rm{SNC}(k,n)$.   
Thus, $f^B_k$ divides $f^X_k$.  
But $P^{-1}X=B$ so the same argument yields $f^X_k$ divides $f^B_k$
and hence for $1\leq k\leq r$,  $f^X_k=u_kf^B_k$ where $u_k\in \bK$ is a unit.
A similar argument shows that if $Y=BQ$ then $f^Y_k=v_kf^B_k$ 
for $1\leq k\leq r$, where $v_k\in \bK$ is a unit.   
Applying these two results to $A=PBQ$ completes the proof.
\end{proof}
\end{lemma}

%{cor:laplacesimple}  simple Laplace.
\begin{lemma}[\bfseries Divisibility and determinantal divisors]
\label{lem:divdetdiv}
\index{determinantal divisors!form divisibility chain}
Let $A\in {\bf M}_{m,n}(\bK)$ and let $(f_0^A, f_1^A, \ldots , f_{r}^A)$ be a maximal sequence of nonzero determinantal divisors of $A$ where $r=\rho(A)$.  
Then $f_k|f_{k+1}$ for $0\leq k < \rho(A)$.
\begin{proof}
Let $A\in {\bf M}_{m,n}(\bK)$, $g\in \rm{SNC}(k+1,m)$, $h\in \rm{SNC}(k+1,n)$  
for $1\leq k < r$.
From corollary~\ref{cor:laplacesimple}, simple Laplace expansion by the $i^{\rm{th}}$
row of  the $(k+1)\times (k+1)$ matrix $A_g^h$  gives
\begin{equation}
\det(A_g^h) = (-1)^i \sum_{j=1}^{k+1} (-1)^j A_g^h(i,j) \det(A_g^h(i\Mid j)).
\end{equation}
Since $A_g^h(i\Mid j)$ is a $k\times k$ matrix, $f^A_k | \det(A_g^h(i\Mid j))$
for $1\leq j \leq k+1$.
Thus, $f^A_k | \det(A_g^h)$ for all choices of $g$ and $h$ and hence $f^A_k | f^A_{k+1}$.
This completes the proof. 
\end{proof}
\end{lemma}

\begin{definition}[\bfseries Sequence of invariant factors]
\label{def:seqinvfac}
\index{determinantal divisors!invariant factors definition}
\index{invariant factors!definition}
Let $A\in {\bf M}_{m,n}(\bK)$ and let $(f_0^A, f_1^A, \ldots , f_{r}^A)$ be a maximal sequence of nonzero determinantal divisors of $A$ where $r=\rho(A)$.  
The sequence $(q^A_1, \ldots, q^A_r)$ where $q^A_k= f^A_k/f^A_{k-1}$, $k=1, \ldots, r$, is 
called a {\em sequence of invariant factors} of the sequence $(f_0^A, f_1^A, \ldots , f_{r}^A)$. From lemma~\ref{lem:detdivequmat}, the sequences of invariant factors of equivalent matrices are the same up to units in $\bK$.
\end{definition}
%THEOREM
\begin{theorem}[\bfseries Smith form and invariant factors]
\label{thm:smiforinvfac}
\index{Smith form!invariant factors}
\index{invariant factors!Smith form}
\index{invariant factors!form divisibility chain}
Let $A\in {\bf M}_{m,n}(\bK)$, $\rho(A)=r>0$. 
Let  
$D=PAQ=\left[
\begin{array}{ccc}
D_r&\vline&\Theta_{r,n-r}\\ \hline
\Theta_{m-r,r}&\vline&\Theta_{m-r,n-r}
\end{array}
\right]$ 
where $D_r=\rm{diag}(d_1, \ldots, d_r)$ and $d_1|d_2|\cdots |d_r$
be a Smith form of $A$  (\ref{thm:smifor}).   
Then   $(q^A_1, \ldots, q^A_r)$ where $q^A_i = d_i$, $i=1, \ldots, r$, is 
a sequence of invariant factors for $A$.  
Thus, the sequence of invariant factors satisfies 
$q^A_1|q^A_2|\cdots |q^A_{r-1}|q^A_r$.
\begin{proof}
Recall lemma~\ref{lem:detdivequmat} which states that equivalent matrices
have the same determinantal divisor  sequences (up to units).
The special structure of $D$ implies that for $k=1, \ldots, r$
\begin{equation}
\label{eq:detdivsmifor}
f^A_k=\gcd \{\prod_{i=1}^k d_{\gamma(i)}\Mid \gamma\in {\rm SNC}(k,r)\}
\end{equation}
where ${\rm SNC}(k,r)$ denotes the strictly increasing functions from 
$\underline{k}\mapsto \underline{r}$.
The fact that $d_1|d_2|\cdots |d_r$ forms a divisibility chain implies that
$\prod_{i=1}^k d_i|\prod_{i=1}^k d_{\gamma(i)}$ for all $\gamma\in {\rm SNC}(k,r).$
Thus, $f^A_k=\prod_{i=1}^k d_i$, $i=1,\ldots, r,$ and 
$q^A_k= f^A_k/f^A_{k-1}=d_k$
define a sequence of determinantal divisors and a corresponding sequence of invariant factors of $A$.
Note 
$d_1|d_2|\cdots |d_{r-1} |d_r$ implies $q^A_1|q^A_2|\cdots |q^A_{r-1}|q^A_r$.
This completes the proof.\\
\end{proof} 
\end{theorem}

%Corollary
\begin{corollary}[\bfseries Invariants with respect to equivalence]
\label{cor:invwitresequ}
\index{equivalence of matrices!three basic invariants}
Let $A, B\in {\bf M}_{m,n}(\bK)$, $\rho(A)=r>0$. 
If  $A$ and $B$ are equivalent matrices (\ref{def:equmat}) then the following sequences are the same up to units:
\begin{description}
\item[$(1)$ Determinantal divisors] $(f_0^A, f_1^A, \ldots , f_{r}^A)$ and $(f_0^B, f_1^B, \ldots , f_{r}^B)$
\item[$(2)$ Invariant factors]  $(q^A_1, \ldots, q^A_r)$ and $(q^B_1, \ldots, q^B_r)$ 
\item[$(3)$ Smith form diagonal entries]  $(d_1^A, \ldots, d_r^A)$ and $(d_1^B, \ldots, d_r^B).$ 
\end{description}
If any one of the sequence pairs $(1)$, $(2)$ or $(3)$ are the same up to units then $A$ and $B$ are equivalent.
\begin{proof}
By theorem~\ref{thm:smiforinvfac}, the Smith form diagonal entries are the same as
the invariant factors up to units.  
Likewise,  by definition~\ref{def:seqinvfac}, the sequence of determinantal divisors determines the sequence of invariant factors and conversly.
By lemma~\ref{lem:detdivequmat}, if $A$ and $B$ are equivalent matrices then they have the same sequences of determinantal divisors up to units.
In particular, (3) implies that $A$ and $B$ are equivalent to the same Smith form (up to units) and are thus equivalent to each other.\\
\end{proof}
\end{corollary}
%Remark
\begin{definition}[\bfseries Elementary divisors] 
\label{def:eldiv}
\index{invariant factors!elementary divisors defined}
\index{elementary divisors!multiset of}
Let $A\in {\bf M}_{m,n}(\bK)$, $\rho(A)=r>0$. 
Let $(q^A_1, \ldots, q^A_r)$,  $i=1, \ldots, r$, be a sequence 
of invariant factors for $A$.
By a set of {\em distinct} primes we mean a set of primes 
$X=\{p_1, p_2, \ldots, p_s\}$, $|X|=s$. 
%chosen from a specified SDR for associates in $\bK$ 
%(e.g., the canonical SDR \ref{def:unitsdr}).  
Assume that $X$ is chosen to be minimal in the sense that $p\in X$ if and only if 
 there exists an invariant factor $q_i^A$,  $i=1, \ldots, r$, such that  $p| q_i^A$.
Factor all of the $q_i^A$ into prime factors as follows:
\begin{equation}
\label{eq:facinvfac}
\begin{array}{cccc}
q_1^A=&p_1^{e_{11}}p_2^{e_{12}}& \cdots& p_s^{e_{1s}}\\
q_2^A=&p_1^{e_{21}}p_2^{e_{22}}& \cdots &p_s^{e_{2s}}\\
\vdots&&\vdots&\\
q_r^A=&p_1^{e_{r1}}p_2^{e_{r2}}& \cdots &p_s^{e_{rs}}.
\end{array}
\end{equation}

A multiset (\ref{rem:notation}) of {\em elementary divisors} of $A$ is  
\[
\{p_i^{e_{ij}}\Mid (i,j)\in \underline{r}\times \underline{s},\, e_{ij}>0\}.
\]
The elementary divisors are determined up to multiplication by units.
\end{definition}
%Remark
\begin{remark}[\bfseries Multisets (\ref{rem:notation}) of elementary divisors]
\label{rem:mulelediv}
\index{elementary divisors!multiset of}
Let $A\in {\bf M}_{m,n}(\bK)$, $\rho(A)=r>0$.
Both sequences $(f^A_0, \ldots, f^A_r)$ 
and $(q^A_1, \ldots, q^A_r)$ form divisibility chains.
Thus, these sequences can be reconstructed from the corresponding multisets
$\{f^A_0, \ldots, f^A_r\}$  and $\{q^A_1, \ldots, q^A_r\}$ by sorting the multisets in order by divisibility.  
The situation for  elementary divisors is similar but requires more discussion.
Referring to equation~\ref{eq:facinvfac}, the divisibility condition,
$q^A_1|q^A_2|\cdots |q^A_{r-1}|q^A_r$, implies that each sequence
of exponents is weakly increasing: 
$e_{1j}\leq  e_{2j} \leq \cdots \leq e_{rj}$, $1\leq j\leq s$.
Thus, knowing the multiset  $\{e_{1j},  e_{2j} , \cdots , e_{rj}\}$
for a particular $j$ is enough to reconstruct the sequence of exponents for that $j$.
Note also that if you know the rank $r$, then knowing the multiset of nonzero
$e_{ij}$ (corresponding to the elementary divisors) is enough to determine the entire sequence
$e_{1j}\leq  e_{2j} \leq \cdots \leq e_{rj}$, $1\leq j\leq s$.
Likewise, the sequences of powers of individual primes are determined by 
their multisets:
%Multisets prime powers
\begin{equation}
\label{eq:colvecpri}
\begin{blockarray}{cccc}
p_1&p_2&&p_s\\
\begin{block}{[c][c]c[c]}
p_1^{e_{11}}&p_2^{e_{12}}& \cdots& p_s^{e_{1s}}\\
p_1^{e_{21}}&p_2^{e_{22}}& \cdots &p_s^{e_{2s}}\\
\vdots&\vdots&\vdots&\vdots\\
p_1^{e_{r1}}&p_2^{e_{r2}}& \cdots &p_s^{e_{rs}}\\
\end{block}
\end{blockarray}.
\end{equation}
Thus, the multiset $\{p_2^{e_{12}},p_2^{e_{22}}, \ldots , p_2^{e_{r2}}\}$
determines the second column in~\ref{eq:colvecpri}.
In fact, if you know $r$, the multiset of all $p_2^{e_{i2}}$ where $e_{i2}>0$
determines the second column.  
Thus, the second column is determined by the elementary divisors of the form
$p_2^{e_{i2}}$.  
Knowing the rank $r$, the multiset of all elementary divisors can be broken down (by determining the primes $p_i$) into the individual multisets corresponding
to the primes $p_i$, $i=1, \ldots, s$. See remark \ref{rem:exaeledivinvfac} for an example.   
\end{remark}
%LEMMA
\begin{lemma}[\bfseries Elementary divisors determine invariant factors]
\label{lem:recfroelediv}
\index{elementary divisors!determine invariant factors}
Let $A\in {\bf M}_{m,n}(\bK)$, $\rho(A)=r>0$. 
Let $(q^A_1, \ldots, q^A_r)$,  $i=1, \ldots, r$, be a sequence 
of invariant factors for $A$.  
Let the elementary divisors be specified as in \ref{eq:facinvfac}.
Then $r$ together with the multiset (\ref{rem:notation}) of  elementary divisors
\[
\{p_i^{e_{ij}}\Mid 1\leq i\leq r\,, 1\leq j\leq s, e_{ij}>0\}
\]
determines the invariant factors and determinantal divisors (up to units).
Thus, two matrices in ${\bf M}_{m,n}(\bK)$are equivalent if and only if they have the same multiset of elementary divisors and same rank (see \ref{cor:invwitresequ}).
\begin{proof}
The idea for the proof is developed in remark~\ref{rem:mulelediv}.
An example is given in remark~\ref{rem:exaeledivinvfac}.
\end{proof}
\end{lemma}

\begin{remark}[\bfseries Example of elementary divisors to invariant factors]
\label{rem:exaeledivinvfac}
\index{elementary divisors!determine invariant factors!example}
Suppose $r=6$ and the multiset of elementary divisors is
\[X=\{2,2,3,3,4,4,5,5,7,7,9,9,9,25,49\}.\]
The multisets $X_2,X_3,X_5,X_7$ associated with the primes $2,3,5,7$ are
\begin{equation}
\label{eq:mulsetpri}
\{2^1,2^1,2^2,2^2\},\;\{3^1,3^1,3^2,3^2,3^2\},\;
\{5^1,5^1,5^2\},\;\{7^1,7^1,7^2\}.
\end{equation}
Sorting each of these sets into numerical order (as they are already listed) and prefixing the number of ones (e.g., in the form $2^0$, $3^0$, $5^0$, $7^0$) needed to make $r=6$ items in the sorted list gives the columns of the following array which represents the factorization of the invariant factors into powers of primes (as in \ref{eq:facinvfac}):
\begin{equation}
\label{eq:facinvfac2}
\begin{array}{cccc}
q_1^A=2^0&3^0& 5^0& 7^0\\
q_2^A=2^0&3^1& 5^0& 7^0\\
q_3^A=2^1&3^1& 5^0& 7^0\\
q_4^A=2^1&3^2& 5^1& 7^1\\
q_5^A=2^2&3^2& 5^1& 7^1\\
q_6^A=2^2&3^2& 5^2& 7^2
\end{array}.
\end{equation}
A recursive approach is to construct the multisets $X_2,X_3,X_5,X_7$ as in
\ref{eq:mulsetpri} and remove the highest powers of each prime to get
$q_6^A=2^23^25^27^2$ and new sets $X'_2,X'_3,X'_5,X'_7$:
\begin{equation}
\label{eq:mulsetpripri}
\{2^1,2^1,2^2\},\;\{3^1,3^1,3^2,3^2\},\;
\{5^1,5^1\},\;\{7^1,7^1\}.
\end{equation}
Proceed recursively from \ref{eq:mulsetpripri} to get $q_5^A, \ldots, q_2^A$.  
The remaining invariant factor, $q_1^A=1$ is determined by knowing the rank $r=6$.\\
\end{remark}
%Section
\section*{Equivalence vs. similarity}
\index{equivalence vs. similarity}
Matrices $X, Y\in {\bf M}_{m,n}(\bK)$ are {\em equivalent}  if there exists
$Q\in {\rm GL}_n(\bK)$ and $P\in {\rm GL}_m(\bK)$ such that $PXQ=Y$ (\ref{def:equmat}).
In this section we take $m=n$ and $\bK=\bF[x]$ where $\bF$ is a field and
$\bF[x]$ is the ring of polynomials with coefficients in $\bF$. 
See remark~\ref{rem:specialringsfields}. 
As previously, for ${\bf M}_{n,n}(\bK)$, we write ${\bf M}_{n}(\bK)$.
Note that ${\bf M}_{n}(\bF)$ is a subring of ${\bf M}_{n}(\bF[x])$.
%The ring of matrices ${\bf M}_{n}(\bF[x])$ is an $\bF$-algebra.
%(definition~\ref{def:vectoralgebra}).  
As an example of an element of ${\bf M}_{n}(\bF[x])$, let
\begin{equation}
\label{eq:matpolent}
P=\left[\begin{array}{cc} x^2/3 & x^3 - x^2/2\\ 2x^3 + 2/5 & 2x - 3 \end{array}\right]
\end{equation}
be in ${\bf M}_{2}(\bF[x])$ where $\bF = \bQ$, the rational numbers.
Note that 
%Expand polynomial matrix
\begin{equation}
\label{eq:exppolmat}
P=P_3 x^3I + P_2 x^2I  + P_1 x^1I + P_0 x^0I
\end{equation}
where $I$ is the $2\times 2$ identity matrix, $x^tI = \rm{diag}(x^t, x^t)$, and 
\[
P_3 = \left[\begin{array}{cc} 0 &1\\ 2& 0 \end{array}\right]\;\;\;
P_2 = \left[\begin{array}{cc} 1/3 &-1/2\\ 0& 0 \end{array}\right]\;\;\;
P_1 = \left[\begin{array}{cc} 0 &0\\ 0& 2 \end{array}\right]\;\;\;
P_0 = \left[\begin{array}{cc} 0 &0\\ 2/5& -3 \end{array}\right].
\]

Recall delta notation, $\delta(\rm{Statement}) = 1$ if ``Statement'' is true, 
$0$ otherwise (\ref{rem:basicsets}).
%The matrix $E_{pq} \in {\bf M}_{2}(\bF[x])$ is defined by $E_{pq}(i,j) = \delta((p,q)=(i,j))$.
The matrix $P_k\in {\bf M}_{2}(\bF)$ is defined by
$P_k(i,j) = {\rm coeff}(x^k, P(i,j))$, $(i,j)\in \underline{2}\times \underline{2}$,
where ${\rm coeff}(x^k, P(i,j))$ denotes the coefficient of $x^k$ in the polynomial $P(i,j)$.\\
%Definition
\begin{definition}[\bfseries Canonical presentation]
\label{def:canprerigeva}
\index{canonical presentation!$P\in {\bf M}_{n}(\bF[x])$}
Let $P\in {\bf M}_{n}(\bF[x])$ and let  $x^kI={\rm diag}(x^k, \ldots, x^k)$ where 
$I$ is the $n\times n$ identity.  
For each $0\leq k\leq m$, define $P_k\in {\bf M}_{n}(\bF)$ by 
$
P_k(i,j) = {\rm coeff}(x^k, P(i,j)), (i,j)\in \underline{n}\times \underline{n},
$
where ${\rm coeff}(x^k, P(i,j))$  denotes the coefficient of $x^k$ in the polynomial $P(i,j)$.
Let $m=\max \{{\rm degree}(P(i,j))\Mid (i,j)\in \underline{n}\times \underline{n}\}$
(see \ref{eq:matpolent} and \ref{eq:exppolmat}).
The  {\em canonical presentation} of $P$ is\\ 
\begin{equation}
\label{eq:canpre}
P=\sum_{k=0}^m P_k x^kI.
\end{equation}
\end{definition} 
%Remark
\begin{remark}[\bfseries Uniqueness of canonical presentation]
\label{rem:unicanpre}
\index{canonical presentation!$P\in {\bf M}_{n}(\bF[x])$!uniqueness}
Note that the canonical presentation of  a matrix $P\in {\bf M}_{n}(\bF[x])$ is unique in the sense that if  $A_k\in {\bf M}_{n}(\bF)$, $0\leq k\leq m$, then (see \ref{eq:canpre})
\begin{equation}
\label{eq:canpreuni}
P=\sum_{k=0}^m P_k x^kI=\sum_{k=0}^m A_k x^kI\;\;\Rightarrow\;\;A_k = P_k,\; 0\leq k\leq m.
\end{equation}
Identity~\ref{eq:canpreuni} is evident when $m=0$ and can be proved by induction on $m$.  Observe that $x^kI$ commutes with every matrix in 
${\bf M}_{n}(\bF[x])$.
% This discussion can be framed in terms of polynomials over rings.
%In this case, $xI$ (or simply $x$) is an {\em indeterminant} or {\em variable}. 
%The ring  ${\bf M}_{n}(\bF[x])$  is identified with the ring of polynomials with coefficients in the ring ${\bf M}_{n}(\bF)$.
Using the concept of a module (\ref{def:module}), the additive group of the ring
${\bf M}_{n}(\bF[x])$  is a left $R$- module  where $R= {\bf M}_{n}(\bF)$.  
Verifying  the axioms of definition~\ref{def:module}  is trivial since $R= {\bf M}_{n}(\bF)$
is a subring of ${\bf M}_{n}(\bF[x])$.  
The set $X=\{x^kI: k=0, 1, \ldots \}$ is a {\em infinite module basis} for this module.
\end{remark}

We next define functions $\rho_A$ and $\lambda_A$ (called {\em right and left evaluation}) from  ${\bf M}_{n}(\bF[x])$ to itself. 
Right evaluation $\rho_A$ turns out to be {\em linear}  for the left $R$-module in the sense 
that $\rho_A(\alpha P + \beta Q) =  \alpha\rho_A(P) + \beta\rho_A(Q)$ for 
$\alpha, \beta\in {\bf M}_{n}(\bF)$ and $P, Q \in  {\bf M}_{n}(\bF[x])$. 
A symmetric result holds for left evaluation $\lambda_A$ which is linear for the right $R$ module
(i.e. $\lambda_A(P\alpha  +Q \beta ) = \lambda_A(P) \alpha + \lambda_A(Q)\beta$).\\
%Definition
\begin{definition}[\bfseries Evaluation by $A\in{\bf M}_{n}(\bF)$]
\label{def:rigeva}
\index{right evaluation!$\rho_A(P)$, $P\in {\bf M}_{n}(\bF[x])$}
Let $A\in{\bf M}_{n}(\bF)$ and let $P\in {\bf M}_{n}(\bF[x])$ have canonical 
presentation $P=\sum_{k=0}^m P_k x^kI.$  
Define functions, {\em right and left evaluation}, $\rho_A, \lambda_A$ 
from ${\bf M}_{n}(\bF[x])$ to ${\bf M}_{n}(\bF)$ by
\[
\rho_A(P) = \sum_{k=0}^m P_k A^k\;\;({\rm right\; eval.})\;\;\;\;
\lambda_A(P)=\sum_{k=0}^m A^kP_k\;\;({\rm left\;eval.}) .
\]
\end{definition}

We derive some basic properties of right evaluation (properties for left evaluation are analogous).\\
%Lemma
\begin{lemma}[\bfseries Properties of right evaluation]
\label{lem:prorigeva}
Let $A, B\in{\bf M}_{n}(\bF)$, $P, Q \in {\bf M}_{n}(\bF[x])$ and $\alpha, \beta \in{\bf M}_{n}(\bF)$.
Then the right evaluation functon $\rho_A$ satisfies
%Equation 1
\begin{equation}
\label{eq:siglin1}
\index{right evaluation!${\bf M}_{n}(\bF)-{\rm linearity}$}
\rho_A(\alpha P + \beta Q) =  \alpha\rho_A(P) + \beta\rho_A(Q)\;\;
({\bf M}_{n}(\bF)-{\rm linearity}).
\end{equation}
%Equation 2
In general,  $\rho_A$ is not multiplicative ($\rho_A(PQ) \neq \rho_A(P)\rho_A(Q)$), but we have
\begin{equation}
\label{eq:siglin2}
\index{right evaluation!quasi-multiplicative ${\bf M}_{n}(\bF[x])$}
\rho_A(PQ)=\sum_kP_k\rho_A(Q)A^k\;\;\;(\rm{quasi-multiplicative\; property}).
\end{equation}
And as special cases %of \ref{eq:siglin2}
%Equations 3
\begin{equation}
\label{eq:siglin3}
\rho_A(Qx^kI)=\rho_A(Q)A^k\;\;\;{\rm and}\;\;\; AB=BA \;\Rightarrow \rho_A(PB)=\rho_A(P)B.
\end{equation}
%Proof
\begin{proof}
%Prove equation 1
{\bf To prove  \ref{eq:siglin1}}:
Let $R=\alpha P + \beta Q$ and let $\sum_k R_kx^kI$ be the canonical presentation of $R$.     
We use \ref{rem:unicanpre} (uniqueness of canonical presentation).
Note that  $\rho_A(R)=\sum_kR_kA^k$ becomes 
\[
\sum_k(\alpha P_k + \beta Q_k) A^k=
\sum_k\alpha P_k A^k + \sum_k \beta  Q_k A^k =
\alpha \rho_A(P) + \beta \rho_A(Q).
\] 
We take the range of values for $k$  to include all
nonzero values of $R_k$, $P_k$ and $Q_k$. 
This proves \ref{eq:siglin1}.\\
%Prove equation 2

{\bfseries To prove~\ref{eq:siglin2}}:  
Note that
$PQ=(\sum_k P_kx^kI)Q=\sum_k P_kQx^kI.$
Linearity~\ref{eq:siglin1} (where the $P_k\in {\bf M}_{n}(\bF)$ play the role of the coefficients
$\alpha,\beta\ldots\in {\bf M}_{n}(\bF)$) implies
\[
\rho_A(PQ)=\rho_A\left(\sum_kP_k(Qx^kI)\right)=\sum_k P_k\rho_A(Qx^kI)=
\sum_kP_k\rho_A(Q)A^k.
\]
The last equality follows from 
\begin{equation}
\rho_A(Qx^kI)=\sum_t Q_t A^{k+t}=\left(\sum_t Q_t A^{t}\right)A^k=\rho_A(Q)A^k
\end{equation}
which, incidentally, proves the first identity of \ref{eq:siglin3}.
To prove the second identity of \ref{eq:siglin3}, take $Q=B$ in \ref{eq:siglin2},
note that $\rho_A(B) = B$ and use $AB=BA$ which implies $BA^k=A^kB$
for all $k$. 
This  completes the proof of the lemma.\\
\end{proof}
\end{lemma}
%Definition
\begin{definition}[\bfseries Characteristic matrix and characteristic polynomial]
\label{def:chamatchapol}
\index{characteristic matrix/polynomial}
Let $A\in{\bf M}_{n}(\bF)$, and let  $I \in{\bf M}_{n}(\bF)$ be the identity.
The matrix $xI-A \in{\bf M}_{n}(\bF[x])$ is called the {\em characteristic matrix}
of $A$.  The polynomial $\det(xI-A)$ is called the {\em characteristic polynomial} of $A$.
\end{definition}
%Theorem	
\begin{theorem} [\bfseries Equivalence implies similarity]
\label{thm:equimpsim}
\index{equivalent $xI-A$ and $xI-B$!implies similar $A$ and $B$}
Let $A, B \in{\bf M}_{n}(\bF)$.  
There exists $P, Q \in \mathrm {GL}_n(\bF[x])$ such that  $xI-A= P(xI-B)Q$ if and only
if there exists $S\in {\rm GL}_n(\bF)$ such that $A=S^{-1}BS$.   
In words, the characteristic matrices of $A$ and $B$ are equivalent if and only
if $A$ and $B$ are similar (\ref{rem:simmatchabas}).
In fact, $S=\rho_A(Q)$ and $S^{-1}=\rho_B(Q^{-1})$.
\begin{proof}
Note that if $A=S^{-1}BS$ then $xI-A=xI-S^{-1}BS=S^{-1}(xI-B)S$.
Thus, we assume $P^{-1}(xI-A)=(xI-B)Q$ or $P^{-1}xI - P^{-1}A = QxI - BQ$.
Apply $\rho_A$  to both sides and use linearity (\ref{eq:siglin1}):
\begin{equation}
\label{eq:appbotsid}
\rho_A(P^{-1}xI) - \rho_A(P^{-1}A) = \rho_A(QxI) - \rho_A(BQ).
\end{equation}
From \ref{eq:siglin3} (first identity)  $\rho_A(P^{-1}xI) = \rho_A(P^{-1})A.$
From \ref{eq:siglin3} (second identity, $B=A$) $\rho_A(P^{-1}A) = \rho_A(P^{-1})A.$
From  \ref{eq:siglin3} (first identity), $\rho_A(QxI) = \rho_A(Q)A.$
From \ref{eq:siglin1} (${\bf M}_{n}(\bF)$- linearity), $\rho_A(BQ) = B\rho_A(Q).$
Substituting these identities into \ref{eq:appbotsid}, we get $\rho_A(Q)A - B\rho_A(Q)=\Theta_n$ or 
$\rho_A(Q)A = B\rho_A(Q).$
To complete the proof we show that $\rho_A(Q) \in {\rm GL}_n(\bF)$ so we can take
$S=\rho_A(Q).$   Let $R= Q^{-1}$ and use \ref{eq:siglin3} (quasi-multiplicative property): 
\begin{equation}
\label{eq:rhoaaacom}
\rho_A(I)=I=\rho_A(RQ)=\sum_kR_k\rho_A(Q)A^k.
\end{equation}
Note that $\rho_A(Q)A = B\rho_A(Q)$ implies by induction that  $\rho_A(Q)A^k = B^k\rho_A(Q).$
Thus, \ref{eq:rhoaaacom} becomes 
$I=\sum_kR_kB^k\rho_A(Q) = \rho_B(R)\rho_A(Q)=\rho_B(Q^{-1})\rho_A(Q).$
This completes the proof.
\end{proof}
\end{theorem}

%Section
\section*{Characteristic matrices and polynomials}

\begin{theorem}[\bfseries Cayley-Hamilton theorem]
\label{thm:cayhamthe}
\index{Cayley-Hamilton theorem}
Let $f_n(x)=\det(xI-A)$
be the characteristic polynomial (\ref{def:chamatchapol}) of $A\in {\bf M}_{n}(\bF)$, $\bF$ a field of characteristic zero.  
Then $f_n(A) = \Theta_n.$
Alternatively stated, $f_n(A) = \Theta_n$ where $f_n$ is the determinantal divisor of $xI-A$ of highest degree.

\begin{proof}
From definition~\ref{def:signcofactor} and corollary~\ref{cor:transposesignedcofactor}
\begin{equation}
\label{eq:adjdet}
\mathrm{adj}(xI-A)\,(xI-A)=I\det(xI-A) = f_n(x)I
\end{equation}
Let $f_n(x)I = \sum_{k=0}^n c_kIx^k$.
Let $P=\mathrm{adj}(xI-A)$, $Q=(xI-A)$  and apply equation~\ref{eq:siglin2} of
lemma~\ref{lem:prorigeva} (quasi-multiplicative property):
\begin{equation}
\label{eq:rhoprod}
f_n(A)=\rho_A(f_n(x)I)= \rho_A(PQ)=\sum_k P_k\rho_A(Q)A^k 
\end{equation}
where $\sum_kP_k x^k$, $P_k\in {\bf M}_{n}(\bF)$, $0\leq k\leq n$, is the canonical
presentation of $P$.
Substituting $\rho_A(Q) = A - A = \Theta_n$ into~\ref{eq:rhoprod} completes the proof.\\
\end{proof}

\end{theorem}

\begin{remark}[\bfseries Reduced Cayley-Hamilton: $q_n(A)=\Theta_n$]
\label{rem:recayham}
\index{Cayley-Hamilton theorem!$q_n(x)$ version}
The Cayley-Hamilton theorem (\ref{thm:cayhamthe}) states that $f_n(A) = \Theta_n$
where $f_n$ is the determinantal divisor of highest degree of $xI-A$ 
(i.e., $f_n(x)=\det(xI-A)$ is the characteristic polynomial of $A$).
In fact, $q_n(A) = \Theta_n$ where $q_n(x)$ is the invariant factor of highest degree of $xI-A$.
 We refer to theorem~\ref{thm:cayhamthe}, equation~\ref{eq:adjdet}
\begin{equation}
\label{eq:adjdetred1}
\mathrm{adj}(xI-A)\,(xI-A)= f_n(x)I.
\end{equation}
By definition, $f_n=q_nf_{n-1}$ so we have
\begin{equation}
\label{eq:adjdetred2}
\mathrm{adj}(xI-A)\,(xI-A)= q_n(x)f_{n-1}(x)I
\end{equation}
Let $P=\mathrm{adj}(xI-A)$ and note (\ref{def:signcofactor}) that 
the set of all entries of $P$ is the set of all signed cofactors of $(xI-A)$ and, by definition,
$f_{n-1}$ is a greatest common divisor of these cofactors. 
Let $\hat{P}=P/f_{n-1}$  so that \ref{eq:adjdetred2}  becomes
\begin{equation}
\label{eq:adjdetred3}
\hat{P}\,(xI-A)= q_n(x)I
\end{equation}
where $\gcd\{\hat{P}(i,j)\Mid (i,j)\in \underline{n}\times \underline{n}\} = 1$ (up to units).
Set $Q=xI-A$ so that, analogous to~\ref{eq:rhoprod}, appling equation~\ref{eq:siglin2} of
lemma~\ref{lem:prorigeva} (quasi-multiplicative property) we obtain 
\begin{equation}
\label{eq:rhoprod2}
q_n(A)=\rho_A(q_n(x)I)= \rho_A(\hat{P}Q)=\sum_k \hat{P}_k\rho_A(Q)A^k.
\end{equation} 
The fact that $\rho_A(Q)=\rho_A(xI-A))=A-A=\Theta_n$ completes the proof.\\
\end{remark}

\begin{definition}[\bfseries Minimal polynomial]
\label{def:minpol}
\index{minimal polynomial!defined}
Let $A\in{\bf M}_n(\bF)$.  
A polynomial $\phi(x)\in\bF[x]$ is called a
{\em minimal polynomial for $A$} if it has minimal degree among all polynomials  $p(x)\in\bF[x]$
such that $p(A)=\Theta_{n,n}\equiv \Theta_n$.  
We call $\phi$ {\em the minimal polynomial} if it is monic 
(i.e., in the canonical SDR for associates for $\bF[x]$ \ref{def:unitsdr}).
\end{definition}

\begin{remark}[\bfseries Minimal polynomial and divisibility]
\label{rem:minpoldiv}
\index{minimal polynomial!divisibility condition}
Let $A\in{\bf M}_n(\bF)$.
Let $p(x)\in \bF[x]$.  It is easily seen that if $\phi(x)$ is minimal for $A$ and $p(A)=\Theta_n$ then $\phi(x)$ divides $p(x)$. 
Otherwise, the remainder $r(x)$ from dividing $p(x)$ by $\phi(x)$ would satisfy
$r(A)=\Theta_n$ and would contradict the minimality of $\phi(x)$.
Thus, $\phi(x)$ is minimal for $A$ if  and only if it divides 
all $p(x)$ such that $p(A)=\Theta_n$.
\end{remark}

\begin{lemma}[\bfseries Minimal polynomial and invariant factors]
\label{lem:minpolinvfac}
\index{minimal polynomial!invariant factor highest degree}
Let $A\in{\bf M}_n(\bF)$ and let $\phi(x)\in\bF[x]$ be the minimal polynomial of $A$.
Then $\phi(x) = q_n(x)$ where $q_n$ is the monic invariant factor of $xI-A$  of highest degree.
\begin{proof}
From remarks~\ref{rem:recayham} and~\ref{rem:minpoldiv} we have that $\phi(x)\Mid q_n(x)$. 
We will show that $q_n(x)\Mid \phi(x)$.
Consider $\phi(x) - \phi(y)\in\bF[x,y]$.
Note that $(x-y)\Mid (\phi(x) - \phi(y))$ (true for any polynomial $\phi(x)$). 
Define $\Phi(x,y)\in \bF[x,y]$ by $\phi(x) - \phi(y) = (x-y)\Phi(x,y)$.
Substituting $x=xI_n$ and $y=A$ preserves this identity since $xI_n$ and $A$ commute.  
Thus, $\phi(xI_n) - \phi(A) = (xI_n-A)\Phi(xI_n,A)$.
Since $\phi(A)=\Theta_n$, we get $\phi(xI_n)=\phi(x)I_n = (xI_n-A)\Phi(xI_n,A)$.
Multiply $\phi(x)I_n = (xI_n-A)\Phi(xI_n,A)$
by $\hat{P}$ from equation~\ref{eq:adjdetred3}.
\begin{equation}
\label{eq:peehatmult}
\hat{P}\phi(x)I_n = \hat{P}(xI_n-A)\Phi(xI_n,A)=q_n(x)\Phi(xI_n,A)
\end{equation}
where the entries of $\hat{P}$ are relatively prime.
Thus, $q_n\Mid \phi(x)$ which was to be shown.
\end{proof}
\end{lemma}

\section*{Rational and Jordan canonical forms.}
\begin{definition}[\bfseries Companion matrix]
\label{def:commat}
\index{rational canonical form!companion matrix}
\index{companion matrix!rational canonical form}
Let $a=(a_0, \ldots, a_{k-1})\in \bF^k$ where $\bF$ is a field of characteristic zero, and let
$p_a (x)=x^k-\sum_{j=0}^{k-1} a_j x^j$.
Define
\begin{equation}
\label{eq:commat}
C(p_a(x)) =
\left[
\begin{array}{cccccc}
0&1&0& \cdots& 0&0 \\
0&0  &1 & \cdots & 0 &0 \\
\vdots&\vdots & \vdots & \vdots & \vdots& \vdots\\
0&0 &0& \cdots & 1&0\\
0&0 &0& \cdots & 0 &1\\
a_0&a_1&a_2&\cdots &a_{k-2} & a_{k-1} 
\end{array}
\right] 
\end{equation} 
to be the {\em companion matrix} of  $p_a(x)$. 
If $a=(a_0)$ then $C(p_a(x)) =(a_0).$\\
\end{definition}

\begin{lemma}[\bfseries Characteristic polynomial of companion matrix]
\label{lem:chapolcommat}
\index{companion matrix!characteristic polynomial}
Let $C=C(p_a(x))\in {\bf M}_{k}(\bF)$ be the companion matrix (\ref{def:commat}) of 
$p_a (x)=x^k-\sum_{j=0}^{k-1} a_j x^j$.
Then $\det(xI-C)=p_a (x)$ and the sequence of nonzero determinantal divisors
of $\,D=xI-C\,$ is   $(f_0^D, f_1^D \ldots , f_{k}^D)=(1,1,\ldots, 1,p_a)$.
\begin{proof}  
First we show that $\det(xI-C)=p_a (x)$.
The proof is by induction.  
The cases $k=1,2$ are easily checked.
Assume that the lemma is true for $k-1$ where $k>2$.
Note that
\begin{equation}
\label{eq:chapolcommat}
D=xI-C =
\left[
\begin{array}{cccccc}
x&-1&0& \cdots& 0&0 \\
0&x &-1 & \cdots & 0 &0 \\
\vdots&\vdots & \vdots & \vdots & \vdots& \vdots\\
0&0 &0& \cdots & -1&0\\
0&0 &0& \cdots & x &-1\\
-a_0&-a_1&-a_2&\cdots &-a_{k-2} &x- a_{k-1} 
\end{array}
\right].
\end{equation} 
Expanding   $\det(xI-C)$ by the first column (Laplace expansion) we get
\begin{equation}
\label{eq:lapexpchamat}
\det(xI-C) = x\det[(xI-C)(1\Mid1)] + (-a_0)(-1)^{k+1}\det[(xI-C)(k\Mid1)].
\end{equation}
Note that by the induction hypothesis,  $\det[(xI-C)(1\Mid1)]=p_{a'}(x)$  where
$a'=(a_1, \ldots, a_{k-1})$.  
Note also that $\det[(xI-C)(k\Mid1)]=(-1)^{k-1}$ and thus 
$(-a_0)(-1)^{k+1}\det[(xI-C)(k\Mid1)] = -a_0$.   
Substituting these results into~\ref{eq:lapexpchamat} proves $\det(xI-C)=p_a (x)$.
It is easily seen in general that $f^D_{k-1}$ is a unit (consider
the submatrix $D(k\Mid 1)$).\\
\end{proof}
\end{lemma}

\begin{lemma}[\bfseries Smith form of a characteristic matrix]
\label{lem:smiforchamat}
\index{Smith form!characteristic matrix}
Let $Q=xI-A$ be the characteristic matrix of $A\in {\bf M}_n(\bF)$. Let
$S={\rm diag}(1, \ldots, 1, q_{k+1}, q_{k+2}, \ldots, q_{n})$
be a Smith form (\ref{def:smifor}, \ref{thm:smifor}) of $Q$ where the $q_{k+j}$,  $j=1,\ldots, n-k$, 
are the non-unit invariant factors. 
Then $S$ is  equivalent to the direct sum
\begin{equation}
\label{eq:equdirsum}
T=
\oplus_{j=1}^{n-k} {\rm diag}(1, \ldots, 1, q_{k+j})
\end{equation}
where ${\rm diag}(1, \ldots, 1, q_{k+j})\in {\bf M}_{d_{k+j}}(\bF[x])$, 
$d_{k+j}=\deg(q_{k+j})$, $j=1,\ldots, n-k$.
\begin{proof}
We  show it is possible to rearrange the diagonal entries of  $S$ to get $T$. 
The theorem requires that ${\rm diag}(1, \ldots, 1, q_{k+j})$ has $d_{k+j} -1$  entries equal to $1$.
Since $\det(S)=\prod_{j=1}^{n-k}q_{k+j}$  we have $\sum_{j=1}^{n-k} d_{k+j} = n$, and hence $\sum_{j=1}^{n-k}( d_{k+j} -1) = n - (n-k) = k$ where $k$ is the number of diagonal elements equal to $1$ in $S$.  
Thus, $T$ is a possible diagonal rearrangement of $S$. Such a rearrangement of $S$ can be achieved by elementary row and column operations and is equivalent to $S$. 
\end{proof}
\end{lemma}
%\hspace*{1 in}\\

\begin{definition}[\bfseries Similarity invariants]
\label{def:siminv}
\index{similarity invariants!defined}
The invariant factors of the characteristic matrix $xI-A$ of $A\in{\bf M}_n(\bF)$ are
called the {\em similarity invariants} of $A$.\\
\end{definition}

\begin{lemma}[\bfseries Companion matrices of non-unit similarity invariants]
\label{lem:commatnonuniinvfac}
\index{companion matrix!of non-unit similarity invariant}
Let $Q=xI-A$ be the characteristic matrix of $A\in {\bf M}_n(\bF)$. Let
\[S={\rm diag}(1, \ldots, 1, q_{k+1}, q_{k+2}, \ldots, q_{n})\]
be a Smith form of $Q$ where the $q_{k+j}$,  $j=1,\ldots, n-k$, 
are the non-unit smilarity invariants of $A$.. 
Then $A$ is similar to $\oplus_{j=1}^{n-k}C(q_{k+j})$ where 
$C(q_{k+j})$ denotes the companion matrix (\ref{def:commat}) of the similarity invariant $q_{k+j}$.
\begin{proof}
We use theorem~\ref{thm:equimpsim} and lemma~\ref{lem:chapolcommat}.
Each term ${\rm diag}(1, \ldots, 1, q_{k+j})$ in 
equation~\ref{eq:equdirsum} is equivalent to the matrix 
$xI_{d_{k+j}} -C(q_{k+j})$  where $d_{k+j}=\deg(q_{k+j})$.
From \ref{eq:equdirsum}, 
we have $Q=xI_n - A$ is equivalent to 
\[\oplus_{j=1}^{n-k}(xI_{d_{k+j}}-C(q_{k+j})) = I_nx - \oplus_{j=1}^{n-k}C(q_{k+j}).\]
By~\ref{thm:equimpsim}, we have $A$ is similar to $\oplus_{j=1}^{n-k}C(q_{k+j})$.\\
\end{proof}
\end{lemma}

\begin{lemma}[\bfseries Companion matrices of elementary divisors]
\label{lem:commatelediv}
\index{companion matrix!of elementary divisors}
We use the notation of lemma~\ref{lem:smiforchamat} and equation \ref{eq:equdirsum}. 
Let $D_t={\rm diag}(1, \ldots, 1, q_{t})$, $k+1\leq t\leq n$,
$q_{t}$  a non-unit invariant factor of $Q=xI-A$,
where ${\rm diag}(1, \ldots, 1, q_{k+j})\in {\bf M}_{d_{k+j}}(\bF[x])$, 
$d_{k+j}=\deg(q_{k+j})$, $j=1,\ldots, n-k$.
Let $q_t=p_1^{e_{t1}}p_2^{e_{t2}} \cdots p_s^{e_{ts}}$
where we assume that $e_{tr}>0$ for $1\leq r\leq s$ so 
all  $p_r^{e_{tr}}$ are elementary divisors of $Q$ (\ref{def:eldiv}).
Then %from \ref{cor:relpridiaent}, 
$D_t$ is equivalent to 
\begin{equation}
\label{eq:commatelediv}
\oplus_{r=1}^{s} {\rm diag}(1, \ldots, 1, p_r^{e_{tr}})
\end{equation}
where ${\rm diag}(1, \ldots, 1, p_r^{e_{tr}})$ is a $\deg(p_r^{e_{tr}})$ square matrix.  
Thus, $C(q_t)$ is similar to $\oplus_{r=1}^{s} C(p_r^{e_{tr}}).$
\begin{proof}
%From lemma~\ref{lem:commatnonuniinvfac} 
We have 
\[
D_t={\rm diag}(1, \ldots, 1, q_{t})={\rm diag}(1, \ldots, 1, p_1^{e_{t1}}p_2^{e_{t2}} \cdots p_s^{e_{ts}}).
\]
Corollary~ \ref{cor:relpridiaent}, a general statement about Smith forms, 
implies that  $D_t$ is a Smith form of
\[
\tilde{D}_t = {\rm diag}(1, \ldots, 1, p_1^{e_{t1}},p_2^{e_{t2}}, \ldots ,p_s^{e_{ts}}).
\]
The hypothesis of corollary~ \ref{cor:relpridiaent} that the 
$p_1^{e_{t1}},p_2^{e_{t2}}, \ldots ,p_s^{e_{ts}}$ are pairwise relatively prime is valid.
Both $\tilde{D}_t$ and  $D_t\in {\bf M}_{n_t}(\bF)$ where $n_t = \deg(q_t)$.
We can now apply the same ``rearranging diagonal entries'' idea used in the proof of lemma~\ref{lem:smiforchamat} 
to show that $\tilde{D}_t$, and hence $D_t$, is equivalent to \ref{eq:commatelediv}.
Analogous to  lemma~\ref{lem:commatnonuniinvfac}, we have that
$C(q_t)$ is similar to $\oplus_{r=1}^{s}C(p_r^{e_{tr}})$.
This completes the proof.\\
\end{proof}
\end{lemma}

\begin{theorem}[\bfseries Frobenius or rational canonical form]
\label{thm:froratcanfor}
\index{Frobenius (rational) canonical form}
Let $xI-A$ be the characteristic matrix of $A\in {\bf M}_n(\bF)$. 
Let ${\mathcal G}$ be the multiset of all non-unit elementary divisors of $xI-A$.
For $g\in {\mathcal G}$ let $\gamma = \deg(g)$.
There exists $P, Q \in {\rm GL}_n(\bF[x])$ and $S\in {\rm GL}_n(F)$ such that 
\begin{equation}
\label{eq:froratcanfor1}
P(xI-A)Q = \bigoplus_{g\in {\mathcal G}} {\rm diag}(1, \ldots , 1, g)
\end{equation}
and
\begin{equation}
\label{eq:froratcanfor2}
S^{-1}AS = \bigoplus_{g\in {\mathcal G}} C(g).
\end{equation}
\begin{proof}
The  proof follows from lemmas~\ref{lem:chapolcommat}, \ref{lem:smiforchamat},  
\ref{lem:commatnonuniinvfac}, and \ref{lem:commatelediv}.\\
\end{proof}
\end{theorem}
\begin{definition}[\bfseries Hypercompanion matrix]
\label{def:hypcommat}
\index{hypercompanion matrix!defined}
Let $\alpha\in \bF$ and 
$p_\alpha (x)= (x-\alpha)^k$.
%Define $H(p_\alpha(x))\in {\bf M}_k(\bF)$ by
Define the {\em hypercompanion matrix} of  $p_\alpha(x)$ by 
\begin{equation}
\label{eq:hypcommat}
H(p_\alpha(x)) =
\left[
\begin{array}{cccccc}
\alpha&1&0& \cdots& 0&0 \\
0&\alpha  &1 & \cdots & 0 &0 \\
\vdots&\vdots & \vdots & \vdots & \vdots& \vdots\\
%0&0 & \cdots &\alpha& 1&0\\
0&0 &0& \cdots & \alpha &1\\
0&0&0&\cdots &0 & \alpha
\end{array}
\right].
\end{equation} 
If $k=1$ then $H(p_\alpha(x)) =(\alpha).$\\
\end{definition}

\begin{remark}[\bfseries Similarity of companion, hypercompanion matrices]
\label{rem:simcomhypmat}
\index{hypercompanion matrix!similar to companion matrix}
Let 
$p_a(x)=x^k-\sum_{j=0}^{k-1} a_j x^j$ where 
$a=(a_0, \ldots, a_{k-1})$ and for $j=0, \ldots, k$,
$-a_j = \binom{k}{j}(-\alpha)^{k-j}x^j$.
Then, by the binomial theorem, $p_a(x)=(x-\alpha)^k$. 
By lemma~\ref{lem:chapolcommat}, $\det(xI-C(p_\alpha))=p_a(x)\equiv(x-\alpha)^k$ and the sequence of nonzero determinantal divisors of $D=xI-C(p_\alpha)$
 is   $(f_0^D, f_1^D \ldots , f_{k}^D)=(1,1,\ldots, 1,p_\alpha)$.
Note that
\begin{equation}
\label{eq:chapolhypcommat}
xI-H(p_\alpha)=
\left[
\begin{array}{cccccc}
x-\alpha&-1&0& \cdots& 0&0 \\
0&x-\alpha &-1 & \cdots & 0 &0 \\
\vdots&\vdots & \vdots & \vdots & \vdots& \vdots\\
0&0 &0& \cdots & -1&0\\
0&0 &0& \cdots & x-\alpha &-1\\
0&0&0&\cdots &0&x-\alpha 
\end{array}
\right].
\end{equation} 
It is easily seen that  $\hat{D}=xI-H(p_\alpha)$  also has its sequence of nonzero
determinantal divisors  
$(f_0^{\hat{D}}, f_1^{\hat{D}}, \ldots , f_{k}^{\hat{D}})=(1,1,\ldots, 1,p_\alpha)$.
Thus, the companion matrix $C(p_\alpha)$ and the hypercompanion matrix $H(p_\alpha)$ have the same similarity invariants and are, therefore, 
similar  matrices  (\ref{thm:equimpsim}).
\end{remark}

\begin{theorem}[\bfseries Jordan canonical form]
\label{thm:jorcanfor}
\index{Jordan canonical form!theorem}
Let $A\in {\bf M}_n(\bF)$.  
Let ${\mathcal G}$ be the multiset of elementary divisors of $xI-A$.
Assume that every  $g\in {\mathcal G}$ is of the form $g=(x-\alpha)^e$, $e>0$.
Then there exists $S\in {\rm GL}_n(F)$ such that 
\begin{equation}
\label{eq:jorcanfor2}
S^{-1}AS = \bigoplus_{g\in {\mathcal G}} H(g)
\end{equation}
where $H(g)$ is the hypercompanion matrix of $g$ (\ref{def:hypcommat}).
\begin{proof}
The proof follows from  theorem~\ref{thm:froratcanfor} (\ref{eq:froratcanfor2}) and
remark~\ref{rem:simcomhypmat}.
\end{proof}
\end{theorem}

\begin{remark}[\bfseries Elementary divisors of  direct sum]
\label{rem:eledivdirsum}
\index{elementary divisors! direct sum}
Let $A=B\oplus C$ where $B\in {\bf M}_b(\bF)$ and $C\in {\bf M}_c(\bF)$ and thus 
$A\in {\bf M}_a(\bF)$ where $a=b+c$.  
Let ${\mathcal G}^B$ and ${\mathcal G}^C$  be the multisets of elementary divisors of $xI_b- B$ and $xI_c - C$.  Let ${\mathcal G} = {\mathcal G}^B \cup {\mathcal G}^C$ be the multiset union of 
${\mathcal G}^B$ and ${\mathcal G}^C$ (\ref{rem:notation}).
From theorem \ref{thm:froratcanfor}
\begin{equation}
\label{eq:equchaB}
P_b(xI_b-B)Q_b = \bigoplus_{g\in {\mathcal G}^B} {\rm diag}(1, \ldots , 1, g)
\end{equation}
and
\begin{equation}
\label{eq:equchaC}
P_c(xI_c-C)Q_c = \bigoplus_{g\in {\mathcal G}^C} {\rm diag}(1, \ldots , 1, g)
\end{equation}
where $P_b,\,\,Q_b \in {\rm GL}_b(\bF[x])$ and $P_c\,\,Q_c \in {\rm GL}_c(\bF[x])$.
We need to prove that, in fact, ${\mathcal G}^A = {\mathcal G}$ where
${\mathcal G}^A$ is the multiset of  elementary divisors of 
$xI_a- A$. 
%Relevent ideas and notation are found in 
%definition~\ref{def:eldiv}, remark~\ref{rem:mulelediv}, lemma~\ref{lem:recfroelediv},
%and remark~\ref{rem:exaeledivinvfac}.

It is clear from equations~\ref{eq:equchaB} and \ref{eq:equchaC} that
\begin{equation}
\label{eq:equchaA}
P_a(xI_a-A)Q_a = \bigoplus_{g\in {\mathcal G}} {\rm diag}(1, \ldots , 1, g)
\end{equation}
where $P_a=P_b \bigoplus P_c$  and $Q_a=Q_b \bigoplus Q_c$. 

In remark~\ref{rem:mulelediv} we saw how to go back and forth between the multiset of elementary divisors and the non-unit invariant factors if we know the rank.  Apply this procedure to the multiset 
${\mathcal G}$ knowing $a$, 
the rank of $A$,  to obtain the associated (with ${\mathcal G}$)
non-unit ``invariant factors'' $q_t=p_1^{e_{t1}}p_2^{e_{t2}} \cdots p_s^{e_{ts}}$. 
At this point, we don't know that ${\mathcal G}$ is  the list of elementary divisors of 
$(xI_a-A)$ so we don't know that this list of ``invariant factors'' is the correct one for $xI_a - A$.
In lemma~\ref{lem:commatelediv} we noted the equivalence of
\[
D_t={\rm diag}(1, \ldots, 1, q_{t})={\rm diag}(1, \ldots, 1, p_1^{e_{t1}}p_2^{e_{t2}} \cdots p_s^{e_{ts}}).
\]
and
\[
\tilde{D}_t = {\rm diag}(1, \ldots, 1, p_1^{e_{t1}},p_2^{e_{t2}}, \ldots ,p_s^{e_{ts}})
\]
by using corollary~ \ref{cor:relpridiaent}, noting  that the hypothesis of that corollary, 
that the $p_i^{e_{ti}}$, $i=1,\ldots s$, are pairwise relatively prime, is valid.
Note that the diagonal matrix of \ref{eq:equchaA} can be rearranged using row and column
matrices so that the 
$\tilde{D}_t = {\rm diag}(1, \ldots, 1, p_1^{e_{t1}},p_2^{e_{t2}}, \ldots ,p_s^{e_{ts}})$
are grouped together.
Using row and column operations, convert these $\tilde{D}_t$ to the equivalent $D_t$.
This shows that, in fact, the $q_t$ are the invariant factors of $xI_a-A$.  
Thus, ${\mathcal G}^A = {\mathcal G}$.
 \end{remark}
\newpage
\mbox{}
\newpage
\hypertarget{index}{ }
\printindex
\newpage
\centerline{NOTES}
%\newpage
%\centerline{NOTES}
%\newpage
%\centerline{NOTES}
\end{document}